A thesis submitted for the degree of
Doctor of Philosophy of
*The Australian National University*

# *Second Maximal Subgroups of the Finite Alternating and Symmetric Groups*

*Alberto Basile*

April 2001

*Al mio caro papà*

## Statement

The results in this thesis are my own except where otherwise stated.

# Acknowledgments

I thank my mum and my sisters Laura, Eveline and Simonetta for their unlimited love.

I thank my supervisor, Laci Kovács, both for his constant and patient assistance and advice throughout the preparation of this thesis, and for all his help and emotional support during my stay here in Australia.

I gratefully acknowledge the financial support of the Australian National University with the ANU PhD Scholarship, and of the Australian Government Department of Employment, Education, Training and Youth Affairs with the International Postgraduate Research Scholarship.

Finally, because of their friendship and of the consequent more or less explicit help, I thank Enrico Saccardo, Jansyl Lagos, Kilmeny Beckering Vinckers, Lydia Ronnenkamp, Giuseppe Jurman and Amanda McPhee.

# Abstract

A subgroup of a finite group $G$ is said to be *second maximal* if it is maximal in every maximal subgroup of $G$ that contains it. A question which has received considerable attention asks: can every positive integer occur as the number of the maximal subgroups that contain a given second maximal subgroup in some finite group $G$? Various reduction arguments are available except when $G$ is almost simple. Following the classification of the finite simple groups, finite almost simple groups fall into three categories: alternating and symmetric groups, almost simple groups of Lie type, sporadic groups and automorphism groups of sporadic groups. This thesis investigates the finite alternating and symmetric groups, and finds that in such groups, except three well known examples, no second maximal subgroup can be contained in more than 3 maximal subgroups.

# Contents

# List of Tables

# List of Figures

CHAPTER 1

# Introduction

An almost simple group is a group having a unique minimal normal subgroup which is non abelian and simple. Up to isomorphism, there are only finitely many almost simple groups containing and normalizing any given finite simple group: they are the subgroups of its automorphism group containing all the inner automorphisms. Following the classification of the finite simple groups (henceforth referred to as CFSG) one of the major areas of research in group theory today is the investigation of the subgroups of the finite almost simple groups, and in particular, the determination of their maximal subgroups. The next natural step in this process and a first important step towards understanding the local structure of the subgroup lattice of all finite groups, is the investigation of the second maximal subgroups of the finite almost simple groups, namely the maximal elements of what remains in the related subgroup lattices after removing the maximal subgroups. According to the CFSG, the finite almost simple groups fall into three classes:

the alternating groups $\mathrm{A}_n$ and the symmetric groups $\mathrm{S}_n$ ($n \geqslant 5$);
the finite almost simple groups of Lie type;
the 26 sporadic simple groups and their automorphisms groups.

In this thesis we consider the second maximal subgroups of the finite alternating and symmetric groups, the maximal subgroups having been discussed in [**LPS87**]. Pálfy [**Pál88**] started things off by dealing with second maximal subgroups of the alternating groups of prime degree. This was extended in the Laurea Thesis [**Bas96**] of the author to the alternating and symmetric groups of prime power degree. In the same work results were obtained for intransitive imprimitive second maximal subgroups in arbitrary degree. However, except in Theorem D below, we do not rely on [**Bas96**] here: while it is available on the Internet[1], it is yet to be published in print. Instead, we rely on the results in [**Pál88**]; on the results of Praeger in [**Pra90**] where inclusions among primitive groups are discussed up to permutational isomorphism; and on the fundamental results of Liebeck, Praeger, Saxl on the factorizations of the almost simple groups in [**LPS90**]. We do not shy away

[1] `http://room.anu.edu.au/~abit/`

from using the Schreier Conjecture[2] a few times because we view it as a byproduct of the CFSG.

The overall picture of the thesis is as follows. In Chapter 2 we present a more detailed setting for the work in this thesis and introduce some specific notation while discussing a few preliminary results.

Chapter 3 contains a collection of results about well known permutation groups. Namely, we provide the arguments used to determine their normalizers (in the parent symmetric group), their parity (whether or not they are contained in the alternating group), and possibly notes about their degree of transitivity/primitivity.

In Chapter 4 we begin our investigation of the second maximal subgroups of the alternating and symmetric groups and we provide a description of all the intransitive second maximal subgroups which are contained in more than two maximal subgroups. The main result of this chapter is the following.

**Theorem A.** *A second maximal subgroup of a finite alternating or symmetric group of degree at least* $5$ *which is contained in more than* $2$ *maximal subgroups is either intransitive with* $3$ *orbits and contained in precisely* $3$ *maximal intransitive subgroups, or it is the stabilizer of a point in* $\mathrm{PSL}(3,2)$ *of degree* $7$ *as in* 2.6.1*, or it is the stabilizer of a point in a* $2$*-primitive maximal subgroup of an alternating group. In the latter case, it is contained in the stabilizer of that point in the alternating group, in another* $2$*-primitive group which is conjugate to the given one in the parent symmetric group (but not in the alternating group), and in nothing else. So, an intransitive second maximal subgroup is never contained in more than* $3$ *maximal subgroups.*

In Chapter 5 we introduce cumbersome yet successful notations to handle non almost simple primitive groups. Successful because we prove the second main result of this thesis:

**Theorem B.** *A primitive second maximal subgroup of a finite alternating or symmetric group which is non almost simple or contained in some non almost simple subgroup is never contained in more than* $2$ *maximal subgroups, unless it is one of the three examples of the kind discovered by Feit and listed in* [**Pál88**, Table II] *and here in* 5.2.1.

Our efforts culminate in Chapter 6 where an exhaustive analysis of the primitive almost simple second maximal subgroups is carried over to yield the third but perhaps most important theorem of this work:

[2]The Schreier Conjecture states that the outer automorphism group of any finite simple group is solvable.

**Theorem C.** *A primitive almost simple second maximal subgroup of a finite alternating or symmetric group which is only contained in almost simple subgroups of which there are at least two, can be recovered at once from Table* 6.B *and is contained in at most* 3 *maximal subgroups.*

There is also an appendix which is not intended to be regarded as part of this thesis. It is a preliminary version of a separate, still unpublished work of the author, to which we refer in Chapter 4. It is bound together with the thesis solely for the convenience of the reader.

As a final remark, note that Chapter 5 and Chapter 6 together provide a complete analysis of the primitive second maximal subgroups of the finite alternating and symmetric groups. In Theorem III.3.1 at p. 62 of [**Bas96**] the author shows that a transitive imprimitive second maximal subgroup of a finite alternating or symmetric group of degree at least 5 is never contained in more than 3 maximal subgroups. Given that a permutation group is either intransitive or transitive imprimitive or primitive, all these results together yield the following conclusive result concerning the second maximal subgroups of finite symmetric and alternating groups.

**Theorem D.** *A second maximal subgroup of a finite alternating or symmetric group of degree at least* 5 *is never contained in more than* 3 *maximal subgroups, unless it is one of the three examples of Feit and Pálfy. Consequently, such a subgroup is never contained in more than* 11 *maximal subgroups.*

A long outstanding problem in Universal Algebra is whether every finite lattice occurs as the congruence lattice of a finite algebra. In the most recent step of a long chain of reductions, Baddeley and Lucchini listed a number of residual questions about almost simple groups [**BL97**]. Theorem D may be viewed as a partial answer to one of these questions. However, elaborating this issue would take us too far from the present context; see Pálfy [**Pál95**] and his featured review MR98j:20022 of [**BL97**] for a more detailed discussion.

## 1.1. General notations

Sets and their elements are written as

$$\mathsf{A}, \mathsf{B}, \mathsf{C}, \dots, \mathsf{X}, \mathsf{Y}, \mathsf{Z}, \qquad \mathsf{a}, \mathsf{b}, \mathsf{c}, \dots, \mathsf{x}, \mathsf{y}, \mathsf{z}.$$

The cardinality of a set $\mathsf{X}$ is denoted by $|\mathsf{X}|$; this is also called the *order* of $\mathsf{X}$ in case $\mathsf{X}$ is finite. A *partition* of a set $\mathsf{X}$ is a set of non empty subsets of $\mathsf{X}$ whose pairwise intersection is empty and whose union equals to $\mathsf{X}$. The cardinality (or order) of a partition is the cardinality of the partition regarded as a set.

Partially ordered sets and hence lattices are usually denoted by

$$\mathcal{A}, \mathcal{B}, \mathcal{C}, \ldots, \mathcal{X}, \mathcal{Y}, \mathcal{Z}.$$

Groups and their elements are written as

$$A, B, C, \ldots, X, Y, Z, \qquad a, b, c, \ldots, x, y, z.$$

Richer algebraic structures and their elements may be shown as

$$\mathfrak{A}, \mathfrak{B}, \mathfrak{C}, \ldots, \mathfrak{X}, \mathfrak{Y}, \mathfrak{Z}, \qquad \mathfrak{a}, \mathfrak{b}, \mathfrak{c}, \ldots, \mathfrak{x}, \mathfrak{y}, \mathfrak{z},$$

but for fields we prefer this special font:

$$\mathbb{F}, \mathbb{G}, \mathbb{L}, \mathbb{P}, \mathbb{F}_q, \ldots$$

Functions are usually written in bold, both latin and greek letters:

$$\boldsymbol{f}, \boldsymbol{g}, \boldsymbol{h}, \ldots \qquad \boldsymbol{\alpha}, \boldsymbol{\beta}, \boldsymbol{\gamma}, \ldots$$

A function $\boldsymbol{f}$ comes together with a domain $\operatorname{Dom} \boldsymbol{f}$ and a codomain $\operatorname{Cod} \boldsymbol{f}$ which is not necessarily equal to the image $\operatorname{Im} \boldsymbol{f}$:

$$\boldsymbol{f} : \operatorname{Dom} \boldsymbol{f} \longrightarrow \operatorname{Cod} \boldsymbol{f}, \qquad \operatorname{Im} \boldsymbol{f} = \left\{ \mathsf{x}\boldsymbol{f} \;\middle|\; \mathsf{x} \in \operatorname{Dom} \boldsymbol{f} \right\}.$$

As we have just seen, we tend to apply functions on the right of their arguments and to compose functions accordingly. For example for $\mathsf{x} \in \operatorname{Dom} \boldsymbol{f}$ and $\operatorname{Cod} \boldsymbol{f} = \operatorname{Dom} \boldsymbol{g}$ we have

$$x(\boldsymbol{fg}) = (x\boldsymbol{f})\boldsymbol{g} = x\boldsymbol{fg}.$$

As common the set of functions from $\mathsf{X}$ to $\mathsf{Y}$ is denoted by $\mathsf{Y}^{\mathsf{X}}$.

We regard non negative integers as sets:

$$0 = \emptyset \qquad \text{the empty set,}$$

$$1 = \{0\}, \quad 2 = \{0, 1\}, \quad \ldots, \quad l = \{0, 1, \ldots, l-1\}, \quad \ldots$$

This has a few unhappy consequences. For example, if we identify $\mathsf{Y} \times \mathsf{Y}$ with $\mathsf{Y}^2$, the set of functions from 2 to $\mathsf{Y}$, then we better write its generic element as $(\mathsf{y}_0, \mathsf{y}_1)$ avoiding a misleading expression $(\mathsf{y}_1, \mathsf{y}_2)$. For similar reasons, we regret we have to regard $\mathrm{S}_n$ and $\mathrm{A}_n$ as the symmetric and alternating group on $n$ rather than on $\{1, 2, \ldots, n\}$.

We may occasionally speak of $l$-tuples or vectors (with $l$ entries) meaning functions whose domain is the positive integer $l$. Then we write

$$\vec{a} = (a_0, a_1, \ldots, a_{l-1}).$$

where we set $a_i = i\vec{a}$ for each $i \in l$. We also write the image of $\vec{a}$ as

$$\tilde{a} = l\vec{a} = \operatorname{Im} \vec{a} = \{a_0, a_1, \ldots, a_{l-1}\}.$$

The set of $l$-tuples with entries in a set $\mathsf{A}$ is $\mathsf{A}^l$. The symmetric groups $\mathrm{S}_l$ and $\operatorname{Sym}\mathsf{A}$ act respectively on the left and on the right of $\mathsf{A}^l$ by composition of functions, so for $\boldsymbol{s} \in \mathrm{S}_l$, $\vec{a} \in \mathsf{A}^l$, $\boldsymbol{\alpha} \in \operatorname{Sym}\mathsf{A}$ one has

$$l \xrightarrow{\boldsymbol{s}} l \xrightarrow{\vec{a}} \mathsf{A} \xrightarrow{\boldsymbol{\alpha}} \mathsf{A}$$
$$l \xrightarrow{\boldsymbol{s}\vec{a}\boldsymbol{\alpha}} \mathsf{A}$$

Since we tend to work with right actions (as we apply functions on the right), we prefer to consider the induced right action of $\mathrm{S}_l$:

$$\vec{a} \cdot \boldsymbol{\sigma} := \boldsymbol{\sigma}^{-1}\vec{a} = (a_{0\boldsymbol{\sigma}^{-1}}, \ldots, a_{(l-1)\boldsymbol{\sigma}^{-1}}).$$

This embeds $\mathrm{S}_l$ and $\operatorname{Sym}\mathsf{A}$ in two commuting subgroups of $\operatorname{Sym}(\mathsf{A}^l)$. The second embedding is also known as diagonal action (see 3.5).

## 1.2. This is "well known"

Throughout this work we use the words "well known" for anything which is not necessarily well known, yet it is well known how to find it in the literature. More precisely, we state here that something is well known if it can be easily found in

- **ATLAS [ATLAS]**,
- **Robinson's *A course in the theory of Groups* [Rob93]**,
- **Beachy's *Abstract Algebra on Line* [Bea96]**
- **Suzuki's *Group Theory I* [Suz82]**,
- **Dixon and Mortimer's *Permutation Groups* [DM96]**.

Clearly, anything stated without a reference or further explanation should be regarded as well known in our sense. For example, we will soon talk about actions of groups on sets or about wreath products in their imprimitive action without any introduction to the topic. This is because not only are they well described in Robinson but you also find "action of a group on a set" and "wreath product" in its index.

As a further advice we remind the reader that information about finite simple groups of a given order is likely to be found in the ATLAS; information about general group theory, construction of new groups is likely to be found both in Robinson and Suzuki; details about abstract algebras like fields, rings and so on (for the little that we will mention to them) are to be found on line at Beachy's site. Finally, facts regarding families of simple groups like the Mathieu simple groups or the projective special linear groups for example, may be sought in this order: the introduction of ATLAS, Dixon and Mortimer, Suzuki, and Robinson.

CHAPTER 2

# Preliminaries

## 2.1. Second maximal subgroups

Let $G$ be a finite group and

$$\mathcal{S} = \{\ H < G \ \mid\ H \text{ is not a maximal subgroup of } G\ \}.$$

The *second maximal* subgroups of $G$ are the maximal elements of $\mathcal{S}$. Clearly, if $N$ is a second maximal subgroup of $G$ and if $N < M < G$, then $N$ is a maximal subgroup of $M$ and $M$ is a maximal subgroup of $G$. As this work is mainly concerned with maximal and second maximal subgroups, we introduce this notation: when $H$ is a maximal subgroup of $G$, write

$$H \lessdot G.$$

Also, we borrow some terminology from Lattice Theory. For subgroups $N$ and $M$ of $G$, the *interval* between $N$ and $M$ is the sublattice of the subgroup lattice of $G$ made of the subgroups of $M$ containing $N$:

$$[N \div M] = \{\ R \leqslant G \ \mid\ N \leqslant R \leqslant M\ \}.$$

If $N$ is a second maximal subgroup of $G$, then the interval $[N \div G]$ has the shape below:

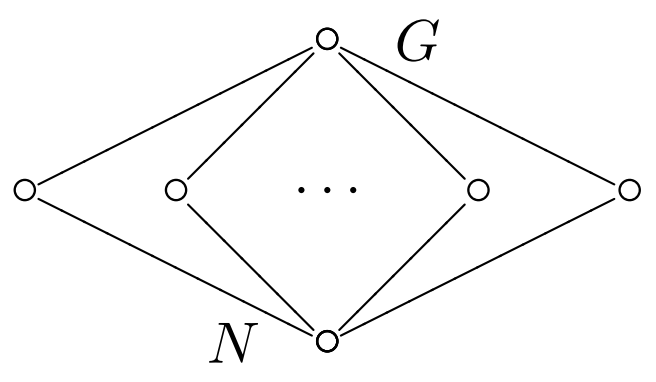

A lattice as above with a smallest, a largest and $r$ pairwise incomparable elements, is denoted by $\mathcal{M}_r$. Therefore $N$ is a second maximal subgroup of $G$ if and only if $[N \div G] = \mathcal{M}_r$ for some positive integer $r$.

We aim to classify the second maximal subgroups of the finite symmetric and alternating groups. Clearly, $S_3$ has just one second maximal subgroup while $A_4$ has exactly three of them (see Figure 2.A). It is easy to check that the second maximal subgroups of $S_4$ are exactly the subgroups of order 3 or 4. Among them, there is the Klein subgroup $V_4$ which is the unique subgroup of order 4 of $A_4$ and is contained in three more subgroups of $S_4$, namely, the three dihedral groups $D_8$ of order 8.

FIGURE 2.A. The subgroup lattice of $A_4$

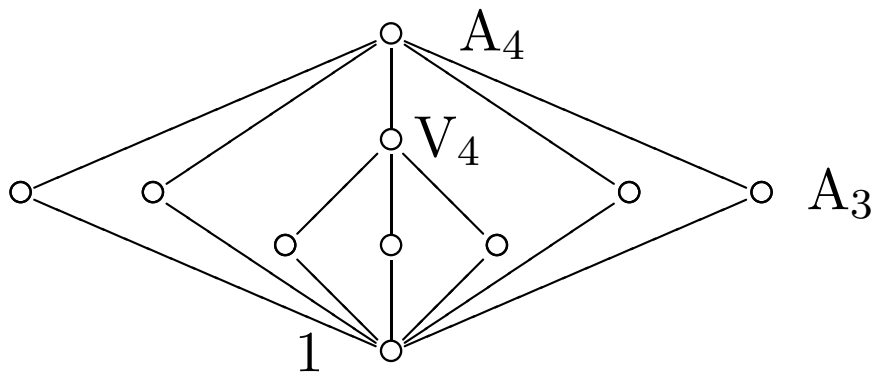

Because of the previous discussion, throughout this thesis we write $\boldsymbol{S}$ for a finite symmetric group of degree at least 5; accordingly, $\boldsymbol{A}$ denotes the related alternating group:

$$\boldsymbol{S} \equiv \operatorname{Sym}\Omega \cong \mathrm{S}_n, \qquad \boldsymbol{A} \equiv \operatorname{Alt}\Omega \cong \mathrm{A}_n, \qquad n \geqslant 5.$$

We seek the second maximal subgroups of $\boldsymbol{S}$ or $\boldsymbol{A}$, especially those contained in a large number of maximal subgroups. Therefore we have little interest in the second maximal subgroups of $\boldsymbol{S}$ which are contained in $\boldsymbol{A}$ because we know already that they are contained in at most one more maximal subgroup of $\boldsymbol{S}$:

2.1.1. PROPOSITION. *Let $H$ be a second maximal subgroup of $\boldsymbol{S}$ and a subgroup of $\boldsymbol{A}$. If $H \lessdot G \neq \boldsymbol{A}$, then $G = \mathrm{N}_{\boldsymbol{S}}\, H$.*

PROOF. Since $H = G \cap \boldsymbol{A}$, $H$ is normal in $G$. Because the degree is at least 5, $H$ is not normal in $\boldsymbol{S}$. Thus $G$ is the normalizer of $H$ in $\boldsymbol{S}$. □

For this reason, we reserve the symbol $\boldsymbol{H}$ either for a second maximal subgroup of $\boldsymbol{A}$, or for a second maximal subgroup of $\boldsymbol{S}$ which is not contained in $\boldsymbol{A}$. Now put $\boldsymbol{U} = \boldsymbol{H}\boldsymbol{A}$ and note that $\boldsymbol{U} \in \{\boldsymbol{S}, \boldsymbol{A}\}$. Then $\boldsymbol{H}$ is a second maximal subgroup of $\boldsymbol{U}$:

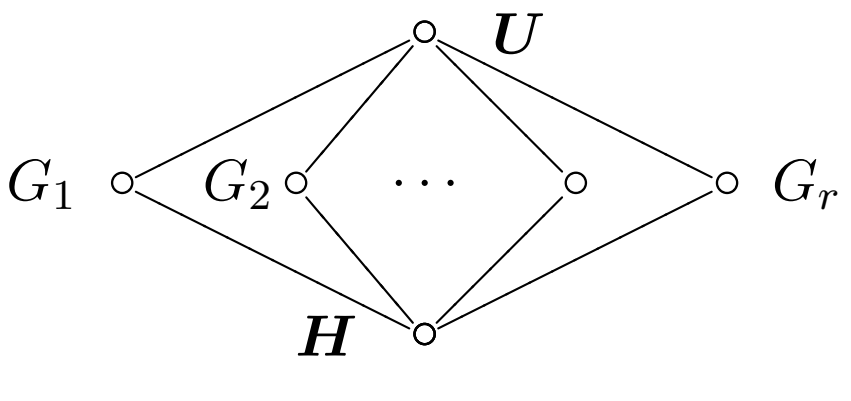

In the following chapters we determine the $\boldsymbol{H}$ as above and possibly the overgroups $G_1, \ldots, G_r$ or at least an upper bound for $r$.

## 2.2. Even, odd, low and high

We separate the subgroups of $\boldsymbol{S}$ into *even* and *odd* subgroups, the even ones being the subgroups of $\boldsymbol{A}$. Accordingly, for an arbitrary permutation representation

$$\varphi : G \longrightarrow \operatorname{Sym} \mathsf{X}$$

we define $G^{\mathrm{e}}$, the even part of $G$, as the complete inverse image by $\varphi$ of $\mathrm{Alt}\,\mathsf{X}$. In particular, if $G$ is a subgroup of $\mathrm{Sym}\,\mathsf{X}$, then we have $G^{\mathrm{e}} = G \cap \mathrm{Alt}\,\mathsf{X}$.

We will also divide the subgroups of $\boldsymbol{S}$ into *low* and *high* subgroups, where we say that a subgroup $G$ of $\boldsymbol{S}$ is low if there are subgroups $K$ and $L$ such that

$$G < K < L < G\boldsymbol{A}.$$

Conversely, a high subgroup $G$ is one such that the interval $[G \div G\boldsymbol{A}]$ is $\mathcal{M}_r$ for some integer $r \geqslant -1$; the trivial cases with $r = 0$, $r = -1$ corresponding to $G$ being maximal in $G\boldsymbol{A}$ and to $G = G\boldsymbol{A}$ respectively. Therefore $G$ is a non trivial high subgroup of $\boldsymbol{S}$ if and only if $G$ is a second maximal subgroup of $G\boldsymbol{A}$.

This terminology may be extended to any almost simple group $S$ with socle $A$. The classification of the non trivial high subgroups of the finite almost simple groups is still a great challenge and may be regarded as a necessary step towards understanding the finite lattices which occur as intervals in subgroup lattices of finite groups; for some history of this long outstanding problem see [**Pál95**].

## 2.3. Stabilizers

Suppose that $G$ acts on $\Omega$. The stabilizer in $G$ of a point $\omega$ is denoted by $\mathrm{St}_G\,\omega$ or by the more common $G_\omega$. Trouble may arise referring to the stabilizer of a subset $\mathsf{X}$ of $\Omega$ without specifying whether it is the setwise stabilizer or the pointwise stabilizer. However, we tend to regard the subscript of $G$ as an element of a set on which there must be a canonic action of $G$. For example, if the subscript is a subset $\mathsf{X}$, the canonic action is the one on the set of all subsets of $\Omega$. Then we want to regard $G_{\mathsf{X}}$ as the setwise stabilizer of $\mathsf{X}$. Similarly, if the subscript is an ordered pair $(\omega_1, \omega_2)$, the canonic action is the diagonal action on $\Omega \times \Omega$. Then we want to regard $G_{(\omega_1,\omega_2)}$ as the stabilizer of both $\omega_1$ and $\omega_2$, that is, $G_{\omega_1} \cap G_{\omega_2}$.

Accordingly, if $\mathsf{X}_1$, $\mathsf{X}_2$ are subsets of $\Omega$, we define

$$G_{(\mathsf{X}_1,\mathsf{X}_2)} \equiv \{\, g \in G \mid (\mathsf{X}_1 g, \mathsf{X}_2 g) = (\mathsf{X}_1, \mathsf{X}_2) \,\} \text{ and}$$
$$G_{\{\mathsf{X}_1,\mathsf{X}_2\}} \equiv \{\, g \in G \mid \{\mathsf{X}_1 g, \mathsf{X}_2 g\} = \{\mathsf{X}_1, \mathsf{X}_2\} \,\}.$$

The first one is the largest subgroup of $G$ which preserves the ordered pair $(\mathsf{X}_1, \mathsf{X}_2)$ and the second is the largest that preserves the set $\{\mathsf{X}_1, \mathsf{X}_2\}$. These definitions extend obviously to any number of subsets. In particular, $G_{\mathsf{X}}$, $G_{(\mathsf{X})}$ and $G_{\{\mathsf{X}\}}$ should all denote the same subgroup of $G$, namely the setwise stabilizer of $\mathsf{X}$. We will only use the first one in this case but we feel obliged to remind the reader that the other two notations may be reserved by other authors for pointwise stabilizers. We refuse to do so as the parentheses convey a meaning which we feel should be preserved. Instead, referring to the pointwise stabilizer of $\mathsf{X}$ we write the longer but unambiguous $\bigcap_{\omega\in\mathsf{X}} G_\omega$.

Clearly, if $\mathsf{X}_1$ and $\mathsf{X}_2$ are subsets of $\Omega$ of different order, then $G_{\{\mathsf{X}_1,\mathsf{X}_2\}} = G_{(\mathsf{X}_1,\mathsf{X}_2)}$. For this reason, in referring to $G_{\{\mathsf{X}_1,\mathsf{X}_2\}}$, we usually assume that $\mathsf{X}_1$ and $\mathsf{X}_2$ have the same order. To deal with cases like this, we call *equipartitions* the partitions in subsets of the same order. There is one and only one equipartition of $\Omega$ of order 1; and there is one and only one equipartition in subsets of order 1. Those two are the *trivial equipartition.*

Stabilizers in $\boldsymbol{S}$ or $\boldsymbol{A}$ of partitions of $\Omega$ in two subsets attract our attention because in general they are maximal subgroups of $\boldsymbol{S}$ or $\boldsymbol{A}$. The same holds for stabilizers of non trivial equipartitions. We explicitly state this in the two propositions below together with facts which are well known or to be found in [**LPS87**].

2.3.1. Proposition. *Let* $\{\mathsf{X}_1,\dots,\mathsf{X}_l\}$ *be a partition of* $\Omega$.

(1) $\boldsymbol{S}_{(\mathsf{X}_1,\dots,\mathsf{X}_l)} \cong (\mathrm{Sym}\,\mathsf{X}_1)\times\cdots\times(\mathrm{Sym}\,\mathsf{X}_l)$.
(2) *If* $l=2$ *and* $|\mathsf{X}_1|\neq|\mathsf{X}_2|$ *then* $\boldsymbol{S}_{(\mathsf{X}_1,\mathsf{X}_2)} \lessdot \boldsymbol{S}$ *and* $\boldsymbol{A}_{(\mathsf{X}_1,\mathsf{X}_2)} \lessdot \boldsymbol{A}$.
(3) *If* $l=2$ *and* $|\mathsf{X}_1|=|\mathsf{X}_2|$, *then* $\boldsymbol{S}_{(\mathsf{X}_1,\mathsf{X}_2)} \lessdot \boldsymbol{S}_{\{\mathsf{X}_1,\mathsf{X}_2\}}$.

2.3.2. Proposition. *Let* $\mathsf{Z}$ *be a non trivial equipartition of* $\Omega$ *of order* $l$ *in subsets of order* $m$.

(1) $\boldsymbol{S}_\mathsf{Z}$ *is a maximal subgroup of* $\boldsymbol{S}$ *isomorphic to* $\mathrm{S}_m \,\mathrm{wr}\, \mathrm{S}_l$.
(2) $\boldsymbol{A}_\mathsf{Z}$ *is transitive.*
(3) $\boldsymbol{A}_\mathsf{Z}$ *is a maximal subgroup of* $\boldsymbol{A}$ *unless* $\boldsymbol{A}=\mathrm{A}_8$, $l=4$.

The isomorphism of the first part of the proposition above is also completely described in Chapter 5. Regarding the last part, one shows that if $\boldsymbol{A}=\mathrm{A}_8$ and $l=4$, then $\boldsymbol{A}_\mathsf{Z}$ is a second maximal subgroup of $\boldsymbol{A}$ and that the maximal subgroups of $\boldsymbol{A}$ containing $\boldsymbol{A}_\mathsf{Z}$ are 2 copies of $\mathrm{AGL}(3,2)$, conjugate in $\boldsymbol{S}$ but not in $\boldsymbol{A}$. For each copy, $\boldsymbol{A}_\mathsf{Z}$ corresponds to the stabilizer of a set of parallel lines.

## 2.4. Blocks and imprimitivity systems

Let $G$ be a transitive permutation group on a set $\Omega$. A subset of $\Omega$ is a *block* for $G$ if its orbit under the action of $G$ on the set of the subsets of $\Omega$ is a partition of $\Omega$. Clearly, the singletons and $\Omega$ itself are blocks. These are the trivial blocks. The *imprimitivity systems* of $G$ are the orbits of the non trivial blocks for $G$. Note that an imprimitivity system is an equipartition of $\Omega$ which is preserved by $G$. As we call trivial the equipartition made of singletons and the equipartition $\{\Omega\}$, the imprimitivity systems of $G$ are the non trivial equipartitions of $\Omega$ preserved by $G$.

We may talk of blocks or imprimitivity systems with an intransitive group. We mean of course that those are blocks or imprimitivity systems related to some bigger transitive group.

The following two facts are well known.

2.4.1. PROPOSITION. *Let $G$ be a transitive permutation group and let $H$ be a normal subgroup of $G$. The orbits of $H$ are blocks for $G$.*

2.4.2. PROPOSITION. *Let $G$ be a transitive permutation group on $\Omega$ and let $\omega$ be a point of $\Omega$. There is a bijection between the set of the equipartitions of $\Omega$ preserved by $G$ and the interval $[G_\omega \div G]$ where an equipartition $\mathsf{Z}$ preserved by $G$ corresponds to the stabilizer in $G$ of the unique element of $\mathsf{Z}$ containing $\omega$. In this correspondence, $G_\omega$ and $G$ correspond to the trivial equipartitions.*

An *imprimitive* group is a transitive group for which there is at least one imprimitivity system. By 2.3.2 a maximal imprimitive subgroup of a symmetric group is a wreath product in imprimitive action of two (non trivial) smaller symmetric groups. In particular, it contains *transpositions*, that is, permutations which move exactly 2 points. Similarly, a maximal imprimitive subgroup of an alternating group contains *double transpositions* , that is, products of two disjoint transpositions.

To better describe the relation between non trivial equipartitions and maximal imprimitive subgroups of the symmetric group we shall require the following notation. Suppose $G$ acts transitively on a set $\mathsf{X}$ by $\varphi : G \longrightarrow \operatorname{Sym}\mathsf{X}$ and that $\mathsf{x} \in \mathsf{X}$. As a lattice, the interval $[G_{\mathsf{x}} \div G]$ does not depend on $\mathsf{x}$, and we write $\mathcal{C}(G, \varphi, \mathsf{X})$ for an unspecified representative of this Lattice isomorphism class. We may write just $\mathcal{C}(G, \mathsf{X})$ or even $\mathcal{C}(G)$ if there are no doubts on the action involved. Two more definitions are needed now: the linear and vertical sum of lattices. Given two lattices $\mathcal{L}_1$ and $\mathcal{L}_2$, we may endow the disjoint union $\mathcal{L}_1 \sqcup \mathcal{L}_2$ with an order which extends the orders of $\mathcal{L}_1$, $\mathcal{L}_2$ and sets all the elements of $\mathcal{L}_1$ less than any element of $\mathcal{L}_2$. With respect to this partial order the disjoint union $\mathcal{L}_1 \sqcup \mathcal{L}_2$ is obviously a lattice called the *linear sum* of the two lattices. The *vertical sum* $\mathcal{L}_1 \bar{\oplus} \mathcal{L}_2$ is the quotient of the linear sum modulo the equivalence relation $\sim$ defined by: $\mathsf{a} \sim \mathsf{b}$ if either

$$\mathsf{a} = \mathsf{b}, \qquad \text{or}$$

$$\mathsf{x}_1 \leqslant \mathsf{a} \leqslant \mathsf{x}_2, \quad \mathsf{x}_1 \leqslant \mathsf{b} \leqslant \mathsf{x}_2 \qquad \text{whenever } \mathsf{x}_1 \in \mathcal{L}_1, \quad \mathsf{x}_2 \in \mathcal{L}_2.$$

Note that all equivalence classes of $\sim$ are singletons except perhaps one: when $\mathcal{L}_1$ has a largest element and $\mathcal{L}_2$ has a smallest element, those two elements form a single equivalence class and become identified in $\mathcal{L}_1 \bar{\oplus} \mathcal{L}_2$. One naturally identifies $\mathcal{L}_1$ and $\mathcal{L}_2$ with their images under $\mathsf{x} \mapsto [\mathsf{x}]_\sim$ so that they both may be regarded as sublattices (indeed as intervals) of their vertical sum.

2.4.3. THEOREM. *If $H$, $K$ act faithfully and transitively on* $\mathsf{X}$, $\mathsf{Y}$ *respectively* (*both* $\mathsf{X}$, $\mathsf{Y}$ *of order at least* 2), *then according to the imprimitive action of* $H \operatorname{wr} K$ *on* $\mathsf{X} \times \mathsf{Y}$ *we have*

$$\mathcal{C}(H \operatorname{wr} K) = \mathcal{C}(H) \bar{\oplus} \mathcal{C}(K).$$

PROOF. To save some writing put $W = H \operatorname{wr} K$; $W$ is transitive on $\mathsf{X} \times \mathsf{Y}$ and acts transitively on $\mathsf{Y}$ via the canonical projection on the top group $\boldsymbol{\pi} : W \longrightarrow K$. Pick $\mathsf{x} \in \mathsf{X}$ and $\mathsf{y} \in \mathsf{Y}$, it is enough to show that

$$\big[W_{(\mathsf{x},\mathsf{y})} \div W\big] = [H_x \div H] \bar{\oplus} [K_{\mathsf{y}} \div K].$$

To this end, note that if we write $H \operatorname{wr} K_{\mathsf{y}}$ for the subgroup of $W$ generated by $H^{\mathsf{Y}-\{\mathsf{y}\}}$ and $K_{\mathsf{y}}$, then

$$W_{\mathsf{y}} = H \times (H \operatorname{wr} K_{\mathsf{y}}) \qquad \text{and} \qquad W_{(\mathsf{x},\mathsf{y})} = H_{\mathsf{x}} \times (H \operatorname{wr} K_{\mathsf{y}}).$$

Since $W_{\mathsf{y}} = \boldsymbol{\pi}^{-1}(K_{\mathsf{y}})$, $[W_{\mathsf{y}} \div W] = [K_{\mathsf{y}} \div K]$. Similarly, consider $\boldsymbol{\pi}_{\mathsf{y}}$, the left canonical projection of $W_{\mathsf{y}} = H \times (H \operatorname{wr} K_{\mathsf{y}})$ onto $H$. Since $W_{(\mathsf{x},\mathsf{y})} = \boldsymbol{\pi}_{\mathsf{y}}^{-1}(H_{\mathsf{x}})$, $\big[W_{(\mathsf{x},\mathsf{y})} \div W_{\mathsf{y}}\big] = [H_{\mathsf{x}} \div H]$. Therefore we only need to show that if

$$W_{(\mathsf{x},\mathsf{y})} < T \leqslant W \qquad \text{and} \qquad T \nleqslant W_{\mathsf{y}},$$

then $W_{\mathsf{y}} \leqslant T$. In fact, say $t \in T$ such that $\mathsf{y}t \neq \mathsf{y}$. Provided $B$ is the base subgroup of $H \operatorname{wr} K_{\mathsf{y}}$,

$$t^{-1}(H \times \{1\})t \leqslant B \leqslant W_{(\mathsf{x},\mathsf{y})} \leqslant T.$$

Thus $H \times \{1\} \leqslant tTt^{-1} = T$ and so $T \geqslant (H \times \{1\})W_{(\mathsf{x},\mathsf{y})} = W_{\mathsf{y}}$, as claimed. □

2.4.4. REMARK. $H$, $\mathsf{X}$, $K$, $\mathsf{Y}$ *as in theorem above. If* $H^e$ *is transitive or* $|\mathsf{Y}| > 2$, *then again*

$$\mathcal{C}\big((H \operatorname{wr} K)^e\big) = \mathcal{C}(H) \bar{\oplus} \mathcal{C}(K).$$

PROOF. Use $W$, $\mathsf{x}$, $\mathsf{y}$, $\boldsymbol{\pi}$, $\boldsymbol{\pi}_{\mathsf{y}}$ as in the proof of the theorem. We may assume that $W^{\mathrm{e}} < W$. Note that $W^{\mathrm{e}}\boldsymbol{\pi} = K$ because if $w \in W$ and $s$ is a transposition of $\mathsf{X} \times \{\mathsf{y}\}$, then either $w$ or $ws$ lies in $W^{\mathrm{e}}$ and $w\boldsymbol{\pi} = (ws)\boldsymbol{\pi}$. Similarly, one shows that $W_{\mathsf{y}}^{\mathrm{e}}\boldsymbol{\pi}_{\mathsf{y}} = H$ (because $(W^{\mathrm{e}})_{\mathsf{y}} = (W_{\mathsf{y}})^{\mathrm{e}}$, we only write $W_{\mathsf{y}}^{\mathrm{e}}$). This proves that $W^{\mathrm{e}}$ is transitive on $\mathsf{X} \times \mathsf{Y}$ and in particular that $W_{\mathsf{y}}^{\mathrm{e}} < W_{\mathsf{y}}$ and $W_{(\mathsf{x},\mathsf{y})}^{\mathrm{e}} < W_{(\mathsf{x},\mathsf{y})}$. Furthermore,

$$\big[W_{\mathsf{y}}^{\mathrm{e}} \div W^{\mathrm{e}}\big] = [K_{\mathsf{y}} \div K] \qquad \text{and} \qquad \Big[W_{(\mathsf{x},\mathsf{y})}^{\mathrm{e}} \div W_{\mathsf{y}}^{\mathrm{e}}\Big] = [H_{\mathsf{x}} \div H].$$

Suppose now

$$W_{(\mathsf{x},\mathsf{y})}^{\mathrm{e}} < T \leqslant W^{\mathrm{e}} \qquad \text{and} \qquad T \nleqslant W_{\mathsf{y}}.$$

Then $T$ contains $R = H^{\mathrm{e}} \times \{1\}$ and $W^{\mathrm{e}}_{(\mathsf{x},\mathsf{y})}$. The order of $RW^{\mathrm{e}}_{(\mathsf{x},\mathsf{y})}$ is

$$\frac{|H^{\mathrm{e}}||W_{(\mathsf{x},\mathsf{y})}|}{2|H^{\mathrm{e}}_{\mathsf{x}}|}.$$

If $H^{\mathrm{e}}$ is transitive, then $|H^{\mathrm{e}} : H^{\mathrm{e}}_{\mathsf{x}}| = |H : H_{\mathsf{x}}|$ and so $T \geqslant W^{\mathrm{e}}_{\mathsf{y}} = RW^{\mathrm{e}}_{(\mathsf{x},\mathsf{y})}$.

Otherwise, if $|\mathsf{Y}| > 2$, we use a different $R$ but the same argument to find that again $T$ contains $W^{\mathrm{e}}_{\mathsf{y}}$. The definition of $R$ is as follows: say $\mathsf{y}_2$, $\mathsf{y}_3$ two distinct elements of $\mathsf{Y} - \{\mathsf{y}\}$ and $t \in T$ such that $\mathsf{y}_2 t = \mathsf{y}$; the largest even subgroup of the base subgroup of $H \operatorname{wr} K_{\mathsf{y}}$, fixing all the points outside $\mathsf{X} \times \{\mathsf{y}_2, \mathsf{y}_3\}$ is $P \cong (H \times H)^{\mathrm{e}}$; put $R = t^{-1}Pt$. □

2.4.5. Remark. *$H$, $\mathsf{X}$ as in theorem above. If $H^e$ is not transitive, then*

$$\mathcal{C}\big((H \operatorname{wr} \mathrm{S}_2)^e\big) \supsetneq \mathcal{C}(H)\bar{\oplus}\mathcal{C}(\mathrm{S}_2).$$

*In fact, if $T \in [(H_{\mathsf{x}} \times H)^e \div (H \operatorname{wr} \mathrm{S}_2)^e]$, then either*

$$T \in [(H_{\mathsf{x}} \times H)^e \div (H \times H)^e] = \mathcal{C}(H),$$

*or $T \geqslant H^e \times H^e$. Note that $[(H^e \times H^e) \div (H \operatorname{wr} \mathrm{S}_2)^e] = \mathcal{M}_3$.*

Proof. First of all, $H_{\mathsf{x}} \leqslant H^{\mathrm{e}}$, so $(H_{\mathsf{x}} \times H)^{\mathrm{e}} = H_{\mathsf{x}} \times H^{\mathrm{e}}$. Furthermore, $|H : H_{\mathsf{x}}|$ is even, so the top group $\mathrm{S}_2$ is even (see (5.5.G) for example). Suppose a transitive $T$ contains $H_{\mathsf{x}} \times H^{\mathrm{e}}$. Then it contains $\{1\} \times H^{\mathrm{e}}$ and so it contains $H^{\mathrm{e}} \times H^{\mathrm{e}}$. The factor group

$$\frac{(H \operatorname{wr} \mathrm{S}_2)^{\mathrm{e}}}{H^{\mathrm{e}} \times H^{\mathrm{e}}}$$

has order 4 and contains at least 2 different subgroups of order 2, namely the image of $(H \times H)^{\mathrm{e}}$ and the image of the top group. Thus it is the elementary abelian group of order 4. □

2.4.6. Corollary. *Let $\mathsf{Z}$ be a non trivial equipartition of $\Omega$ and let $\boldsymbol{A}$ be the alternating group on $\Omega$. If $\Omega$ has order at least 6, and if $\mathsf{Y}$ is an imprimitivity system of $\boldsymbol{A}_{\mathsf{Z}}$, then $\mathsf{Y} = \mathsf{Z}$.*

2.4.7. Corollary. *If $\Omega$ has order 4 and $\boldsymbol{A}$, $\mathsf{Z}$ are as above, then $\boldsymbol{A}_{\mathsf{Z}}$ has three distinct imprimitivity systems.*

2.4.8. Corollary. *Let $\boldsymbol{S}$ be a finite symmetric group and let $G$ be a transitive subgroup of $\boldsymbol{S}$. There is a bijection between the set of the imprimitivity systems of $G$ and the set of the maximal imprimitive subgroups of $\boldsymbol{S}$ containing $G$.*

## 2.5. *Toccata e fuga* on primitive groups

***Toccata.*** A *primitive* group is a transitive group without imprimitivity systems. The following is an immediate consequence of 2.4.2.

2.5.1. Proposition. *A transitive group is primitive if and only if its point stabilizers are maximal subgroups.*

However, that is not the only way to prove that a group is primitive! For example, one easily proves that the symmetric and alternating groups in their natural action are primitive[1] by showing that any block is a trivial block. The other primitive groups, that is, the primitive groups not containing the symmetric or alternating group of their degree, are referred as *proper* primitive groups.

2.5.2. Proposition. *The even part of proper primitive groups is transitive.*

Proof. $G^{\mathrm{e}}$ is a normal subgroup of $G$ and its orbits are blocks for $G$. Since $G$ is primitive, the blocks are trivial and either $G^{\mathrm{e}}$ is transitive or it is the trivial group $\{1\}$. But $G^{\mathrm{e}} = \{1\}$ forces $G = \mathrm{S}_2$. □

Unfortunately, there are transitive groups whose even part is not transitive. Easy examples are the transitive cyclic subgroups of even degree. Also, there are primitive groups whose even part is not primitive, the smallest degree being 21 ([**LPS87**, §1 Rem.3], see also 3.4.1).

***...e fuga.*** Each element of a proper primitive group of large degree moves a lot of points. Results in this direction were published as early as 1871 by Jordan. We just recall two striking consequences here.

2.5.3. Proposition. *No proper primitive group contains transpositions or cycles of length* 3 ([**DM96**, § 7.4])*. A proper primitive group with a double transposition has degree less than* 9 ([**Wie64**, p. 43])*.*

Since maximal imprimitive subgroups of the symmetric or alternating groups contain transpositions or double transpositions, proposition above shows why the maximal imprimitive subgroups of degree at least 9 are indeed maximal subgroups.

## 2.6. Exercise: put all together

Second maximal subgroups, even or odd, stabilizers of points, intransitive groups, imprimitive groups, primitive groups. What if we try to put all together?

---

[1]with the exception of $\mathrm{A}_2 = 1$ which is not even transitive!

2.6.1. Lemma. *Use $\Omega$, $\boldsymbol{S}$, $\boldsymbol{A}$, $\boldsymbol{H}$ and $\boldsymbol{U}$ as in* 2.1. *Suppose that* $\{\mathsf{X},\mathsf{Y}\}$ *is a partition of* $\Omega$ *and that* $T$ *is a maximal imprimitive subgroup of* $\operatorname{Sym}\mathsf{Y}$. *If* $\boldsymbol{H} \leqslant (\operatorname{Sym}\mathsf{X}) \times T$ *and if there is a proper primitive group* $G$ *such that* $\boldsymbol{H} < G$, *then* $\boldsymbol{H}$ *is even,* $\boldsymbol{A} = \mathrm{A}_7$, $\left|\mathsf{X}\right| = 1$, $G \cong \mathrm{PSL}(3,2)$ *and* $\boldsymbol{H}$ *is the stabilizer of* $\mathsf{X}$ *in* $G$.

Proof. Note that $(\operatorname{Sym}\mathsf{X}) \times T < (\operatorname{Sym}\mathsf{X}) \times (\operatorname{Sym}\mathsf{Y}) < \boldsymbol{S}$. If $\boldsymbol{H}$ is a second maximal subgroup of $\boldsymbol{S}$ then $\boldsymbol{H} = (\operatorname{Sym}\mathsf{X}) \times T$. But $T$, being maximal imprimitive, is isomorphic to the wreath product of 2 symmetric groups in imprimitive action and contains transpositions. Then $G$, containing $H$ which contains $T$, contains transpositions. By the first part of 2.5.3, this is not possible.

It follows that $\boldsymbol{H}$ may only be a second maximal subgroup of $\boldsymbol{A}$. In particular $\boldsymbol{H}$ is even. Then

$$\boldsymbol{H} \leqslant \big((\operatorname{Sym}\mathsf{X}) \times T\big)^{\mathrm{e}} < \big((\operatorname{Sym}\mathsf{X}) \times (\operatorname{Sym}\mathsf{Y})\big)^{\mathrm{e}} < \boldsymbol{A}.$$

This forces $\boldsymbol{H} = \big((\operatorname{Sym}\mathsf{X}) \times T\big)^{\mathrm{e}}$. In particular $\boldsymbol{H}$ contains $\operatorname{Alt}\mathsf{X}$ and by the first part of 2.5.3 again, $\left|\mathsf{X}\right| < 3$. Also, $\boldsymbol{H}$ contains $T^{\mathrm{e}}$ which contains double transpositions. By the second part of 2.5.3, the order of $\Omega$ is less than 9 and hence the order of $\mathsf{Y}$ is either 4 or 6. Then $\boldsymbol{H}$ may only be one of the 4 following groups:

$$\big(\mathrm{S}_2 \times (\mathrm{S}_2 \operatorname{wr} \mathrm{S}_2)\big)^{\mathrm{e}} < \mathrm{A}_6 = \boldsymbol{A};$$
$$(\mathrm{S}_2 \operatorname{wr} \mathrm{S}_2)^{\mathrm{e}} < \mathrm{A}_4 < \mathrm{A}_5 = \boldsymbol{A};$$
$$\big(\mathrm{S}_2 \times (\mathrm{S}_2 \operatorname{wr} \mathrm{S}_3)\big)^{\mathrm{e}} < \mathrm{A}_8 = \boldsymbol{A};$$
$$(\mathrm{S}_2 \operatorname{wr} \mathrm{S}_3)^{\mathrm{e}} < \mathrm{A}_6 < \mathrm{A}_7 = \boldsymbol{A}.$$

It is easy to see that in the first two cases there are no primitive groups between $\boldsymbol{H}$ and $\boldsymbol{A}$. In the third case $\boldsymbol{H}$ is low because it is contained in $(\mathrm{S}_2 \operatorname{wr} \mathrm{S}_4)^{\mathrm{e}}$ which is contained in a copy of $\mathrm{AGL}(3,2)$. However, in the last case $\boldsymbol{H}$ is contained in two copies of $\mathrm{PSL}(3,2)$, conjugate in $\boldsymbol{S}$ but not in $\boldsymbol{A}$. $\boldsymbol{H}$ is the stabilizer of a point in each of them. In fact, that point may be regarded as a point of the Fano plane where a $\mathrm{PSL}(3,2)$ acts naturally. Then the blocks stabilized by $\boldsymbol{H}$ are exactly the remaining points of the lines passing through the stabilized point.

There are no other primitive groups containing such $\boldsymbol{H}$, so $[\boldsymbol{H} \div \boldsymbol{A}]$ is an $\mathcal{M}_3$:

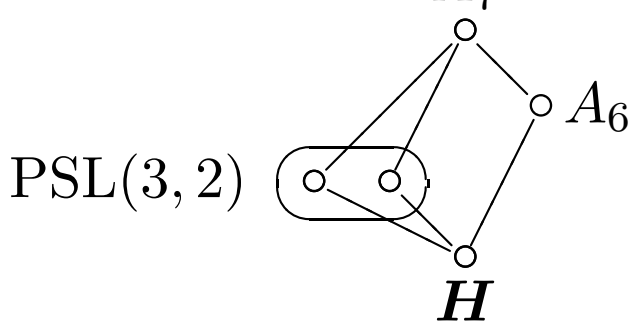

This completes the proof. □

## 2.7. Isomorphic maximal subgroups

Isomorphic maximal subgroups of the symmetric or alternating groups are usually but not always conjugate in the symmetric group. There are well known exceptions in degree 6 where both $\mathrm{PSL}_2(5)$ and $\mathrm{PGL}_2(5)$ acting on the projective line of 6 points are maximal primitive and isomorphic to the stabilizers of a point in $\mathrm{A}_6$, $\mathrm{S}_6$ respectively. Since conjugation in $\mathrm{S}_6$ preserves transitivity they are not conjugate in $\mathrm{S}_6$. Also, maximal imprimitive (but transitive) subgroups of $\mathrm{S}_6$, $\mathrm{A}_6$ with blocks of order 2, are isomorphic to maximal intransitive subgroups with an orbit of order 2 in view of

$$\mathrm{S}_2 \operatorname{wr} \mathrm{S}_3 = (\mathrm{S}_2)^3 \rtimes \mathrm{S}_3 = \Delta \times (M \rtimes \mathrm{S}_3) \cong \mathrm{S}_2 \times \mathrm{S}_4,$$

where $\Delta$ is the diagonal of $(\mathrm{S}_2)^3$ and $M$ is its *deleted permutation module* (as described in [**KL90**, p. 185]).

A similar phenomenon which we shall find useful later is described below. It is proved with numerical arguments in [**DM96**, §5.2] that with only a few exceptions, any subgroup of $\mathrm{A}_n$, $\mathrm{S}_n$ with "small" index must be intransitive. Then with the support of the ATLAS to rule out the exceptions, one checks the following

2.7.1. Proposition. *Suppose $n \neq 6$. If a maximal subgroup of $\mathrm{S}_n$ or $\mathrm{A}_n$ has the same order as a maximal intransitive subgroup, then these two subgroups are conjugate by one element of $\mathrm{A}_n$.*

## 2.8. The homogeneous mark of two subgroups

If $K$ and $L$ are subgroups of a finite group $G$, their *homogenous mark* in $G$, which we write as $\mathrm{hm}_G(K, L)$, is the number of $G$-conjugates of $K$ either containing $L$ or contained in $L$. That is, the order of

$$\{\ M \ \mid \ M = K^g \text{ for some } g \in G \text{ and } (M \leqslant L \text{ or } M \geqslant L)\ \}.$$

We tend to omit the subscript $G$ when no confusion arises. Clearly $\mathrm{hm}_G(G, K) = 1$ and $\mathrm{hm}_G(K, G) = |G : \mathrm{N}_G\, K|$. In fact,

$$\mathrm{hm}_G(K, L) \leqslant |G : \mathrm{N}_G\, K| \text{ for any } L \leqslant G.$$

In order to count the subgroups containing a given subgroup $K$ we will often use the following lemma which we name after Pálfy as we saw it first in [**Pál88**].

2.8.1. Lemma (Pálfy). *Let $K$, $L$ be subgroups of $G$, then*

$$|G : \mathrm{N}(L)|\, \mathrm{hm}(K, L) = |G : \mathrm{N}(K)|\, \mathrm{hm}(L, K).$$

*As a consequence, if* $K \leqslant L$, *then*

$$\mathrm{hm}(L, K) = \frac{\left|\mathrm{N}(K) : K\right|}{\left|\mathrm{N}(L) : L\right|} \frac{\mathrm{hm}(K, L)}{\left|L : K\right|}.$$

For the sake of completeness, we remind the reader that if $K$ is a subgroup of $L$, then the number

$$\left|\mathrm{N}_G\, L : L\right| \mathrm{hm}_G(L, K)$$

is called the *mark* of $K$ in $G$ on $L$ and corresponds to the number of fixed points of $K$ acting on the coset space of the right cosets of $L$ in $G$. Given a finite group $G$, one may choose representative subgroups for the conjugacy classes of its subgroups and compute the marks of any given couple of such representatives. These marks may then be conveniently stored in a lower triangular matrix called the table of marks. Quoting Burnside [**Bur11**, §185],

> *The number of distinct ways in which a given group can be represented as a permutation group of given degree is determined at once by its table of marks .*

Nowadays, there are tables of marks available for many finite groups and from them one could easily compute the related homogeneous marks.

## 2.9. Maximal subgroups of direct products

Given a group-isomorphism $\boldsymbol{h} : \frac{A}{M} \longrightarrow \frac{B}{N}$, the *Goursat subgroup* of $A \times B$ corresponding to $\boldsymbol{h}$ is defined as follows:

$$\mathbf{G}_{\boldsymbol{h}} \equiv \{(a, b) \mid (Ma)\boldsymbol{h} = Nb\}.$$

It is proved that all subgroups of a direct product $L \times R$ are Goursat subgroups:

2.9.1. Theorem. *Let* $H$ *be a subgroup of* $L \times R$, *and say* $\boldsymbol{\lambda}$, $\boldsymbol{\rho}$ *the left and right canonical projections respectively. Put*

$$M = \big(H \cap (L \times \{1\})\big)\boldsymbol{\lambda} \qquad \textit{and} \qquad N = \big(H \cap (\{1\} \times R)\big)\boldsymbol{\rho},$$

*then*

1. $M \trianglelefteq H\boldsymbol{\lambda}$ *and* $N \trianglelefteq H\boldsymbol{\rho}$,
2. *there is an isomorphism* $\boldsymbol{h} : \frac{H\boldsymbol{\lambda}}{M} \longrightarrow \frac{H\boldsymbol{\rho}}{N}$ *such that* $H = \mathbf{G}_{\boldsymbol{h}}$.

For example, it is straightforward to see that $(S_m \times S_n)^{\mathrm{e}} = \mathbf{G}_{\boldsymbol{a}}$ where $\boldsymbol{a}$ is the unique isomorphism between $\frac{S_m}{A_m}$ and $\frac{S_n}{A_n}$.

We now state a few results about maximal subgroups of direct products. The complete underlying theory and proofs may be found in the appendix.

2.9.2. PROPOSITION. *If $G \lessdot L \times R$ then one and only one of the following three cases holds:*

(1) $G = T \times R$ *for some* $T \lessdot L$,
(2) $G = L \times T$ *for some* $T \lessdot R$,
(3) $G = \mathbf{G}_{\boldsymbol{a}}$ *for some isomorphism* $\boldsymbol{a} : \frac{L}{M} \longrightarrow \frac{R}{N}$ *with* $M$ *maximal normal subgroup of* $L$.

2.9.3. PROPOSITION. *If* $G \lessdot \mathbf{G}_{\boldsymbol{a}}$ *where* $\boldsymbol{a} : \frac{L}{M} \longrightarrow \frac{R}{N}$ *is an isomorphism with* $\frac{L}{M}$ *simple and different from* $\{1\}$, *then one and only one of the following three cases holds:*

(1) $G = \mathbf{G}_{\boldsymbol{b}}$ *where* $\boldsymbol{b} = \boldsymbol{ca}$ *and* $\boldsymbol{c} : \frac{T}{T \cap M} \longrightarrow \frac{L}{M}$ *is the canonic isomorphism for some* $T \lessdot L$ *with* $M$ *not contained in* $T$,
(2) $G = \mathbf{G}_{\boldsymbol{b}}$ *where* $\boldsymbol{b} = \boldsymbol{ac}$ *and* $\boldsymbol{c} : \frac{R}{N} \longrightarrow \frac{T}{T \cap N}$ *is the canonic isomorphism for some* $T \lessdot R$ *with* $N$ *not contained in* $T$,
(3) $G = \mathbf{G}_{\boldsymbol{b}}$ *where* $\boldsymbol{b}$ *is the restriction of* $\boldsymbol{a}$ *to* $\frac{T}{M}$ *for some maximal subgroup* $T$ *of* $L$ *containing* $M$.

In particular, it is now very easy to check that

2.9.4. COROLLARY. *The maximal subgroups of* $S_a \times S_b$ *are*

- *the* $T \times S_b$ *with* $T \lessdot S_a$,
- *the* $S_a \times T$ *with* $T \lessdot S_b$, *and*
- $(S_a \times S_b)^e$.

*The maximal subgroups of* $(S_a \times S_b)^e$ *are*

- *the* $(T \times S_b)^e$ *with* $T \lessdot S_a$ *and* $T \neq A_a$,
- *the* $(S_a \times T)^e$ *with* $T \lessdot S_b$ *and* $T \neq A_b$, *and*
- $A_a \times A_b$.

CHAPTER 3

# Permutation representations

## 3.1. Permutation representations, actions and permutation groups

Recall that a group $G$ is said to *act* (on the right) of a set $\mathsf{X}$ (or that there is a *right action* of $G$ on $\mathsf{X}$) when there exists a function

$$\cdot : \mathsf{X} \times G \longrightarrow \mathsf{X}, \qquad (\mathsf{x}, g) \mapsto \mathsf{x} \cdot g$$

such that $(\mathsf{x} \cdot h) \cdot g = \mathsf{x} \cdot (hg)$ and $\mathsf{x} \cdot 1 = \mathsf{x}$ whenever $\mathsf{x} \in \mathsf{X}$, $g, h \in G$.

Following [**Rob93**, §1.6] we use the languages of actions and permutation representations interchangeably. Furthermore, we tend to identify permutation representations with their images so that the language used for permutation groups is sometime carried on actions and permutation representations as well. For example, we may say that two actions commute when the images of the related permutation representations are commutating subgroups of a same permutation group.

Two permutation representations

$$\boldsymbol{\varphi} : G \longrightarrow \operatorname{Sym} \mathsf{X} \qquad \text{and} \qquad \boldsymbol{\psi} : G \longrightarrow \operatorname{Sym} \mathsf{Y}$$

are said to be *equivalent* if there exists a bijection $\boldsymbol{b} : \mathsf{X} \longrightarrow \mathsf{Y}$ such that for each $g \in G$, $g\boldsymbol{\psi} = \boldsymbol{b}^{-1}(g\boldsymbol{\varphi})\boldsymbol{b}$, that is,

$$\big(\mathsf{x} \cdot (g\boldsymbol{\varphi})\big)\boldsymbol{b} = (\mathsf{x}\boldsymbol{b}) \cdot (g\boldsymbol{\psi}) \qquad \forall \mathsf{x} \in \mathsf{X},\ g \in G.$$

The map $\boldsymbol{b}$ is then said *equivalence* of the permutation representations.

In [**Asc93**, §4] two permutation representations

$$\boldsymbol{\varphi} : G \longrightarrow \operatorname{Sym} \mathsf{X} \qquad \text{and} \qquad \boldsymbol{\psi} : L \longrightarrow \operatorname{Sym} \mathsf{Y}$$

are said to be *quasiequivalent* if there exists a group isomorphism $\boldsymbol{a} : G \longrightarrow L$ and a bijection $\boldsymbol{b} : \mathsf{X} \longrightarrow \mathsf{Y}$ such that for each $g \in G$, $g\boldsymbol{a}\boldsymbol{\psi} = \boldsymbol{b}^{-1}(g\boldsymbol{\varphi})\boldsymbol{b}$, that is,

$$\big(\mathsf{x} \cdot (g\boldsymbol{\varphi})\big)\boldsymbol{b} = (\mathsf{x}\boldsymbol{b}) \cdot (g\boldsymbol{a}\boldsymbol{\psi}) \qquad \forall \mathsf{x} \in \mathsf{X},\ g \in G.$$

We call the ordered pair $(\boldsymbol{a}, \boldsymbol{b})$ (or just the bijection $\boldsymbol{b}$ if $\boldsymbol{a}$ is the identity) *permutational isomorphism*. In fact, we say in this case that $\boldsymbol{\varphi}$, $\boldsymbol{\psi}$, or $G\boldsymbol{\varphi}$, $L\boldsymbol{\psi}$, or $G$ on $\mathsf{X}$, $L$ on $\mathsf{Y}$ are *permutationally isomorphic*, as opposed to $G$, $L$ being just abstractly isomorphic.

We usually distinguish actions of a same group up to permutational isomorphism or more exceptionally up to equivalence. The following lemma is extremely useful in discerning between different transitive permutation representations of a group.

3.1.1. LEMMA. *Let* $\boldsymbol{\varphi} : G \longrightarrow \operatorname{Sym} \mathsf{X}$ *and* $\boldsymbol{\psi} : G \longrightarrow \operatorname{Sym} \mathsf{Y}$ *be transitive permutation representations (not necessarily faithful),* $\mathsf{x} \in \mathsf{X}$*,* $\mathsf{y} \in \mathsf{Y}$*.*

(1) $\boldsymbol{\varphi}$*,* $\boldsymbol{\psi}$ *are equivalent if and only if there is* $g \in G$ *such that* $g^{-1} G_{\mathsf{x}} g = G_{\mathsf{y}}$*.*
(2) $G$ *on* $\mathsf{X}$ *and* $G$ *on* $\mathsf{Y}$ *are permutationally isomorphic if and only if there is* $\boldsymbol{a} \in \operatorname{Aut} G$ *such that* $G_{\mathsf{x}} \boldsymbol{a} = G_{\mathsf{y}}$*.*

PROOF. This is generally known. The first part for example corresponds exactly to [**Asc93**, (5.9)(2)]. The second part is left to the reader □

In this chapter we discuss a few very well known actions, occasionally filling in a few details that have not yet been written or at least collected in a unique place. We are especially interested in the degree of transitivity, full normalizers and parity which we describe separately in the three sections below.

## 3.2. Degree of transitivity, primitivity

A transitive permutation group $G$ on a non empty set $\mathsf{X}$ is 2-transitive if for some $\mathsf{x} \in \mathsf{X}$, $G_{\mathsf{x}}$ is transitive on $\mathsf{X} - \{\mathsf{x}\}$. Then one readily checks that $G_{\mathsf{y}}$ is transitive on $\mathsf{X} - \{\mathsf{y}\}$ for each $\mathsf{y} \in \mathsf{X}$. Inductively, for a positive integer $m > 1$, $G$ is $(m+1)$-transitive if for some $\mathsf{x} \in \mathsf{X}$, $G_{\mathsf{x}}$ is $m$-transitive on $\mathsf{X} - \{\mathsf{x}\}$.

Similarly, a primitive permutation group $G$ on $\mathsf{X}$ is 2-primitive if for some $\mathsf{x} \in \mathsf{X}$, $G_{\mathsf{x}}$ is primitive on $\mathsf{X} - \{\mathsf{x}\}$. Also, for a positive integer $m > 1$, $G$ is $(m+1)$-primitive if for some $\mathsf{x} \in \mathsf{X}$, $G_{\mathsf{x}}$ is $m$-primitive on $\mathsf{X} - \{\mathsf{x}\}$.

The *degree of transitivity* (*primitivity*) of $G$ is then the largest integer $m$ such that $G$ is $m$-transitive ($m$-primitive).

For example, provided $m > 2$, the degrees of transitivity of $\mathrm{S}_m$ and $\mathrm{A}_m$ are $m$ and $m-2$ respectively. Also, it is easy to check that a 2-transitive group must be primitive; consequently, if a group is not $m$-primitive, then it may not be $(m+1)$-transitive either.

The research on the 2-transitive permutation groups culminated in their classification which in view of the well known theorem below is a by-product of the CFSG.

3.2.1. THEOREM ([**Bur11**, §154, Th. XIII]). *The socle of a finite* 2*-transitive permutation group is either a regular elementary abelian* $p$*-group, or a non regular non abelian simple group.*

The complete list of the socles of the 2-transitive permutation groups of the latter type is given in Cameron Table 3.A ([**Cam81**], see also [**DM96**] for example). Column "deg-t" shows the maximum degree of transitivity of a group with socle appearing in first column.

TABLE 3.A. Socles of the 2-transitive almost simple groups [**Cam81**]

| Soc | deg | deg-t | *notes* |
|---|---|---|---|
| $A_n$, $n \geqslant 5$ | $n$ | $n$ | 2 rep. if $n = 6$ |
| $PSL_d(q)$, $d \geqslant 2$ | $(q^d - 1)/(q - 1)$ | 3 if $d = 2$ | $(d, q) \neq (2, 2), (2, 3)$ |
| | | 2 if $d > 2$ | 2 rep. se $d > 2$ |
| $PSU_3(q)$ | $q^3 + 1$ | 2 | $q > 2$ |
| $Sz(q)$ | $q^2 + 1$ | 2 | $q = 2^{2a+1} > 2$ |
| ${}^2G_2(q)$ | $q^3 + 1$ | 2 | $q = 3^{2a+1} > 3$ |
| $PSp_{2d}(2)$ | $2^{d-1}(2^d \pm 1)$ | 2 | $d > 2$ |
| $PSL_2(11)$ | 11 | 2 | 2 rep. |
| $PSL_2(8)$ | 28 | 2 | |
| $A_7$ | 15 | 2 | 2 rep. |
| $M_{11}$ | 11 | 4 | |
| $M_{11}$ | 12 | 3 | |
| $M_{12}$ | 12 | 5 | 2 rep. |
| $M_{22}$ | 22 | 3 | |
| $M_{23}$ | 23 | 4 | |
| $M_{24}$ | 24 | 5 | |
| HS | 176 | 2 | 2 rep. |
| $Co_3$ | 276 | 2 | |

The socle is always deg-t transitive, except for

- $A_n$ of degree $n$ which is $(n-2)$-transitive,
- $PSL_2(8)$of degree 28 which is only transitive primitive, and
- $PSL_2(q)$, $q$ odd, of degree $1 + q$ which is 2-transitive.

There are no 6-transitive proper primitive groups. Also, checking that in Table 3.A the abstractly isomorphic socles having the same degree are indeed permutationally isomorphic, yields to the result below.

3.2.2. PROPOSITION. *If two maximal 2-transitive subgroups of $S_n$ are isomorphic, then they are conjugate in $S_n$.*

## 3.3. Full normalizer

Let $L$ be a transitive subgroup of $\boldsymbol{S} = \operatorname{Sym}\Omega$ and $T$ be the stabilizer in $L$ of a point $\omega$. By the *full normalizer* of $L$ we mean its normalizer in $\boldsymbol{S}$. It is shown

in [**DM96**, Ch.4] that

$$\mathrm{C}_{\boldsymbol{S}}\, L \cong \frac{\mathrm{N}_L\, T}{T}$$

and that the conjugation in $\mathrm{N}_{\boldsymbol{S}}\, L$ induces a group homomorphism

$$\Psi : \mathrm{N}_{\boldsymbol{S}}\, L \longrightarrow \operatorname{Aut} L$$

whose kernel is $\mathrm{C}_{\boldsymbol{S}}\, L$. Now, $\operatorname{Aut} L$ acts naturally on the subgroup lattice of $L$ and hence on its power-set. Well, the image of $\Psi$ is the group of the automorphisms of $L$ which preserve the conjugacy class of $T$ in $L$. Following [**KL90**, §3.2], we denote this class by

$$[T]_L = \{\, l^{-1} T\, l \;\mid\; l \in L \,\},$$

so $\operatorname{Im}\Psi = \mathrm{St}_{\operatorname{Aut} L}[T]_L$. By the Frattini Argument (see 6.1)

$$\operatorname{Im}\Psi = (\mathrm{St}_{\operatorname{Aut} L}\, T)(\operatorname{Inn} L) \tag{3.3.A}$$

and this determines completely $\mathrm{N}_{\boldsymbol{S}}\, L$. In fact, if $\boldsymbol{\beta}\Psi = \boldsymbol{a}$ with $T^{\boldsymbol{a}} = T$, then there is some $m \in \mathrm{N}_L\, T$ such that $\boldsymbol{\beta}$ is defined by

$$\boldsymbol{\beta} : (\omega \cdot l) \mapsto \omega \cdot (m l^{\boldsymbol{a}}), \qquad l \in L. \tag{3.3.B}$$

Suppose now that $L$ is the socle of a primitive group $G$ of type almost simple. We claim that in this case $\Psi$ is a permutational isomorphism between $\mathrm{N}_{\boldsymbol{S}}\, L$ and the group $\operatorname{Im}\Psi$ acting on $[T]_L$.

Note that by our assumptions $L$ is non abelian simple and $T \neq \{1\}$ (the last assertion is a consequence of the Schreier Conjecture). Also, $\operatorname{Inn} L$ is the unique minimal normal subgroup of $\operatorname{Im}\Psi$.

First of all, the action of $\operatorname{Im}\Psi$ on $[T]_L$ is faithful. If not, the kernel of the action would contain all the inner automorphisms of $L$ which is impossible because $\operatorname{Inn} L$ is transitive on $[T]_L$.

Then we prove that $T$ is self-normalizing in $L$ ($L$ might not be primitive). If not, $T = L \cap G_\omega$ is normal in $G_\omega$. Thus $T$ is normal in $\langle \mathrm{N}_L\, T, G_\omega \rangle$ which is equal to $G$ by primitivity. Then $T$ would be also normal in $L$ but it is not.

As a consequence, the equation $\mathrm{C}_{\boldsymbol{S}}\, L = \{1\}$ holds. Moreover,

$$\big|[T]_L\big| = \big|L : T\big|.$$

In particular, $\operatorname{Im}\Psi$ and $\mathrm{N}_{\boldsymbol{S}}\, L$ are abstractly isomorphic via $\Psi$ and have at least the same degree as permutation groups.

To prove that indeed $\Psi$ is a permutational isomorphism, we just show that if $\boldsymbol{s} \in \mathrm{N}_{\boldsymbol{S}}\, L$ and $\omega\boldsymbol{s} = \omega$, then $T^{\boldsymbol{s}\Psi} = T$. In fact, stabilizers of points with respect to the two actions must have the same order and if the image of one is contained in the other, then they are one the image of the other, proving in this way that the isomorphism is a permutational isomorphism.

Provided $T = \boldsymbol{S}_\omega \cap L$,

$$T^{\boldsymbol{s}\Psi} = \boldsymbol{s}^{-1} T\, \boldsymbol{s} = \boldsymbol{s}^{-1}(\boldsymbol{S}_\omega \cap L)\, \boldsymbol{s} = \boldsymbol{S}_{\omega \boldsymbol{s}} \cap L = \boldsymbol{S}_\omega \cap L = T.$$

Clearly, the full normalizer of $G$ is contained in the full normalizer of $L$. Since the image of $G$ by $\Psi$ is denoted by $\mathrm{Aut}_G\, L$ (see [**AS85**]), we have:

3.3.1. PROPOSITION. *Let $\boldsymbol{S}$ be the parent symmetric group of a primitive group $G$ of type almost simple, $L$ the socle of $G$ and $T$ the stabilizer of one point in $L$. Then $\mathrm{N}_{\boldsymbol{S}}\, G$ and $(\mathrm{N}_{\mathrm{Aut}\, L}\, \mathrm{Aut}_G\, L) \cap \big((\mathrm{St}_{Aut L}\, T)(\mathrm{Inn}\, L)\big)$ on $[T]_L$ are permutationally isomorphic.*

In view of 3.3.1, the determination of the full normalizer of $L$ is straightforward when $L$ is one of the simple groups listed in ATLAS. As an example, consider $L$ as an abstract group and suppose that we want to determine all the transitive permutation representations of $L$ of degree $n$ such that the related full normalizers are primitive. Of course, we are only interested in a set of representatives for the conjugacy classes in $\mathrm{S}_n$ of these permutation representations. To save some writing, put $\boldsymbol{S} = \mathrm{S}_n$.

To find the representations such that $L$ is also primitive, consider the conjugacy classes (in $L$) of the maximal subgroups of $L$ of index $n$. There might be more than one, either fused in some automorphic extension (i.e. joined in the maximal subgroups table with a vertical line at the column corresponding to that automorphic extension) or not. The actions of $L$ on the sets of the right cosets of these maximal subgroups yield at least one representative of each conjugacy class of primitive subgroups of $\boldsymbol{S}$ abstractly isomorphic to $L$. Since classes fused by automorphic extensions yield conjugate subgroups of $\boldsymbol{S}$ (and hence conjugate normalizers), it is enough to determine the full normalizer related to only one representation for each set of fused classes.

The normalizer in $\boldsymbol{S}$ of $L$ acting on one of its maximal subgroups of index $n$, is the largest automorphic extension of $L$ which preserves the conjugation class of that maximal subgroup (that is, which has a semicolon in the entry of the table whose row is the given conjugation class and whose column is the chosen automorphic extension). This is enough to find all the maximal primitive almost simple subgroups of $\boldsymbol{S}$ whose socle is primitive and abstractly isomorphic to $L$.

To find the representations such that $L$ is not primitive, observe that the intersection of $L$ with a maximal subgroup of index $n$ of any automorphic extension $E$ must be properly contained in a maximal subgroup of $L$; let us denote by $\mathsf{C}_0$ its conjugacy class. Denote by $\mathsf{C}$ the orbit of $\mathsf{C}_0$ under the action by conjugation of $E$; the elements of $\mathsf{C}$ are the conjugacy classes of maximal subgroups of $L$ which are joined to $\mathsf{C}_0$ in the column correspondent to $E$ (there might be only one of such classes as in a maximal primitive $\mathrm{S}_5$ of $\mathrm{M}_{12}.2$). Any automorphic extension of $L$ containing $E$ which preserves $\mathsf{C}$ (i.e. which behaves as $E$ does on the elements of $\mathsf{C}$) is isomorphic to a subgroup of $\boldsymbol{S}$. The largest of such extensions is then isomorphic to the full normalizer of $L$. In this way one determines up to permutational isomorphisms all maximal primitive almost simple subgroups of $\boldsymbol{S}$ which have an imprimitive socle abstractly isomorphic to $L$.

We close this section with an immediate corollary of [**Kle87**, 1.3.2] which shows that the maximal subgroups of a simple group $L$ which *lift* to maximal subgroups of an automorphic extension $G$ are exactly the ones for which $G$ is contained in the full normalizers of the related permutation representations of $L$.

3.3.2. Lemma. *Let $L$ be a simple normal subgroup of $G$ and let $T$ be a non trivial proper subgroup of $L$. Then $G = L\,\mathrm{N}_G T$ if and only if $[T]_G = [T]_L$. Furthermore, if $T \lessdot L$, then $[T]_G = [T]_L$ if and only if $\mathrm{N}_G T \lessdot G$.*

## 3.4. Parity

As far as my little knowledge is concerned, discussing the *parity* of a permutation group, that is, whether it is a subgroup of the alternating group, is a quite intriguing and at times rather challenging task. Some help comes from the fact that the even subgroup of an odd permutation group is normal, indeed of index 2. If a group has no such subgroups, like all groups of odd order or the non abelian simple groups for example, then all its permutation representations must be even. Yet, this is not enough to unravel most cases where I know no other technique than brute force:

> given a set of convenient generators, count modulo 2 the transpositions that they may be composed of in your permutation representation.

It is well known that if $\boldsymbol{s}$ of $\mathrm{S}_n$ has an orbit of length $d$, then $\boldsymbol{s}$ induces on that orbit a permutation which is composition of $d-1$ transpositions. Of course, if $\boldsymbol{s}$ has more than one orbit, then $\boldsymbol{s}$ is composition of the transpositions induced on each orbit. In particular, only the orbits of even length contribute to the count modulo 2. Summarizing:

$$\mathrm{Par}\,\boldsymbol{s} \equiv \sum_{2|e|f} \frac{1}{e} \odot_{\boldsymbol{s}}^{e} \mod 2, \tag{3.4.A}$$

where

$$\operatorname{Par} \boldsymbol{s} := \begin{cases} 0 & \text{if } \boldsymbol{s} \in \mathrm{A}_n, \\ 1 & \text{otherwise} \end{cases}$$

and $\odot_{\boldsymbol{s}}^d$ is the number of points belonging to $\boldsymbol{s}$-orbits of length $d$. Observe now that if we denote by $\operatorname{Fix} \boldsymbol{s}^d$ the number of fixed points of $\boldsymbol{s}^d$, then

$$\operatorname{Fix}(\boldsymbol{s}^e) = \sum_{d|e} \odot_{\boldsymbol{s}}^d.$$

By Möbius inversion [**Cam94**, §12.7,Rem.3] one has

$$\odot_{\boldsymbol{s}}^e = \sum_{d|e} \mu(d) \operatorname{Fix}(\boldsymbol{s}^{e/d})$$

where the Möbius function $\mu$ is defined by

$$\mu(d) := \begin{cases} (-1)^l & \text{if } d \text{ is the product of } l \text{ distinct primes}, \\ 0 & \text{otherwise.} \end{cases}$$

Therefore (3.4.A) becomes

$$\operatorname{Par} \boldsymbol{s} \equiv \sum_{2|e|f} \frac{1}{e} \sum_{2\nmid d|e} \mu(d) \big(\operatorname{Fix}(\boldsymbol{s}^{e/d}) - \operatorname{Fix}(\boldsymbol{s}^{e/2d})\big) \mod 2. \tag{3.4.B}$$

Of course, if $\boldsymbol{s}$ has order 2, then as expected

$$\operatorname{Par} \boldsymbol{s} \equiv \frac{1}{2}\big(n - \operatorname{Fix}(\boldsymbol{s})\big) \mod 2.$$

Character Theory can be of great help because if $\chi$ is the permutation character of a given permutation representation, then $\operatorname{Fix}(\boldsymbol{s}) = \chi(\boldsymbol{s})$. Suppose for example that $H$ is the stabilizer of a point in a transitive permutation group $G$. In [**Isa76**, 5.14] it is shown that the permutation character of the action is $\chi = (1_H)^G$ and from equation [**Isa76**, 5.1] we get

$$\operatorname{Fix}(g) = \chi(g) = \big|\mathrm{C}_G\, g\big| \frac{\big|H \cap \mathrm{Cl}(g)\big|}{\big|H\big|}\ .$$

We use the formula above to show that the primitive permutation representation of degree 21 of $\mathrm{PGL}_2(7)$ is odd.

3.4.1. Proposition. *There is a primitive subgroup $G$ of $\mathrm{S}_{21}$ whose even part is not primitive.*

Proof. Consider page 3 of the ATLAS. Call $G = \mathrm{PGL}_2(7)$ an automorphic extension of $S = \mathrm{PSL}_2(7)$. $G$ has a maximal subgroup $H$ of index 21 such that $H \cong \mathrm{D}_{16}$ and $H \cap S$ is contained in a maximal subgroup $M \cong \mathrm{S}_4$ of $S$. Since $S$ is not contained in $H$, $\operatorname{Core}_G H = \{1\}$ and there is a faithful action of $G$ of degree 21. Note that $G$ is primitive because $H \lessdot G$. To prove that $G$ is an odd subgroup, it is enough to show that if $g \in \mathtt{2B}$ (the conjugation class of the ATLAS character table), then $g$ is odd. Now,

$$\operatorname{Fix}(g) = \left|\mathrm{C}_G\, g\right| \frac{\left|H \cap \mathtt{2B}\right|}{\left|H\right|} = \frac{12}{16}\left|H \cap \mathtt{2B}\right|.$$

The conjugacy classes of $H = \langle\, \rho, \tau \;\mid\; \rho^2 = \tau^2 = \tau\rho\tau\rho = 1 \,\rangle$ are the $\mathrm{Cl}(h)$ as below:

| $\mathrm{Cl}(h)$ | 2a | 2b | 4a | 8a | 8b |
|---|---|---|---|---|---|
| $\left\|\mathrm{C}_H\, h\right\|$ | 2 | 16 | 8 | 8 | 8 |
| $h$ | $\tau$ | $\rho^4$ | $\rho^2$ | $\rho$ | $\rho^3$ |

Of course, $H \cap S$ has index 2 in $H$. Thus, it is a subgroup of order 8 of $M$. It follows that $H \cap S \cong \mathrm{D}_8$ and so $H \cap S = \langle \rho^2, \tau \rangle$. The elements $\rho^4$, $\tau$, $\tau\rho^2$, $\tau\rho^4$, $\tau\rho^6$ must lie in $S$ and hence in the class 2A of $G$. However, the other 4 elements of 2a, namely $\tau\rho$, $\tau\rho^3$, $\tau\rho^5$, $\tau\rho^7$ cannot lie in 2A because each of them together with 2A generates $H$. Therefore $\left|H \cap \mathtt{2B}\right| = 4$ and $\operatorname{Fix}(g) = 3$. This implies that

$$\operatorname{Par} g \underset{(2)}{\equiv} \frac{1}{2}(21 - 3) = 9$$

is odd and that $G^{\mathrm{e}} = S$. From $H \cap S < M < S$ it follows that $G^{\mathrm{e}}$ is not primitive. □

### 3.5. Diagonal action on a cartesian power

For each $l > 1$ the *diagonal action* of $\operatorname{Sym} \mathsf{X}$ on $\mathsf{X}^l$ is

$$\boldsymbol{\vartheta}_l : \operatorname{Sym} \mathsf{X} \longrightarrow \operatorname{Sym}(\mathsf{X}^l)$$

where for all $\mathsf{x}_1, \ldots, \mathsf{x}_l \in \mathsf{X}$, $\boldsymbol{s} \in \operatorname{Sym} \mathsf{X}$

$$(\mathsf{x}_1, \ldots, \mathsf{x}_l)(\boldsymbol{s}\boldsymbol{\vartheta}_l) := (\mathsf{x}_1\boldsymbol{s}, \ldots, \mathsf{x}_l\boldsymbol{s}).$$

3.5.1. Proposition. *For all $\boldsymbol{s} \in \operatorname{Sym} \mathsf{X}$*

$$\operatorname{Par}(\boldsymbol{s}\boldsymbol{\vartheta}_l) \equiv l\left|\mathsf{X}\right| \operatorname{Par} \boldsymbol{s} \mod 2.$$

Proof. The composition of $\boldsymbol{\vartheta}_l$ with Par is a group homomorphism from $\operatorname{Sym} \mathsf{X}$ to the group of order 2. Therefore it is enough to show that the proposition holds when $\boldsymbol{s}$ is a transposition.

Assume that $\boldsymbol{s}$ is a transposition and put $n = \left|\mathsf{X}\right|$. Then the number of fixed points of $\boldsymbol{s}\boldsymbol{\vartheta}_l$ is $(n-2)^l$ and so $\boldsymbol{s}\boldsymbol{\vartheta}_l$ moves $n^l - (n-2)^l = 2\sum_{i=0}^{l-1} n^i (n-2)^{l-i-1}$ in orbits of length 2. Therefore $\boldsymbol{s}\boldsymbol{\vartheta}_l$ is product of transpositions whose number is

$$\sum_{i=0}^{l-1} n^i (n-2)^{l-i-1} \equiv l n^{l-1} \mod 2.$$

Since $l > 1$, $n^{l-1} \equiv n$ thus $\operatorname{Par} \boldsymbol{s}\boldsymbol{\vartheta}_l \equiv ln \mod 2$. □

## 3.6. Power set action

The *power set* of a set $\mathsf{X}$ is

$$\mathcal{P}(\mathsf{X}) := \{ \mathsf{A} \mid \mathsf{A} \subseteq \mathsf{X} \}.$$

For each positive integer $l < |\mathsf{X}|$, the $l$-homogeneous component of $\mathcal{P}(\mathsf{X})$ is

$$\binom{\mathsf{X}}{l} := \{ \mathsf{A} \subseteq \mathsf{X} \mid \mathsf{A} \text{ has order } l \}.$$

The *power set action* of $\operatorname{Sym}\mathsf{X}$ on $\mathcal{P}(\mathsf{X})$ is given by

$$\mathsf{A} \cdot \boldsymbol{s} := \{ \mathsf{a}\boldsymbol{s} \mid \mathsf{a} \in \mathsf{A} \} \qquad \text{for all } \mathsf{A} \subseteq \mathsf{X},\ \boldsymbol{s} \in \operatorname{Sym}\mathsf{X}.$$

This clearly induces an action on each homogeneous component of $\mathcal{P}(\mathsf{X})$. If $\boldsymbol{s} \in \operatorname{Sym}\mathsf{X}$, the permutation induced on $\mathcal{P}(\mathsf{X})$ is denoted by $\mathcal{P}(\boldsymbol{s})$, while the permutation induced on the $l$-homogeneous component $\binom{\mathsf{X}}{l}$ is denoted by $\binom{\boldsymbol{s}}{l}$.

We observe that $\operatorname{Sym}\mathsf{X}$ is primitive on the $l$-homogeneous component of $\mathcal{P}(\mathsf{X})$ provided that $|\mathsf{X}| \neq 2l$ because the stabilizer of a subset $\mathsf{A}$ of order $l$ coincides with $(\operatorname{Sym}\mathsf{A}) \times \operatorname{Sym}(\mathsf{X} - \mathsf{A})$ which is maximal in $\operatorname{Sym}\mathsf{X}$ if and only if $|\mathsf{A}| \neq |\mathsf{X} - \mathsf{A}|$ (see 2.3.1). However, if $1 < l < |\mathsf{X}| - 1$ then $\operatorname{Sym}\mathsf{X}$ is not 2-transitive on the $l$-homogeneous component (of course it is $|\mathsf{X}|$-transitive on the 1-homogeneous component, and also on the $(|\mathsf{X}| - 1)$-homogeneous component).

3.6.1. Proposition. *Let $\mathsf{X}$ be a set of order $n$, $0 < l < n$. For all $\boldsymbol{s} \in \operatorname{Sym}\mathsf{X}$*

$$\operatorname{Par}\binom{\boldsymbol{s}}{l} \equiv \binom{n-2}{l-1}\operatorname{Par}\boldsymbol{s} \mod 2.$$

Proof. It is enough to prove the result for a transposition $\boldsymbol{s}$, say $\boldsymbol{s} = (\mathsf{a}\ \mathsf{b})$. Then $\binom{\boldsymbol{s}}{l}$ has all non trivial orbits of length 2 and moves exactly the sets which intersect $\{\mathsf{a}, \mathsf{b}\}$ exactly in one point. Their number is twice the number of subsets of order $l - 1$ of $\mathsf{X} - \{\mathsf{a}, \mathsf{b}\}$. Hence $\binom{\boldsymbol{s}}{l}$ is made of $\binom{n-2}{l-1}$ transpositions. □

## 3.7. Field automorphisms as permutations

Let $\mathbb{F}$ be a field of order $q = p^f$ for some prime $p$. The group $G$ of the field automorphisms of $\mathbb{F}$ is a cyclic group of order $f$ generated for example by the Frobenius automorphism

$$\boldsymbol{\pi} : \mathbb{F} \longrightarrow \mathbb{F}, \qquad \mathfrak{a} \mapsto \mathfrak{a}^p.$$

One may check that for any positive divisor $d$ of $f$ the set of fixed points of $\boldsymbol{\pi}^d$ is the subfield of order $p^d$. This shows at once that $G$ has many orbits on $\mathbb{F}$, in fact at least $p$ and not necessarily of the same order. For example, if $f = 3$, then $\boldsymbol{\pi}$ is product of $p(p-1)(p+1)/3$ cycles of length 3 and so the number of orbits is $p + p(p-1)(p+1)/3$.

We are interested in the parity of $\boldsymbol{\pi}$. Of course, $\boldsymbol{\pi}$ is even if $f$ is odd. Unfortunately, when $f$ is even, the parity of $\boldsymbol{\pi}$ is somewhat less obvious. We shall appeal to the following easy lemma.

3.7.1. Lemma. *For all* $b \geqslant 1$, $(4k \pm 1)^{(2^b)} \equiv 1 \mod 2^{b+2}$.

Proof. In fact,

$$(4k \pm 1)^{(2^b)} = (\pm 1)^{(2^b)} + \sum_{i \geqslant 1} \binom{2^b}{i} (\pm 4k)^i.$$

But $\binom{2^b}{i} = \frac{2^b}{i}\binom{2^b-1}{i-1}$ and the highest power of 2 dividing $i$ (recall $i \geqslant 1$) is always dividing $4^{i-1}$. □

3.7.2. Theorem. *On a field of order* $p^f$, *p prime, the Frobenius automorphism is even unless* $p^f = 4$ *or* $f$ *is even and* $p \equiv 3 \mod 4$.

Proof. We may assume that $f = 2^a k$ with $a \geqslant 1$ and $k$ is odd. Set $q = p^k$, there is a chain of fields

$$\mathrm{GF}(q) < \mathrm{GF}\big(q^2\big) < \ldots < \mathrm{GF}\big(q^{(2^a)}\big).$$

For $b = 1, \ldots, a$, the elements of $\mathrm{GF}\big(q^{(2^b)}\big) - \mathrm{GF}\big(q^{(2^{b-1})}\big)$ are permuted by $\boldsymbol{\pi}^k$ in orbits of length $2^b$; so in $\frac{1}{2^b}\big(q^{(2^b)} - q^{(2^{b-1})}\big)$ orbits. The issue is therefore the parity of

$$\sum_{b=1}^{a} \frac{1}{2^b}\left(q^{(2^b)} - q^{(2^{b-1})}\right).$$

When $p = 2$ the sum is even unless $q = 2$ and $a = 1$ (note that for $q = 2$ the first two terms of the sum are odd).

Suppose $q$ odd. By the lemma above

$$q^{(2^b)} \equiv q^{(2^{b-1})} \equiv 1 \mod 2^{b+1}$$

whenever $b \geqslant 2$, so the summands for $b = 2, \ldots, a$ are even. While at $b = 1$ we have

$$\frac{1}{2}(q^2 - q)$$

which is even if and only if $q \equiv 1 \mod 4$, that is, $p \equiv 1 \mod 4$. □

## 3.8. Affine groups on the related vector space

Let $\mathsf{V}$ be a vector space. For $\mathsf{v} \in \mathsf{V}$ we call

$$/\mathsf{v}/ : \mathsf{V} \longrightarrow \mathrm{Sym}\,\mathsf{V}, \qquad \mathsf{u} \mapsto \mathsf{u} + \mathsf{v}$$

the map that sends $\mathsf{v}$ to the translation by $\mathsf{v}$, that is, $/\star/$ is the right regular action of $\mathsf{V}$ on itself[1]. The group of the translations on $\mathsf{V}$ and the general linear group $\mathrm{GL}(\mathsf{V})$ regarded as subgroups of $\mathrm{Sym}\,\mathsf{V}$ generate the *affine group* on $\mathsf{V}$

$$\mathrm{AGL}(\mathsf{V}) = \big\langle /\mathsf{V}/, \mathrm{GL}(\mathsf{V}) \big\rangle.$$

In fact, for all $\boldsymbol{g} \in \mathrm{GL}(\mathsf{V})$, $\mathsf{v} \in \mathsf{V}$ we have

$$\boldsymbol{g}^{-1}/\mathsf{v}/\boldsymbol{g} = /\mathsf{v}\boldsymbol{g}/$$

and so $\mathrm{AGL}(\mathsf{V}) = /\mathsf{V}/ \rtimes \mathrm{GL}(\mathsf{V})$.

As usual we write $\mathrm{AGL}_n(q)$ for the affine group on a vector space of dimension $n$ over a field of order $q$.

**Degree of transitivity, primitivity.** Such degrees depend on the parameters $n$, $q$ as showed in Table 3.B. Rows 1, 2 and 5 follow from the permutational

TABLE 3.B. Degree of transitivity, primitivity of $\mathrm{AGL}_n(q)$

| $n$ | $q$ | degree of transitivity | degree of primitivity | row |
|---|---|---|---|---|
| 1 | 2 | 2 | 2 | 1 |
| 1 | 3 | 3 | 3 | 2 |
| 1 | $q-1$ is a Mersenne prime | 2 | 2 | 3 |
| 1 | all other $q$ | 2 | 1 | 4 |
| 2 | 2 | 4 | 4 | 5 |
| $>2$ | 2 | 3 | $\geqslant 2$ | 6 |
| $>1$ | $\geqslant 3$ | 2 | 1 | 7 |

isomorphisms $\mathrm{AGL}_1(2) \cong \mathrm{S}_2$, $\mathrm{AGL}_1(3) \cong \mathrm{S}_3$ and $\mathrm{AGL}_2(2) \cong \mathrm{S}_4$ (which together with $\mathrm{AGL}_1(4) \cong \mathrm{A}_4$ provide the complete list of permutational isomorphisms with symmetric or alternating groups). Affine groups may not send 3 collinear points to 3 non collinear points, it follows that $\mathrm{AGL}_n(q)$ is not 3-transitive provided $n \geqslant 2$ and $q \geqslant 3$. However, $\mathrm{AGL}_n(q)$ is 2-transitive because $\mathrm{GL}_n(q)$ is clearly transitive on the non zero points. Furthermore, the lines through the origin deprived of the

[1] see §5.10 for an explanation of this notation.

origin are non trivial blocks for $\mathrm{GL}_n(q)$ so that $\mathrm{AGL}_n(q)$ is not 2-primitive. This justifies row 7. Similarly, $\mathrm{AGL}_n(2)$ may not send a plane to 4 non complanar points so that if $n > 2$, then $\mathrm{AGL}_n(2)$ is not 4-transitive. However, if $n \geqslant 2$, then $\mathrm{GL}_n(2)$ is 2-transitive on the non zero points thus $\mathrm{AGL}_n(2)$ is 3-transitive as reported in row 6. We remain with rows 3 and 4, that is $\mathrm{AGL}_1(q)$ when $q > 3$. The stabilizer of one point is a cyclic regular subgroup on the remaining $q-1$ points. Therefore it is transitive but not 2-transitive; and it is primitive if and only if $q-1$ is prime, that is, $q$ is a power of 2 and so $q-1$ is a *Mersenne prime*[2].

**Parity.** The parity of affine groups is ruled by the parity of $q$ with only two exceptions.

3.8.1. Proposition. $\mathrm{AGL}_n(q)$ *is even if and only if* $q$ *is even except for* $(n, q) = (1, 2),\ (2, 2)$.

Proof. if $(n, q) \neq (1, 2)$, then all translations are even. In fact, each proper translation has order $p$, where $q = p^r$, and moves all points. Therefore it is a product of $q^n/p$ cycles of length $p$. If $p$ is odd, then cycles of length $p$ are even. If $p = 2$ but $(n, q) \neq (1, 2)$, then $4 \mid q^n$, hence $q^n/p$ is even.

If $n = 1$ we observe that $\mathrm{GL}_1(q)$ is generated by a cycle of length $q-1$ therefore $\mathrm{AGL}_1(q)$ is even if and only if $q$ is even (and $q \neq 2$).

Assume now $n > 1$ and consider $\boldsymbol{u} \in \mathrm{GL}_n(q)$ such that $e_0\boldsymbol{g} = \mathfrak{u}e_0$, $e_i\boldsymbol{g} = e_i$ for all $i = 1 \dots n-1$, where $\vec{e} = (e_0, \dots, e_{n-1})$ is a basis of the vector space and $\mathfrak{u}$ generates the multiplicative group of the field. This $\boldsymbol{u}$ is product of cycles of length $(q-1)$ and fixes exactly the points whose $e_0$-component is zero: there are $q^{n-1}$ such points. Hence $\boldsymbol{u}$ is the product of $q^{n-1}$ cycles, each of length $q-1$. As a consequence, if $q$ is odd, then $\boldsymbol{u}$ is odd.

It remains to prove that $\mathrm{GL}_n(q)$ is even when $n > 1$, $q$ is even and $(n, q) \neq (2, 2)$. To do so, we observe that $\mathrm{GL}_n(q) = \langle \mathrm{SL}_n(q), \boldsymbol{u} \rangle$ and we prove that $\mathrm{SL}_n(q)$ is even. By [**Rob93**, 3.2.10] $\mathrm{SL}_n(q)$ is generated by *transvections* $\boldsymbol{T}_{ij}(\mathfrak{f})$ where $\mathfrak{f}$ is a non zero scalar of the field and

$$e_h\boldsymbol{T}_{ij}(\mathfrak{f}) = e_h + \delta_{hj}\mathfrak{f}e_i.$$

Such a $\boldsymbol{T}_{ij}(\mathfrak{f})$ fixes exactly the points whose $e_j$-component is zero: there are $q^{n-1}$ of these points. Moreover, the non trivial orbits of $\boldsymbol{T}_{ij}(\mathfrak{f})$ all have length 2. Therefore $\boldsymbol{T}_{ij}(\mathfrak{f})$ is product of $(q-1)q^{n-1}/2$ transpositions and is even because $4 \mid q^{n-1}$. □

[2]`http://www.utm.edu/research/primes/mersenne.shtml`

**Full normalizer of the translation subgroup.** If $q$ is a power of the prime $p$, then $\mathsf{V}$ may be regarded as a vector space, say $\mathsf{V}_p$, over the prime subfield $\mathbb{F}_p$. Since $/\mathsf{V}/ = /\mathsf{V}_p/$ and $\mathrm{GL}(\mathsf{V})$ is a subgroup of $\mathrm{GL}(\mathsf{V}_p)$, we have that $\mathrm{AGL}(\mathsf{V}) \leqslant \mathrm{AGL}(\mathsf{V}_p)$. This shows that the normalizer of $/\mathsf{V}/$ is at least as large as $\mathrm{AGL}(\mathsf{V}_p)$. And it may not be any larger (see [**DM96**, Cor. 4.2B] for example).

## 3.9. Projective groups on their projective space

Let $\mathsf{V}$ be a vector space of dimension at least 2. The set of subspaces of $\mathsf{V}$ of dimension $l$ is known as the *Grassmannian* $\mathrm{G}_l(\mathsf{V})$. It is clear that each $\boldsymbol{g} \in \mathrm{GL}(\mathsf{V})$ induces a canonical bijection

$$\mathrm{G}_l(\boldsymbol{g}) : \mathrm{G}_l(\mathsf{V}) \longrightarrow \mathrm{G}_l(\mathsf{V}).$$

The *projective space* on $\mathsf{V}$, $\mathrm{PG}(\mathsf{V})$, may be identified with $\mathrm{G}_1(\mathsf{V})$. We write $\mathrm{PG}_d(q)$ for the projective space on a vector space of dimension $d$ over a field of order $q$. The number of points of $\mathrm{PG}_d(q)$ is

$$\frac{q^d - 1}{q - 1} = 1 + q + \cdots + q^{d-1}.$$

This number may exceptionally be a prime power. Here we are only interested to the power of 2 case.

3.9.1. Lemma. *Suppose that $1 + q + \cdots + q^{n-1}$ is a power of 2 for some positive integer $q \neq 1$, then $n \leqslant 2$.*

Proof. If $n > 1$ then $q$ must be odd, hence $n$ must be even. Say $n = 2e$ so that we have

$$2^l = \frac{q^e - 1}{q - 1}(q^e + 1).$$

If $e = 1$ then $n = 2$ as wanted. Otherwise $1 + q + \cdots + q^{e-1}$ is a positive power of 2 and $e$ is even, say $e = 2f$. In particular $q^e + 1 = (q^2)^f + 1$ is a power of 2 and also congruent to 2 modulo 4. This forces $q^e + 1 = 2$ and hence $n = e = 0$, a contradiction. □

Since we are only interested to the case when $d > 1$, we see that if the number of points of $\mathrm{PG}_d(q)$ is a power of 2, then $d = 2$ and so $q$ is a power of 2 minus 1. In particular, if $q$ is prime, then $q$ is a Mersenne prime.

**The action.** The Projective Special Linear group $\mathrm{PSL}_d(q)$ has a faithful 2-transitive action on $\mathrm{PG}_d(q)$. Apart from the cases when $(d, q) = (2, 2)$ or $(2, 3)$, $\mathrm{PSL}_d(q)$ is a non abelian simple group.

To exhibit suitable generators for the normalizer of $\mathrm{PSL}_d(q)$ in $\mathrm{Sym}\,\mathrm{PG}_d(q)$, say $(e_0, \ldots, e_{d-1})$ a basis of the underlying vector space, $\mathfrak{u}$ a generator of the multiplicative group of the field and $\boldsymbol{\pi}$ the Frobenius automorphism of the field: $\mathfrak{a} \mapsto \mathfrak{a}^p$, where $p$ is the characteristic of the field.

The Projective General Linear group $\mathrm{PGL}_d(q)$ is generated by $\mathrm{PSL}_d(q)$ together with the map $\mathrm{G}_1(\boldsymbol{u})$ where $\boldsymbol{u}$ is defined by

$$e_0\boldsymbol{u} = \mathfrak{u}e_0, \qquad e_i\boldsymbol{u} = e_i \quad \text{for all } i = 1 \ldots d-1.$$

The normalizer of $\mathrm{PSL}_d(q)$ in $\mathrm{Sym}\,\mathrm{PG}_d(q)$ is $\mathrm{P\Gamma L}_d(q)$ which is generated by $\mathrm{PGL}_d(q)$ together with $\mathrm{PG}_d(q, \boldsymbol{\pi})$ where for each automorphism $\boldsymbol{\alpha}$ of the field, $\mathrm{PG}_d(q, \boldsymbol{\alpha})$ is the permutation induced on $\mathrm{PG}_d(q)$ by $\boldsymbol{\alpha}$. The subgroup generated by $\mathrm{PSL}_d(q)$ and $\mathrm{PG}_d(q, \boldsymbol{\pi})$ is denoted by $\mathrm{P\Sigma L}_d(q)$.

Observe that the stabilizer in $\mathrm{P\Gamma L}_d(q)$ of a projective point must preserve the system made of the projective lines through that point. Since projective lines have at least 3 points, $\mathrm{P\Gamma L}_d(q)$ may not be 2-primitive for $d > 2$. However, $\mathrm{PGL}_2(q)$ is 3-transitive on $\mathrm{PG}_2(q)$ and so $\mathrm{P\Gamma L}_2(q)$ is also 3-transitive. Of course, they are both 2-primitive. The order of $\mathrm{PGL}_2(q)$ is $(q+1)q(q-1)$ therefore[3] no subgroup of $\mathrm{PGL}_2(q)$ is 3-transitive. We remind that $\mathrm{PSL}_2(q) = \mathrm{PGL}_2(q)$ if $q$ is even but that $\mathrm{PSL}_2(q)$ has index 2 in $\mathrm{PGL}_2(q)$ when $q$ is odd.

3.9.2. Proposition. *There is one and only one conjugacy class of transitive subgroups of* $\mathrm{S}_{q+1}$ *isomorphic to* $\mathrm{PSL}_2(q)$

Proof. Say $p$ the characteristic of the field. Observe first that a $p$-Sylow subgroups of $\mathrm{GL}_d(q)$ is conjugate to the subgroup of the upper unitriangular matrices (which has a faithful image in $\mathrm{PSL}_d(q)$) and has $d$ orbits on $\mathrm{PG}_d(q)$ of order 1, $q$, …, $q^{d-1}$. Therefore a $p$-Sylow subgroup of $\mathrm{PSL}_2(q)$ has exactly 2 orbits of order 1 and $q$ respectively.

A transitive subgroup of $\mathrm{S}_{q+1}$ isomorphic to $\mathrm{PSL}_2(q)$ has the stabilizer of a point isomorphic to a subgroup of index $q+1$ of $\mathrm{PSL}_2(q)$. In view of 3.1.1, it is more than enough to show that there is only one conjugacy class of these subgroups in $\mathrm{PSL}_2(q)$.

A subgroup of index $q+1$ of $\mathrm{PSL}_2(q)$ has order dividing $q(q-1)$ and contains a $p$-Sylow subgroup. It is not transitive on $\mathrm{PG}_2(q)$ because $q+1$ may not divide $q-1$. Therefore it has on it the same orbits of its $p$-Sylow subgroups. In particular, it is the stabilizer of a projective point in $\mathrm{PSL}_2(q)$. But the stabilizers of projective points form a unique conjugacy class. □

[3] A $k$-transitive group of degree $n$ and order $n(n-1)\cdots(n-k+1)$ is said *sharply* $k$-transitive.

**Parity.**

3.9.3. PROPOSITION. $\mathrm{PGL}_d(q)$ *is even if and only if* $(d+1)q$ *is even and larger than* 6.

PROOF. Note first that $\mathrm{PGL}_2(2)$ and $\mathrm{PGL}_2(3)$ are odd because they are permutationally isomorphic to $\mathrm{S}_3$ and $\mathrm{S}_4$ respectively. Suppose now that $(d,q)$ is not equal to $(2,2)$ or $(2,3)$. Then $\mathrm{PSL}_d(q)$ is simple, hence even. Since $\mathrm{PGL}_d(q)$ is generated by $\mathrm{PSL}_d(q)$ together with $\mathrm{G}_1(\boldsymbol{u})$, it is enough to determine the sign of $\mathrm{G}_1(\boldsymbol{u})$.

It is easy to check that $\mathrm{G}_1(\boldsymbol{u})$ fixes exactly the points with homogeneous coordinates

$$[0, \mathfrak{f}_1, \ldots, \mathfrak{f}_{d-1}], \qquad [1, 0, \ldots, 0].$$

Their number is $1 + \frac{q^{d-1}-1}{q-1}$ so that $\mathrm{G}_1(\boldsymbol{u})$ moves $q^{d-1}-1$ points in orbits of length $q-1$. Thus $\mathrm{G}_1(\boldsymbol{u})$ is made of cycles of length $(q-1)$ in number of $\frac{q^{d-1}-1}{q-1}$.

If $q$ is even, then $q-1$ is odd and $\mathrm{G}_1(\boldsymbol{u})$ is even.

If $q$ is odd, then $\operatorname{Par} \mathrm{G}_1(\boldsymbol{u}) \equiv 1 + q + \ldots + q^{d-2} \mod 2$. This is congruent to $1 + (d-2)q$ which in turn is congruent to $1 + d$.

Therefore $\mathrm{G}_1(\boldsymbol{u})$ is even if and only if $q$ is even or $d$ is odd, if and only if $(d+1)q$ is even. □

3.9.4. PROPOSITION. *Let* $\boldsymbol{\alpha}$ *be an automorphism of the field of order* $q$. *Then*

$$\operatorname{Par}\big(\mathrm{PG}_d(q,\boldsymbol{\alpha})\big) \underset{(2)}{\equiv} \left(1 + q\frac{(d+1)(d-2)}{2}\right) \operatorname{Par} \boldsymbol{\alpha}.$$

PROOF. We regard $\mathrm{PG}_d(q)$ as the disjoint union of the following $\mathrm{PG}_d(q,\boldsymbol{\alpha})$ invariant sets

$$\{1\} \times \mathbb{F}_q^{d-1}, \{(0,1)\} \times \mathbb{F}_q^{d-2}, \ldots, \{(\underbrace{0,\ldots,0}_{d-2},1)\} \times \mathbb{F}_q, \{(\underbrace{0,\ldots,0}_{d-1},1)\}$$

where $\mathbb{F}_q$ denotes the field of order $q$. Then $\mathrm{PG}_d(q,\boldsymbol{\alpha})$ may be regarded as the external product

$$\boldsymbol{\alpha\vartheta}_{d-1} \times \boldsymbol{\alpha\vartheta}_{d-2} \times \cdots \times \boldsymbol{\alpha\vartheta}_2 \times \mathrm{id}$$

where $\boldsymbol{\vartheta}_i$ is the diagonal action defined in 3.5. hence its sign is congruent modulo 2 to

$$(1 + q\sum_{i=2}^{d-1} i) \operatorname{Par} \boldsymbol{\alpha} \underset{(2)}{\equiv} \left(1 + q\frac{(d+1)(d-2)}{2}\right) \operatorname{Par} \boldsymbol{\alpha}.$$

□

In particular, $\mathrm{PG}_d(q,\boldsymbol{\pi})$ is even whenever $\boldsymbol{\pi}$ is even (see 3.7.2). However, one can readily check that $\mathrm{PG}_d(4,\boldsymbol{\pi})$ is odd. The other cases are summarized as follows:

3.9.5. COROLLARY. *Let $q = p^f$ where $p$ is prime, $p \equiv 3 \mod 4$ and $f$ is even. $\mathrm{PG}_d(q, \boldsymbol{\pi})$ is even if and only if $d \equiv 0, 1 \mod 4$.*

From all these results one gets the even part of $\mathrm{P\Gamma L}_d(q)$ as in table 3.C.

TABLE 3.C. Even part of $\mathrm{P\Gamma L}_d(p^f)$, $(2,2) \neq (d, p^f) \neq (2,3)$

| $p$ | $f$ | $d$ | $\mathrm{P\Gamma L}_d(p^f)^{\mathrm{e}}$ |
|---|---|---|---|
| 2 | 2 | all | PGL |
| 2 | $\neq 2$ | all | PΓL |
| $1 \mod 4$ | all | even | PΣL |
| | | odd | PΓL |
| $3 \mod 4$ | odd | even | PΣL |
| | | odd | PΓL |
| $3 \mod 4$ | even | $0 \mod 4$ | PΣL |
| | | $1 \mod 4$ | PΓL |
| | | $2 \mod 4$ | $\left\langle \mathrm{PSL}, \mathrm{PG}(\boldsymbol{\pi}^2), \mathrm{G}_1(\boldsymbol{u})\,\mathrm{PG}(\boldsymbol{\pi}) \right\rangle$ |
| | | $3 \mod 4$ | $\left\langle \mathrm{PGL}, \mathrm{PG}(\boldsymbol{\pi}^2) \right\rangle$ |

### 3.10. Doubly transitive Suzuki groups

We follow [**Suz62**, §13] and construct the *Suzuki groups* as subgroups of the symplectic groups $\mathrm{Sp}_4(q)$ where $q = 2^{2a+1}$ for some positive integer $a$. Here we identify $\mathrm{Sp}_4(q)$ with the $4 \times 4$ invertible matrices $\mathsf{M}$ with entries in $\mathbb{F} = \mathbb{F}_q$ such that $\mathsf{M}^t \mathsf{T} \mathsf{M} = \mathsf{T}$, where $\mathsf{T}$ is the matrix associated to a standard non degenerate bilinear symplectic form:

$$\mathsf{T} = \begin{bmatrix} 0 & 0 & 0 & 1 \\ 0 & 0 & 1 & 0 \\ 0 & 1 & 0 & 0 \\ 1 & 0 & 0 & 0 \end{bmatrix}.$$

Note that in view of the isomorphism $\mathrm{Sp}_4(q) \cong \mathrm{PSp}_4(q)$ we adopt language and notation of projective geometry and surround matrices with square brakets. If $r^2 = 2q$, the mapping $\boldsymbol{\theta} : \mathfrak{a} \mapsto \mathfrak{a}^r$ is an automorphism of $\mathbb{F}$ whose square is the Frobenius automorphism: $\boldsymbol{\pi} : \mathfrak{a} \mapsto \mathfrak{a}^2$. Let $\mathfrak{a}$ and $\mathfrak{b}$ be arbitrary elements of $\mathbb{F}$ and let $(\mathfrak{a}, \mathfrak{b})$ denote the matrix

$$(\mathfrak{a}, \mathfrak{b}) = \begin{bmatrix} 1 & 0 & 0 & 0 \\ \mathfrak{a} & 1 & 0 & 0 \\ \mathfrak{a}^{1+\boldsymbol{\theta}} + \mathfrak{b} & \mathfrak{a}^{\boldsymbol{\theta}} & 1 & 0 \\ \mathfrak{a}^{2+\boldsymbol{\theta}} + \mathfrak{a}\mathfrak{b} + \mathfrak{b}^{\boldsymbol{\theta}} & \mathfrak{b} & \mathfrak{a} & 1 \end{bmatrix},$$

where for each integer $i$ we denote $\mathfrak{a}^i\mathfrak{a}^{\boldsymbol{\theta}}$ by $\mathfrak{a}^{i+\boldsymbol{\theta}}$. Matrix multiplication gives

$$(\mathfrak{a},\mathfrak{b})(\mathfrak{c},\mathfrak{d}) = (\mathfrak{a}+\mathfrak{c}, \mathfrak{a}\mathfrak{c}^{\boldsymbol{\theta}}+\mathfrak{b}+\mathfrak{d}).$$

The totality $Q(q)$ of these matrices $(\mathfrak{a},\mathfrak{b})$ forms a subgroup of order $q^2$ and exponent 4 of $\mathrm{Sp}_4(q)$ where $(\mathfrak{a},\mathfrak{b})$ has order 2 if anf only if $\mathfrak{a} = 0 \neq \mathfrak{b}$. We call $D(q)$ the subgroup of $\mathsf{Q}(q)$ made of its $q-1$ involutions.

To each non zero element $\mathfrak{u}$ of the field we associate a diagonal matrix, denoted by $(\mathfrak{u})$, with $\mathfrak{s}_1$, $\mathfrak{s}_2$, $\mathfrak{s}_3$ and $\mathfrak{s}_4$ in the main diagonal, where $\mathfrak{s}_1^{\boldsymbol{\theta}} = \mathfrak{u}^{1+\boldsymbol{\theta}}$, $\mathfrak{s}_2^{\boldsymbol{\theta}} = \mathfrak{u}$, $\mathfrak{s}_3 = \mathfrak{s}_2^{-1}$ and $\mathfrak{s}_4 = \mathfrak{s}_1^{-1}$. The totality $K(q)$ of the matrices $(\mathfrak{u})$ forms a cyclic subgroup of order $q-1$ of $\mathrm{Sp}_4(q)$ isomorphic to $\mathbb{F}^\times$. Matrix multiplication yelds that

$$(\mathfrak{u})^{-1}(\mathfrak{a},\mathfrak{b})(\mathfrak{u}) = (\mathfrak{a}\mathfrak{u}, \mathfrak{b}\mathfrak{u}^{1+\boldsymbol{\theta}}).$$

Therefore the group $H(q)$ generated by $Q(q)$ and $K(q)$ is a group of order $q^2(q-1)$. Let $\tau$ denote the matrix $\mathsf{T}$; clearly $\tau$ lies in $\mathrm{Sp}_4(q)$. The Suzuki group $\mathrm{Sz}(q)$ is the subgroup of $\mathrm{Sp}_4(q)$ generated by $H(q)$ and $\tau$. There is just one conjugacy class of the subgroups of $\mathrm{Sp}_4(q)$ which are isomorphic to $\mathrm{Sz}(q)$.

3.10.1. Proposition. *If* $q \geqslant 8$, $\mathrm{Sz}(q)$ *is a simple group of order* $q^2(q-1)(q^2+1)$ *acting doubly transitively on the Suzuki ovoid*

$$\mathcal{O} = \left\{ \begin{bmatrix} 1 \\ \mathfrak{a} \\ \mathfrak{a}^{1+\boldsymbol{\theta}} + \mathfrak{b} \\ \mathfrak{a}^{2+\boldsymbol{\theta}} + \mathfrak{a}\mathfrak{b} + \mathfrak{b}^{\boldsymbol{\theta}} \end{bmatrix} \;\middle|\; \mathfrak{a},\mathfrak{b} \in \mathbb{F} \right\} \cup \left\{ \begin{bmatrix} 0 \\ 0 \\ 0 \\ 1 \end{bmatrix} \right\}$$

*and transitively on the other points of* $\mathrm{PG}_4(q)$.

(See [**Lie87**, Lemma 2.11] for the last assertion).

To save some writing, put $S = \mathrm{Sz}(q)$, $o = [1\ 0\ 0\ 0]^t$ and $\infty = [0\ 0\ 0\ 1]^t$. One readily checks that $Q(q)$ is a 2-Sylow subgroup of $S$ which acts regularly on $\mathcal{O}-\{\infty\}$, that $S_\infty = H(q)$ and that $S_{(\infty,o)} = K(q)$. Provided

$$S_{(\infty,o)} < D(q)K(q) < S_\infty,$$

$S$ is not doubly primitive on $\mathcal{O}$. In particular, it is not 3-transitive either.

The field automorphisms of $\mathrm{Sp}_4(q)$ normalize $S$ because they commute with $\boldsymbol{\theta}$ and it is shown in [**Suz62**] that $S$ has no other outer automorphisms. In particular $\mathrm{Aut}\, S$ is the extension of $S$ by a cyclic group of order $2a+1$. Note that the automorphisms act on $\mathcal{O}$ and so

$$\mathrm{Aut}\, S \hookrightarrow \mathrm{Sym}\, \mathcal{O} \cong \mathrm{S}_{q^2+1}.$$

In fact, $\mathrm{Aut}\, S \hookrightarrow \mathrm{Alt}\, \mathcal{O}$ because $2a+1$ is odd.

One proves in a completely similar way to 3.9.2 that there is one and only one conjugacy class of subgroups of $S$ of index $q^2+1$. This shows that there is one and only one conjugacy class of transitive subgroups of $\operatorname{Sym}\mathcal{O}$ isomorphic to $S$. There is also a primitive action of $S$ on the 2-homogeneous component of $\mathcal{P}(\mathcal{O})$, that is, the set of the subsets of order 2 of $\mathcal{O}$, which we denote by $\binom{\mathcal{O}}{2}$. In fact, [**Suz62**, §15] shows that $S$ has one and only one conjugacy class of maximal subgroups of order $2(q-1)$, therefore there is one and only one conjugacy class of primitive subgroups of $\operatorname{Sym}\binom{\mathcal{O}}{2}$ isomorphic to $S$.

CHAPTER 4

# The intransitive case

Recall the notation established in 2.1 where $\boldsymbol{S} \equiv \operatorname{Sym}\Omega$ is a finite symmetric group of degree $n \geqslant 5$, $\boldsymbol{A}$ is the related alternating group, $\boldsymbol{H}$ is a second maximal subgroup of $\boldsymbol{S}$ or $\boldsymbol{A}$, and $\boldsymbol{U} = \boldsymbol{H}\boldsymbol{A}$:

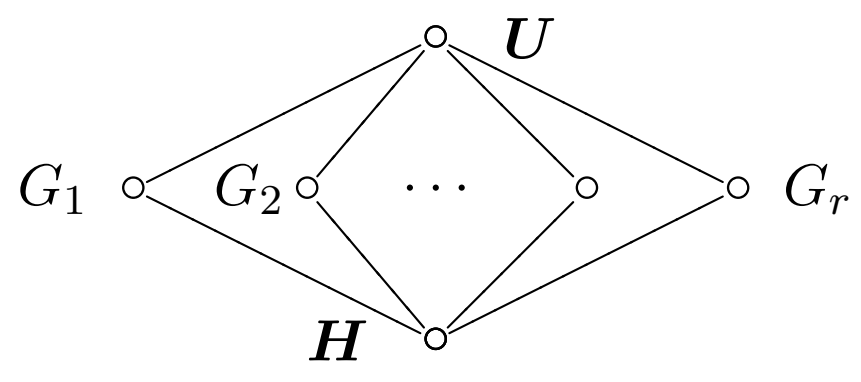

Throughout this chapter we make the further assumption that $\boldsymbol{H}$ is an intransitive subgroup of $\boldsymbol{S}$ and we prove the following result.

**Theorem A.** *If $r > 2$, then either $\boldsymbol{H}$ is intransitive with 3 orbits and contained in precisely 3 maximal intransitive subgroups, or $\boldsymbol{U} = \boldsymbol{A}$ and either $\boldsymbol{H}$ is the stabilizer of a point in* $\mathrm{PSL}(3,2)$ *of degree 7 as in 2.6.1, or $\boldsymbol{H}$ is the stabilizer of a point in a 2-primitive maximal subgroup. In the latter case, it is contained in the stabilizer of that point in $\boldsymbol{A}$, in another 2-primitive group which is conjugate to the given one in $\boldsymbol{S}$ (but not in $\boldsymbol{A}$), and in nothing else. So, in any case $r = 3$.*

We follow the outline given in [**Pál88**] for the intransitive second maximal subgroups of the alternating groups of prime degree because those ideas apply successfully to any degree and to the intransitive second maximal subgroups of the symmetric groups as well. However, more care is needed here because we may not assume that each $G_i$ is primitive.

### 4.1. There are at most three orbits

Assume that $\boldsymbol{H}$ has at least three orbits and choose two orbits $\mathsf{X}_1$, $\mathsf{X}_2$ in a way that

$$|\mathsf{X}_1| \leqslant |\mathsf{X}_2| \leqslant |\mathsf{X}_3|$$

where $\mathsf{X}_3$ is the union of the other orbits. It is enough to show that $\mathsf{X}_3$ is an orbit. Surely, $\mathsf{X}_3$ is an orbit of $\boldsymbol{U}_{(\mathsf{X}_1,\mathsf{X}_2,\mathsf{X}_3)}$ because $n > 4$ and so $\mathsf{X}_2$ has order at least 2. Thus we show that $\boldsymbol{H} = \boldsymbol{U}_{(\mathsf{X}_1,\mathsf{X}_2,\mathsf{X}_3)}$.

Pick $\omega_1 \in \mathsf{X}_1$, $\omega_2 \in \mathsf{X}_2$ and $\omega_3, \omega_3' \in \mathsf{X}_3$ with $\omega_3 \neq \omega_3'$. It is easy to check that

$$(\omega_1\ \omega_2\ \omega_3) \in \boldsymbol{U} - \boldsymbol{U}_{(\mathsf{X}_1,\mathsf{X}_2\cup\mathsf{X}_3)},$$
$$(\omega_2\ \omega_3\ \omega_3') \in \boldsymbol{U}_{(\mathsf{X}_1,\mathsf{X}_2\cup\mathsf{X}_3)} - \boldsymbol{U}_{(\mathsf{X}_1,\mathsf{X}_2,\mathsf{X}_3)}.$$

This shows that

$$\boldsymbol{H} \leqslant \boldsymbol{U}_{(\mathsf{X}_1,\mathsf{X}_2,\mathsf{X}_3)} < \boldsymbol{U}_{(\mathsf{X}_1,\mathsf{X}_2\cup\mathsf{X}_3)} < \boldsymbol{U};$$

but $\boldsymbol{H}$ is second maximal in $\boldsymbol{U}$ and hence $\boldsymbol{H} = \boldsymbol{U}_{(\mathsf{X}_1,\mathsf{X}_2,\mathsf{X}_3)}$.

## 4.2. Second maximal subgroups with three orbits

In this section we assume that $\boldsymbol{H}$ has precisely three orbits: $\mathsf{X}_1$ of order $n_1$, $\mathsf{X}_2$ of order $n_2$ and $\mathsf{X}_3$ of order $n_3$, with $n_1 \leqslant n_2 \leqslant n_3$. We have just seen that $\boldsymbol{H} = \boldsymbol{U}_{(\mathsf{X}_1,\mathsf{X}_2,\mathsf{X}_3)}$ in this case. We may also assume that $n_1 + n_2 \neq n_3$, otherwise

$$\boldsymbol{H} < \boldsymbol{U}_{(\mathsf{X}_1\cup\mathsf{X}_2,\mathsf{X}_3)} < \boldsymbol{U}_{\{\mathsf{X}_1\cup\mathsf{X}_2,\mathsf{X}_3\}} < \boldsymbol{U}.$$

Also, if $n_i = n_j > 1$, then

$$\boldsymbol{H} < \boldsymbol{U}_{\{\mathsf{X}_i,\mathsf{X}_j\}} < \boldsymbol{U}_{\mathsf{X}_i\cup\mathsf{X}_j} < \boldsymbol{U}$$

which cannot happen. Thus, either $n_1 = n_2 = 1$ or $n_1 < n_2 < n_3$. In any case $2n_3 \neq n < 3n_3$ and $n_3 > 2$.

4.2.1. Proposition. *If $\boldsymbol{H}$ has three orbits then $[\boldsymbol{H} \div \boldsymbol{U}] = \mathcal{M}_3$.*

Proof. The subgroup $\boldsymbol{H} = \boldsymbol{U}_{(\mathsf{X}_1,\mathsf{X}_2,\mathsf{X}_3)}$ must contain 3-cycles on $\mathsf{X}_3$. By 2.5.3 there are not proper primitive groups containing $\boldsymbol{H}$ and we claim that all the proper subgroups of $\boldsymbol{U}$ containing $\boldsymbol{H}$ are intransitive as they cannot be transitive imprimitive. It is enough to show that any congruence preserved by $\boldsymbol{H}$ is trivial.

Since $2n_3 \neq n < 3n_3$, $\mathsf{X}_3$ is not a block of $\boldsymbol{H}$. Also, if $\mathsf{B}$ is a nontrivial block of $\boldsymbol{H}$ containing $\mathsf{X}_3$, then $\Omega - \mathsf{B}$ must be a block too. Assume now that $\mathsf{Z}$ is a congruence preserved by $\boldsymbol{H}$. Consider

$$\mathsf{Y} = \mathsf{Z} \cap \mathsf{X}_3 := \{\ \mathsf{B} \cap \mathsf{X}_3 \ \mid\ \mathsf{B} \in \mathsf{Z}\ \}$$

which is a partition of $\mathsf{X}_3$ preserved by $\boldsymbol{H}$ acting on $\mathsf{X}_3$. Since $\boldsymbol{H}$ is primitive on $\mathsf{X}_3$, $\mathsf{Y}$ must be trivial. Then either $\mathsf{X}_3 \subset \mathsf{B}$ for some $\mathsf{B} \in \mathsf{Z}$, or $|\mathsf{B} \cap \mathsf{X}_3| \leqslant 1$ for all $\mathsf{B} \in \mathsf{Z}$.

In the first case $\mathsf{B} \cap \mathsf{X}_i$ is not empty for some $i \neq 3$ and $\Omega - \mathsf{B}$ is a block which must intersect both $\mathsf{X}_1$ and $\mathsf{X}_2$ (as they alone are too small). Then pick $\omega_{i,1} \in \mathsf{X}_i \cap \mathsf{B}$, $\omega_{i,2} \in \mathsf{X}_i - \mathsf{B}$, $\omega_{3,1}, \omega_{3,2} \in \mathsf{X}_3$. The double transposition $(\omega_{i,1}\ \omega_{i,2})(\omega_{3,1}\ \omega_{3,2})$ lies in $\boldsymbol{H}$ and splits $\mathsf{B}$, against our assumptions that $\mathsf{B}$ is a block of $\boldsymbol{H}$.

Therefore we remain with the second case: $\left|\mathsf{B} \cap \mathsf{X}_3\right| \leqslant 1$ for all $\mathsf{B} \in \mathsf{Z}$. In particular, $\left|\mathsf{Z}\right| \geqslant n_3$ and there are at least three distinct blocks $\mathsf{B}_1$, $\mathsf{B}_2$, $\mathsf{B}_3$ in $\mathsf{Z}$ intersecting $\mathsf{X}_3$. Assume that they are not trivial, that is, they have at least order 2. Either $\mathsf{X}_1$ or $\mathsf{X}_2$ intersects two of them, say that $\mathsf{X}_1$ intersects $\mathsf{B}_1$, $\mathsf{B}_2$ for example. Pick $\omega_1 \in \mathsf{X}_1 \cap \mathsf{B}_1$, $\omega_2 \in \mathsf{X}_1 \cap \mathsf{B}_2$ and a transposition $\boldsymbol{s}$ of $\mathsf{X}_3 - \mathsf{B}_1$. Again, the double transposition $(\omega_1\ \omega_2)\boldsymbol{s}$ lies in $\boldsymbol{H}$ and splits $\mathsf{B}_1$. This is impossible therefore all the $\mathsf{B}_i$ must have order 1, which is what we wanted to prove.

As a consequence, all the maximal subgroups of $\boldsymbol{U}$ containing $\boldsymbol{H}$ are intransitive and their orbits must be union of orbits of $\boldsymbol{H}$. There are only three groups satisfying these conditions, namely

$$\boldsymbol{U}_{(\mathsf{X}_1 \cup \mathsf{X}_2, \mathsf{X}_3)}, \quad \boldsymbol{U}_{(\mathsf{X}_2 \cup \mathsf{X}_3, \mathsf{X}_1)}, \quad \boldsymbol{U}_{(\mathsf{X}_1 \cup \mathsf{X}_3, \mathsf{X}_2)}.$$

Hence $r = 3$. □

## 4.3. Second maximal subgroups with two orbits

In this section we assume that $\boldsymbol{H}$ has precisely two orbits: $\mathsf{X}_1$ of order $n_1$ and $\mathsf{X}_2$ of order $n_2$, with $n_1 \leqslant n_2$. Let us deal first with the case where $n_1 = n_2$.

4.3.1. Proposition. *If* $n_1 = n_2$ *then* $[\boldsymbol{H} \div \boldsymbol{U}] = \mathcal{M}_1$.

Proof. Since $\boldsymbol{H} \leqslant \boldsymbol{U}_{(\mathsf{X}_1, \mathsf{X}_2)} < \boldsymbol{U}_{\{\mathsf{X}_1, \mathsf{X}_2\}} < \boldsymbol{U}$ and since $\boldsymbol{H}$ is high, we conclude

$$\boldsymbol{H} = \boldsymbol{U}_{(\mathsf{X}_1, \mathsf{X}_2)}.$$

Provided $n \geqslant 6$, $n_i \geqslant 3$ and $\boldsymbol{H}$ contains 3-cycles. In particular, there are no primitive maximal subgroups of $\boldsymbol{U}$ containing $\boldsymbol{H}$. Moreover, $\boldsymbol{H}$ acts primitively on both $\mathsf{X}_1$ and $\mathsf{X}_2$. If $\mathsf{B}$ is a non trivial block for $\boldsymbol{H}$ properly containing one of the $\mathsf{X}_i$, then $n \geqslant 2\left|\mathsf{B}\right| > 2n_i = n$, a contradiction. If $\mathsf{B}$ is a non trivial block for $\boldsymbol{H}$ intersecting both $\mathsf{X}_1$ and $\mathsf{X}_2$, then $\left|\mathsf{B} \cap \mathsf{X}_i\right| = 1$ for each $i$ and there are at least two other blocks $\mathsf{B}_1$, $\mathsf{B}_2$ (in the same imprimitivity system of some transitive groups containing $\boldsymbol{H}$) intersecting both $\mathsf{X}_1$ and $\mathsf{X}_2$. The product of a transposition that swaps $\mathsf{B} \cap \mathsf{X}_1$ with $\mathsf{B}_1 \cap \mathsf{X}_1$ and a transposition that swaps $\mathsf{B} \cap \mathsf{X}_2$ with $\mathsf{B}_2 \cap \mathsf{X}_2$ lies in $\boldsymbol{H}$ and splits $\mathsf{B}$, again a contradiction.

Therefore the only non trivial blocks for $\boldsymbol{H}$ are $\mathsf{X}_1$ and $\mathsf{X}_2$ themselves. This proves that $[\boldsymbol{H} \div \boldsymbol{U}] = \mathcal{M}_1$. □

Now assume $n_1 < n_2$. We show that either $\boldsymbol{H}$ is contained in a transitive imprimitive group and $[\boldsymbol{H} \div \boldsymbol{U}] = \mathcal{M}_2$ or all the transitive groups containing $\boldsymbol{H}$ are primitive.

4.3.2. Proposition. *If* $\boldsymbol{H} < \boldsymbol{S}_{\mathsf{Z}}$ *for some non trivial equipartition* $\mathsf{Z}$ *of* $\Omega$, *then*

(1) $\mathsf{X}_1 \in \mathsf{Z}$ *and* $\mathsf{X}_2$ *is the union of all the other blocks of* $\mathsf{Z}$.
(2) $\boldsymbol{H} = (\operatorname{Sym}\mathsf{X}_1 \times T) \cap \boldsymbol{U}$ *where* $T$ *is the maximal imprimitive subgroup of* $\operatorname{Sym}\mathsf{X}_2$ *preserving* $\mathsf{Z} - \{\mathsf{X}_1\}$.
(3) $[\boldsymbol{H} \div \boldsymbol{U}] = \mathcal{M}_2$.

PROOF. We prove that $\boldsymbol{H}$ acts on $\mathsf{X}_2$ as a maximal imprimitive subgroup of $\operatorname{Sym}\mathsf{X}_2$. To this end, consider

$$\mathsf{Z}_i = \mathsf{Z} \cap \mathsf{X}_i := \left\{ \mathsf{B} \cap \mathsf{X}_i \;\middle|\; \mathsf{B} \in \mathsf{Z},\ \mathsf{B} \cap \mathsf{X}_i \neq \emptyset \right\}, \qquad i = 1, 2$$

which by transitivity of $\boldsymbol{H}$ on $\mathsf{X}_i$ is an equipartition of $\mathsf{X}_i$. If the elements of $\mathsf{Z}_2$ are singletons, then

$$n_1 \geqslant \left|\mathsf{Z}_1\right| \geqslant \left|\mathsf{Z}_2\right| = n_2$$

which contradicts $n_1 < n_2$. If $\mathsf{Z}_2 = \{\mathsf{X}_2\}$ then $n \geqslant 2n_2$ which again contradicts $n_1 < n_2$. Therefore $\mathsf{Z}_2$ is a non trivial congruence for $\boldsymbol{H}$ acting on $\mathsf{X}_2$. Call $T$ the maximal imprimitive subgroup of $\operatorname{Sym}\mathsf{X}_2$ preserving $\mathsf{Z}_2$ (there is only one by 2.4.3). Then we must have

$$\boldsymbol{H} \leqslant (\operatorname{Sym}\mathsf{X}_1 \times T) \cap \boldsymbol{U} < \boldsymbol{U}_{(\mathsf{X}_1, \mathsf{X}_2)} < \boldsymbol{U}$$

which implies $\boldsymbol{H} = (\operatorname{Sym}\mathsf{X}_1 \times T) \cap \boldsymbol{U}$.

To prove (1) and (2), we just need to show that $\mathsf{X}_1 \in \mathsf{Z}$. First of all, $n_1 > 1$ otherwise there is $\mathsf{B} \in \mathsf{Z}$ such that $\left|\mathsf{B} \cap \mathsf{X}_2\right| = m - 1$ whilst $\left|\mathsf{B}' \cap \mathsf{X}_2\right| = m$ for all the other $\mathsf{B}' \in \mathsf{Z}$. But this is impossible, $\mathsf{Z}_2$ being an equipartition. Said this, observe that $\mathsf{Z}_1$ is a trivial equipartition because it is preserved by $\boldsymbol{H}$ which acts on $\mathsf{X}_1$ as the full symmetric group on $\mathsf{X}_1$. If $\mathsf{Z}_1$ is made of singletons (so there are more than 1) call $\mathsf{B}_1$, $\mathsf{B}_2$ two blocks of $\mathsf{Z}$ which intersect $\mathsf{X}_1$. There is a double transposition $\boldsymbol{s}$ of $\boldsymbol{H}$ which swaps $\mathsf{B}_1 \cap \mathsf{X}_1$ and $\mathsf{B}_2 \cap \mathsf{X}_1$ and swaps two points of $\mathsf{B}_1 \cap \mathsf{X}_2$. Then $\boldsymbol{s}$ splits $\mathsf{B}_1$ against our assumptions.

It follows that $\mathsf{Z}_1 = \{\mathsf{X}_1\}$ and this yields $\mathsf{X}_1 \in \mathsf{Z}$, or otherwise, $\mathsf{Z}_2$ may not be an equipartition of $\mathsf{X}_2$.

To prove (3), recall 2.6.1. If $\boldsymbol{H}$ is contained in proper primitive groups, then $\boldsymbol{U} = A_7$ and $\mathsf{X}_1$ has order 1. But then $\mathsf{X}_1$ may not lie in $\mathsf{Z}$. A contradiction.

Therefore we conclude that all the transitive proper subgroups of $\boldsymbol{U}$ containing $\boldsymbol{H}$ are imprimitive. But if $\boldsymbol{H} \leqslant \boldsymbol{S}_{\mathsf{Z}'}$, then following this same proof with $\mathsf{Z}'$ in place of $\mathsf{Z}$, we may show that $\mathsf{X}_1 \in \mathsf{Z}'$. Also, as $\boldsymbol{H}$ acts on $\mathsf{X}_2$ exactly as $T$ does, $T$ must preserve $\mathsf{Y}' = \mathsf{Z}' - \{\mathsf{X}_1\}$. But $T$ preserves only one congruence (2.4.3), hence $\mathsf{Y}' = \mathsf{Y}$ which implies $\mathsf{Z}' = \mathsf{Z}$.

That proves that the only imprimitive subgroup of $\boldsymbol{U}$ containing $\boldsymbol{H}$ is $\boldsymbol{U}_{\mathsf{Z}}$. Since the unique intransitive subgroup containing $\boldsymbol{H}$ is $\boldsymbol{U}_{(\mathsf{X}_1,\mathsf{X}_2)}$, we have $[\boldsymbol{H} \div \boldsymbol{U}] = \mathcal{M}_2$. □

In fact, with only one exception the primitive maximal subgroups of $\boldsymbol{U}$ containing $\boldsymbol{H}$ are 2-primitive.

4.3.3. THEOREM. *If $\boldsymbol{H}$ is contained in a primitive maximal subgroup of $\boldsymbol{U}$, then either $\boldsymbol{U} = A_7$ and $r = 3$ as in* 2.6.1 *or all the maximal transitive subgroups of $\boldsymbol{U}$ containing $\boldsymbol{H}$ are* 2*-primitive, of the same order and not conjugate in $\boldsymbol{U}$. Furthermore, $\boldsymbol{H}$ is the stabilizer of a point in each of them.*

PROOF. Let us call $G_2, \ldots, G_r$ the transitive proper subgroups of $\boldsymbol{U}$ containing $\boldsymbol{H}$. By the previous proposition they are all primitive and so $\boldsymbol{H}$ may not contain transpositions or 3-cycles (2.5.3).

Assume first that $\boldsymbol{H}$ is odd, so that $\boldsymbol{U} = \boldsymbol{S}$. Our $\boldsymbol{H}$ must be a maximal subgroup of $\operatorname{Sym}\mathsf{X}_1 \times \operatorname{Sym}\mathsf{X}_2$. Recall 2.9.4 and accordingly distinguish three cases.

(1) $\boldsymbol{H} = (\operatorname{Sym}\mathsf{X}_1 \times \operatorname{Sym}\mathsf{X}_2)^{\mathrm{e}}$. But this $\boldsymbol{H}$ is even against our assumptions.
(2) $\boldsymbol{H} = T \times \operatorname{Sym}\mathsf{X}_2$ with $T \lessdot \operatorname{Sym}\mathsf{X}_1$. But this $\boldsymbol{H}$ contains transpositions because $n_2 > 1$.
(3) $\boldsymbol{H} = \operatorname{Sym}\mathsf{X}_1 \times T$ with $T \lessdot \operatorname{Sym}\mathsf{X}_2$. But this $\boldsymbol{H}$ contains transpositions on $\mathsf{X}_1$ unless $n_1 = 1$.

Therefore $\boldsymbol{H} = T$ with $T \lessdot \operatorname{Sym}\mathsf{X}_2$ and $\mathsf{X}_1 = \{\omega\}$ for some $\omega \in \Omega$. It follows that $\boldsymbol{H} = (G_i)_\omega$ for all $i$ from 2 to $r$. As a consequence, $G_i$ is 2-transitive and $|G_i| = n|\boldsymbol{H}|$. Furthermore, $T$ may not be imprimitive (otherwise $T$ and hence $\boldsymbol{H}$ would contain transpositions) so that $G_i$ must be 2-primitive.

We now apply Pálfy Lemma 2.8.1 to show that no two $G_i$ are conjugate. As stabilizer of $\omega$, $\boldsymbol{H}$ is not normal in $G_i$ for any $i > 1$. Nor $\boldsymbol{H}$ is normal in $\operatorname{Sym}\mathsf{X}_2$, because we are assuming $\boldsymbol{H}$ odd. Thus $\boldsymbol{H}$ is selfnormalizing and so are the $G_i$.

Assume now that $\boldsymbol{H}^{\mathsf{s}} < G$ where $\mathsf{s} \in \boldsymbol{S}$ and $G = G_2$ for example. Then

$$\boldsymbol{H}^{\mathsf{s}} = (G_\omega)^{\mathsf{s}} = (G \cap \boldsymbol{S}_\omega)^{\mathsf{s}} = G^s \cap \boldsymbol{S}_{\omega\mathsf{s}} \leqslant \boldsymbol{S}_{\omega\mathsf{s}}.$$

This shows that indeed $\boldsymbol{H}^{\mathsf{s}} = G_{\omega\mathsf{s}}$ as they have the same order. But $\omega\mathsf{s} = \omega\mathsf{g}$ for some $\mathsf{g} \in G$, thus

$$\boldsymbol{H}^{\mathsf{s}} = G_{\omega\mathsf{s}} = G_{\omega\mathsf{g}} = \boldsymbol{H}^{\mathsf{g}}.$$

Hence $\operatorname{hm}(\boldsymbol{H}, G) = |G : \boldsymbol{H}|$ and $\operatorname{hm}(G, \boldsymbol{H}) = 1$ as we claimed.

If $\boldsymbol{H}$ is even, then the proof is somewhat similar. As before, we observe that $\boldsymbol{H}$ is a maximal subgroup of $\boldsymbol{U}_{(\mathsf{X}_1,\mathsf{X}_2)} = (\operatorname{Sym}\mathsf{X}_1 \times \operatorname{Sym}\mathsf{X}_2)^{\mathrm{e}}$ and we distinguish three cases.

(1) $\boldsymbol{H} = \operatorname{Alt} \mathsf{X}_1 \times \operatorname{Alt} \mathsf{X}_2$. But this $\boldsymbol{H}$ contains 3-cycles on $\mathsf{X}_2$.

(2) $\boldsymbol{H} = (T \times \operatorname{Sym} \mathsf{X}_2)^{\mathrm{e}}$ with $T \lessdot \operatorname{Sym} \mathsf{X}_1$. Again, $\boldsymbol{H}$ contains 3-cycles in this case which we have seen before may not happen.

(3) $\boldsymbol{H} = (\operatorname{Sym} \mathsf{X}_1 \times T)^{\mathrm{e}}$ with $T \lessdot \operatorname{Sym} \mathsf{X}_2$, $T \neq \operatorname{Alt} \mathsf{X}_2$. Since $\boldsymbol{H}$ has only two orbits, $T$ must be transitive. In fact, we may assume that $T$ is primitive otherwise we apply 2.6.1 to find that $\boldsymbol{U} = A_7$ and $r = 3$.

Since $\boldsymbol{H}$ may not contain 3-cycles on $\mathsf{X}_1$, we must have $n_1 \leqslant 2$. Suppose that $\mathsf{X}_1 = \{\mathsf{a}, \mathsf{b}\}$, $\mathsf{c} \in \mathsf{X}_2$ and $G = G_2$ for example. We observe that $G_{\mathsf{a}}$ contains $T^{\mathrm{e}}$ which is transitive on $\mathsf{X}_2$ (2.5.2). Therefore the orbits of $G_{\mathsf{a}}$ are $\{\mathsf{a}\}$, $\{\mathsf{b}\}$, $\mathsf{X}_2$. Beside, $\{\mathsf{a}, \mathsf{b}\}$ is an orbit of $G_{\mathsf{c}}$ because $T_{\mathsf{c}}$ is not contained in $\boldsymbol{A}$ and hence the transposition $(a\ b)$ multiplied by an odd permutation of $T_{\mathsf{c}}$ lies in $\boldsymbol{H}$. But $G_{\mathsf{c}}$, which is conjugate to $G_{\mathsf{a}}$ in $G$, must have exactly three orbits of order 1, 1, $n_2$, so, $n_2 = 2$. This contradicts $n_1 < n_2$ and so $n_1 = 1$.

Therefore we may now assume that $\boldsymbol{H}$ is a maximal primitive subgroup of $\operatorname{Alt} \mathsf{X}_2$ and that $n_1 = 1$.

As, $\boldsymbol{H}$ is primitive on $\mathsf{X}_2$ and it is the stabilizer of a point for each of the subgroups $G_2, \ldots, G_r$, they are all 2-primitive and of the same order.

Now $\boldsymbol{H}$ has index at most 2 in its normalizer in $\boldsymbol{S}$. As in the odd case, one shows that $\operatorname{hm}_{\boldsymbol{S}}(\boldsymbol{H}, G) = \left|G : \boldsymbol{H}\right|$ for any maximal transitive subgroup $G$ of $\boldsymbol{A}$ containing $\boldsymbol{H}$. Hence $\operatorname{hm}_{\boldsymbol{S}}(G, \boldsymbol{H}) \leqslant 2$ and $\operatorname{hm}_{\boldsymbol{A}}(G, \boldsymbol{H}) = 1$ which concludes the proof. □

Call as in the proof of the theorem above

$$G_2, \ldots, G_r,$$

the transitive maximal subgroups of $\boldsymbol{U}$ containing $\boldsymbol{H}$, and suppose $r > 3$. Then they are all 2-primitive and $\boldsymbol{H}$ is the stabilizer of a point in each of them. A 2-transitive group is either affine or almost simple with details given in Table 3.A. We claim that none of the $G_i$ is affine.

If not so, suppose that$G = G_j$ is an affine group. By maximality, $G$ is $\mathrm{AGL}(d, p) \cap \boldsymbol{U}$ and $\Omega$ is a vector space of dimension $d$ on a field of prime order $p$. Since $G$ must also be 2-primitive (and maximal in $\boldsymbol{U}$), we see from Table 3.B, p. 29 that $d > 2$ and $p = 2$. In particular $\mathrm{AGL}(d, p)$ is even and $\boldsymbol{U} = \boldsymbol{A}$. All the maximal affine subgroups of $\boldsymbol{A}$ are conjugate in $\boldsymbol{S}$ and they split in precisely two $\boldsymbol{A}$-conjugacy classes. Because we assume $r > 3$ and because the $G_i$ are not conjugate in $\boldsymbol{A}$, there must be a non affine 2-transitive $G_k$, say $L = G_k$. In particular $L$ is an almost simple group and $\operatorname{Soc} L$ has a subgroup of index $2^d$. Here is where we first use the Guralnick Theorem.

4.3.4. Theorem ([**Gur83**]). *Let $G$ be a nonabelian simple group with $H < G$ and $\left|G : H\right| = p^a$, $p$ prime. One of the following holds.*

(1) $G = \mathrm{A}_n$ *and* $H \cong \mathrm{A}_{n-1}$ *with* $n = p^a$.
(2) $G = \mathrm{PSL}(n, q)$ *and $H$ is the stabilizer of a (projective) point or hyperplane (Note that* $\frac{q^n-1}{q-1} = p^a$, *and that for this to happen, $n$ must be prime).*
(3) $G = \mathrm{PSL}(2, 11)$ *and* $H \cong \mathrm{A}_5$ *of index* 11.
(4) $G = \mathrm{M}_{23}$ *and* $H \cong \mathrm{M}_{22}$ *of index* 23.
(5) $G = \mathrm{M}_{11}$ *and* $H \cong \mathrm{M}_{10}$ *of index* 11.
(6) $G = \mathrm{PSp}(4, 3)$ *and $H$ is the parabolic subgroup of index* 27.

Then $L$ is a projective group $\mathrm{P\Gamma L}(n, q)^{\mathrm{e}}$ where $\frac{q^n-1}{q-1} = 2^d$. By 3.9.1, $n = 2$ and $2^d = 1 + q$. Thus $q$ is odd and so is $\mathrm{P\Gamma L}(2, q)$. The order of $\mathrm{P\Gamma L}(2, q)^{\mathrm{e}}$ is

$$\frac{1}{2}(q+1)q(q-1)f = \frac{1}{2}2^d(2^d-1)(2^d-2)f \tag{4.3.A}$$

where $q = p_1^f$ and $p_1$ is prime. However, the order of $\mathrm{AGL}(d, 2)$ is

$$2^d(2^d-1)(2^d-2)\cdots(2^d-2^{d-1}). \tag{4.3.B}$$

By 4.3.3 $G$ and $L$ must have the same order. One checks that (4.3.A) and (4.3.B) are not equal for $d = 3$ or 5 (for $d = 4$, $2^d - 1$ is not a prime power). Therefore we may assume that $d > 5$ and $q > 31$ and then

$$\frac{1}{2}f = (p_1^f - 3)(p_1^f - 7) \cdot \ldots > p_1^f - 3 \geqslant 3^f - 3.$$

This forces $f = 1$ and hence $\frac{1}{2}$ to be product of integers. The contradiction arose by assuming that $G_j$ is an affine group. Therefore provided $r > 3$, none of the $G_i$ is affine.

4.3.5. Proposition. *If $[\boldsymbol{H} \div \boldsymbol{U}] = \mathcal{M}_r$ with $r > 3$, then $\boldsymbol{H}$ has degree bigger than* 276 *and all the transitive maximal subgroups of $\boldsymbol{U}$ containing $\boldsymbol{H}$ are almost simple.*

Proof. The transitive maximal subgroups of $\boldsymbol{U}$ containing $\boldsymbol{H}$ are almost simple because they cannot be imprimitive by 4.3.2 and they cannot be primitive of affine type for the discussion above. Therefore they are 2-transitive and almost simple. There is precisely one conjugacy class of each of them in $\boldsymbol{S}$, hence at most two $\boldsymbol{A}$-conjugacy classes when they are contained in $\boldsymbol{A}$. Unless they have non isomorphic socles, there are at most 2 of them containing $\boldsymbol{H}$ because they are not conjugate in $\boldsymbol{U}$. As that yields $r \leqslant 3$, we only remain with the degrees such that there are

2-transitive almost simple groups having non isomorphic socles. Up to degree 276 those degrees are:

- **11:** with socles $\mathrm{PSL}_2(11)$, $\mathrm{M}_{11}$. They are both selfnormalizing, both even of order 660 and 7,920 respectively.
- **12:** with socles $\mathrm{PSL}_2(11)$, $\mathrm{M}_{11}$ and $\mathrm{M}_{12}$. $\mathrm{PSL}_2(11)$ has index 2 in its normalizer which is odd and of order 1,320, while the other two are both selfnormalizing, both even of order 7,920 and 95,040 respectively.
- **24:** with socles $\mathrm{PSL}_2(23)$, $\mathrm{M}_{24}$. $\mathrm{PSL}_2(23)$ has index 2 in its normalizer which is odd and of order 12,144, while $\mathrm{M}_{24}$ is selfnormalizing, even and of much bigger order.
- **28:** with socles $\mathrm{Sp}_6(2)$, $\mathrm{PSU}_3(3)$ and $\mathrm{PSL}_2(8)$, The first one is selfnormalizing, $\mathrm{PSU}_3(3)$ has index 2 in its normalizer, while $\mathrm{PSL}_2(8)$ has index 3. The normalizers of the first and the third are even of order 1,512 and 1,451,520 respectively. The order of the second is 6,048 and the order of its full normalizer is 12,096.

As the orders of the normalizers in $\boldsymbol{U}$ are different, None of them has a common point stabilizer and so they may not contain the same $\boldsymbol{H}$. □

4.3.6. THEOREM. *If* $[\boldsymbol{H} \div \boldsymbol{U}] = \mathcal{M}_r$, *then* $r \leqslant 3$.

PROOF. The argument is similar to the one of proposition above but we have to deal with the infinite parametrised families of 2-transitive simple groups. Note however that because the transitive maximal subgroups of $\boldsymbol{U}$ containing $\boldsymbol{H}$ must be 2-primitive, of all the projective linear groups one only needs to consider the projective linear groups on projective lines, that is, the ones having socle $\mathrm{PSL}_2(q)$ for some $q$. We show for example that above a same $\boldsymbol{H}$ there cannot be almost simple groups $G$, $L$ having socles isomorphic to $\mathrm{PSL}_2(q_1)$ and $\mathrm{Sz}(q_2)$. If so,

$$q_1 = q_2^2 = 2^{2f}$$

for some odd integer $f$, thus $G$ is even of order $(q_2^2+1)q_2^2(q_2^2-1)2f$. On the other hand $L$ is even of order $(q_2^2+1)q_2^2(q_2-1)f$. But this contradicts the fact that $G$ and $L$ should have the same order. As a further example, suppose $G$ has socle isomorphic to $\mathrm{PSL}_2(q)$ and $L$ is almost simple with socle $\mathrm{Sp}_{2d}(2)$ and degree $2^{d-1}(2^d \pm 1)$. Now $L = \mathrm{Sp}_{2d}(2)$, therefore $L$ is even and so $G$ is even too. Comparing the degrees

$$1 + q = 2^{d-1}(2^d \pm 1)$$

one finds that $q$ is odd hence $G$ has index 2 in its normalizer. It follows that $G$ has order $\frac{1}{2}(q+1)q(q-1)f$ where $q = p^f$ and $p$ is prime, while $L$ has order $2^{(d^2)}k$ where $k$ is an odd integer. The highest power of 2 dividing the order of $G$ is $2^{d-1}a$ where $a$ is the highest power of 2 dividing $f$. Clearly $a < f < 2d$, so $2^{d-1}a < 2^d d < 2^{(d^2)}$. It follows that $G$ and $L$ have not the same order.

In entirely similar ways one shows that any two transitive maximal subgroups of $\boldsymbol{U}$ containing $\boldsymbol{H}$ must have isomorphic socle. Hence there are at most two of them. □

CHAPTER 5

# Primitive second maximal subgroups

We intend to supplement the qualitative results of [**Pra90**] on the inclusion problem for finite primitive groups with a quantitative analysis restricted to the case where the smallest subgroup is a second maximal subgroup of the parent symmetric or alternating group. This amounts to counting the inclusions that in [**Pra90**] are given up to permutational isomorphisms, under the assumption that the smallest group is fixed and second maximal. In other words, adopting the notation $\boldsymbol{S}$, $\boldsymbol{A}$, $\boldsymbol{H}$, $\boldsymbol{U}$ of 2.1, we investigate the maximal subgroups of $\boldsymbol{U}$ containing $\boldsymbol{H}$ in the assumption that $\boldsymbol{H}$ is a given primitive group.

$\boldsymbol{S} \equiv \operatorname{Sym}\Omega \cong \mathrm{S}_n, \quad n \geqslant 5$

$\boldsymbol{A} \equiv \operatorname{Alt}\Omega \cong \mathrm{A}_n$

$\boldsymbol{U} = \boldsymbol{S}$ or $\boldsymbol{A}$

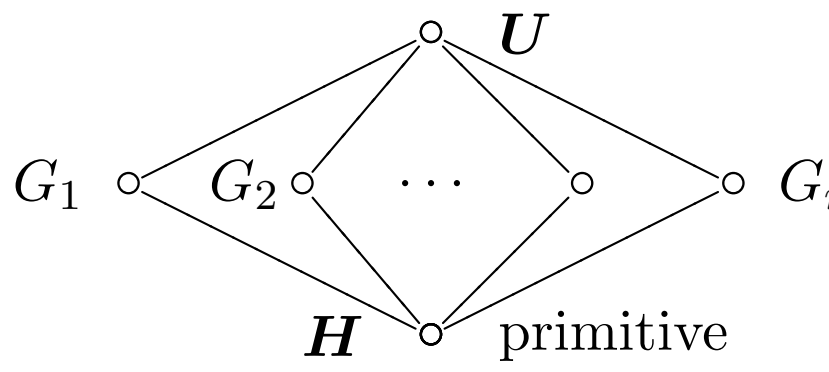

Primitive groups are subdivided by the O'Nan-Scott Theorem into types according to the permutational isomorphism type of their socles. Since these socles are powers of finite simple groups, the O'Nan-Scott Theorem is used to reduce questions concerning primitive groups to questions involving either much smaller groups or, ultimately, finite simple groups. The advantage in dealing with questions of the latter kind is that there is a classification of the finite simple groups. Then one may obtain the sought answers by an exhaustive analysis of the involved simple groups.

Accordingly, we decided to split the primitive case in two parts leaving the exhaustive analysis to next and last chapter of this thesis. In this chapter we investigate all cases where at least one of the groups to discuss is non almost simple and we prove the following theorem.

**Theorem B.** *If $\boldsymbol{H}$ or at least one of the $G_i$ is non almost simple, then $r \leqslant 2$, unless it is one of the three examples of the kind discovered by Feit and listed in* [**Pál88**, Table II] *or here in* 5.2.1.

We now proceed with a brief discussion of the O'Nan-Scott Theorem followed by some introductory material. In reading these first sections, one should not worry yet about terms like "blow-ups" or "compound diagonal", as their precise meaning will be given in due course.

## 5.1. On the O'Nan-Scott Theorem

At the *Santa Cruz Conference on Finite Groups* in 1979, M. O'Nan and L. L. Scott discussed a classification scheme for primitive groups. It became "a theorem of O'Nan and Scott" in [**Cam81**] and finally "the O'Nan-Scott Theorem" in [**LPS88**]. Nowadays, any satisfactory treatment of the O'Nan-Scott Theorem comprises two main parts:

(1) a list of examples or types of primitive permutation groups;
(2) an argument proving that each primitive group is permutationally isomorphic, or has the full normalizer of its socle permutationally isomorphic, to a unique group of the list.

This subdivides primitive permutation groups in types according to the permutational isomorphism type of their socle. However, there is not yet full agreement on the hierarchy and names of the different types. Some authors prefer a hierarchy with names reflecting the structure of the argument given in (2); others may prefer the hierarchy induced by inclusions among subgroups; others may decide for a flat hierarchy with symbolic names reminding the structure of the groups.

Here we follow [**Pra90**] where primitive groups are subdivided in a flat hierarchy of eight types:

- affine, or Holomorph of Abelian, or HA ;
- Almost Simple, or AS;
- Holomorph of non-abelian Simple, or HS;
- Simple Diagonal, or SD;
- Twisted Wreath, or TW

and the blow-ups of types AS, HS, SD, namely

- Product Action, or PA;
- Holomorph of Composite, or HC;
- Compound Diagonal, or CD.

The primitive groups of affine type are the ones having a non trivial abelian normal subgroup. This is the only minimal normal subgroup in such a primitive group. Primitive groups of other types may have one or two minimal normal subgroups. If there are two minimal normal subgroups, then these are the centralizers of each other (in the full symmetric group), they are regular and isomorphic. The primitive groups are of type HS or HC depending on whether these two regular subgroups are simple or not. Primitive groups of the other five types have only one minimal normal subgroup. This is either simple, when the primitive groups are of type AS, or a power $T^l$ of a non abelian simple group $T$ ($l > 1$). In the latter case, the stabilizer of a point is either $\{1\}$ or isomorphic to $T$ or isomorphic to $R^l$ where $\{1\} < R < T$, or else, $l = ab$ with $a, b > 1$ and the stabilizer is isomorphic to $T^b$. Then the primitive groups are of type TW, SD, PA or CD respectively. Surprisingly, types HA, HS, HC and AS are recognisable by the abstract isomorphism type of the socle alone. The permutational isomorphism type of the socle is required to distinguish among types PA, SD, CD, TW. This is summarized in Figure 5.A.

It was shown in [**Pra90**] that primitive groups of certain types are always contained in primitive groups of certain other types. For example, each primitive group is contained in one of type AS (such as the parent symmetric group). The diagram of all such inevitable inclusions is given in Figure 5.B. In particular, it shows that all primitive groups of type TW or HC are low.

## 5.2. Outline

The O'Nan-Scott Theorem provides qualitative information: an insight into the structure of the primitive groups. In [**Pra90**] the O'Nan-Scott Theorem is used to investigate the general qualitative "inclusion problem": is a primitive group subgroup of another? The answer is given up to permutational isomorphism, that is, the answer is "yes" whenever the first is permutationally isomorphic to a subgroup of the other. However, we are dealing with a slightly different task: *counting* the subgroups containing a given primitive group, in fact, a second maximal primitive group. Quoting [**Kov89b**]:

> counting the number of mathematical objects of a certain kind is often undertaken as a test problem ("if you can't count them, you don't really know them").

To solve this kind of problems, one needs a language which support the distinction among the mathematical objects of the same kind: in our case, one should have different names for different permutationally isomorphic permutation groups.

This is why we improvise languages supporting distinction among permutationally isomorphic permutation groups for types PA, HS, SD, HC, CD. Our approach, endowing "standard" primitive groups with a permutational isomorphism whose images are the required groups, is naïve and unsatisfactory especially because it introduces a cumbersome notation and it fails in providing unique names: with our notations, a same group may be addressed in too many ways. However, this will eventually produce the sought result.

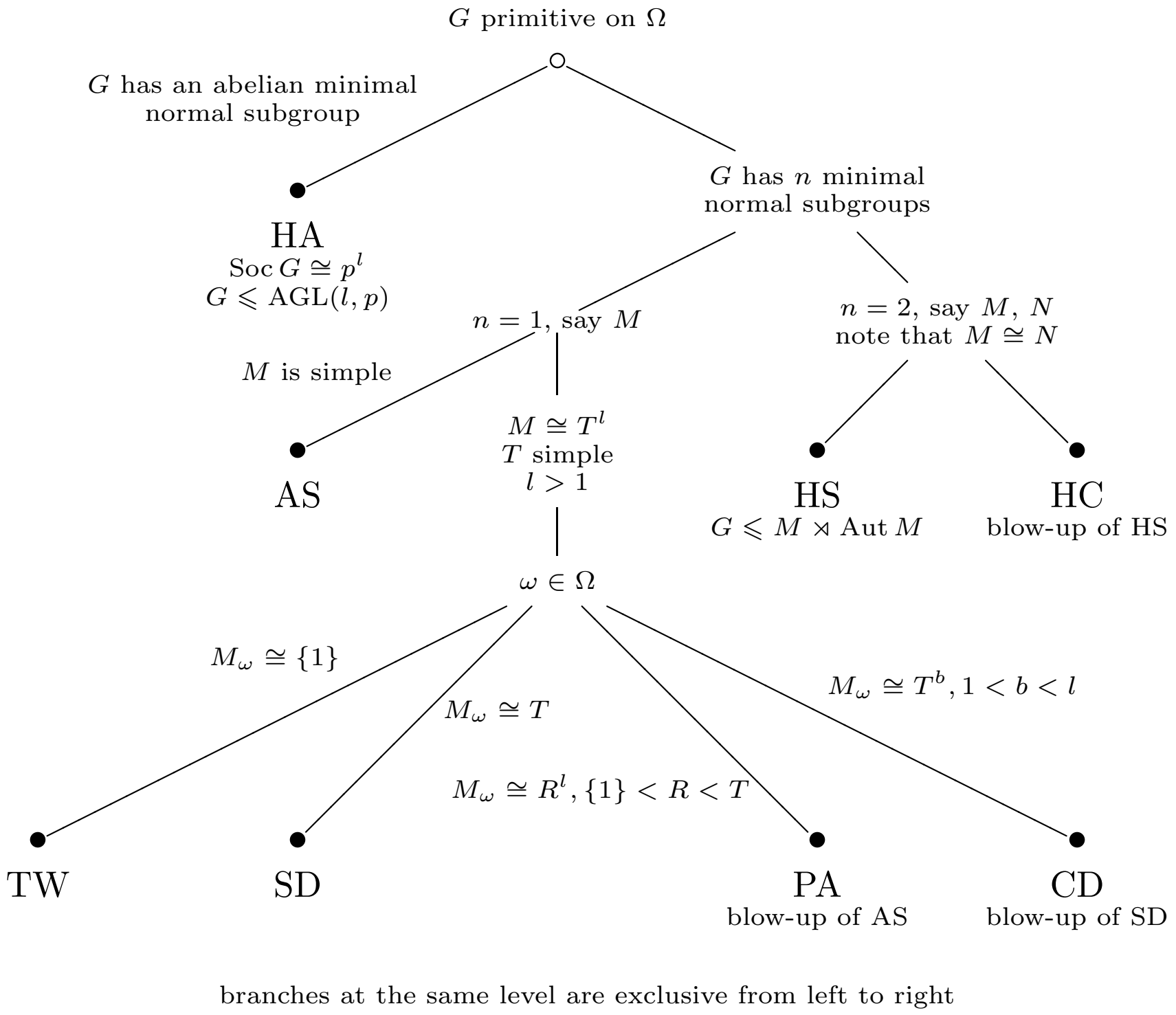

FIGURE 5.A. The eight types of primitive groups

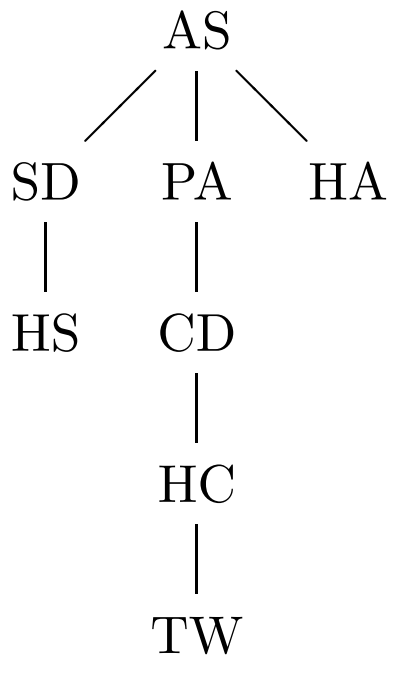

FIGURE 5.B. Inevitable inclusions among primitive groups

Whenever needed we translate [**Pra90**] results into our languages. In doing so we tend to quote (and change) as little as possible and that little is often over-weakened in the translation. The reason behind this "weakening" relies on the fact that we only need these results for the very special case in which the smallest subgroup is second maximal. It should be possible to follow our arguments without having to read [**Pra90**] first. Of course, our arguments do rely on [**Pra90**] but then, they also rely on the CFSG. Yet, this does not mean that a reader should master all of it before reading our work! There is a little price to pay: sometime we have to prove again facts that would clearly emerge from a careful reading of [**Pra90**]. For example, when it is there stated that $(G, G_1)$ is a blow-up of an inclusion $(L, L_1)$, we report that $G$ is a blow-up of $L$, $G_1$ is a blow-up of $L_1$ and $G \leqslant G_1$; deliberately omitting the details that make of this inclusion among blow-ups, the blow-up of an inclusion. Instead we provide the missing details with non strictly necessary proofs. Then, some minor refinement is required to show that a primitive second maximal subgroup of type PA, HS, SD or CD is never contained in more than two maximal subgroups.

About the second maximal subgroups of type HA (affine case), we reduce to the case where the second maximal subgroup is even and of prime degree. Then we appeal to the classification of the second maximal subgroups of the alternating groups of prime degree [**Pál88**]. We only quote the relevant part of this classification.

5.2.1. PROPOSITION. *Let $p$ be a prime and let $\boldsymbol{H}$ be a transitive subgroup of $\mathrm{A}_p$ such that $[\boldsymbol{H} \div \mathrm{A}_p] \cong \mathcal{M}_r$ for some $r > 2$. Then $\boldsymbol{H}$ is the normalizer of a cycle of length $p$ in $r-1$ almost simple groups having socle isomorphic to $S$, where $p$, $r$, $S$ are as follows.*

| $p$ | $\boldsymbol{H}$ | $S$ | $r$ |
|---|---|---|---|
| 13 | 13.3 | $\mathrm{PSL}_3(3)$ | *5* |
| 31 | 31.5 | $\mathrm{PSL}_5(2)$ | *7* |
| 31 | 31.3 | $\mathrm{PSL}_3(5)$ | *11* |

Finally, we deal with primitive almost simple groups which are contained in some non almost simple group. Again, we have the support of a qualitative result: a classification of such inclusions up to permutational isomorphisms [**LPS87**]. Then we only need to show that such inclusions involve low subgroups, unless the degree is 8 where the situation is under control.

We remain with the case where all the subgroups containing the second maximal primitive group are almost simple; this is left to next Chapter.

## 5.3. Preliminaries

The *socle* of a group $G$ is the subgroup generated by its minimal normal subgroups and it is commonly denoted by $\mathrm{Soc}\, G$. The socle of a primitive group $G$ is always isomorphic to a power of some simple group. If $\mathrm{Soc}\, G \cong T^l$ with $T$ simple but non abelian, then $\mathrm{Soc}\, G$ has exactly $l$ minimal normal subgroups: these are the *component*[1] of $G$ and we denote them by

$$(\sqrt[l]{G})_0\,, \quad \ldots, \quad (\sqrt[l]{G})_{l-1}\ .$$

Also, set $\sqrt[l]{G} = (\sqrt[l]{G})_0$ and note that $\sqrt[l]{G} \cong \sqrt[l]{\operatorname{Soc} G}$ so that one may write $(\sqrt[l]{\operatorname{Soc} G})^l \cong \operatorname{Soc} G$. In fact, $\operatorname{Soc} G$ is the internal direct product of the $(\sqrt[l]{G})_i$, and for $i \in l$ we denote by $\boldsymbol{\rho}_i^G\colon \operatorname{Soc} G \longrightarrow (\sqrt[l]{G})_i$ the canonical projection which restricts to the identity on $(\sqrt[l]{G})_i$ and has kernel $\langle\, (\sqrt[l]{G})_j \mid i \neq j \in l\,\rangle$.

Sometime $\operatorname{Soc} G$ is not big enough and then we use the *shoe* of $G$, that is, the subgroup of $G$ generated by the intersection of the normalizers in $G$ of the minimal normal subgroups of $\operatorname{Soc} G$. Of course, $\operatorname{Soc} G \leqslant \operatorname{Shoe} G$. For example, for $A$ almost simple,

$$\operatorname{Soc}(A^2 \rtimes \mathrm{S}_2) = (\operatorname{Soc} A)^2, \quad \operatorname{Shoe}(A^2 \rtimes \mathrm{S}_2) = A^2\ .$$

The next two lemmas are due to Praeger. They are of fundamental importance for analysing inclusions among primitive groups.

5.3.1. Lemma ([**Pra90**, 4.5]). *Suppose that $H < G$ with $H$ primitive and that neither $H$ nor $G$ is of affine type. Then $H \cap \operatorname{Soc} G \neq \{1\}$ and one of the following holds:*

(1) $\operatorname{Soc} H = \operatorname{Soc} G$*;*
(2) $\operatorname{Soc} H < \operatorname{Soc} G$*, $H$ is of type TW and $G$ is of type HC;*
(3) *$H$ is a core-free subgroup of $G$ (that is, $\bigcap_{g\in G} H^g = \{1\}$) and $G$ is not of type TW.*

When $H$ is maximal in $G$, the last point of lemma above is refined as follows.

5.3.2. Lemma ([**Pra90**, 4.6]). *Suppose that $H \lessdot G$ with $H$ primitive and $H$ core-free in $G$. Suppose also that $G$ has $m$ components while $H$ has $l$ components. If neither* $\operatorname{Soc} H$*, nor* $\operatorname{Soc} G$ *is abelian, then one of the following is true:*

(1) *$H \cap \operatorname{Soc} G = \langle\, (H \cap \operatorname{Soc} G)\boldsymbol{\rho}_i^G \mid i \in m\,\rangle$, $H$ and $G$ are transitive by conjugation on $\{\, (\sqrt[m]{G})_i \mid i \in m\,\}$. Moreover, $l = am$ for some $a \geqslant 1$, $\operatorname{Soc} H < \operatorname{Soc} G$ and there is a partition $\{\mathsf{E}_0, \ldots, \mathsf{E}_{m-1}\}$ of order $m$ of $l$ such that for all $i \in m$*

$$\big(\operatorname{Soc} H\big)\boldsymbol{\rho}_i^G = \langle\, (\sqrt[l]{H})_j \mid j \in \mathsf{E}_i\,\rangle$$

*In this case $G$ is not of type HS, HC or TW.*

[1] The components of a group are its subnormal quasisimple subgroups.

(2) *$H \cap \operatorname{Soc} G = \operatorname{Soc} H$, $\sqrt[m]{G} \cong \sqrt[l]{H}$ and $m = al$ for some $a \geqslant 2$. Further, there is a partition $\{\mathsf{E}_0, \dots, \mathsf{E}_{l-1}\}$ of order $l$ of $m$ such that for all $j \in l$*

$$(\sqrt[l]{H})_j \boldsymbol{\rho}_i^G = \begin{cases} (\sqrt[m]{G})_i & \text{if } i \in \mathsf{E}_j \\ \{1\} & \text{otherwise.} \end{cases}$$

*In this case $G$ is not of type TW or AS.*

## 5.4. The maximal subgroups table argument

Let $\boldsymbol{H}$ be a second maximal subgroup of $\boldsymbol{S}$ or $\boldsymbol{A}$ as described in 2.1 and suppose that $\boldsymbol{H} < G < \boldsymbol{U} = \boldsymbol{AH}$ for some primitive almost simple subgroup $G$ with socle $L$. Suppose also that $L \not< \boldsymbol{H}$.

We claim that the number of conjugate in $\boldsymbol{S}$ of $G$ containing $\boldsymbol{H}$, that is $\operatorname{hm}(G, \boldsymbol{H})$, is bounded by $k\left|\mathrm{N}(\operatorname{Soc} \boldsymbol{H}) : \boldsymbol{H}\right|$ where $k$ is the number of conjugacy classes in $G$ of maximal subgroups of $G$ isomorphic to $\boldsymbol{H}$ and the normalizer is taken in $\boldsymbol{S}$. In fact, it is bounded by $\frac{1}{2}k\left|\mathrm{N}(\boldsymbol{H}) : \boldsymbol{H}\right|$ when $G$ is not selfnormalizing in $\boldsymbol{S}$.

This is called *maximal subgroups table argument* because the determination of $k$ is straightforward when maximal subgroups tables like the ones in the ATLAS are available for $G$.

The proof of the claim goes as follows. From $L \not< \boldsymbol{H}$ we get $\boldsymbol{H} \ntrianglelefteq G$; thus $\left|G : \boldsymbol{H}\right|$ is the number of conjugate subgroups of $\boldsymbol{H}$ in $G$. Suppose now that $\boldsymbol{H}^s < G$ for some $s \in \boldsymbol{S}$. Then $\boldsymbol{H} < G^{s^{-1}}$ which together with the assumption that $\boldsymbol{H}$ is high implies $\boldsymbol{H}^s \lessdot G$. Therefore both $\boldsymbol{H}$ and $\boldsymbol{H}^s$ must be in the column of the maximal subgroups of $G$. This provides an immediate bound for $\operatorname{hm}(G, \boldsymbol{H})$:

$$k\left|G : \boldsymbol{H}\right|,$$

where $k$ is the number of conjugacy classes of maximal subgroups of $G$ isomorphic to $\boldsymbol{H}$. By Pálfy Lemma 2.8.1, $\operatorname{hm}(G, \boldsymbol{H}) \leqslant k\left|\mathrm{N}(\boldsymbol{H}) : \boldsymbol{H}\right|$ and if $G$ is not self-normalizing in $\boldsymbol{S}$, then $\operatorname{hm}(G, \boldsymbol{H}) \leqslant \frac{1}{2}k\left|\mathrm{N}(\boldsymbol{H}) : \boldsymbol{H}\right|$. The claim follows because $\mathrm{N}(\boldsymbol{H}) \leqslant \mathrm{N}(\operatorname{Soc} \boldsymbol{H})$.

It goes without saying that if $\boldsymbol{H}$ is primitive of type almost simple and its socle is listed in the ATLAS, then $\mathrm{N}(\operatorname{Soc} \boldsymbol{H})$ may be found as explained at the end of 3.3.

## 5.5. Wreath products

We are only interested in finite wreath products $A \operatorname{wr} B$ where $B$ is a permutation group (not necessarily transitive) of degree at least 2, say for example that $B \leqslant \mathrm{S}_l$ for some $l > 1$. The wreath product $A \operatorname{wr} B$ is

$$W = B \ltimes A^l$$

where for $\vec{a} \in A^l$, $\boldsymbol{\beta} \in B$ we set

$$\boldsymbol{\beta}^{-1}\vec{a}\boldsymbol{\beta} = (a_{0\boldsymbol{\beta}^{-1}}, \dots, a_{(l-1)\boldsymbol{\beta}^{-1}}).$$

Accordingly, any element $w$ of $W$ is uniquely written as

$$w = \boldsymbol{\beta}\vec{a}, \qquad \boldsymbol{\beta} \in B, \quad \vec{a} \in A^l.$$

The *top subgroup* of $W$ is $B$ and the *top projection* is

$$\boldsymbol{\rho} : W \longrightarrow B; \qquad \boldsymbol{\beta}\vec{a} \mapsto \boldsymbol{\beta}. \tag{5.5.A}$$

This is a group homomorphism whose kernel is $A^l$, the *base subgroup* of $W$. In fact, $\boldsymbol{\rho}$ is an action of $W$ on $l$ and for each $i \in l$ we denote by $W_i$ the stabilizer in $W$ of $i$ according to this action. One sees that

$$W_i \cong (B_i \ltimes A^{l-1}) \times A.$$

The right canonical projection onto $A$ according to this factorization is denoted by

$$\boldsymbol{\rho}_i : W_i \longrightarrow A; \qquad \boldsymbol{\beta}\vec{a} \mapsto a_i. \tag{5.5.B}$$

In fact, $\boldsymbol{\rho}_i$ may be regarded as a function from $W$ which restricts to a group homomorphism from $W_i$. We also call

$$\boldsymbol{\kappa}_i : A \longrightarrow W, \qquad i \in l \tag{5.5.C}$$

the $l$ canonical inclusions with images in $A^l$. Then $\boldsymbol{\kappa}_i\boldsymbol{\rho}_i = \mathrm{id}_A$ and, for $i \neq j \in l$, $\boldsymbol{\kappa}_i\boldsymbol{\rho}_j$ is the trivial endomorphism of $A$ with kernel $A$. The image of $\boldsymbol{\kappa}_i$ is called $i$th (coordinate) factor of the base subgroup. $W$ intertwines the factors of its base subgroup, hence the name wreath product:

$$\boldsymbol{\beta}^{-1}(A\boldsymbol{\kappa}_i)\boldsymbol{\beta} = A\boldsymbol{\kappa}_{i\boldsymbol{\beta}}, \qquad i \in l, \ \boldsymbol{\beta} \in \mathrm{S}_l. \tag{5.5.D}$$

A *large subgroup* of $W = A \,\mathrm{wr}\, B$ is a subgroup $G \leqslant W$ such that $G\boldsymbol{\rho}$ is transitive on $l$ and

$$(G \cap W_i)\boldsymbol{\rho}_i = A \qquad \forall i \in l.$$

The set of large subgroups of $W$ is closed under conjugation in $W$. But if $C < A$ and $G$ is a large subgroup of $C \,\mathrm{wr}\, \mathrm{S}_l$ regarded as a subgroup of $W$, then there may be conjugates $G^w$ of $G$ in $W$ which are not large subgroups of any related wreath product. For example, say $C$ a subgroup of order 2 of $A = \mathrm{S}_3$, $w \in A \,\mathrm{wr}\, \mathrm{S}_2$ such that $w\boldsymbol{\rho} = 1$, $w\boldsymbol{\rho}_0 = 1$, $w\boldsymbol{\rho}_1 = (0\ 1\ 2)$. Then $G = w^{-1}(C \,\mathrm{wr}\, \mathrm{S}_2)w$ is not a large subgroup of anything because $(G \cap W_0)\boldsymbol{\rho}_0 = C$ but $(G \cap W_1)\boldsymbol{\rho}_1 \neq C$. As a consequence some care is needed in handling objects defined in terms of large subgroups (see the blow-ups in 5.9 for example).

If $A$ acts on a set $\mathsf{X}$, the wreath product $W = A \operatorname{wr} B$ may be endowed with the *imprimitive action* on the coproduct

$$\mathsf{X} \times l = \coprod_{i=0}^{l-1} \mathsf{X}$$

or with the *product action* on the product

$$\mathsf{X}^l = \prod_{i=0}^{l-1} \mathsf{X}.$$

We now describe these two actions.

**5.5.1. Imprimitive action.** The imprimitive action of $W$ is defined by

$$\lfloor \mathsf{X} \times l : \star \rfloor : W \longrightarrow \operatorname{Sym}(\mathsf{X} \times l), \qquad w \mapsto \lfloor \mathsf{X} \times l : w \rfloor;$$

(regard the symbol $\lfloor\ \rfloor$ as a broken $\coprod$) where for $\boldsymbol{\beta} \in B$, $\vec{a} \in A^l$

$$(\mathsf{x}, i)\lfloor \mathsf{X} \times l : \boldsymbol{\beta} \rfloor := (\mathsf{x}, i\boldsymbol{\beta}), \tag{5.5.E}$$

$$(\mathsf{x}, i)\lfloor \mathsf{X} \times l : \vec{a} \rfloor := (\mathsf{x} \cdot a_i, i). \tag{5.5.F}$$

However, if $w \in W$, we may write just $\lfloor w \rfloor$ in place of the most correct $\lfloor \mathsf{X} \times l : w \rfloor$, leaving the missing details to the context. Accordingly, the image of the imprimitive action in $\operatorname{Sym}(\mathsf{X} \times l)$ is denoted by $\lfloor \mathsf{X} \times l : A \operatorname{wr} B \rfloor$ or just $\lfloor A \operatorname{wr} B \rfloor$. It follows immediately from (5.5.E) and (5.5.F) that whenever $\boldsymbol{\beta} \in B$ and $a \in A$,

$$\operatorname{Par}\lfloor \boldsymbol{\beta} \rfloor \equiv \big|\mathsf{X}\big| \operatorname{Par} \boldsymbol{\beta} \mod 2, \tag{5.5.G}$$

$$\operatorname{Par}\lfloor a\boldsymbol{\kappa}_i \rfloor = \operatorname{Par} a. \tag{5.5.H}$$

We denote the $l$ canonical inclusions of $\mathsf{X}$ in $\mathsf{X} \times l$ by

$$\boldsymbol{\varepsilon}_i : \mathsf{X} \longrightarrow \mathsf{X} \times l; \qquad \mathsf{x} \mapsto (\mathsf{x}, i) \tag{5.5.I}$$

and we denote by $\vec{\boldsymbol{\varepsilon}}$ the vector $(\boldsymbol{\varepsilon}_0, \ldots, \boldsymbol{\varepsilon}_{l-1})$. For $w \in W$, set

$$\vec{\boldsymbol{\varepsilon}} \cdot \lfloor w \rfloor := (\boldsymbol{\varepsilon}_0 \lfloor w \rfloor, \ldots, \boldsymbol{\varepsilon}_{l-1} \lfloor w \rfloor).$$

One readily checks that

$$\vec{\boldsymbol{\varepsilon}} \cdot \lfloor \boldsymbol{\beta} \rfloor = (\boldsymbol{\varepsilon}_{0\boldsymbol{\beta}}, \ldots, \boldsymbol{\varepsilon}_{(l-1)\boldsymbol{\beta}}),$$
$$\vec{\boldsymbol{\varepsilon}} \cdot \lfloor \vec{a} \rfloor = (a_0\boldsymbol{\varepsilon}_0, \ldots, a_{l-1}\boldsymbol{\varepsilon}_{l-1}).$$

It follows that $\lfloor A \operatorname{wr} B \rfloor$ preserves the equipartition of $\mathsf{X} \times l$ made of the images of the $\boldsymbol{\varepsilon}_i$

$$\operatorname{Is} \vec{\boldsymbol{\varepsilon}} = \big\{ \operatorname{Im} \boldsymbol{\varepsilon}_i \;\big|\; i \in l \big\}.$$

In fact, we shall see that the stabilizer in $\mathrm{Sym}(\mathsf{X}\times l)$ of $\mathrm{Is}\,\vec{\boldsymbol{\varepsilon}}$ is $\lfloor \mathrm{Sym}\,\mathsf{X} \,\mathrm{wr}\, \mathrm{S}_l \rfloor$. Observe that

$$\mathrm{Im}\,\boldsymbol{\varepsilon}_{i(w\boldsymbol{\rho})} = (\mathrm{Im}\,\boldsymbol{\varepsilon}_i)\lfloor w \rfloor, \qquad i \in l, \quad w \in W.$$

Therefore the actions of this stabilizer on $l$ and on $\mathrm{Is}\,\vec{\boldsymbol{\varepsilon}}$ are equivalent.

**5.5.2. Product action.** The product action of $W$ is defined by

$$\lceil \mathsf{X}^l : \star \rceil : W \longrightarrow \mathrm{Sym}(\mathsf{X}^l), \qquad w \mapsto \lceil \mathsf{X}^l : w \rceil;$$

(regard the symbol $\lceil\,\rceil$ as a broken $\prod$) where for $\boldsymbol{\beta} \in B$, $\vec{a} \in A^l$

$$\vec{\mathsf{x}}\lceil \mathsf{X}^l : \boldsymbol{\beta} \rceil := \boldsymbol{\beta}^{-1}\vec{\mathsf{x}} = (\mathsf{x}_{0\boldsymbol{\beta}^{-1}}, \dots, \mathsf{x}_{(l-1)\boldsymbol{\beta}^{-1}}),$$
$$\vec{\mathsf{x}}\lceil \mathsf{X}^l : \vec{a} \rceil := (\mathsf{x}_0 a_0, \dots, \mathsf{x}_{l-1} a_{l-1}).$$

However, if $w \in W$, we may write just $\lceil w \rceil$ in place of the most correct $\lceil \mathsf{X}^l : w \rceil$, leaving the missing details to the context. Accordingly, the image of the product action in $\mathrm{Sym}(\mathsf{X}^l)$ is denoted by $\lceil \mathsf{X}^l : A \,\mathrm{wr}\, B \rceil$ or just $\lceil A \,\mathrm{wr}\, B \rceil$. The parity issue is now settled as follows. To save some writing, say $m = |\mathsf{X}|$. Then

$$\mathrm{Par}\lceil \boldsymbol{\beta} \rceil \equiv m^{l-1}\frac{m-1}{2}\,\mathrm{Par}\,\boldsymbol{\beta} \mod 2, \tag{5.5.J}$$
$$\mathrm{Par}\lceil a\boldsymbol{\kappa}_i \rceil \equiv m^{l-1}\,\mathrm{Par}\,a \mod 2, \tag{5.5.K}$$

whenever $\boldsymbol{\beta} \in B$ and $a \in A$.

We denote the canonical projections of $\mathsf{X}^l$ onto $\mathsf{X}$ by

$$\boldsymbol{\pi}_i : \mathsf{X}^l \longrightarrow \mathsf{X}; \qquad \vec{\mathsf{x}} \mapsto \mathsf{x}_i \tag{5.5.L}$$

and we denote by $\vec{\boldsymbol{\pi}}$ the vector $(\boldsymbol{\pi}_0, \dots, \boldsymbol{\pi}_{l-1})$. For $w \in W$, set

$$\lceil w \rceil \cdot \vec{\boldsymbol{\pi}} := (\lceil w \rceil \boldsymbol{\pi}_0, \dots, \lceil w \rceil \boldsymbol{\pi}_{l-1}).$$

Then

$$\lceil \boldsymbol{\beta} \rceil \cdot \vec{\boldsymbol{\pi}} = (\boldsymbol{\pi}_{0\boldsymbol{\beta}^{-1}}, \dots, \boldsymbol{\pi}_{(l-1)\boldsymbol{\beta}^{-1}}), \tag{5.5.M}$$
$$\lceil \vec{a} \rceil \cdot \vec{\boldsymbol{\pi}} = (\boldsymbol{\pi}_0 a_0, \dots, \boldsymbol{\pi}_{l-1} a_{l-1}). \tag{5.5.N}$$

By a *hyperplane* of $\mathsf{X}^l$ we mean a fiber of any one of the projections $\boldsymbol{\pi}_i$:

$$\mathrm{Hyp}(i, x) := \mathrm{Fb}(\boldsymbol{\pi}_i, \mathsf{x}) = \{\, \vec{\mathsf{x}} \mid \vec{\mathsf{x}}\boldsymbol{\pi}_i = \mathsf{x} \,\}.$$

By a *codirection* we mean the set consisting of all the fibers of one of the $\boldsymbol{\pi}_i$, that is, the quotient set of that $\boldsymbol{\pi}_i$:

$$\mathrm{Cd}(i) := \mathrm{Qt}\,\boldsymbol{\pi}_i = \{\, \mathrm{Fb}(\boldsymbol{\pi}_i, \mathsf{x}) \mid \mathsf{x} \in \mathsf{X} \,\}.$$

Clearly $\mathrm{Hyp}(i,x)\lceil\boldsymbol{\beta}\rceil = \mathrm{Hyp}(i\boldsymbol{\beta},x)$ and $\mathrm{Hyp}(i,x)\lceil\vec{a}\rceil = \mathrm{Hyp}(i,xa_i)$. This implies $\mathrm{Cd}(i)\lceil\boldsymbol{\beta}\rceil = \mathrm{Cd}(i\boldsymbol{\beta})$ and $\mathrm{Cd}(i)\lceil\vec{a}\rceil = \mathrm{Cd}(i)$. It follows that $\lceil A \,\mathrm{wr}\, B\rceil$ preserves the set

$$\mathsf{H} := \left\{ \mathrm{Hyp}(i,x) \;\middle|\; i \in l, \quad \mathsf{x} \in \mathsf{X} \right\}$$

of hyperplanes of $\mathsf{X}^l$ and also preserves the set of its codirections

$$\mathrm{Qs}\,\vec{\boldsymbol{\pi}} = \left\{ \mathrm{Cd}(i) \;\middle|\; i \in l \right\}$$

which is in itself a non trivial equipartition of $\mathsf{H}$. In fact, we shall see that $\lceil \mathrm{Sym}\,\mathsf{X} \,\mathrm{wr}\, \mathrm{S}_l\rceil$ is the largest subgroup of $\mathrm{Sym}(\mathsf{X}^l)$ which preserves $\mathsf{H}$ and this subgroup also preserves $\mathrm{Qs}\,\vec{\boldsymbol{\pi}}$. Because of

$$\mathrm{Cd}(i(w\boldsymbol{\rho})) = \mathrm{Cd}(i)\lceil w\rceil, \qquad i \in l, \quad w \in W, \tag{5.5.O}$$

the actions of this stabilizer on $l$ and on $\mathrm{Qs}\,\vec{\boldsymbol{\pi}}$ are equivalent. Moreover, note that each codirection is an imprimitivity system of the base subgroup of the wreath product.

## 5.6. Inner wreath products

The aim of this section is to provide some language for working with the inner wreath products. These are the permutation groups which are permutationally isomorphic to wreath products either in imprimitive or in product action. Consequently there are inner wreath products of two kinds, namely inner wreath product on coproduct decompositions and inner wreath product on product decompositions. In [**Kov89b**] there is a conclusive detailed description of such subgroups in terms of convenient decompositions of the subgroups themselves. Here we prefer a more naive approach.

Given a set $\mathsf{X}$ of order at least 2 and an integer $l > 1$ we will first define equicoproduct decompositions of a set $\mathsf{Y}$ as bijections $\boldsymbol{f} : \mathsf{X} \times l \longrightarrow \mathsf{Y}$. In this way, we may compose canonical inclusions of $\mathsf{X}$ into $\mathsf{X} \times l$ with these bijections to yield inclusions into $\mathsf{Y}$. Similarly, we will define equiproduct decompositions as bijections $\boldsymbol{p} : \mathsf{Y} \longrightarrow \mathsf{X}^l$ (opposite direction to enable the composition with the canonical projections of $\mathsf{X}^l$ onto $\mathsf{X}$). Then we simply write $\lfloor \boldsymbol{f} : \boldsymbol{s}\rfloor$ or $\lceil \boldsymbol{p} : \boldsymbol{s}\rceil$ for the permutation of $\mathsf{Y}$ corresponding via $\boldsymbol{f}$ or $\boldsymbol{p}$ to the permutation $\boldsymbol{s}$ of $\mathsf{X} \times l$ or $\mathsf{X}^l$:

$$\lfloor \boldsymbol{f} : \boldsymbol{s}\rfloor := \boldsymbol{f}^{-1}\boldsymbol{s}\boldsymbol{f},$$
$$\lceil \boldsymbol{p} : \boldsymbol{s}\rceil := \boldsymbol{p}\boldsymbol{s}\boldsymbol{p}^{-1}.$$

The reader should not worry to remember whether it is a conjugation by $\boldsymbol{f}$ or by $\boldsymbol{f}^{-1}$ in the first case. There is one and only one sensible way to compose $\boldsymbol{f}$, $\boldsymbol{s}$, $\boldsymbol{f}^{-1}$ together once we agree that $\boldsymbol{f}$ goes from $\mathsf{X} \times l$ to $\mathsf{Y}$ and we agree to this because this alone enables composition of $\boldsymbol{f}$ with the canonical inclusions. Similarly, in the second case, there is one and only one way to define $\boldsymbol{p}$ in a way that enables composition with the canonical projections of $\mathsf{X}^l$ onto $\mathsf{X}$. Then the definition of $\lceil \boldsymbol{p} : \boldsymbol{s} \rceil$ may not be any different.

The details of all this follow.

**Equicoproduct decompositions.** A *coproduct* is a family of functions

$$\boldsymbol{f}_\lambda : \mathsf{X}_\lambda \longrightarrow \mathsf{Y}$$

with a common codomain $\mathsf{Y}$ satisfying the Coproduct Universal Property.

CUP. For each family of functions $\boldsymbol{g}_\lambda : \mathsf{X}_\lambda \longrightarrow \mathsf{Z}$ with a common codomain $\mathsf{Z}$, there is a unique $\boldsymbol{g} : \mathsf{Y} \longrightarrow \mathsf{Z}$ such that $\boldsymbol{g}_\lambda = \boldsymbol{f}_\lambda \boldsymbol{g}$ for all $\lambda$.

In this case each $\boldsymbol{f}_\lambda$ is one-to-one, and the $\boldsymbol{f}_\lambda$ are called canonical inclusions of the coproduct. Informally, one may also say that $\mathsf{Y}$ is the coproduct of the $\mathsf{X}_\lambda$. For example, $\mathsf{X} \times l$ together with the canonical inclusions (5.5.I) is the coproduct of $l$ copies of $\mathsf{X}$. Since $\operatorname{Dom} \boldsymbol{\varepsilon}_i = \mathsf{X}$ for all $i$, we say that $\mathsf{X} \times l$ is an *equicoproduct*. We are only interested in equicoproducts with finitely many canonical inclusions, and we collect the latter in a vector $\vec{\boldsymbol{f}} \in (\mathsf{Y}^{\mathsf{X}})^l$ where $l$ is the number of canonical inclusions, $\mathsf{X}$ is their domain and $\mathsf{Y}$ is their codomain. It is easily seen that if $\boldsymbol{f} : \mathsf{X} \times l \longrightarrow \mathsf{Y}$ is a bijection, then $\mathsf{Y}$ together with

$$\vec{\boldsymbol{f}} := \vec{\boldsymbol{\varepsilon}} \cdot \boldsymbol{f} = (\boldsymbol{\varepsilon}_0 \boldsymbol{f}, \dots, \boldsymbol{\varepsilon}_{l-1} \boldsymbol{f})$$

is an equicoproduct. What makes the CUP so important is that the converse is also true.

5.6.1. THEOREM. *A set $\mathsf{Y}$ together with $\vec{\boldsymbol{f}} \in (\mathsf{Y}^{\mathsf{X}})^l$ is an equicoproduct if and only if there is a bijection $\boldsymbol{f} : \mathsf{X} \times l \longrightarrow \mathsf{Y}$ such that $\vec{\boldsymbol{f}} = \vec{\boldsymbol{\varepsilon}} \cdot \boldsymbol{f}$. Furthermore, in this case, the bijection is uniquely determined.*

Given sets $\mathsf{Y}$, $\mathsf{X}$ of order at least 2 such that $|\mathsf{Y}| = l|\mathsf{X}|$ for some integer $l > 1$, we define the non empty set of the equicoproduct decompositions

$$\operatorname{ECDec}(\mathsf{Y}, \mathsf{X}, l) := \{\, \boldsymbol{f} : \mathsf{X} \times l \longrightarrow \mathsf{Y} \mid \boldsymbol{f} \text{ is a bijection} \,\}.$$

To save some writing we may occasionally omit the third argument $l$ which is uniquely determined by the order of $\mathsf{Y}$ and $\mathsf{X}$. Also, by definition, in writing $\operatorname{ECDec}(\mathsf{Y}, \mathsf{X}, l)$ or $\operatorname{ECDec}(\mathsf{Y}, \mathsf{X})$ we implicitly assume that $\mathsf{X}$ and $\mathsf{Y}$ have suitable orders.

Now, let $A$ act on $X$, let $B$ be a subgroup of $\mathrm{S}_l$ and let $\boldsymbol{f} \in \mathrm{ECDec}(\mathsf{Y}, \mathsf{X}, l)$. All the elements of $\mathrm{Sym}\,\mathsf{Y}$ may be written as

$$\lfloor \boldsymbol{f} : \boldsymbol{s} \rfloor := \boldsymbol{f}^{-1}\boldsymbol{s}\boldsymbol{f}$$

with $\boldsymbol{s} \in \mathrm{Sym}(\mathsf{X} \times l)$. But to avoid redundancy, when $w \in A \,\mathrm{wr}\, B$, set

$$\lfloor \boldsymbol{f} : w \rfloor := \lfloor \boldsymbol{f} : \lfloor \mathsf{X} \times l : w \rfloor \rfloor.$$

The inner wreath products on equicoproduct decompositions are the permutation groups $\lfloor \boldsymbol{f} : A \,\mathrm{wr}\, B \rfloor$ with $\boldsymbol{f} \in \mathrm{ECDec}(\mathsf{Y}, \mathsf{X}, l)$, $A \leqslant \mathrm{Sym}\,\mathsf{X}$ and $B \leqslant \mathrm{S}_l$.

5.6.2. Proposition. *If* $\boldsymbol{f}$, $\boldsymbol{g} \in \mathrm{ECDec}(\mathsf{Y}, \mathsf{X})$ *then the two inner wreath products* $\lfloor \boldsymbol{f} : A \,\mathrm{wr}\, B \rfloor$ *and* $\lfloor \boldsymbol{g} : A \,\mathrm{wr}\, B \rfloor$ *are conjugate in* $\mathrm{Sym}\,\mathsf{Y}$.

Proof. $\mathrm{Sym}\,\mathsf{Y}$ with the action given by composition of functions is transitive on $\mathrm{ECDec}(\mathsf{Y}, \mathsf{X})$ and so there exists $\boldsymbol{s} \in \mathrm{Sym}\,\mathsf{Y}$ such that $\boldsymbol{g} = \boldsymbol{f}\boldsymbol{s}$. One checks that whenever $\boldsymbol{t} \in \mathrm{Sym}(\mathsf{X} \times l)$,

$$\boldsymbol{s}^{-1}\lfloor \boldsymbol{f} : \boldsymbol{t} \rfloor \boldsymbol{s} = \boldsymbol{s}^{-1}\boldsymbol{f}^{-1}\boldsymbol{t}\boldsymbol{f}\boldsymbol{s} = \lfloor \boldsymbol{f}\boldsymbol{s} : \boldsymbol{t} \rfloor = \lfloor \boldsymbol{g} : \boldsymbol{t} \rfloor.$$

Therefore $\boldsymbol{s}^{-1}\lfloor \boldsymbol{f} : A \,\mathrm{wr}\, B \rfloor \boldsymbol{s} = \lfloor \boldsymbol{g} : A \,\mathrm{wr}\, B \rfloor$. □

Clearly, $\mathrm{Sym}\,\mathsf{Y}$ acts regularly by composition of functions on the right of $\mathrm{ECDec}(\mathsf{Y}, \mathsf{X})$. It follows that any permutation of $\mathsf{Y}$ is uniquely determined by its action on any given $\boldsymbol{f} \in \mathrm{ECDec}(\mathsf{Y}, \mathsf{X})$. For example, when $\boldsymbol{\beta} \in \mathrm{S}_l$, $\lfloor \boldsymbol{f} : \boldsymbol{\beta} \rfloor$ is the unique permutation $\boldsymbol{s}$ of $\mathsf{Y}$ such that

$$(\boldsymbol{f}_0\boldsymbol{s}, \ldots, \boldsymbol{f}_{(l-1)}\boldsymbol{s}) = (\boldsymbol{f}_{0\boldsymbol{\beta}}, \ldots, \boldsymbol{f}_{(l-1)\boldsymbol{\beta}}). \tag{5.6.A}$$

Similarly, when $\vec{\boldsymbol{a}} \in (\mathrm{Sym}\,\mathsf{X})^l$, $\lfloor \boldsymbol{f} : \vec{\boldsymbol{a}} \rfloor$ is the unique permutation $\boldsymbol{s}$ of $\mathsf{Y}$ such that

$$(\boldsymbol{f}_0\boldsymbol{s}, \ldots, \boldsymbol{f}_{(l-1)}\boldsymbol{s}) = (\boldsymbol{a}_0\boldsymbol{f}_0, \ldots, \boldsymbol{a}_{l-1}\boldsymbol{f}_{l-1}). \tag{5.6.B}$$

Suppose now that $\boldsymbol{f} \in \mathrm{ECDec}(\mathsf{Y}, \mathsf{X}, l)$. Our goal is to show that the inner wreath product $\lfloor \boldsymbol{f} : (\mathrm{Sym}\,\mathsf{X}) \,\mathrm{wr}\, \mathrm{S}_l \rfloor$ is the stabilizer in $\mathrm{Sym}\,\mathsf{Y}$ of the equipartition made of the images of the canonical inclusions collected in $\vec{\boldsymbol{f}}$. We begin with introducing some new terms.

Given a subset $\mathsf{A}$ of $\mathsf{Y}^{\mathsf{X}}$ we define the set of images from $\mathsf{A}$

$$\mathrm{Is}\,\mathsf{A} := \left\{ \mathrm{Im}\,\boldsymbol{a} \;\middle|\; \boldsymbol{a} \in \mathsf{A} \right\}$$

and the left closure

$$\mathrm{LC}\,\mathsf{A} = (\mathrm{Sym}\,\mathsf{X})\mathsf{A} = \left\{ \boldsymbol{s}\boldsymbol{a} \;\middle|\; \boldsymbol{s} \in \mathrm{Sym}\,\mathsf{X},\, \boldsymbol{a} \in \mathsf{A} \right\}.$$

It follows at once from the definitions that

$$\mathsf{A} \subseteq \mathrm{LC}\,\mathsf{A}, \qquad \mathrm{LC}\,\mathrm{LC}\,\mathsf{A} = \mathrm{LC}\,\mathsf{A}, \qquad \mathrm{Is}\,\mathsf{A} = \mathrm{Is}\,\mathrm{LC}\,\mathsf{A}.$$

Furthermore, if $\boldsymbol{s} \in \operatorname{Sym} \mathsf{Y}$, then

$$(\operatorname{Is} \mathsf{A})\boldsymbol{s} = \operatorname{Is}(\mathsf{A}\boldsymbol{s}), \qquad (\operatorname{LC} \mathsf{A})\boldsymbol{s} = \operatorname{LC}(\mathsf{A}\boldsymbol{s}).$$

Finally, we recall from 1.1 that if $\vec{\boldsymbol{f}} \in (\mathsf{Y}^{\mathsf{X}})^l$, then we denote $\{\boldsymbol{f}_0, \dots, \boldsymbol{f}_{l-1}\}$ by $\tilde{\boldsymbol{f}}$.

5.6.3. THEOREM. *Let* $\boldsymbol{f}, \boldsymbol{g} : \mathsf{X} \times l \longrightarrow \mathsf{Y}$*, with* $\boldsymbol{f}$ *bijective. Set* $\vec{\boldsymbol{g}} = \vec{\boldsymbol{\varepsilon}} \cdot \boldsymbol{g}$*,* $\vec{\boldsymbol{f}} = \vec{\boldsymbol{\varepsilon}} \cdot \boldsymbol{f}$*. The following are equivalent:*

(a) $\operatorname{Is} \tilde{\boldsymbol{g}} = \operatorname{Is} \tilde{\boldsymbol{f}}$.
(b) $\boldsymbol{g} \in \boldsymbol{f} \lfloor \boldsymbol{f} : (\operatorname{Sym} \mathsf{X}) \operatorname{wr} \mathrm{S}_l \rfloor$.
(c) $\operatorname{LC} \tilde{\boldsymbol{g}} = \operatorname{LC} \tilde{\boldsymbol{f}}$.
(d) $\boldsymbol{g}$ *is a bijection and* $\tilde{\boldsymbol{g}} \subset \operatorname{LC} \tilde{\boldsymbol{f}}$.

*In particular* $\lfloor \boldsymbol{f} : (\operatorname{Sym} \mathsf{X}) \operatorname{wr} \mathrm{S}_l \rfloor = \operatorname{St} \operatorname{Is} \tilde{\boldsymbol{f}} = \operatorname{St} \operatorname{LC} \tilde{\boldsymbol{f}}$.

PROOF. Assume (a). Then there is a function $\boldsymbol{\sigma} : l \longrightarrow l$ such that

$$\operatorname{Im} \boldsymbol{f}_i = \operatorname{Im} \boldsymbol{g}_{i\boldsymbol{\sigma}} \qquad \forall i \in l.$$

But there are $l$ different $\operatorname{Im} \boldsymbol{f}_i$'s therefore $\boldsymbol{\sigma}$ must be surjective and hence bijective, so $\boldsymbol{\sigma} \in \mathrm{S}_l$. For each $i \in l$, $\boldsymbol{f}_{i\boldsymbol{\sigma}^{-1}}$ is injective and we set $\boldsymbol{t}_i = \boldsymbol{g}_i \boldsymbol{f}_{i\boldsymbol{\sigma}^{-1}}$ so that $\boldsymbol{t}_i \in \operatorname{Sym} \mathsf{X}$ and $\boldsymbol{g}_i = \boldsymbol{t}_i \boldsymbol{f}_{i\boldsymbol{\sigma}^{-1}}$. Then $\vec{\boldsymbol{g}} = \vec{\boldsymbol{f}} \lfloor \boldsymbol{f} : \boldsymbol{\sigma} \rfloor \lfloor \boldsymbol{f} : \vec{\boldsymbol{t}} \rfloor$ which is equivalent to $\boldsymbol{g} = \boldsymbol{f} \lfloor \boldsymbol{f} : \boldsymbol{\sigma} \rfloor \lfloor \boldsymbol{f} : \vec{\boldsymbol{t}} \rfloor$. This proves (c) and also (b).

Clearly (b) and (c) together imply (d) so that (d) follows from (a).

On the other hand (a) follows at once from (b).

Now assume (c), then $\operatorname{Is} \tilde{\boldsymbol{g}} = \operatorname{Is} \operatorname{LC} \tilde{\boldsymbol{g}} = \operatorname{Is} \operatorname{LC} \tilde{\boldsymbol{f}} = \operatorname{Is} \tilde{\boldsymbol{f}}$ which is (a).

If by any chance (d) holds, then $\operatorname{LC} \tilde{\boldsymbol{g}} \subseteq \operatorname{LC} \operatorname{LC} \tilde{\boldsymbol{f}} = \operatorname{LC} \tilde{\boldsymbol{f}}$ thus $\operatorname{Is} \tilde{\boldsymbol{g}} \subseteq \operatorname{Is} \tilde{\boldsymbol{f}}$. But $\operatorname{Is} \tilde{\boldsymbol{g}}$ has $l$ elements because $\boldsymbol{g} \in \operatorname{ECDec}(\mathsf{Y}, \mathsf{X}, l)$. Then $\operatorname{Is} \tilde{\boldsymbol{g}} = \operatorname{Is} \tilde{\boldsymbol{f}}$ which is(a).

This proves the equivalence. Assume now that $\boldsymbol{s} \in \operatorname{Sym} \mathsf{Y}$ and say $\boldsymbol{g} = \boldsymbol{f}\boldsymbol{s}$. Then $(\operatorname{Is} \tilde{\boldsymbol{f}})\boldsymbol{s} = \operatorname{Is} \tilde{\boldsymbol{g}}$ and $(\operatorname{LC} \tilde{\boldsymbol{f}})\boldsymbol{s} = \operatorname{LC} \tilde{\boldsymbol{g}}$. Therefore $\boldsymbol{s}$ preserves $\operatorname{Is} \tilde{\boldsymbol{f}}$ if and only if $\operatorname{Is} \tilde{\boldsymbol{f}} = \operatorname{Is} \tilde{\boldsymbol{g}}$, if and only if $\operatorname{LC} \tilde{\boldsymbol{f}} = \operatorname{LC} \tilde{\boldsymbol{g}}$, if and only if $\boldsymbol{s} \in \lfloor \boldsymbol{f} : (\operatorname{Sym} \mathsf{X}) \operatorname{wr} \mathrm{S}_l \rfloor$, if and only if $\boldsymbol{s}$ preserves $\operatorname{LC} \tilde{\boldsymbol{f}}$. □

**Equiproduct decompositions.** Products are dual to coproducts in a categorical sense. Consequently, this section follows the pattern set in the previous one and a few facts are stated without proof. Any further comment on the underlying duality is left to the reader.

A *product* is a family of functions

$$\boldsymbol{p}_\lambda : \mathsf{Y} \longrightarrow \mathsf{X}_\lambda$$

with a common domain $\mathsf{Y}$ satisfying the Product Universal Property.

PUP. For each family of functions $\boldsymbol{g}_\lambda : \mathsf{Z} \longrightarrow \mathsf{X}_\lambda$ with a common domain $\mathsf{Z}$, there is a unique $\boldsymbol{g} : \mathsf{Z} \longrightarrow \mathsf{Y}$ such that $\boldsymbol{g}_\lambda = \boldsymbol{g}\boldsymbol{p}_\lambda$ for all $\lambda$.

In this case each $\boldsymbol{p}_\lambda$ is surjective, and the $\boldsymbol{p}_\lambda$ are called canonical projections of the product. Informally, $\mathsf{Y}$ is also said to be the product of the $\mathsf{X}_\lambda$. For example $\mathsf{X}^l$ together with the canonical projections (5.5.L) is the product of $l$ copies of $\mathsf{X}$. Since $\operatorname{Cod} \boldsymbol{\pi}_i = \mathsf{X}$ for all $i$, we say that $\mathsf{X}^l$ is an *equiproduct*. We are only interested in equiproducts with finitely many canonical projections, and we collect the latter in a vector $\vec{\boldsymbol{p}} \in (\mathsf{X}^{\mathsf{Y}})^l$ where $l$ is the number of canonical projections, $\mathsf{X}$ is their codomain and $\mathsf{Y}$ is their domain. It is easily seen that if $\boldsymbol{p} : \mathsf{Y} \longrightarrow \mathsf{X}^l$ is a bijection, then $\mathsf{Y}$ together with

$$\vec{\boldsymbol{p}} := \boldsymbol{p} \cdot \vec{\boldsymbol{\pi}} = (\boldsymbol{p}\boldsymbol{\pi}_0, \ldots, \boldsymbol{p}\boldsymbol{\pi}_{l-1})$$

is an equiproduct. What makes the PUP so important is that the converse is also true.

5.6.4. THEOREM. *A set* $\mathsf{Y}$ *together with* $\vec{\boldsymbol{p}} \in (\mathsf{X}^{\mathsf{Y}})^l$ *is an equiproduct if and only if there is a bijection* $\boldsymbol{p} : \mathsf{Y} \longrightarrow \mathsf{X}^l$ *such that* $\vec{\boldsymbol{p}} = \boldsymbol{p} \cdot \vec{\boldsymbol{\pi}}$*. Furthermore, in this case, the bijection is uniquely determined.*

Given sets $\mathsf{Y}$, $\mathsf{X}$ of order at least 2 such that $|\mathsf{Y}| = |\mathsf{X}|^l$ for some integer $l > 1$, we define the non-empty set of the equiproduct decompositions

$$\operatorname{EPDec}(\mathsf{Y}, \mathsf{X}, l) := \{\, \boldsymbol{p} : \mathsf{Y} \longrightarrow \mathsf{X}^l \;\mid\; \boldsymbol{p} \text{ is a bijection} \,\}.$$

To save some writing we may occasionally omit the third argument $l$ which is uniquely determined by the order of $\mathsf{Y}$ and $\mathsf{X}$. Also, by definition, in writing $\operatorname{EPDec}(\mathsf{Y}, \mathsf{X}, l)$ or $\operatorname{EPDec}(\mathsf{Y}, \mathsf{X})$ we implicitly assume that $\mathsf{X}$ and $\mathsf{Y}$ have suitable orders.

Now, let $A$ act on $X$, let $B$ be a subgroup of $\mathrm{S}_l$ and let $\boldsymbol{p} \in \operatorname{EPDec}(\mathsf{Y}, \mathsf{X}, l)$. All the elements of $\operatorname{Sym} \mathsf{Y}$ may be written as

$$\lceil \boldsymbol{p} : \boldsymbol{s} \rceil := \boldsymbol{p}\boldsymbol{s}\boldsymbol{p}^{-1}$$

with $\boldsymbol{s} \in \operatorname{Sym} \mathsf{X}^l$. But to avoid redundancy, when $w \in A \operatorname{wr} B$, set

$$\lceil \boldsymbol{p} : w \rceil := \lceil \boldsymbol{p} : \lceil \mathsf{X}^l : w \rceil \rceil.$$

The inner wreath products on equiproduct decompositions are the permutation groups $\lceil \boldsymbol{p} : A \operatorname{wr} B \rceil$ with $\boldsymbol{p} \in \operatorname{EPDec}(\mathsf{Y}, \mathsf{X}, l)$, $A \leqslant \operatorname{Sym} \mathsf{X}$ and $B \leqslant \mathrm{S}_l$.

5.6.5. PROPOSITION. *If* $\boldsymbol{p}$, $\boldsymbol{g} \in \operatorname{EPDec}(\mathsf{Y}, \mathsf{X})$ *then the two inner wreath products* $\lceil \boldsymbol{p} : A \operatorname{wr} B \rceil$, $\lceil \boldsymbol{g} : A \operatorname{wr} B \rceil$ *are conjugate in* $\operatorname{Sym} \mathsf{Y}$.

Clearly, $\operatorname{Sym}\mathsf{Y}$ acts regularly on the left of $\operatorname{EPDec}(\mathsf{Y},\mathsf{X})$. It follows that any permutation of $\mathsf{Y}$ is uniquely determined by its action on any given $\boldsymbol{p} \in \operatorname{EPDec}(\mathsf{Y},\mathsf{X})$. For example, when $\boldsymbol{\beta} \in \mathrm{S}_l$, $\lceil \boldsymbol{p} : \boldsymbol{\beta} \rceil$ is the unique permutation $\boldsymbol{s}$ of $\mathsf{Y}$ such that

$$(\boldsymbol{s}\boldsymbol{p}_0, \ldots, \boldsymbol{s}\boldsymbol{p}_{l-1}) = (\boldsymbol{p}_{0\boldsymbol{\beta}^{-1}}, \ldots, \boldsymbol{p}_{(l-1)\boldsymbol{\beta}^{-1}}). \tag{5.6.C}$$

Similarly, when $\vec{\boldsymbol{a}} \in (\operatorname{Sym}\mathsf{X})^l$, $\lceil \boldsymbol{p} : \vec{\boldsymbol{a}} \rceil$ is the unique permutation $\boldsymbol{s}$ of $\mathsf{Y}$ such that

$$(\boldsymbol{s}\boldsymbol{p}_0, \ldots, \boldsymbol{s}\boldsymbol{p}_{l-1}) = (\boldsymbol{p}_0\boldsymbol{a}_0, \ldots, \boldsymbol{p}_{l-1}\boldsymbol{a}_{l-1}). \tag{5.6.D}$$

Suppose now that $\boldsymbol{p} \in \operatorname{EPDec}(\mathsf{Y},\mathsf{X},l)$. Our goal is to show that the inner wreath product $\lceil \boldsymbol{p} : (\operatorname{Sym}\mathsf{X}) \operatorname{wr} \mathrm{S}_l \rceil$ is the stabilizer in $\operatorname{Sym}\mathsf{Y}$ of the set of hyperplanes related to the equiproduct decomposition $\boldsymbol{p}$. As in the case of the coproduct we begin with introducing some terms.

Let $\boldsymbol{f} : \mathsf{U} \longrightarrow \mathsf{V}$ be a function. As common, we denote by $\mathsf{A}\boldsymbol{f}^{-1}$ the pre-image of a subset $\mathsf{A}$ of $\mathsf{V}$:

$$\mathsf{A}\boldsymbol{f}^{-1} := \{\ \mathsf{u} \in \mathsf{U} \ \mid\ \mathsf{u}\boldsymbol{f} \in \mathsf{A}\ \}.$$

The nonempty pre-images of singletons are called *fibers*:

$$\mathsf{v} \in \operatorname{Im}\boldsymbol{f}, \qquad \operatorname{Fb}(\boldsymbol{f},\mathsf{v}) := \{\mathsf{v}\}\boldsymbol{f}^{-1}.$$

The set of fibers of a function $\boldsymbol{f}$ is a partition of $\operatorname{Dom}\boldsymbol{f}$ called *quotient set* of $\boldsymbol{f}$ which we denote by:

$$\operatorname{Qt}\boldsymbol{f} := \{\ \operatorname{Fb}(\boldsymbol{f},\mathsf{v}) \ \mid\ \mathsf{v} \in \operatorname{Im}\boldsymbol{f}\ \}.$$

Given a subset $\mathsf{A}$ of $\mathsf{X}^{\mathsf{Y}}$ we define the set of fibers from $\mathsf{A}$:

$$\operatorname{Fs}\mathsf{A} := \bigcup_{\boldsymbol{a}\in\mathsf{A}} \operatorname{Qt}\boldsymbol{a} = \{\ \operatorname{Fb}(\boldsymbol{a},\mathsf{y}) \ \mid\ \boldsymbol{a} \in \mathsf{A},\ \mathsf{y} \in \operatorname{Im}\boldsymbol{a}\ \},$$

the set of quotients from $\mathsf{A}$:

$$\operatorname{Qs}\mathsf{A} = \{\ \operatorname{Qt}\boldsymbol{a} \ \mid\ \boldsymbol{a} \in \mathsf{A}\ \}$$

and the right closure of $\mathsf{A}$:

$$\operatorname{RC}\mathsf{A} = \mathsf{A}(\operatorname{Sym}\mathsf{X}) = \{\ \boldsymbol{a}\boldsymbol{s} \ \mid\ \boldsymbol{s} \in \operatorname{Sym}\mathsf{X},\ \boldsymbol{a} \in \mathsf{A}\ \}.$$

It follows at once from these definitions that

$$\mathsf{A} \subseteq \operatorname{RC}\mathsf{A}, \qquad \operatorname{RC}\operatorname{RC}\mathsf{A} = \operatorname{RC}\mathsf{A},$$

$$\operatorname{Fs}\mathsf{A} = \operatorname{Fs}\operatorname{RC}\mathsf{A}, \qquad \operatorname{Qs}\mathsf{A} = \operatorname{Qs}\operatorname{RC}\mathsf{A},$$

$$\operatorname{Fs}\mathsf{A} = \bigcup \operatorname{Qs}\mathsf{A}. \tag{5.6.E}$$

Furthermore, if $\boldsymbol{s} \in \operatorname{Sym}\mathsf{Y}$ then

$$(\operatorname{Fs}\mathsf{A})\boldsymbol{s} = \operatorname{Fs}(\boldsymbol{s}^{-1}\mathsf{A}), \qquad (\operatorname{Qs}\mathsf{A})\boldsymbol{s} = \operatorname{Qs}(\boldsymbol{s}^{-1}\mathsf{A}), \qquad \boldsymbol{s}(\operatorname{RC}\mathsf{A}) = \operatorname{RC}(\boldsymbol{s}\mathsf{A}).$$

If $\boldsymbol{p} \in \mathrm{EPDec}(\mathsf{Y}, \mathsf{X}, l)$, then the elements of $\mathrm{Fs}\,\tilde{\boldsymbol{p}}$ are called hyperplanes of the product decomposition $\boldsymbol{p}$ and the elements of $\mathrm{Qs}\,\tilde{\boldsymbol{p}}$ are called codirections of the product decomposition $\boldsymbol{p}$. Codirections are equipartitions of $\mathsf{Y}$ and the set of the codirections is an equipartition of the set of the hyperplanes. It follows from the PUP that each codirection is in bijection with $\mathsf{X}$ and that codirections are the maximal sets of pairwise disjoint hyperplanes. In fact, for each $\mathsf{y} \in \mathsf{Y}$ we have

$$\{\mathsf{y}\} = \bigcap_{i \in l} \mathrm{Fb}(\boldsymbol{p}_i, \mathsf{y}\boldsymbol{p}_i)$$

which follows from the validity in $\mathsf{X}^l$ of

$$\{(\mathsf{x}_0, \ldots, \mathsf{x}_{l-1})\} = \bigcap_{i \in l} \mathrm{Hyp}(i, \mathsf{x}_i)$$

and shows that any two hyperplanes of different codirections must intersect in a non empty set. In particular, any permutation of $\mathsf{Y}$ which sends hyperplanes to hyperplanes must also send codirections to codirections. On the other hand, if a permutation of $\mathsf{Y}$ sends codirections to codirections then in particular it sends hyperplanes to hyperplanes. Therefore $\mathrm{St}\,\mathrm{Fs}\,\tilde{\boldsymbol{p}} = \mathrm{St}\,\mathrm{Qs}\,\tilde{\boldsymbol{p}}$ and we now show that this stabilizer is $\lceil \boldsymbol{p} : (\mathrm{Sym}\,\mathsf{X}) \,\mathrm{wr}\, \mathrm{S}_l \rceil$.

5.6.6. Theorem. *Let $\boldsymbol{p}, \boldsymbol{q} : \mathsf{Y} \longrightarrow \mathsf{X}^l$, with $\boldsymbol{p}$ bijective. Set $\vec{\boldsymbol{q}} = \boldsymbol{q} \cdot \vec{\boldsymbol{\pi}}$, $\vec{\boldsymbol{p}} = \boldsymbol{p} \cdot \vec{\boldsymbol{\pi}}$. The following are equivalent:*

(a) $\mathrm{Fs}\,\tilde{\boldsymbol{q}} = \mathrm{Fs}\,\tilde{\boldsymbol{p}}$.
(b) $\mathrm{Qs}\,\tilde{\boldsymbol{q}} = \mathrm{Qs}\,\tilde{\boldsymbol{p}}$.
(c) $\boldsymbol{q} \in \lceil \boldsymbol{p} : (\mathrm{Sym}\,\mathsf{X}) \,\mathrm{wr}\, \mathrm{S}_l \rceil \boldsymbol{p}$.
(d) $\mathrm{RC}\,\tilde{\boldsymbol{q}} = \mathrm{RC}\,\tilde{\boldsymbol{p}}$.
(e) *$\boldsymbol{q}$ is a bijection and* $\tilde{\boldsymbol{q}} \subset \mathrm{RC}\,\tilde{\boldsymbol{p}}$.

*In particular* $\lceil \boldsymbol{p} : (\mathrm{Sym}\,\mathsf{X}) \,\mathrm{wr}\, \mathrm{S}_l \rceil = \mathrm{St}\,\mathrm{Fs}\,\tilde{\boldsymbol{p}} = \mathrm{St}\,\mathrm{Qs}\,\tilde{\boldsymbol{p}} = \mathrm{St}\,\mathrm{RC}\,\tilde{\boldsymbol{p}}$.

Proof. We show that (c) $\Rightarrow$ (d) $\Rightarrow$ (a), (c) $\Rightarrow$ (e) $\Rightarrow$ (a), (b) $\Rightarrow$ (a) and, finally, (a) $\Rightarrow$ ((b) and (c)).

From (c), (5.6.C) and (5.6.D), we get (d) which implies (a) because

$$\mathrm{Fs}\,\tilde{\boldsymbol{q}} = \mathrm{Fs}\,\mathrm{RC}\,\tilde{\boldsymbol{q}} = \mathrm{Fs}\,\mathrm{RC}\,\tilde{\boldsymbol{p}} = \mathrm{Fs}\,\tilde{\boldsymbol{p}}.$$

Of course, (c) and (d) imply (e) which implies (a) because $\mathrm{Fs}\,\tilde{\boldsymbol{q}} \subseteq \mathrm{Fs}\,\mathrm{RC}\,\tilde{\boldsymbol{p}} = \mathrm{Fs}\,\tilde{\boldsymbol{p}}$ and the two sets have the same number of elements.

Clearly (b) implies (a) because $\mathrm{Fs}\,\tilde{\boldsymbol{p}} = \bigcup \mathrm{Qs}\,\tilde{\boldsymbol{p}}$ and $\mathrm{Fs}\,\tilde{\boldsymbol{q}} = \bigcup \mathrm{Qs}\,\tilde{\boldsymbol{q}}$ which are applications of (5.6.E).

Now assume (a). Any $\mathrm{Qt}\,\boldsymbol{q}_i$ is a partition of $\mathsf{Y}$ and hence a maximal set of pairwise disjoint element of $\mathrm{Fs}\,\tilde{\boldsymbol{p}}$, that is, a codirection of $\boldsymbol{p}$. Thus there is a function

$\boldsymbol{s} : l \longrightarrow l$ such that

$$\mathrm{Qt}\,\boldsymbol{q}_i = \mathrm{Qt}\,\boldsymbol{p}_{i\boldsymbol{s}} \qquad \forall i \in l.$$

In particular $\mathrm{Qt}\,\boldsymbol{q}_i$ has the same order of $\mathsf{X}$, say $m$, and $\mathrm{Qs}\,\tilde{\boldsymbol{q}} \subseteq \mathrm{Qs}\,\tilde{\boldsymbol{p}}$. Next argument shows that in fact, they are equal, hence proving (b). Since $\bigcup_{i\in l} \mathrm{Qt}\,\boldsymbol{q}_i = \mathrm{Fs}\,\tilde{\boldsymbol{q}} = \mathrm{Fs}\,\tilde{\boldsymbol{p}}$ has $ml$ elements, the $\mathrm{Qt}\,\boldsymbol{q}_i$'s are all different and their number is $l$. Thus $\mathrm{Qs}\,\tilde{\boldsymbol{q}} = \mathrm{Qs}\,\tilde{\boldsymbol{p}}$ and so $\boldsymbol{s} \in \mathrm{S}_l$. Define $\vec{\boldsymbol{x}} \in (\mathsf{X}^{\mathsf{X}})^l$ by

$$\mathrm{Fb}(\boldsymbol{q}_i, \mathsf{x}) = \mathrm{Fb}(\boldsymbol{p}_{i\boldsymbol{s}}, \mathsf{x}\boldsymbol{x}_i).$$

Since $\boldsymbol{q}_i$ has $m$ fibers, $\boldsymbol{x}_i$ is injective for each $i$ and hence $\vec{\boldsymbol{x}} \in (\mathrm{Sym}\,\mathsf{X})^l$. Then

$$\mathrm{Fb}(\boldsymbol{q}_i, \mathsf{x}) = \mathrm{Fb}(\boldsymbol{p}_{i\boldsymbol{s}}\boldsymbol{x}_i^{-1}, \mathsf{x}) \qquad \forall i \in l$$

which implies $\vec{\boldsymbol{q}} = \lfloor \boldsymbol{p} : \vec{\boldsymbol{x}} \rfloor \lfloor \boldsymbol{p} : \boldsymbol{s}^{-1} \rfloor \vec{\boldsymbol{p}}$ and hence (c).

This proves the equivalence and the rest follows from an argument similar to that at the end of the proof of 5.6.3 □

Assume now that $\boldsymbol{p} \in \mathrm{EPDec}(\mathsf{Y}, \mathsf{X}, l)$. Set $\mathsf{H} = \mathrm{Fs}\,\tilde{\boldsymbol{p}}$ and say

$$\boldsymbol{f} : \mathsf{X} \times l \longrightarrow \mathsf{H}; \qquad (\mathsf{x}, i) \mapsto \mathrm{Fb}(\boldsymbol{p}_i, \mathsf{x}).$$

This $\boldsymbol{f}$ is bijective because it is a surjective function between two sets with the same order. Since $\lceil \boldsymbol{p} : (\mathrm{Sym}\,\mathsf{X}) \,\mathrm{wr}\, \mathrm{S}_l \rceil$ is the stabilizer in $\mathrm{Sym}\,\mathsf{H}$ of the non trivial equipartition $\mathrm{Qs}\,\tilde{\boldsymbol{p}}$, one should be able to view the group as an inner wreath product on coproduct decomposition. In fact,

5.6.7. Theorem. *The map* $\boldsymbol{f}$ *so defined lies in* $\mathrm{ECDec}(\mathsf{H}, \mathsf{X}, l)$*, and the action of* $\mathrm{Sym}\,\mathsf{Y}$ *on subsets of* $\mathsf{Y}$ *induces a faithful action* $\varphi$ *of* $\lceil \boldsymbol{p} : (\mathrm{Sym}\,\mathsf{X}) \,\mathrm{wr}\, \mathrm{S}_l \rceil$ *on* $\mathsf{H}$ *such that for each* $w \in (\mathrm{Sym}\,\mathsf{X}) \,\mathrm{wr}\, \mathrm{S}_l$

$$\lceil \boldsymbol{p} : w \rceil \varphi = \lfloor \boldsymbol{f} : w \rfloor$$

*In particular,* $\varphi$ *identifies* $\lceil \boldsymbol{p} : (\mathrm{Sym}\,\mathsf{X}) \,\mathrm{wr}\, \mathrm{S}_l \rceil$ *with* $\lfloor \boldsymbol{f} : (\mathrm{Sym}\,\mathsf{X}) \,\mathrm{wr}\, \mathrm{S}_l \rfloor$.

Proof. It is enough to show that for each $\boldsymbol{\beta} \in \mathrm{S}_l$, $\vec{\boldsymbol{x}} \in (\mathrm{Sym}\,\mathsf{X})^l$ we have

$$(*) \qquad \lceil \boldsymbol{p} : \boldsymbol{\beta} \rceil \varphi = \lfloor \boldsymbol{f} : \boldsymbol{\beta} \rfloor \quad \text{and} \quad \lceil \boldsymbol{p} : \vec{\boldsymbol{x}} \rceil \varphi = \lfloor \boldsymbol{f} : \vec{\boldsymbol{x}} \rfloor.$$

This is done using (5.6.C), (5.6.A), (5.6.D) and (5.6.B) as follows: for each $i \in l$, (5.6.C) says that $\lceil \boldsymbol{p} : \boldsymbol{\beta}^{-1} \rceil \boldsymbol{p}_i = \boldsymbol{p}_{i\boldsymbol{\beta}}$. Then for each $\mathsf{x} \in \mathsf{X}$

$$\mathsf{x}\boldsymbol{f}_i(\lceil \boldsymbol{p} : \boldsymbol{\beta} \rceil \varphi) = \mathrm{Fb}(\boldsymbol{p}_i, \mathsf{x}) \lceil \boldsymbol{p} : \boldsymbol{\beta} \rceil = \mathrm{Fb}(\lceil \boldsymbol{p} : \boldsymbol{\beta}^{-1} \rceil \boldsymbol{p}_i, \mathsf{x}) =$$
$$= \mathrm{Fb}(\boldsymbol{p}_{i\boldsymbol{\beta}}, \mathsf{x}) = \mathsf{x}\boldsymbol{f}_{i\boldsymbol{\beta}}.$$

Therefore $\lceil \boldsymbol{p} : \boldsymbol{\beta} \rceil \varphi = \lfloor \boldsymbol{f} : \boldsymbol{\beta} \rfloor$ follows from (5.6.A). This proves the first equality of $(*)$ and the second is proved in a similar way, using (5.6.D) and (5.6.B) in place of (5.6.C) and (5.6.A). □

## 5.7. On the base subgroup

It is fairly easy to show that if $A$ is a transitive subgroup of $\operatorname{Sym}\mathsf{X}$, then the even part of the base subgroup of an inner wreath product $\lfloor \boldsymbol{f} : A \operatorname{wr} \mathrm{S}_l \rfloor$ is intransitive and the set of its orbits equals to $\operatorname{Is}\tilde{\boldsymbol{f}}$. In other words, the base subgroup of an inner imprimitive wreath product (or maybe even something a little bit smaller) determines the fundamental part of the related equicoproduct decomposition.

Dually, when $A$ is primitive (but not regular!), the base subgroup of an inner wreath product $\lceil \boldsymbol{p} : A \operatorname{wr} \mathrm{S}_l \rceil$ is transitive imprimitive and the set of its imprimitivity systems (having order equal to the order of $\mathsf{X}$) is equal to $\operatorname{Qs}\tilde{\boldsymbol{p}}$, the set of the codirections of $\boldsymbol{p}$. However, the duality breaks down and exceptions arise when we consider only the even part of the base subgroup. A detailed statement with proof follows.

5.7.1. Theorem. *Let $\boldsymbol{p} \in \operatorname{EPDec}(\mathsf{Y}, m, l)$. If $A$ is a primitive but not regular subgroup of $\mathrm{S}_m$, then $\lceil \boldsymbol{p} : A^l \rceil^{\mathrm{e}}$ is transitive imprimitive on $\mathsf{Y}$ with many imprimitivity systems. However, it has precisely $l$ imprimitivity systems of order $m$, namely the $l$ codirections of $\boldsymbol{p}$, except when $m$ is a prime $\equiv 3 \mod 4$, $A$ is the dihedral group $\mathrm{D}_{2m}$, and $l = 2$ (in which case the set of the imprimitivity systems of order $m$ of $\lceil \boldsymbol{p} : A^l \rceil^{\mathrm{e}}$ has order $m+1$ and contains the $l = 2$ codirections).*

Proof. We may assume that $\mathsf{Y} = m^l$ and that $\boldsymbol{p}$ is the identity. To save some writing, put $L := \lceil A^l \rceil^{\mathrm{e}}$. In order to prove that $L$ is transitive, it is enough to show that $A^{\mathrm{e}}$ is transitive. This follows from 2.5.2 because $A$ is non regular and so $m > 2$.

The part regarding the imprimitivity systems preserved by $L$ is less obvious. Pick a point $\mathsf{x} \in m$ and put $\mathsf{y} = (\mathsf{x}, \dots, \mathsf{x}) \in m^l$. By 2.4.2 there is a bijection $\boldsymbol{\varphi}$ taking the set of the imprimitivity systems of order $m$ of $L$ to the set of the subgroups of $L$ of index $m$ containing $L_{\mathsf{y}}$. It takes no effort to check that $L_{\mathsf{y}} = \lceil (A_{\mathsf{x}})^l \rceil^{\mathrm{e}}$ and that $(\operatorname{Qs}\tilde{\boldsymbol{p}})\boldsymbol{\varphi} = \{\, T_i \mid i \in l \,\}$, where

$$T_i = \langle\, A_{\mathsf{x}}\boldsymbol{\kappa}_i, A\boldsymbol{\kappa}_j \mid i \neq j \in l \,\rangle^{\mathrm{e}}.$$

Therefore we show that the subgroups of $L$ containing $L_{\mathsf{y}}$ are precisely the $T_i$. It is a good idea to rule out the exceptions first.

Suppose that $A^{\mathrm{e}}$ is regular. Then $A^{\mathrm{e}}$ has index 2, order at least 3, and it is a minimal normal subgroup of $A$ (because non trivial normal subgroups of $A$ are transitive). In fact, $A^{\mathrm{e}}$ is the unique minimal normal subgroup of $A$, so it is the socle of $A$. Moreover, $A_{\mathsf{x}}$ has order 2, so $A_{\mathsf{x}}$ is a solvable group. By [**DM96**, Th. 4.7B(iii)] the primitive groups with non abelian regular socle have non solvable point stabilizers. Therefore $A^{\mathrm{e}}$ is elementary abelian of order $p^d$ for some prime $p$; $A$ is a primitive subgroup of an affine group; $A_{\mathsf{x}}$ is an irreducible linear group of order 2. This forces $d = 1$ and $m = p$ (which then is odd), so $A \cong \mathrm{D}_{2m}$. As $A > A^{\mathrm{e}}$, we must have $m \equiv 3 \mod 4$. Now assume that $L_{\mathsf{y}} < T < L$ with $T$ of index $m$ in $L$. As $G := \lceil (A^{\mathrm{e}})^l \rceil$ is transitive, $L = L_{\mathsf{y}} G$, so $L = TG$ and hence $T \cap G$ has index $m$ in $G$. Since $m$ is prime $T \cap G$ is maximal in $G$, moreover, $T \cap G$ is normal both in $T$ and in $G$ which is abelian, thus it is normal in $L$. By Dedekind's law

$$T = T \cap L = T \cap (L_{\mathsf{y}} G) = L_{\mathsf{y}}(T \cap G).$$

Therefore determining the $T$ such that $L_{\mathsf{y}} < T < L$, amounts to determine the $L$-invariant maximal subgroups of the elementary abelian $p$-group $G$. We leave it to the reader to verify that they are exactly the $l$ yielding the $T_i$ when $l > 2$ and they are $m + 1$, 2 of which yielding the $T_i$, when $l = 2$.

We may now assume that $A^{\mathrm{e}}$ is non regular, so $(A^{\mathrm{e}})_{\mathsf{x}} > \{1\}$. As $(A^{\mathrm{e}})_{\mathsf{x}} = (A_{\mathsf{x}})^{\mathrm{e}}$, for the reminder of this proof we simply write $A^{\mathrm{e}}_{\mathsf{x}}$. We shall need that the normal closure $N$ of $A^{\mathrm{e}}_{\mathsf{x}}$ in $A$ is $A^{\mathrm{e}}$, yielding that a transitive subgroup of $A$ containing $A^{\mathrm{e}}_{\mathsf{x}}$ contains $A^{\mathrm{e}}$ as well. Indeed, as a non trivial normal subgroup in a primitive group, $N$ is transitive, so $A_{\mathsf{x}} N = A$, and then

$$\bigl|A : N\bigr| = \bigl|A_{\mathsf{x}} : A_{\mathsf{x}} \cap N\bigr| \leqslant \bigl|A_{\mathsf{x}} : A^{\mathrm{e}}_{\mathsf{x}}\bigr| \leqslant 2;$$

but by its definition, $N$ lies in $A^{\mathrm{e}}$. As above, say $L_{\mathsf{y}} < T < L$ with $T$ of index $m$ in $L$, we need to show that $T = T_i$ for some $i$. For each $i \in l$ the projection $T\boldsymbol{\rho}_i$ is either $A_x$ or $A$ because by primitivity $A_x \lessdot A$ and $L_{\mathsf{y}}\boldsymbol{\rho}_i = A_x$. The last equality holds because $(a, a, 1, \ldots, 1) \in L_{\mathsf{y}}$ for example, even if $a \in A_x$ is odd. Of course,

$$T \leqslant \bigl\langle\, T\boldsymbol{\rho}_i\boldsymbol{\kappa}_i \;\bigm|\; i \in l \,\bigr\rangle^{\mathrm{e}}.$$

If $T\boldsymbol{\rho}_i = A_x$ for at least one $i$, say $j$, then $T = T_j$ because the first is contained in the second and they have both the same index in $L$. Otherwise, $T$ is a subcartesian subgroup of $\lceil A^l \rceil$ and we aim to a contradiction. For each $i \in l$

$$\{1\} < A^{\mathrm{e}}_{\mathsf{x}}\boldsymbol{\kappa}_i \leqslant T \cap (A\boldsymbol{\kappa}_i) \trianglelefteq A\boldsymbol{\kappa}_i.$$

Since $A$ is primitive, $T \cap (A\boldsymbol{\kappa}_i)$ corresponds to a transitive subgroup of $A$ containing $A^{\mathrm{e}}_{\mathsf{x}}$, thus containing $N$ which we showed above is equal to $A^{\mathrm{e}}$. It follows that

$$T \geqslant \bigl\langle \lceil (A^{\mathrm{e}})^l \rceil, L_{\mathsf{y}} \bigr\rangle.$$

Again, $A^{\mathrm{e}}$ is transitive so that $A = A^{\mathrm{e}} A_x$. Then each element of $L$ is of the form $\vec{a}\vec{b}$ with $\vec{a} \in \lceil (A^{\mathrm{e}})^l \rceil$, $\vec{b} \in L_{\mathsf{y}}$ ($\vec{b}$ must be even because $\vec{a}\vec{b}$ and $\vec{a}$ are even). This shows that $L = \bigl\langle \lceil (A^{\mathrm{e}})^l \rceil, L_{\mathsf{y}} \bigr\rangle$, so $T = L$ against our assumptions. This completes the proof. □

## 5.8. Associative Law for wreath products

The canonical identification of $(\mathsf{A}\times\mathsf{B})\times\mathsf{C}$ with $\mathsf{A}\times(\mathsf{B}\times\mathsf{C})$ leads to the associative law (up to permutational isomorphisms)

$$\lfloor\lfloor \mathrm{S}_a \operatorname{wr} \mathrm{S}_b\rfloor \operatorname{wr} \mathrm{S}_c\rfloor \cong \lfloor \mathrm{S}_a \operatorname{wr}\lfloor \mathrm{S}_b \operatorname{wr} \mathrm{S}_c\rfloor\rfloor.$$

Similarly, the canonical identification of $(\mathsf{A}^{\mathsf{B}})^{\mathsf{C}}$ with $\mathsf{A}^{\mathsf{B}\times\mathsf{C}}$ leads to the associative law (up to permutational isomorphisms)

$$\lceil\lceil \mathrm{S}_a \operatorname{wr} \mathrm{S}_b\rceil \operatorname{wr} \mathrm{S}_c\rceil \cong \lceil \mathrm{S}_a \operatorname{wr}\lfloor \mathrm{S}_b \operatorname{wr} \mathrm{S}_c\rfloor\rceil.$$

We have the language to show this in great detail. However, as it is a language mostly composed of cumbersome terminology, we can only expect fairly tedious proofs.

5.8.1. Proposition. *Let $a, b, c$ be integers greater than* 1. *Let*

$$\boldsymbol{\zeta} : b \times c \longrightarrow d = bc$$

*be any bijection and define*

$$\begin{aligned} \boldsymbol{\varphi} &: (a \times b) \times c \longrightarrow a \times d; & ((h,i),j) &\mapsto (h,(i,j)\boldsymbol{\zeta}) \\ \boldsymbol{\psi} &: a^d \longrightarrow (a^b)^c; & \boldsymbol{f} &\mapsto (j \mapsto (i \mapsto (i,j)\boldsymbol{\zeta}\boldsymbol{f})). \end{aligned}$$

*Then*

$$\lfloor\boldsymbol{\varphi} : \lfloor \mathrm{S}_a \operatorname{wr} \mathrm{S}_b\rfloor \operatorname{wr} \mathrm{S}_c\rfloor = \lfloor \mathrm{S}_a \operatorname{wr}\lfloor\boldsymbol{\zeta} : \mathrm{S}_b \operatorname{wr} \mathrm{S}_c\rfloor\rfloor \tag{5.8.A}$$

*and*

$$\lceil\boldsymbol{\psi} : \lceil \mathrm{S}_a \operatorname{wr} \mathrm{S}_b\rceil \operatorname{wr} \mathrm{S}_c\rceil = \lceil \mathrm{S}_a \operatorname{wr}\lfloor\boldsymbol{\zeta} : \mathrm{S}_b \operatorname{wr} \mathrm{S}_c\rfloor\rceil. \tag{5.8.B}$$

Proof. We only prove (5.8.B) here. To avoid confusion, let us make the convention that for this proof we always have $h \in a$, $i \in b$, $j \in c$, $k \in d$. Call

$$\begin{aligned} \boldsymbol{\pi}^1_k &: a^d \longrightarrow a, & k &\in d, \\ \boldsymbol{\pi}^2_j &: (a^b)^c \longrightarrow a^b, & j &\in c, \\ \boldsymbol{\pi}^3_i &: a^b \longrightarrow a, & i &\in b \end{aligned}$$

the canonical projections (5.5.L), and

$$\begin{aligned} \boldsymbol{\kappa}^1_k &: \mathrm{S}_a \longrightarrow \mathrm{S}_a \operatorname{wr} \mathrm{S}_d, & k &\in d, \\ \boldsymbol{\kappa}^2_j &: \mathrm{S}_{a^b} \longrightarrow \mathrm{S}_{a^b} \operatorname{wr} \mathrm{S}_c, & j &\in c, \\ \boldsymbol{\kappa}^3_i &: \mathrm{S}_a \longrightarrow \mathrm{S}_a \operatorname{wr} \mathrm{S}_b, & i &\in b \\ \boldsymbol{\kappa}^4_j &: \mathrm{S}_b \longrightarrow \mathrm{S}_b \operatorname{wr} \mathrm{S}_c, & j &\in c \end{aligned}$$

the canonical inclusions (5.5.C). Because of our convention, we may sometime omit the superscript.

Note that the $\boldsymbol{\pi}^1_k$, $\boldsymbol{\pi}^2_j$, $\boldsymbol{\pi}^3_i$ satisfy

$$\boldsymbol{\pi}^1_{(i,j)\zeta} = \boldsymbol{\psi}\, \boldsymbol{\pi}^2_j\, \boldsymbol{\pi}^3_i \qquad \forall\, i \in b,\ j \in c.$$

Equation (5.8.B) is proved by computing $\lceil \boldsymbol{\psi} : \boldsymbol{s} \rceil$ for suitable $\boldsymbol{s}$ which together generate $\lceil \lceil \mathrm{S}_a \operatorname{wr} \mathrm{S}_b \rceil \operatorname{wr} \mathrm{S}_c \rceil$. This is done by computing first their composition with a generic $\boldsymbol{\pi}^1_k$, say $k = (u, v)\boldsymbol{\zeta}$ for example.

Start with $\boldsymbol{s} = \boldsymbol{a}\boldsymbol{\kappa}^3_i \boldsymbol{\kappa}^2_j$ where $\boldsymbol{a} \in \mathrm{S}_a$.

$$\boldsymbol{\psi}\lceil (a^b)^c : \boldsymbol{s} \rceil \boldsymbol{\psi}^{-1} \boldsymbol{\pi}^1_k = \boldsymbol{\psi} \lceil \boldsymbol{s} \rceil \boldsymbol{\pi}^2_v \boldsymbol{\pi}^3_u = \boldsymbol{\psi} \lceil (\boldsymbol{a}\boldsymbol{\kappa}^3_i)\boldsymbol{\kappa}^2_j \rceil \boldsymbol{\pi}^2_v \boldsymbol{\pi}^3_u.$$

Now, by (5.6.D)

$$\lceil (\boldsymbol{a}\boldsymbol{\kappa}^3_i)\boldsymbol{\kappa}^2_j \rceil \boldsymbol{\pi}^2_v = \begin{cases} \boldsymbol{\pi}^2_v \lceil a^b : \boldsymbol{a}\boldsymbol{\kappa}^3_i \rceil & \text{if } j = v, \\ \boldsymbol{\pi}^2_v & \text{otherwise.} \end{cases}$$

Similarly,

$$\lceil \boldsymbol{a}\boldsymbol{\kappa}^3_i \rceil \boldsymbol{\pi}^3_u = \begin{cases} \boldsymbol{\pi}^3_u \boldsymbol{a} & \text{if } i = u, \\ \boldsymbol{\pi}^3_u & \text{otherwise.} \end{cases}$$

Whence,

$$(\boldsymbol{\psi}\lceil (a^b)^c : \boldsymbol{s} \rceil \boldsymbol{\psi}^{-1}) \boldsymbol{\pi}^1_k = \begin{cases} \boldsymbol{\pi}^1_k \boldsymbol{a} & \text{if } k = (i,j)\boldsymbol{\zeta}, \\ \boldsymbol{\pi}^1_k & \text{otherwise.} \end{cases} \tag{5.8.C}$$

This together with (5.6.D) imply

$$\lceil \boldsymbol{\psi} : \boldsymbol{a}\boldsymbol{\kappa}_i \boldsymbol{\kappa}_j \rceil = \lceil a^d : \boldsymbol{a}\boldsymbol{\kappa}_{(i,j)\zeta} \rceil. \tag{5.8.D}$$

In a similar way if $\boldsymbol{b} \in \mathrm{S}_b$, then

$$(\boldsymbol{\psi}\lceil (a^b)^c : \boldsymbol{b}\boldsymbol{\kappa}^2_j \rceil \boldsymbol{\psi}^{-1}) \boldsymbol{\pi}^1_k = \begin{cases} \boldsymbol{\pi}^1_{(u\boldsymbol{b}^{-1}, v)\zeta} & \text{if } j = v, \\ \boldsymbol{\pi}^1_k & \text{otherwise.} \end{cases} \tag{5.8.E}$$

On the other hand, by (5.5.F)

$$\lfloor \boldsymbol{\zeta} : \boldsymbol{b}^{-1}\boldsymbol{\kappa}^4_j \rfloor : k \mapsto \begin{cases} (u\boldsymbol{b}^{-1}, v)\boldsymbol{\zeta} & \text{if } j = v, \\ k & \text{otherwise.} \end{cases}$$

so that

$$\lceil a^d : \lfloor \boldsymbol{\zeta} : \boldsymbol{b}\boldsymbol{\kappa}^4_j \rfloor \rceil \boldsymbol{\pi}^1_k = \begin{cases} \boldsymbol{\pi}^1_{(u\boldsymbol{b}^{-1}, v)\zeta} & \text{if } j = v, \\ \boldsymbol{\pi}^1_k & \text{otherwise.} \end{cases} \tag{5.8.F}$$

Then (5.8.E), (5.8.F) and (5.6.C) imply

$$\lceil \boldsymbol{\psi} : \boldsymbol{b}\boldsymbol{\kappa}_j^2 \rceil = \lceil a^d : \lfloor \boldsymbol{\zeta} : \boldsymbol{b}\boldsymbol{\kappa}_j^4 \rfloor \rceil. \tag{5.8.G}$$

Finally, if $\boldsymbol{c} \in \mathrm{S}_c$,

$$(\boldsymbol{\psi} \lceil (a^b)^c : \boldsymbol{c} \rceil \boldsymbol{\psi}^{-1}) \boldsymbol{\pi}_k^1 = \boldsymbol{\pi}^1_{(u, v\boldsymbol{c}^{-1})\boldsymbol{\zeta}}, \tag{5.8.H}$$

On the other hand, by (5.5.E)

$$\lfloor \boldsymbol{\zeta} : \boldsymbol{c}^{-1} \rfloor : k \mapsto (u, v\boldsymbol{c}^{-1})\boldsymbol{\zeta},$$

so that

$$\lceil a^d : \lfloor \boldsymbol{\zeta} : \boldsymbol{c} \rfloor \rceil \boldsymbol{\pi}_k^1 = \boldsymbol{\pi}^1_{k \lfloor \boldsymbol{\zeta} : \boldsymbol{c}^{-1} \rfloor} = \boldsymbol{\pi}^1_{(u, v\boldsymbol{c}^{-1})\boldsymbol{\zeta}}. \tag{5.8.I}$$

Then (5.8.H), (5.8.I) and (5.6.C) imply

$$\lceil \boldsymbol{\psi} : \boldsymbol{c} \rceil = \lceil a^d : \lfloor \boldsymbol{\zeta} : \boldsymbol{c} \rfloor \rceil.$$

Which together with (5.8.D) and (5.8.G) imply (5.8.B). □

An application of the proof of the the associative law (5.8.B) follows. For $a \geqslant 5$ and $b, c > 1$, say $H = \lceil \lceil \mathrm{S}_a \operatorname{wr} \mathrm{S}_b \rceil \operatorname{wr} \mathrm{S}_c \rceil$ and call $\boldsymbol{\sigma}$ the action of $H$ by conjugation on the set $\mathsf{E}$ of the minimal normal subgroups of its socle $\operatorname{Soc} H$. As expected ($\mathsf{E}$ being in bijection with the set of coordinate factors of the base subgroup of $\mathrm{S}_a \operatorname{wr} \mathrm{S}_{bc}$), we show that $\boldsymbol{\sigma}$ is equivalent to the action induced by restriction to $\lceil \boldsymbol{\psi} : H \rceil$ of the top projection $\boldsymbol{\rho} : \mathrm{S}_a \operatorname{wr} \mathrm{S}_{bc} \longrightarrow \mathrm{S}_{bc}$.

5.8.2. Proposition. *$H\boldsymbol{\sigma}$ on $\mathsf{E}$ and $\lceil \boldsymbol{\psi} : H \rceil \boldsymbol{\rho}$ on $bc$ are permutationally isomorphic. In particular, $H^e$ (not to confuse with $(H\boldsymbol{\sigma})^e$) has a unique imprimitivity system on $\mathsf{E}$ with respect to $\boldsymbol{\sigma}$, whose blocks are the*

$$\mathsf{E}_j = \{ \mathrm{A}_a \boldsymbol{\kappa}_i^3 \boldsymbol{\kappa}_j^2 \mid i \in b \}, \qquad j \in c.$$

*In fact, up to the order of the minimal normal subgroups of $(\mathrm{A}_{a^b})^c$, one may assume that for each $j \in c$*

$$\mathsf{E}_j = \{ (\sqrt[bc]{H})_k \mid (\sqrt[bc]{H})_k \leqslant (\sqrt[c]{\lceil \mathrm{S}_{a^b} \operatorname{wr} \mathrm{S}_c \rceil})_j \}.$$

Proof. Call $\boldsymbol{\sigma}'$ the action by conjugation of $\lceil \boldsymbol{\psi} : H \rceil$ on the set $\mathsf{E}'$ of the minimal normal subgroups of its socle. Of course, $H\boldsymbol{\sigma}$ on $\mathsf{E}$ and $\lceil \boldsymbol{\psi} : H \rceil$ on $\mathsf{E}'$ are permutationally isomorphic. Note that

$$\mathsf{E}' = \{ \mathrm{A}_a \boldsymbol{\kappa}_k^1 \mid k \in bc \}.$$

Then (5.5.D) shows that $\lceil \boldsymbol{\psi} : H \rceil \boldsymbol{\sigma}'$ on $\mathsf{E}'$ and $\lceil \boldsymbol{\psi} : H \rceil \boldsymbol{\rho}$ on $bc$ are permutationally isomorphic, the isomorphism being induced by

$$\mathrm{A}_a \boldsymbol{\kappa}_k^1 \longleftrightarrow k, \qquad k \in bc.$$

By the associative law (5.8.B), $\lceil \boldsymbol{\psi} : H \rceil \boldsymbol{\rho} = \lfloor \boldsymbol{\zeta} : \mathrm{S}_b \operatorname{wr} \mathrm{S}_c \rfloor$. This, together with 2.4.6, shows that at least $H\boldsymbol{\sigma}$ has a unique imprimitivity system on $\mathsf{E}$.

Now, if the top subgroup $\lceil a^{bc} : \mathrm{S}_{bc} \rceil$ is even then $\lceil \boldsymbol{\psi} : H^{\mathrm{e}} \rceil = \lceil \boldsymbol{\psi} : H \rceil^{\mathrm{e}}$ contains $\lceil a^{bc} : \lfloor \boldsymbol{\zeta} : \mathrm{S}_b \operatorname{wr} \mathrm{S}_c \rfloor \rceil$. Otherwise, $a$ is odd by (5.5.J) and therefore, by (5.5.K), the base subgroup $\lceil a^{bc} : {\mathrm{S}_a}^{bc} \rceil$ is odd. It follows that the image by $\boldsymbol{\rho}$ of $\lceil \boldsymbol{\psi} : H^{\mathrm{e}} \rceil$ is again $\lceil a^{bc} : \lfloor \boldsymbol{\zeta} : \mathrm{S}_b \operatorname{wr} \mathrm{S}_c \rfloor \rceil$. Therefore $H^{\mathrm{e}}\boldsymbol{\sigma} = H\boldsymbol{\sigma}$ even though this is not necessarily equal to $(H\boldsymbol{\sigma})^{\mathrm{e}}$.

To finish, note that the blocks of the imprimitivity system on $\mathsf{E}'$ are the

$$\mathsf{E}'_j = \{\ \mathrm{A}_a \boldsymbol{\kappa}^1_{(i,j)\zeta} \ \mid\ i \in b\ \}, \qquad\qquad j \in c$$

and that by (5.8.D), $\boldsymbol{\psi}^{-1}(\mathrm{A}_a \boldsymbol{\kappa}^1_{(i,j)\zeta})\boldsymbol{\psi} = \mathrm{A}_a \boldsymbol{\kappa}_i^3 \boldsymbol{\kappa}_j^2$. Therefore $\mathsf{E}_j$, which is equal to $\boldsymbol{\psi}^{-1}\mathsf{E}'_j\boldsymbol{\psi}$, is a block for $H\boldsymbol{\sigma}$ as claimed. □

## 5.9. Inner Blow-ups

In this section we refer to the blow-up construction introduced by Kovács in [**Kov89a**]. Given a primitive group $A$ with socle $M$, an integer $l > 1$ and a large subgroup $P$ of the wreath product $(A/M) \operatorname{wr} \mathrm{S}_l$, the blow-up of $A$ by $P$, $G = A \uparrow P$, is defined as the complete inverse image of $P$ under the natural homomorphism of $A \operatorname{wr} \mathrm{S}_l$ onto $(A/M) \operatorname{wr} \mathrm{S}_l$. Sometime there is no need to specify the second argument $P$, and so we just write $G = A \uparrow \star_l$ where the subscript is there to remind that the $\star$ plays the role of a large subgroup of $(A/M) \operatorname{wr} \mathrm{S}_l$.

We usually regard $A$ as a subgroup of some $\mathrm{S}_m$, and so we also regard a blow-up $A \uparrow \star_l$ as a subgroup of $W = \mathrm{S}_m \operatorname{wr} \mathrm{S}_l$. We observe that if $W_0 = \mathrm{St}_W\, 0$ with respect to the action induced by the top projection $\boldsymbol{\rho}$ of $W$, then a blow-up $G$ in $W$ must be a subgroup of $(G \cap W_0)\boldsymbol{\rho}_0 \operatorname{wr} \mathrm{S}_l$, in fact, a large subgroup.

The main problem with blow-ups is that the set of blow-ups contained in $W$, that is, of the subgroups of $W$ which are blow-ups, is not closed under conjugation in $W$. For example, say $A$ a subgroup of order 5 of $\mathrm{S}_5$ and $w \in \mathrm{S}_5 \operatorname{wr} \mathrm{S}_2$ such that $w\boldsymbol{\rho} = 1$, $w\boldsymbol{\rho}_0 = 1$ and $w\boldsymbol{\rho}_1 \notin \mathrm{N}(A)$. Then $G = w^{-1}(A \operatorname{wr} \mathrm{S}_2)w$ is not a blow-up of $W$ because $(G \cap W_0)\boldsymbol{\rho}_0 = A$ but $G \nleq A \operatorname{wr} \mathrm{S}_2$.

As a consequence almost anything involving blow-ups makes heavy use of permutational isomorphisms (which unfortunately are usually kept implicit). Then we decided to equip blow-ups with explicit permutational isomorphisms and to call the images of such isomorphisms inner blow-ups (not to confuse with internal blow-up decompositions of permutation groups as introduced in [**Kov89a**]). This yields classes of subgroups closed under conjugation and allows the encoding of more information in statements involving blow-ups. We also lose some kind of uniqueness though because an inner blow-up is image of many different permutational isomorphisms. However, this will not interfere with our specific needs.

Blow-ups were originally thought of as subgroups of a wreath product $W$ in product action on some $\mathsf{X}^l$. We have inner wreath products on product decompositions of a set $\mathsf{Y}$ as images of $W$ via permutational isomorphisms:

$$\lceil \boldsymbol{p} : W \rceil = \boldsymbol{p} W \boldsymbol{p}^{-1}, \qquad \boldsymbol{p} \in \mathrm{EPDec}(\mathsf{Y}, \mathsf{X}, l).$$

Coherently, for $A \leqslant \mathrm{Sym}\,\mathsf{X}$, $P$ large subgroup of $(A/\operatorname{Soc} A) \operatorname{wr} \mathrm{S}_l$, we define the inner blow-up of $A$ by $P$ related to $\boldsymbol{p}$ as

$$\lceil \boldsymbol{p} : A \uparrow P \rceil = \boldsymbol{p}(A \uparrow P)\boldsymbol{p}^{-1}.$$

Now both the set of blow-ups contained in $\mathrm{Sym}\,\mathsf{Y}$ and the set of blow-ups contained in $\lceil \boldsymbol{p} : W \rceil$ are closed under conjugation because

$$\boldsymbol{y}^{-1} \lceil \boldsymbol{p} : A \uparrow P \rceil \boldsymbol{y} = \lceil \boldsymbol{y}^{-1}\boldsymbol{p} : A \uparrow P \rceil \qquad \forall \boldsymbol{y} \in \mathrm{Sym}\,\mathsf{Y}.$$

We are especially interested in the inner blow-ups which are second maximal in their parent symmetric or alternating group.

5.9.1. LEMMA. *Let $\boldsymbol{p} \in \mathrm{EPDec}(\mathsf{Y}, m, l)$, let $A$ be a subgroup of $\mathrm{S}_m$ and let $H = \lceil \boldsymbol{p} : A \uparrow \star_l \rceil$ be an inner blow-up. If $H$ is second maximal in $\mathrm{Sym}\,\mathsf{Y}$ or $\mathrm{Alt}\,\mathsf{Y}$, then either $A = \mathrm{S}_m$ or $\mathrm{A}_m \neq A \lessdot \mathrm{S}_m$ and $H^e = \lceil \boldsymbol{p} : A \operatorname{wr} \mathrm{S}_l \rceil^e$.*

PROOF. We may assume that $\mathsf{Y} = \mathsf{X}^l$ with $\mathsf{X} = m$ and that $\boldsymbol{p}$ is the identity. To prove that $A$ is maximal when it is a proper subgroup of $\mathrm{S}_m$, assume $A < B \leqslant \mathrm{S}_m$. Of course

$$\left|\lceil B \operatorname{wr} \mathrm{S}_l \rceil^{\mathrm{e}}\right| \geqslant \frac{1}{2}\left|B\right|^l l! \qquad \text{and} \qquad \left|\lceil A \operatorname{wr} \mathrm{S}_l \rceil^{\mathrm{e}}\right| \leqslant \left|A\right|^l l!.$$

Thus the index of $\lceil A \operatorname{wr} \mathrm{S}_l \rceil^{\mathrm{e}}$ in $\lceil B \operatorname{wr} \mathrm{S}_l \rceil^{\mathrm{e}}$ is at least

$$\frac{1}{2}[B : A]^l \geqslant 2^{l-1} \geqslant 2$$

so that $\lceil A \operatorname{wr} \mathrm{S}_l \rceil^{\mathrm{e}}$ must be a proper subgroup of $\lceil B \operatorname{wr} \mathrm{S}_l \rceil^{\mathrm{e}}$. In particular

$$H^{\mathrm{e}} \leqslant \lceil A \operatorname{wr} \mathrm{S}_l \rceil^{\mathrm{e}} < \lceil B \operatorname{wr} \mathrm{S}_l \rceil^{\mathrm{e}} \leqslant \lceil \mathrm{S}_m \operatorname{wr} \mathrm{S}_l \rceil^{\mathrm{e}} < \mathrm{Alt}\,\mathsf{Y}$$

and of course

$$H \leqslant \lceil A \operatorname{wr} \mathrm{S}_l \rceil < \lceil B \operatorname{wr} \mathrm{S}_l \rceil \leqslant \lceil \mathrm{S}_m \operatorname{wr} \mathrm{S}_l \rceil < \mathrm{Sym}\,\mathsf{Y}.$$

Since we assume $H$ high, $B$ must be equal to $\mathrm{S}_m$ and $A$ is maximal. Furthermore, either $H$ is odd and equal to $\lceil A \operatorname{wr} \mathrm{S}_l \rceil$ or $H$ is even and equal to $\lceil A \operatorname{wr} \mathrm{S}_l \rceil^{\mathrm{e}}$. But in any case $H^{\mathrm{e}} = \lceil A \operatorname{wr} \mathrm{S}_l \rceil^{\mathrm{e}}$ as stated. We only need to prove that $A \neq \mathrm{A}_m$ but this comes from the fact that neither $\lceil \mathrm{A}_m \operatorname{wr} \mathrm{S}_l \rceil$ nor $\lceil \mathrm{A}_m \operatorname{wr} \mathrm{S}_l \rceil^{\mathrm{e}}$ is high. □

We are now in a position to take full advantage of the following which is a weak version of [**Pra90**, 7.1]. Recall that by $\star_l$ we mean a suitable large subgroup of some wreath product having $\mathrm{S}_l$ as top group. However, a $\star_l$ appearing twice or more in a same statement, is not to be interpreted as the identifier of a same group; we would use a more common identifier otherwise.

5.9.2. Proposition. *Let $H = \lceil \boldsymbol{p} : A \uparrow \star_l \rceil$ be an inner blow-up for some $\boldsymbol{p} \in \mathrm{EPDec}(\mathsf{Y}, m, l)$ and some almost simple primitive subgroup $A$ of $\mathrm{S}_m$. If $H \lessdot G \leqslant \operatorname{Sym} \mathsf{Y}$, then one of the following is true.*

(a) *$A \cong \mathrm{PSL}(2,7)$ of degree 8 and there is a transitive subgroup $L$ of $\mathrm{S}_l$ such that $H = \lceil \boldsymbol{p} : A \operatorname{wr} L \rceil$ and $G = \lceil \boldsymbol{p} : B \operatorname{wr} L \rceil$ where $A < B \cong \mathrm{AGL}_3(2)$ of degree 8. Note that $H$ is low.*

(b) *There is a primitive almost simple subgroup $B$ of $\mathrm{S}_m$ such that $A \leqslant B$ and $G = \lceil \boldsymbol{q} : B \uparrow \star_l \rceil$ for some $\boldsymbol{q} \in \mathrm{EPDec}(\mathsf{Y}, m, l)$.*

(c) *$l = ab > b$. There are $\boldsymbol{r}, \boldsymbol{q} \in \mathrm{EPDec}(\mathsf{Y}, m^a, b)$ and $\boldsymbol{p}' \in \mathrm{EPDec}(m^a, m, a)$ such that $H = \lceil \boldsymbol{r} : B \uparrow \star_b \rceil$ with $B = \lceil \boldsymbol{p}' : A \uparrow \star_a \rceil$ and $G = \lceil \boldsymbol{q} : C \uparrow \star_b \rceil$ with $\mathrm{A}_{m^a} \leqslant C \leqslant \mathrm{S}_{m^a}$.*

(d) *$m$ is a square and there is $\boldsymbol{q} \in \mathrm{EPDec}(\mathsf{Y}, m, l)$ and a proper subgroup $B$ of $\mathrm{S}_m$ such that $G = \lceil \boldsymbol{q} : B \uparrow \star_l \rceil$. Note that $H$ is low.*

In fact, there is much more in [**Pra90**, 7.1]. For example, in case (b), one expects that there is an intimate connection between $\boldsymbol{p}$ and $\boldsymbol{q}$. This connection is made evident in [**Pra90**] but here we content ourselves to cut and paste from it according to our needs. The next proposition is stated following the notation set in 2.1 where $\boldsymbol{H}$ is a non trivial high subgroup of $\boldsymbol{S} = \operatorname{Sym} \Omega$, that is, $\boldsymbol{H}$ is a second maximal subgroup of $\boldsymbol{U} = (\operatorname{Alt} \Omega)\boldsymbol{H}$. Recall that a primitive subgroup of $\operatorname{Sym} \mathsf{Y}$ of type PA is an inner blow-up $\lceil \boldsymbol{p} : A \uparrow \star_l \rceil$ for some $\boldsymbol{p} \in \mathrm{EPDec}(\mathsf{Y}, m, l)$ and almost simple primitive subgroup $A$ of $\mathrm{S}_m$.

5.9.3. Theorem. *Suppose that the high subgroup $\boldsymbol{H}$ of $\boldsymbol{S}$ is primitive of type* PA. *For example, say $\boldsymbol{H} = \lceil \boldsymbol{p} : A \uparrow \star_l \rceil$ with $\boldsymbol{p} \in \mathrm{EPDec}(\Omega, m, l)$ and $A$ primitive almost simple subgroup of $\mathrm{S}_m$. Suppose also that $\boldsymbol{H} < G < \boldsymbol{U}$. Then one or two of the following is true.*

(1) *$A = \mathrm{S}_m$ and $G = \lceil \boldsymbol{p} : \mathrm{S}_m \operatorname{wr} \mathrm{S}_l \rceil \cap \boldsymbol{U}$.*

(2) *$\mathrm{A}_m \neq A \lessdot \mathrm{S}_m$, $H = \lceil \boldsymbol{p} : A \operatorname{wr} \mathrm{S}_l \rceil \cap \boldsymbol{U}$, $G = \lceil \boldsymbol{p} : \mathrm{S}_m \operatorname{wr} \mathrm{S}_l \rceil \cap \boldsymbol{U}$.*

(3) $l = ab$ *for* $a, b > 1$. *There are* $\boldsymbol{r} \in \mathrm{EPDec}(\Omega, m^a, b)$, $\boldsymbol{p}' \in \mathrm{EPDec}(m^a, m, a)$ *such that* $H = \lceil \boldsymbol{r} : B \,\mathrm{wr}\, \mathrm{S}_b \rceil \cap \boldsymbol{U}$ *with* $B = \lceil \boldsymbol{p}' : \mathrm{S}_m \,\mathrm{wr}\, \mathrm{S}_a \rceil$. *In this case* $\mathrm{Qs}\, \tilde{\boldsymbol{r}}$ *is determined by* $\boldsymbol{H}$ *and so are* $b = \left|\mathrm{Qs}\, \tilde{\boldsymbol{r}}\right|$ *and* $G = \mathrm{St}_{\boldsymbol{U}} \mathrm{Qs}\, \tilde{\boldsymbol{r}}$.

*In particular* $[\boldsymbol{H} \div \boldsymbol{U}] = \mathcal{M}_k$ *with* $k \leqslant 2$.

Proof. This is really a corollary of 5.9.2 applied to $H = \boldsymbol{H}$ where the further assumption that $\boldsymbol{H}$ is high is used to eliminate those cases (a), (d) and to refine that (b) in these (1), (2) and that (c) in this (3). Also, note that the socle of a group $\boldsymbol{H}$ as in (2) cannot possibly be isomorphic to the socle of a group $\boldsymbol{H}$ as in (1), (3), thus proving the fact that at most two of the three cases hold.

Consider 5.9.2(b) first. $G = \lceil \boldsymbol{q} : \mathrm{S}_m \,\mathrm{wr}\, \mathrm{S}_l \rceil \cap \boldsymbol{U}$ because of the maximality of $G$ in $\boldsymbol{U}$. Also, being $\boldsymbol{H}$ high, 5.9.1 says that either $A = \mathrm{S}_m$ or $\mathrm{A}_m \neq A \lessdot \mathrm{S}_m$ and $\boldsymbol{H} = \lceil \boldsymbol{p} : A \,\mathrm{wr}\, \mathrm{S}_l \rceil \cap \boldsymbol{U}$. Thus, proving (1) and (2) amounts to show that really $G = \lceil \boldsymbol{p} : \mathrm{S}_m \,\mathrm{wr}\, \mathrm{S}_l \rceil \cap \boldsymbol{U}$. This is done by proving that the codirections of $\boldsymbol{p}$ are the codirections of $\boldsymbol{q}$, that is, $\mathrm{Qs}\, \tilde{\boldsymbol{p}} = \mathrm{Qs}\, \tilde{\boldsymbol{q}}$ because by 5.6.6 $\lceil \boldsymbol{p} : \mathrm{S}_m \,\mathrm{wr}\, \mathrm{S}_l \rceil = \mathrm{St}\, \mathrm{Qs}\, \tilde{\boldsymbol{p}}$ and $\lceil \boldsymbol{q} : \mathrm{S}_m \,\mathrm{wr}\, \mathrm{S}_l \rceil = \mathrm{St}\, \mathrm{Qs}\, \tilde{\boldsymbol{q}}$.

If $\mathrm{Soc}\, H = \mathrm{Soc}\, G$, then $\mathrm{Qs}\, \tilde{\boldsymbol{q}} = \mathrm{Qs}\, \tilde{\boldsymbol{p}}$ follows from 5.7.1. Otherwise, $A < \mathrm{S}_m$ and $\boldsymbol{H} = \lceil \boldsymbol{p} : A \,\mathrm{wr}\, \mathrm{S}_l \rceil \cap \boldsymbol{U}$. With notations 5.3 we may assume that $(\sqrt[l]{G})_i = \lceil \boldsymbol{q} : \mathrm{A}_m \boldsymbol{\kappa}_i \rceil$ for all $i$, where $\boldsymbol{\kappa}_i : \mathrm{A}_m \longrightarrow \mathrm{A}_m \,\mathrm{wr}\, \mathrm{S}_l$ is defined as in (5.5.C). By Lemmas 5.3.1 and 5.3.2, $\mathrm{Soc}\, H < \mathrm{Soc}\, G$ and we may assume

$$\big(\mathrm{Soc}\, H\big) \boldsymbol{\rho}_i^G = (\sqrt[l]{H})_i \leqslant (\sqrt[l]{G})_i \qquad \forall i \in l.$$

By 5.5.D,

$$\mathrm{N}_{\boldsymbol{H}}\, (\sqrt[l]{H})_i = \boldsymbol{H} \cap \mathrm{N}_G\, (\sqrt[l]{G})_i \qquad \forall\, i \in l$$

so that $\mathrm{Shoe}\, \boldsymbol{H} \leqslant \mathrm{Shoe}\, G$. Now, $A$ is primitive of type almost simple thus $\mathrm{Soc}\, A$ is non regular. This shows that $A^{\mathrm{e}}$, containing $\mathrm{Soc}\, A$, is non regular. By 5.7.1, $(\mathrm{Shoe}\, \boldsymbol{H})^{\mathrm{e}} = \lceil \boldsymbol{p} : A^l \rceil^{\mathrm{e}}$ has exactly $l$ imprimitivity systems of order $m$, the codirections of $\boldsymbol{p}$. These must be the $l$ imprimitivity systems of order $m$ of $\mathrm{Shoe}\, G = \lceil \boldsymbol{q} : (\mathrm{S}_m)^l \rceil$, which are indeed the codirections of $\boldsymbol{q}$, as claimed.

Next, assume 5.9.2(c). Since $\boldsymbol{H}$ is high,

$$\begin{aligned}
&B \lessdot \mathrm{S}_{m^a},\ \text{which implies } A = \mathrm{S}_m \text{ and } B = \lceil \boldsymbol{p}' : A \,\mathrm{wr}\, \mathrm{S}_a \rceil,\\
&\boldsymbol{H} = \lceil \boldsymbol{r} : B \,\mathrm{wr}\, \mathrm{S}_b \rceil \cap \boldsymbol{U},\\
&G = \lceil \boldsymbol{q} : \mathrm{S}_{m^a} \,\mathrm{wr}\, \mathrm{S}_b \rceil \cap \boldsymbol{U},\ \text{note that } G = \mathrm{St}_{\boldsymbol{U}}\, \mathrm{Qs}\, \tilde{\boldsymbol{q}}.
\end{aligned}$$

First, let us prove that the codirections of $\boldsymbol{r}$ are the codirections of $\boldsymbol{q}$. To save some writing, call $M = \mathrm{Soc}\, \boldsymbol{H}$ and $N = \mathrm{Soc}\, G$. Note that $M = \lceil \boldsymbol{p} : (\mathrm{A}_m)^l \rceil$ and $N = \lceil \boldsymbol{q} : (\mathrm{A}_{m^a})^b \rceil$. By Lemmas 5.3.1 and 5.3.2, there is an equipartition

$\mathsf{E} = \{\mathsf{E}_0, \dots, \mathsf{E}_{b-1}\}$ of $l$ such that for each $j \in b$

$$M\boldsymbol{\rho}_j^G = \big\langle\ (\sqrt[l]{H})_k \ \big|\ k \in \mathsf{E}_j\ \big\rangle.$$

Since $G$ preserves the set of the minimal normal subgroups of its socle, this equipartition is preserved by $\boldsymbol{H}$. But the action of $\boldsymbol{H}$ on $l$ is known (see 5.8.2), $H$ has only one imprimitivity system on $l$ and so one may write $\mathsf{E} = \big\{\ \mathsf{E}'_j \ \big|\ j \in b\ \big\}$ where

$$\mathsf{E}'_j = \big\{\ k \in l \ \big|\ (\sqrt[l]{H})_k \leqslant (\sqrt[b]{\lceil \boldsymbol{r} : \mathrm{S}_{m^a} \operatorname{wr} \mathrm{S}_b \rceil})_j\ \big\}, \qquad j \in b.$$

Of course, if $k \in \mathsf{E}'_j$ then $(\sqrt[l]{H})_k \boldsymbol{\rho}_j^G = (\sqrt[l]{H})_k$ thus $\mathsf{E}'_j \subseteq \mathsf{E}_j$ and hence $\mathsf{E}'_j = \mathsf{E}_j$ for all $j \in b$. Clearly,

$$\lceil \boldsymbol{r} : B^b \rceil^{\mathrm{e}} \leqslant \mathrm{St}_{\boldsymbol{H}}\, (\sqrt[b]{\lceil \boldsymbol{r} : \mathrm{S}_{m^a} \operatorname{wr} \mathrm{S}_b \rceil})_j = \mathrm{St}_{\boldsymbol{H}}\, \mathsf{E}'_j.$$

On the other hand, $\mathrm{St}_{\boldsymbol{H}}\, \mathsf{E}_j = \boldsymbol{H} \cap \mathrm{N}_G\, (\sqrt[b]{G})_j$ because $(\sqrt[b]{G})_j$ is the image of $\boldsymbol{\rho}_j^G$. Therefore

$$\lceil \boldsymbol{r} : B^b \rceil^{\mathrm{e}} \leqslant \operatorname{Shoe} G.$$

But $\operatorname{Shoe} G = \lceil \boldsymbol{q} : (\mathrm{S}_{m^a})^b \rceil \cap \boldsymbol{U}$ has exactly $b$ imprimitivity systems of order $m^a$ which are the codirections of $\boldsymbol{q}$. These are imprimitivity systems of order $m^a$ of $\lceil \boldsymbol{r} : B^b \rceil^{\mathrm{e}}$ which are the codirections of $\boldsymbol{r}$ as wanted.

With an entirely similar argument one shows that if $\boldsymbol{H}$ is also equal to some $\lceil \boldsymbol{r}' : B_1 \operatorname{wr} \mathrm{S}_{b_1} \rceil$, then $\operatorname{Qs} \tilde{\boldsymbol{r}'} = \operatorname{Qs} \tilde{\boldsymbol{q}} = \operatorname{Qs} \tilde{\boldsymbol{r}}$. Therefore $\boldsymbol{H}$ determines $\operatorname{Qs} \tilde{\boldsymbol{r}}$ and hence $b = \left|\operatorname{Qs} \tilde{\boldsymbol{r}}\right|$ and $G = \mathrm{St}_{\boldsymbol{U}} \operatorname{Qs} \tilde{\boldsymbol{r}}$. This completely proves our claims. □

## 5.10. The Holomorph

Given a group $G$, the *right regular action* of $G$ on itself is the homomorphism $g \mapsto /g/$ where $/g/ \in \operatorname{Sym} G$ is defined by

$$x/g/ := xg, \qquad x, g \in G.$$

The choice of diagonal bars is not arbitrary but was suggested to us while dealing with primitive groups of simple diagonal type, of which some holomorphs are a special case. Here the diagonal bars are leaning on the right reminding in this way that the action is induced by a multiplication on the right. Accordingly, we use the following notations for the *left regular action* of $G$ and the *conjugation action* of $G$:

$$x\backslash g\backslash := g^{-1}x, \qquad x \in G,\ g \in G,$$
$$x\backslash g/ := g^{-1}xg, \qquad x \in G,\ g \in G.$$

There is also a *natural action* of $\operatorname{Aut} G$ on $G$ and for the purpose of obtaining a homogeneous notation at a later stage, we set

$$\backslash \boldsymbol{a}/ := \boldsymbol{a}, \qquad \boldsymbol{a} \in \operatorname{Aut} G.$$

Thus $\backslash G/ = \backslash \mathrm{Inn}\, G/$ and for each $g \in G$, $\backslash g/ = \backslash g\backslash /g/ = /g/\backslash g\backslash$. In particular, $\backslash G\backslash \leqslant \backslash \mathrm{Aut}\, G//G/$ and $/G/ \leqslant \backslash G\backslash\backslash \mathrm{Aut}\, G/$.

The *holomorph* of $G$, which is denoted by $\mathrm{Hol}\, G$, is the normalizer in $\mathrm{Sym}\, G$ of $/G/$. It is well known that both $/G/$ and $\backslash G\backslash$ are regular and one the centralizer of the other in $\mathrm{Sym}\, G$. Also, $\backslash \mathrm{Aut}\, G/$ fixes $1 \in G$, therefore

$$\backslash \mathrm{Aut}\, G/ \cap /G/ = \{1\} = \backslash G\backslash \cap \backslash \mathrm{Aut}\, G/.$$

The validity of

$$\backslash \boldsymbol{a}/^{-1}/g/\backslash \boldsymbol{a}/ = /g\boldsymbol{a}/, \qquad \text{and} \qquad \backslash \boldsymbol{a}/^{-1}\backslash g\backslash\backslash \boldsymbol{a}/ = \backslash g\boldsymbol{a}\backslash,$$

shows that both $/G/$ and $\backslash G\backslash$ are normalized by $\backslash \mathrm{Aut}\, G/$. Therefore $\mathrm{Hol}\, G$ contains $\backslash \mathrm{Aut}\, G//G/ = \backslash G\backslash\backslash \mathrm{Aut}\, G/$. On the other hand (see §3.3), there is a homomorphism $\Psi : \mathrm{Hol}\, G \longrightarrow \mathrm{Aut}\, G$ whose kernel is $\backslash G\backslash$. Thus

$$\mathrm{Hol}\, G = \backslash \mathrm{Aut}\, G/ \ltimes /G/ = \backslash G\backslash \rtimes \backslash \mathrm{Aut}\, G/.$$

Coherently with what we did with wreath products, if $\boldsymbol{f} : G \longrightarrow \mathsf{Y}$ is a bijection and $\boldsymbol{s} \in \mathrm{Sym}\, G$, we set

$$\backslash \boldsymbol{f} : \boldsymbol{s}/ := \boldsymbol{f}^{-1}\boldsymbol{s}\boldsymbol{f}.$$

Also, to simplify notations, we put

$$\begin{aligned} \backslash \boldsymbol{f} : g\backslash &:= \boldsymbol{f}^{-1}\backslash g\backslash \boldsymbol{f}, \qquad & /\boldsymbol{f} : g/ &:= \boldsymbol{f}^{-1}/g/\boldsymbol{f}, \qquad & g &\in G, \\ \backslash \boldsymbol{f} : g/ &:= \boldsymbol{f}^{-1}\backslash g/\boldsymbol{f}, \qquad & \backslash \boldsymbol{f} : \boldsymbol{a}/ &:= \boldsymbol{f}^{-1}\backslash \boldsymbol{a}/\boldsymbol{f}, \qquad & \boldsymbol{a} &\in \mathrm{Aut}\, G. \end{aligned}$$

The *inner holomorph* in $\mathsf{Y}$ with respect to a bijection $\boldsymbol{f} : G \longrightarrow \mathsf{Y}$ is $\mathrm{Hol}(\mathsf{Y}, \boldsymbol{f}) = \backslash \boldsymbol{f} : \mathrm{Hol}\, G/$. Of course, for each $\boldsymbol{s} \in \mathrm{Sym}\, \mathsf{Y}$,

$$\boldsymbol{s}^{-1}\, \mathrm{Hol}(\mathsf{Y}, \boldsymbol{f})\boldsymbol{s} = \mathrm{Hol}(\mathsf{Y}, \boldsymbol{f}\boldsymbol{s}),$$

therefore the set of the inner holomorphs in $\mathsf{Y}$ is closed under conjugation.

## 5.11. The Generalized Holomorph

If $K$ is a subgroup of $G$, the coset space of the right cosets of $K$ in $G$ is usually denoted by $(G : K)$. We regret that in this thesis we denote this coset space by $\binom{G}{K\diamond}$. This is because we need a compact yet meaningful notation for the actions on this coset space, which is consistent with the notation introduced for the imprimitive and product actions of wreath products.

For the reminder of this section suppose that $K$ is a core-free subgroup of $G$ and put $N = \mathrm{N}_G\, K$ and $A = \mathrm{St}_{\mathrm{Aut}\, G}\, K$. The actions modulo $K$ of $G$, $N$, $A$ on $\binom{G}{K\diamond}$ are denoted by $/\tbinom{G}{K_\diamond}:\star/$, $\backslash\tbinom{G}{K_\diamond}:\star\backslash$ and $\backslash\tbinom{G}{K_\diamond}:\star/$ respectively, and are defined as follows:

$$\begin{aligned} (Kx)/\tbinom{G}{K_\diamond}:g/ &:= Kxg, && x\in G,\ g\in G,\\ (Kx)\backslash\tbinom{G}{K_\diamond}:n\backslash &:= n^{-1}Kx = Kn^{-1}x, && x\in G,\ n\in N,\\ (Kx)\backslash\tbinom{G}{K_\diamond}:\boldsymbol{a}/ &:= K(x\boldsymbol{a}), && x\in G,\ \boldsymbol{a}\in A. \end{aligned}$$

Note that $\backslash\tbinom{G}{K_\diamond}:N\backslash$ is the centralizer of $/\tbinom{G}{K_\diamond}:G/$ in $\mathrm{Sym}(G:K)$. Also,

$$(*) \qquad \backslash\tbinom{G}{K_\diamond}:N\backslash \leqslant \backslash\tbinom{G}{K_\diamond}:A//\tbinom{G}{K_\diamond}:N/$$

because for each $n\in N$, $\boldsymbol{a}:=\backslash n/\in A$ and $\backslash\tbinom{G}{K_\diamond}:\boldsymbol{a}/=\backslash\tbinom{G}{K_\diamond}:n\backslash/\tbinom{G}{K_\diamond}:n/$. As $\mathrm{Hol}\,G$ is the normalizer of $/G/$ in $\mathrm{Sym}\,G$, so we define $\mathrm{GHol}(G,K)$, the *generalized holomorph* of $G$ modulo $K$, as the normalizer of $/\tbinom{G}{K_\diamond}:G/$ in $\mathrm{Sym}\,\tbinom{G}{K_\diamond}$.

The validity of

$$\backslash\tbinom{G}{K_\diamond}:\boldsymbol{a}/^{-1}/\tbinom{G}{K_\diamond}:g/\backslash\tbinom{G}{K_\diamond}:\boldsymbol{a}/ = /\tbinom{G}{K_\diamond}:g\boldsymbol{a}/, \qquad \boldsymbol{a}\in A,$$

shows that $/\tbinom{G}{K_\diamond}:G/$ is normalized by $\backslash\tbinom{G}{K_\diamond}:A/$. Therefore $\mathrm{GHol}(G,K)$ contains $\backslash\tbinom{G}{K_\diamond}:A//\tbinom{G}{K_\diamond}:G/$. On the other hand (see (3.3.A)) there is a homomorphism $\Psi:\mathrm{GHol}(G,K)\longrightarrow(\mathrm{Inn}\,G)A$ whose kernel is $\backslash\tbinom{G}{K_\diamond}:N\backslash$. Since $\mathrm{Inn}\,G=/\tbinom{G}{K_\diamond}:G/\Psi$ and $A=\backslash\tbinom{G}{K_\diamond}:A/\Psi$, it follows from $(*)$ that

$$\mathrm{GHol}(G,K)=\backslash\tbinom{G}{K_\diamond}:A//\tbinom{G}{K_\diamond}:G/.$$

Again, as we did in the previous section, if $\boldsymbol{f}:\tbinom{G}{K_\diamond}\longrightarrow\mathsf{Y}$ is a bijection and $\boldsymbol{s}\in\mathrm{Sym}\,\tbinom{G}{K_\diamond}$, we define

$$\backslash\boldsymbol{f}:\boldsymbol{s}/:=\boldsymbol{f}^{-1}\boldsymbol{s}\boldsymbol{f}.$$

This yields naturally to the definition of $\mathrm{GHol}(\mathsf{Y},\boldsymbol{f}):=\backslash\boldsymbol{f}:\mathrm{GHol}(G,K)/$, the *inner generalized holomorph* in $\mathrm{Sym}\,\mathsf{Y}$ with respect to $\boldsymbol{f}$.

## 5.12. Primitive groups of Simple Diagonal type

The primitive groups of simple diagonal type are permutationally isomorphic to convenient subgroups of a generalized holomorph $\mathrm{GHol}(G,\Delta)$ where $G=T^l$ is a power of a non abelian simple group $T$ and $\Delta$ is the diagonal of $G$

$$\Delta=\mathtt{diag}(T,l)=\{\ \underbrace{(t,\ldots,t)}_{l}\ \mid\ t\in T\ \}.$$

Observe that $N=\mathrm{N}_G\,\Delta=\Delta$ in this case. Also, due to the fact that $T$ is non abelian simple, $\mathrm{Aut}\,G=(\mathrm{Aut}\,T)\,\mathrm{wr}\,\mathrm{S}_l$ and

$$A=\mathrm{St}_{\mathrm{Aut}\,G}\,\Delta=\mathrm{S}_l\ltimes\mathtt{diag}(\mathrm{Aut}\,T,l)\cong\mathrm{S}_l\times\mathrm{Aut}\,T.$$

For this reason, we set

$$\begin{aligned} \backslash\boldsymbol{a}/ &: \Delta\vec{t} \mapsto \Delta(t_0\boldsymbol{a}, \ldots, t_{l-1}\boldsymbol{a}), && \boldsymbol{a} \in \operatorname{Aut} T, \\ \backslash\boldsymbol{\alpha}/ &: \Delta\vec{t} \mapsto \Delta(t_{0\boldsymbol{\alpha}^{-1}}, \ldots, t_{(l-1)\boldsymbol{\alpha}^{-1}}), && \boldsymbol{\alpha} \in \mathrm{S}_l \end{aligned}$$

and of course for $\boldsymbol{b} = (\boldsymbol{\alpha}, \boldsymbol{a}) \in \mathrm{S}_l \times \operatorname{Aut} T$, $\backslash\boldsymbol{b}/ = \backslash\boldsymbol{\alpha}/\backslash\boldsymbol{a}/$. In this notation

$$\operatorname{GHol}(T^l, \Delta) = \backslash\mathrm{S}_l \times \operatorname{Aut} T//\tbinom{T^l}{\Delta\diamond} : T^l/.$$

Note that $\backslash\mathrm{S}_l \times \operatorname{Aut} T/$ is the stabilizer of $\Delta$ in $\operatorname{GHol}(T^l, \Delta)$ and that

$$\backslash\mathrm{S}_l \times \operatorname{Aut} T/ \cap /\tbinom{T^l}{\Delta\diamond} : T^l/ = \backslash\operatorname{Inn} T/ = /\tbinom{T^l}{\Delta\diamond} : \Delta/ = \backslash\tbinom{T^l}{\Delta\diamond} : \Delta\backslash.$$

As a consequence,

$$\frac{\operatorname{GHol}(T^l, \Delta)}{/\tbinom{T^l}{\Delta\diamond} : T^l/} \cong \frac{\backslash\mathrm{S}_l \times \operatorname{Aut} T/}{\backslash\operatorname{Inn} T/} \cong \mathrm{S}_l \times \operatorname{Out} T.$$

Call $\backslash{\downarrow}/ : \operatorname{GHol}(T^l, \Delta) \longrightarrow \mathrm{S}_l \times \operatorname{Out} T$ the homomorphism induced by the canonical projection of $\operatorname{GHol}(T^l, \Delta)$ onto its factor group over $/\tbinom{T^l}{\Delta\diamond} : T^l/$. Then $\backslash{\downarrow}/$ induces a lattice isomorphism between the interval $\Big[/\tbinom{T^l}{\Delta\diamond} : T^l/ \div \operatorname{GHol}(T^l, \Delta)\Big]$ and the subgroup lattice of $\mathrm{S}_l \times \operatorname{Out} T$. The complete inverse image by $\backslash{\downarrow}/$ of a subgroup $P$ of $\mathrm{S}_l \times \operatorname{Out} T$ is denoted by $\backslash{\uparrow} P/$. For example

$$\begin{aligned} \backslash{\uparrow} (\mathrm{S}_l \times \operatorname{Out} T)/ &= \operatorname{GHol}(T^l, \Delta), \\ \backslash{\uparrow} \{1\}/ &= /\tbinom{T^l}{\Delta\diamond} : T^l/ \text{ and} \\ \backslash{\uparrow} \mathrm{S}_l/ &\cong T \operatorname{wr} \mathrm{S}_l. \end{aligned}$$

It is now time to view the holomorph of the non-abelian simple group $T$ as a subgroup of the generalized holomorph of $T^2$ modulo $\Delta$.

5.12.1. Proposition. *Let $T$ be a non-abelian simple group and set*

$$\boldsymbol{f} : T \longrightarrow \tbinom{T^2}{\Delta\diamond}, \qquad t \mapsto \Delta(1, t).$$

*Then $\boldsymbol{f}$ is a bijection and*

$$\begin{aligned} /\boldsymbol{f} : t/ &= /\tbinom{T^2}{\Delta\diamond} : (1, t)/, && t \in T, \\ \backslash\boldsymbol{f} : t\backslash &= /\tbinom{T^2}{\Delta\diamond} : (t, 1)/, && t \in T, \\ \backslash\boldsymbol{f} : \boldsymbol{a}/ &= \backslash\boldsymbol{a}/, && \boldsymbol{a} \in \operatorname{Aut} T. \end{aligned}$$

*In particular,* $\operatorname{Hol}\big(\tbinom{T^2}{\Delta\diamond}, \boldsymbol{f}\big) < \operatorname{GHol}(T^2, \Delta)$*, likewise,* $\operatorname{Hol} T$ *is permutationally isomorphic to* $\backslash{\uparrow} \operatorname{Out} T/$.

This is the reason of the dichotomy in the following theorem.

5.12.2. Theorem ([**DM96**, Th 4.5A]). *Say* $\boldsymbol{\sigma} : \mathrm{S}_l \times \mathrm{Out}\, T \longrightarrow \mathrm{S}_l$ *the canonical projection. Then* $H = \backslash\!\uparrow P/$ *is primitive if and only if*

(1) $l = 2$ *and* $P\boldsymbol{\sigma} = \{1\}$ *(so* $H$ *is like a subgroup of* $\mathrm{Hol}\, T$*), or*
(2) $P\boldsymbol{\sigma}$ *is primitive.*

In fact, in [**Pra90**] it is shown that $\backslash\!\uparrow /$ is surjective on the set of proper primitive subgroups of $\mathrm{Sym}\binom{T^l}{\Delta_\diamond}$ containing $\backslash\!\uparrow \{1\}/$. Now, a primitive group is of type HS (Holomorph of non-abelian Simple) if it is permutationally isomorphic to a group $H$ as in (1). It is of type Simple Diagonal or more simply of type SD if it is permutationally isomorphic to a group $H$ as in (2). Therefore we have:

5.12.3. Theorem. *Assume that the high subgroup* $\boldsymbol{H}$ *of* $\boldsymbol{S}$ *is primitive of type* HS *or* SD*. Say for example that* $\boldsymbol{H} = \boldsymbol{f}^{-1}\backslash\!\uparrow P/\boldsymbol{f}$ *with* $\boldsymbol{f} : \binom{T^l}{\Delta_\diamond} \longrightarrow \Omega$ *a bijection and* $P \leqslant \mathrm{S}_l \times \mathrm{Out}\, T$*. Assume also that* $\boldsymbol{H} < G < \boldsymbol{U}$*. Then* $G = \mathrm{GHol}(\Omega, \boldsymbol{f}) \cap \boldsymbol{U}$ *and* $[\boldsymbol{H} \div \boldsymbol{U}] = \mathcal{M}_1$*.*

## 5.13. Primitive groups of Compound Diagonal type

The primitive groups of Compound Diagonal type are the inner blow-ups of primitive groups of Simple Diagonal type. The inclusions of a primitive group of Compound Diagonal type in other groups are described below (recall that the symbol $\star$ may be used in a same statement to name different unspecified objects).

5.13.1. Proposition ([**Pra90**, 8.1]). *Suppose that* $H < G < \mathrm{Sym}\, \Omega$ *with* $H$ *primitive of type* CD *and* $G \neq \mathrm{Alt}\, \Omega$*. Say* $H = \lceil \boldsymbol{p} : A \uparrow \star_l \rceil$ *with* $A \leqslant \mathrm{GHol}(T^c, \Delta)$ *of degree* $m = |T|^{c-1}$ *and* $\boldsymbol{p} \in \mathrm{EPDec}(\Omega, \binom{T^c}{\Delta_\diamond}, l)$ *for example. Then one of the following is true:*

(1) $\mathrm{Soc}\, G = \mathrm{Soc}\, H$ *and there is* $\boldsymbol{q} \in \mathrm{EPDec}(\Omega, \binom{T^c}{\Delta_\diamond}, l)$ *such that* $G$ *is contained in* $\lceil \boldsymbol{q} : \mathrm{GHol}(T^c, \Delta) \,\mathrm{wr}\, \mathrm{S}_l \rceil$*. Note that* $H$ *is low.*
(2) $l = ab > a$*. There are* $\boldsymbol{r}, \boldsymbol{q} \in \mathrm{EPDec}(\Omega, m^a, b)$ *such that* $H = \lceil \boldsymbol{r} : B \uparrow \star_b \rceil$ *and* $G = \lceil \boldsymbol{q} : C \uparrow \star_b \rceil$ *with* $\mathrm{A}_{m^a} \leqslant C \leqslant \mathrm{S}_{m^a}$*. Either* $a = 1$ *and* $A$ *is permutationally isomorphic to* $B$*; or* $H$ *is low,* $a > 1$ *and there is* $\boldsymbol{p}_1 \in \mathrm{EPDec}(m^a, \binom{T^c}{\Delta_\diamond}, a)$ *such that* $B = \lceil \boldsymbol{p}_1 : A \uparrow \star_a \rceil$*.*

5.13.2. Theorem. *Suppose that the non trivial high subgroup* $\boldsymbol{H}$ *of* $\boldsymbol{S}$ *is primitive of type* CD*. Say* $\boldsymbol{H} = \lceil \boldsymbol{p} : A \uparrow \star_l \rceil$ *with* $A \leqslant \mathrm{GHol}(T^c, \Delta)$ *of degree* $m = |T|^{c-1}$ *and* $\boldsymbol{p} \in \mathrm{EPDec}(\Omega, \binom{T^c}{\Delta_\diamond}, l)$ *for example. Suppose also that* $\boldsymbol{H} < G < \boldsymbol{U} = (\mathrm{Alt}\, \Omega)\boldsymbol{H}$*. Then* $\boldsymbol{H} = \lceil \boldsymbol{p} : A \,\mathrm{wr}\, \mathrm{S}_l \rceil \cap \boldsymbol{U}$*,* $A = \mathrm{GHol}(T^c, \Delta) \lessdot \mathrm{S}_m$ *and* $G = \mathrm{St}_{\boldsymbol{U}}\, \mathrm{Qs}\, \tilde{\boldsymbol{p}}$*. In particular* $[\boldsymbol{H} \div \boldsymbol{U}] = \mathcal{M}_1$*.*

PROOF. $\boldsymbol{H}$ is high therefore by 5.9.1 $A = \mathrm{GHol}(T^c, \Delta)$. Moreover, $A$ must be a maximal subgroup of $\mathrm{S}_m$ and $\boldsymbol{H} = \lceil \boldsymbol{p} : A \,\mathrm{wr}\, \mathrm{S}_l \rceil \cap \boldsymbol{U}$. Point 2 of 5.13.1 says that $G = \lceil \boldsymbol{q} : \mathrm{S}_m \,\mathrm{wr}\, \mathrm{S}_l \rceil \cap \boldsymbol{U}$ (remember that $G \lessdot \boldsymbol{U}$ and $\lceil \boldsymbol{q} : \mathrm{A}_m \,\mathrm{wr}\, \mathrm{S}_l \rceil \cap \boldsymbol{U}$ is not). Thus it is enough to show that the codirections of $\boldsymbol{p}$ and $\boldsymbol{q}$ are the same. In fact, provided $\mathrm{Qs}\, \tilde{\boldsymbol{p}} = \mathrm{Qs}\, \tilde{\boldsymbol{q}}$,

$$G = \mathrm{St}_{\boldsymbol{U}}\, \mathrm{Qs}\, \tilde{\boldsymbol{q}} = \mathrm{St}_{\boldsymbol{U}}\, \mathrm{Qs}\, \tilde{\boldsymbol{p}} \,.$$

Apply 5.7.1 to $\mathrm{Soc}\, G = \lceil \boldsymbol{q} : (\mathrm{A}_m)^l \rceil$ to find that $\mathrm{Qs}\, \tilde{\boldsymbol{q}}$ is the set of the imprimitivity systems of order $m$ of $\mathrm{Soc}\, G$. These must be $l$ different imprimitivity systems of order $m$ of $\mathrm{Soc}\, H$ because 5.3.1 and 5.3.2 imply $\mathrm{Soc}\, H \leqslant \mathrm{Soc}\, G$. But $\mathrm{Soc}\, H = \lceil \boldsymbol{p} : (\mathrm{Soc}\, A)^l \rceil$ has exactly $l$ imprimitivity systems of order $m$ which are the codirections of $\boldsymbol{p}$. Thus the codirections of $\boldsymbol{q}$ are exactly the codirections of $\boldsymbol{p}$, as we claimed. □

## 5.14. The affine case

Again, we start quoting a proposition from [**Pra90**].

5.14.1. PROPOSITION ([**Pra90**, 5.1]). *Suppose that $H < G$ with $H$ primitive of affine type and $H$ of degree $n = r^l$, $l \geqslant 1$, $r$ a prime. If $\mathrm{A}_n \neq G \neq \mathrm{S}_n$ then either $\mathrm{Soc}\, H = \mathrm{Soc}\, G$ and so $G$ is of affine type or $\mathrm{Soc}\, H \leqslant \mathrm{Soc}\, G$ and one of the following holds:*

1. *$n = 11$, $\mathrm{Soc}\, G = \mathrm{PSL}(2, 11)$.*
2. *$n = 11$, $\mathrm{Soc}\, G = \mathrm{M}_{11}$.*
3. *$n = 23$, $\mathrm{Soc}\, G = \mathrm{M}_{23}$.*
4. *$n = 27$, $\mathrm{Soc}\, G = \mathrm{PSp}(4, 3)$, $H \cap \mathrm{Soc}\, G = 3^3 \rtimes \mathrm{S}_4$.*
5. *$l = 1$ and $r = 1 + q + \cdots + q^{d-1}$ for some $q = p^f > 1$, $p$ prime. Here $G$ is almost simple with socle isomorphic to $\mathrm{PSL}_d(q)$. Note that $d$ must be a prime and that $f = d^a$ for some $a \geqslant 0$. In particular, if $d$ is even, then $d = 2$ and $r = 1 + 2^{2^a}$ is a Fermat prime.*
6. *$l = ab > a$; there is $\boldsymbol{r} \in \mathrm{EPDec}(n, r^a, b)$ and a subgroup $A$ of $\mathrm{S}_{r^a}$ such that $G = \lceil \boldsymbol{r} : A \uparrow \star_b \rceil$. Note that if $A < \mathrm{S}_{r^a}$, then $H$ is low.*

**Degrees 11 and 23.** An affine primitive group of degree $r$ prime must contain a cycle of length $r$. With the ATLAS and with Pálfy Lemma one finds that the interval of the subgroups of $\mathrm{S}_{11}$ containing a cycle of length 11 and the interval of the subgroups of $\mathrm{S}_{23}$ containing a cycle of length 23 are as in Figure 5.C.

**Degree 27.** The discussion below is based on the maximal subgroups Table 5.A. Since $\mathrm{PGSp}_4(3)$ has only one conjugacy class of subgroups of index 27, the transitive subgroups of $\mathrm{S}_{27}$ isomorphic to $\mathrm{PGSp}_4(3)$ are all conjugate in $\mathrm{S}_{27}$. Say $G$ one of these subgroups of $\mathrm{S}_{27}$, $S = \mathrm{Soc}\, G$ and $T$ the stabilizer of one point in $G$. Since $T \cap S \lessdot S$, $S$ and $G$ are primitive. Call $H$ a maximal subgroup of $G$ isomorphic to $3^3 \rtimes (\mathrm{S}_4 \times 2)$. Note that $H \cap S \cong 3^3 \rtimes \mathrm{S}_4$. $H$ is a subgroup of order $3^3 \cdot 48$ of $\mathrm{N}(3^3) \cong \mathrm{AGL}_3(3)$. Since all maximal subgroups of order 48 of

$$\frac{\mathrm{AGL}_3(3)}{3^3} \cong \mathrm{GL}_3(3) = C_2 \times \mathrm{SL}_3(3) \cong C_2 \times \mathrm{PSL}_3(3)$$

are conjugate, $H$ must be permutationally isomorphic to $\mathrm{AGL}(1,3) \,\mathrm{wr}\, \mathrm{S}_3 = \mathrm{S}_3 \,\mathrm{wr}\, \mathrm{S}_3$ in product action. (Note that $\mathrm{S}_m \,\mathrm{wr}\, \mathrm{S}_l$ is maximal in $\mathrm{S}_{m^l}$ only for $m \geqslant 5$). By 5.5.J, $H$ is an odd subgroup of $\mathrm{S}_{27}$, and so is $G$. Using the maximal subgroup table argument, one finds that

$$[H \cap S \div \mathrm{A}_{27}] \cong \mathcal{M}_2 \cong [H \div \mathrm{S}_{27}].$$

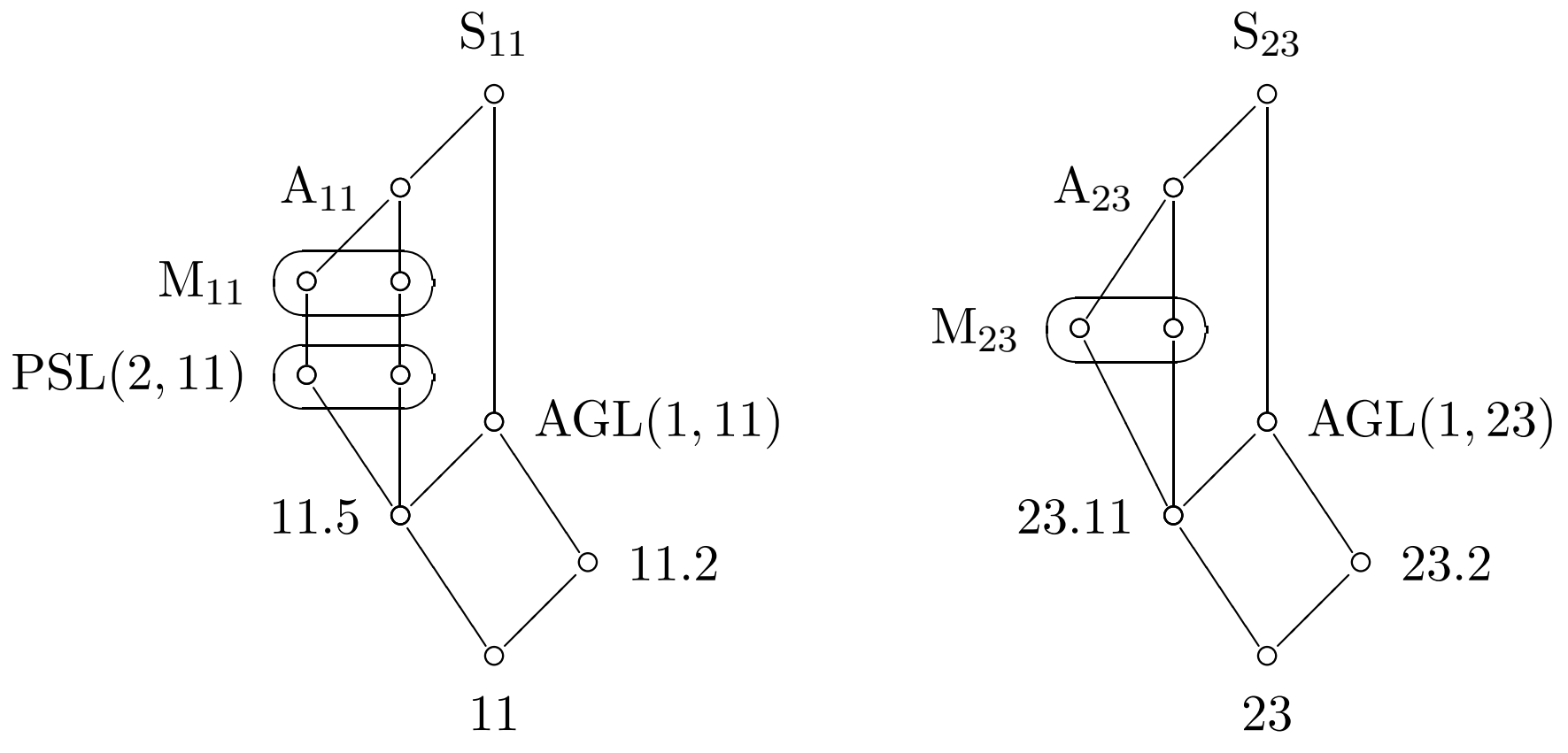

FIGURE 5.C. Subgroups containing a $p$-Sylow of $\mathrm{S}_p$; $p = 11, 23$

TABLE 5.A. Maximal subgroups of $\mathrm{PSp}_4(3)$

| Order | Index | Structure | in $\mathrm{PGSp}_4(3)$ | Characters |
|---|---|---|---|---|
| 960 | 27 | $264 \rtimes \mathrm{A}_5$ | $: 2^4 \rtimes \mathrm{S}_5$ | $1a + 6a + 20a$ |
| 720 | 36 | $\mathrm{S}_6$ | $: \mathrm{S}_6 \times 2$ | $1a + 15b + 20a$ |
| 648 | 40 | $3_+^{1+2} \rtimes 2\mathrm{A}_4$ | $: 3_+^{1+2} \rtimes 2\mathrm{S}_4$ | $1a + 15a + 24a$ |
| 648 | 40 | $3^3 \rtimes \mathrm{S}_4$ | $: 3^3 \rtimes (\mathrm{S}_4 \times 2)$ | $1a + 15b + 24a$ |
| 576 | 45 | $2 \cdot (\mathrm{A}_4 \times \mathrm{A}_4).2$ | $: \ast.2$ | $1a + 20a + 24a$ |

**Projective spaces with a prime number of projective points.** We show that the inclusion $H < G$ with $H$ primitive of affine type and $G$ proper primitive group of almost simple type and socle isomorphic to $\mathrm{PSL}_d(q)$ as in (5) do appear but that $G < \mathrm{A}_r$ in this case (recall that $r$ is prime). Then one may appeal to a paper from Pálfy where all second maximal subgroups of $\mathrm{A}_r$ have been classified.

Say $\mathbb{K}$ a field of order $q^d$ and $\mathbb{F}$ its subfield of order $q$. Regard $\mathbb{K}$ as a vector space $\mathsf{V}$ over $\mathbb{F}$. It is well known that $\mathrm{GL}(\mathsf{V})$ and $\Gamma\mathrm{L}(\mathsf{V})$ are the centralizer and the normalizer respectively of $\backslash\mathbb{F}\backslash$ in $\mathrm{Aut}(\mathsf{V},+)$, where $(\mathsf{V},+)$ is the abelian group supporting the structure of vector space of $\mathsf{V}$ and, for each $\mathfrak{a} \in \mathbb{K}$, $\backslash\mathfrak{a}\backslash$ is the multiplication by the scalar $\mathfrak{a}$. Say $\mathfrak{u}$ a generator of the multiplicative group $\mathbb{K}^\times$ of $\mathbb{K}$. of course, $\backslash\mathfrak{u}\backslash \in \mathrm{GL}(\mathsf{V})$ and its image in $\mathrm{PGL}(\mathsf{V})$ generates a cyclic group which is transitive on $\mathrm{PG}(\mathsf{V}) = \mathrm{PG}_d(q)$. Provided that $r$ is prime, this is a primitive subgroup of affine type of $\mathrm{PGL}(\mathsf{V})$. We may assume $r > 5$ because $\mathrm{PGL}_2(4) = \mathrm{A}_5$, against our assumptions. Then $q(d+1) > 6$ and by 3.9.3 $\mathrm{PGL}_d(q)$ is even. As a consequence (see 3.7.2), $\mathrm{P\Gamma L}_d(q)$ is even when $f$ is odd. However, if $f$ is even, then $d = 2$ but $q \neq 4$ (otherwise $r = 5$) so that $\mathrm{P\Gamma L}_d(q)$ is even again.

**Conclusions.**

5.14.2. Theorem. *Suppose that the high subgroup $\boldsymbol{H}$ of $\boldsymbol{S} = \mathrm{Sym}\,\Omega$ is primitive of type* HA *and degree $n = r^l$, where $r$ is prime. Suppose also that $\boldsymbol{H} < G < \boldsymbol{U}$ where $\boldsymbol{U} = (\mathrm{Alt}\,\Omega)\boldsymbol{H}$. Then one of the following is true:*

(1) *$G = \mathrm{N}_{\boldsymbol{U}}\,\mathrm{Soc}\,\boldsymbol{H}$ is of affine type;*
(2) *$n = 11$, $G \cong \mathrm{AGL}(1,11)$, $\boldsymbol{H} \cong 11.2$, $[\boldsymbol{H} \div \mathrm{S}_{11}] = \mathcal{M}_1$;*
(3) *$n = 23$, $\mathrm{Soc}\,G = \mathrm{M}_3$, $[\boldsymbol{H} \div \mathrm{A}_{23}] = \mathcal{M}_2$;*
(4) *$n = 27$, $\mathrm{Soc}\,G = \mathrm{PSp}_4(3)$, $[\boldsymbol{H} \div \boldsymbol{U}] = \mathcal{M}_2$;*
(5) *$l = 1$ and $r = \frac{q^d-1}{q-1}$ for some $q = p^f > 1$, $p$ prime. $G$ has socle isomorphic to $\mathrm{PSL}_d(q)$ and $\boldsymbol{U}$ is the finite alternating group of prime degree $r$. Pálfy classifies all such cases and shows the if $[\boldsymbol{H} \div \mathrm{A}_r] = \mathcal{M}_k$ with $k > 2$, then $k = 5, 7, 11$ (see 5.2.1);*
(6) *$l = ab > a$. There is a vector space $\mathsf{V}$ of dimension $a$ over the field of order $r$ and there is $\boldsymbol{p} \in \mathrm{EPDec}(\Omega, \mathsf{V}, b)$ such that $\boldsymbol{H} = \lceil \boldsymbol{p} : \mathrm{AGL}(\mathsf{V})\,\mathrm{wr}\,\mathrm{S}_b \rceil \cap \boldsymbol{U}$. In this case $\boldsymbol{H}$ determines $\mathrm{Qs}\,\tilde{\boldsymbol{p}}$, $b = \left|\mathrm{Qs}\,\tilde{\boldsymbol{p}}\right|$ and $G = \mathrm{St}_{\boldsymbol{U}}\,\mathrm{Qs}\,\tilde{\boldsymbol{p}}$.*

*In particular, if $[\boldsymbol{H} \div \boldsymbol{U}] \cong \mathcal{M}_k$, then $k \in \{1, 2, 5, 7, 11\}$.*

Proof. It is enough to show that the last part of the theorem follows from the last part of 5.14.1. First of all, by maximality of $G$, $G = \lceil \boldsymbol{r} : \mathrm{S}_{r^a}\,\mathrm{wr}\,\mathrm{S}_b \rceil \cap \boldsymbol{U}$ and $r^a \geqslant 5$. In fact, we may assume that $r^a$ is a vector space, say $\mathsf{V}$, of dimension $a$ over $\mathbb{F}_r$. So

$$G = \lceil \boldsymbol{r} : (\mathrm{Sym}\,\mathsf{V})\,\mathrm{wr}\,\mathrm{S}_b \rceil \cap \boldsymbol{U},$$

$$\mathrm{N}(G) = \lceil \boldsymbol{r} : (\mathrm{Sym}\,\mathsf{V})\,\mathrm{wr}\,\mathrm{S}_b \rceil$$

and

$$\operatorname{Soc} G = \lceil \boldsymbol{r} : (\operatorname{Alt} \mathsf{V})^b \rceil.$$

Up to the order of the components of $G$, we may also assume that for each $i \in b$ $(\sqrt[b]{G})_i = \lceil \boldsymbol{r} : (\operatorname{Alt} \mathsf{V})\boldsymbol{\kappa}_i \rceil$. Note that $(\operatorname{Soc} \boldsymbol{H})\boldsymbol{\rho}_i^G$ must be an elementary abelian subgroup of $(\sqrt[b]{G})_i$ which is transitive on $\mathsf{V}$. Thus $(\operatorname{Soc} \boldsymbol{H})\boldsymbol{\rho}_i^G$, being an elementary abelian regular subgroup of $\operatorname{Alt} \mathsf{V}$ ([**Rob93**, 1.6.3]), is conjugate in $\operatorname{Sym} \mathsf{V}$ to $/\mathsf{V}/$, the group of the translations of $\mathsf{V}$. It follows that $\operatorname{Soc} \boldsymbol{H}$ is conjugate in $\operatorname{Shoe} \mathrm{N}(G)$ to $\lceil \boldsymbol{r} : /\mathsf{V}/^b \rceil$. Since $\boldsymbol{H} = \mathrm{N}_G \operatorname{Soc} \boldsymbol{H}$, this shows that $\boldsymbol{H}$ is conjugate in $\mathrm{N}(G)$ to $\lceil \boldsymbol{r} : \mathrm{AGL}(\mathsf{V}) \operatorname{wr} \mathrm{S}_b \rceil \cap \boldsymbol{U}$. In particular, Lemma 5.14.3 below shows that $b$ is determined by $\boldsymbol{H}$ and hence that all overgroups of $\boldsymbol{H}$ of type PA are conjugate in $\operatorname{Sym} \Omega$. Their number is given by Pálfy Lemma:

$$\operatorname{hm}(G, \boldsymbol{H}) = \frac{\left|\mathrm{N}(\boldsymbol{H}) : \boldsymbol{H}\right|}{\left|\mathrm{N}(G) : G\right|} \frac{\operatorname{hm}(\boldsymbol{H}, G)}{\left|G : \boldsymbol{H}\right|}$$

To compute $\operatorname{hm}(\boldsymbol{H}, G)$, assume that for some $\boldsymbol{s} \in \operatorname{Sym} \Omega$, $\boldsymbol{H}^{\boldsymbol{s}} \leqslant G$; then $(\operatorname{Soc} \boldsymbol{H})^{\boldsymbol{s}} \leqslant G$. Then $(\operatorname{Soc} \boldsymbol{H})^{\boldsymbol{s}}$ is conjugate in $\mathrm{N}(G)$ to $\lceil \boldsymbol{r} : /\mathsf{V}/^b \rceil$, that is, $(\operatorname{Soc} \boldsymbol{H})^{\boldsymbol{s}} = (\operatorname{Soc} \boldsymbol{H})^{\boldsymbol{g}}$ for some $\boldsymbol{g} \in \mathrm{N}(G)$. From $\boldsymbol{H}^{\boldsymbol{s}} \leqslant \mathrm{N}\big((\operatorname{Soc} \boldsymbol{H})^{\boldsymbol{s}}\big)$, from $\boldsymbol{H}^{\boldsymbol{s}} \leqslant G$ and from $[\boldsymbol{H}^{\boldsymbol{s}} \div \boldsymbol{U}] \cong \mathcal{M}_k$, follows

$$\begin{aligned} \boldsymbol{H}^{\boldsymbol{s}} &= \mathrm{N}\big((\operatorname{Soc} \boldsymbol{H})^{\boldsymbol{s}}\big) \cap G = \mathrm{N}\big((\operatorname{Soc} \boldsymbol{H})^{\boldsymbol{g}}\big) \cap G = \big(\mathrm{N}(\operatorname{Soc} \boldsymbol{H}) \cap G\big)^{\boldsymbol{g}} \\ &= \boldsymbol{H}^{\boldsymbol{g}}. \end{aligned}$$

Note that $\mathrm{N}(\boldsymbol{H}) \leqslant \mathrm{N}(G)$ because

$$\mathrm{N}\big((\mathrm{AGL}(\mathsf{V}) \operatorname{wr} \mathrm{S}_b) \cap \operatorname{Alt} \mathsf{V}^b\big) = \mathrm{AGL}(\mathsf{V}) \operatorname{wr} \mathrm{S}_b$$

Therefore $\operatorname{hm}(\boldsymbol{H}, G) = \left|\mathrm{N}(G) : \mathrm{N}(\boldsymbol{H})\right|$ and

$$\operatorname{hm}(G, \boldsymbol{H}) = \frac{\left|\mathrm{N}(\boldsymbol{H}) : \boldsymbol{H}\right|\left|\mathrm{N}(G) : \mathrm{N}(\boldsymbol{H})\right|}{\left|\mathrm{N}(G) : G\right|\left|G : \boldsymbol{H}\right|} = \frac{\left|\mathrm{N}(G) : \boldsymbol{H}\right|}{\left|\mathrm{N}(G) : \boldsymbol{H}\right|} = 1 \ .$$

This shows that $G$ is determined by $\boldsymbol{H}$ and so are the imprimitivity systems of order $r^a$ of $\operatorname{Soc} G$, that is, the codirections of $\boldsymbol{r}$. □

5.14.3. Lemma. *Suppose* $\left|\mathrm{AGL}(a, p) \operatorname{wr} \mathrm{S}_b\right| = \left|\mathrm{AGL}(a', p) \operatorname{wr} \mathrm{S}_{b'}\right|$ *with* $ab = a'b'$, *then* $a = a'$ *and* $b = b'$.

Proof. Clearly, $p^{ab}$ divides both the orders so that

$$\left|\mathrm{GL}(a, p) \operatorname{wr} \mathrm{S}_b\right| = \left|\mathrm{GL}(a', p) \operatorname{wr} \mathrm{S}_{b'}\right|.$$

Call $\lfloor \log_p x! \rfloor$ the exponent of the highest power of $p$ dividing $x!$ (6.4.A); then consider the highest powers of $p$ dividing the two members to obtain

$$b\binom{a}{2} + \lfloor \log_p b! \rfloor = b'\binom{a'}{2} + \lfloor \log_p b'! \rfloor \ .$$

Assume now that $b < b'$; It is enough to show that this leads to contradictions. The substitution of $ab$ with $a'b'$ and the hypothesis $b < b'$ gives

$$0 < a'b'(a - a') = 2[\lfloor \log_p b'! \rfloor - \lfloor \log_p b! \rfloor] < 2b' \ .$$

Thus, $a'(a - a') < 2$ which forces $a' = 1$, $a = 2$ and $b' = 2b$. Then

$$\left|\mathrm{GL}(2, p) \,\mathrm{wr}\, \mathrm{S}_b\right| = \left|\mathrm{GL}(1, p) \,\mathrm{wr}\, \mathrm{S}_{2b}\right|$$

which implies

$$p^b (p + 1)^b = 2b \cdot (2b - 1) \cdots (b + 1) \ .$$

In particular, $b \neq 1$ and $b \neq 2$. If $p = 2$, then 5 may not divide the right member and so $2b < 5$. But this implies $b \leqslant 2$. Therefore $p \neq 2$ and $p+2$ is coprime to $p$ and $p + 1$. Since there exist divisors of $p + 2$ which do not divide $p^b (p + 1)^b$, $p + 2 > 2b$ and so

$$2b \cdot (2b - 1) \cdots (b + 1) > (2b - 2)^b (2b - 1)^b \ .$$

But this implies $2b > (2b - 2)^b (2b - 1)$ and hence $b = 1$ which we have seen above yields to contradiction. □

## 5.15. Almost simple in non almost simple

Assume that $\boldsymbol{H}$ is primitive of type almost simple and that $\boldsymbol{H} < G < \boldsymbol{U}$ for some non almost simple group $G$. A classification of the non maximal primitive almost simple groups is given in [**LPS87**]. From it, one extracts a table describing the primitive almost simple groups $K$ which are contained in some non almost simple group $G$:

| $\mathrm{Soc}\, K$ | $\mathrm{N}(\mathrm{Soc}\, G)$ | degree | comments |
|---|---|---|---|
| $\mathrm{PSL}(2, 7)$ | $\mathrm{AGL}_3(2)$ | 8 | |
| $\mathrm{A}_6$ | $\mathrm{S}_6 \,\mathrm{wr}\, \mathrm{S}_2$ | 36 | $G$ in product action, $K$ contains an outer automorphism of $\mathrm{S}_6$, $\mathrm{M}_{12}$, $\mathrm{P\Gamma Sp}_4(q)$ respectively. |
| $\mathrm{M}_{12}$ | $\mathrm{S}_{12} \,\mathrm{wr}\, \mathrm{S}_2$ | 144 | |
| $\mathrm{PSp}_4(q)$, $q$ even, $q > 2$ | $\mathrm{S}_m \,\mathrm{wr}\, \mathrm{S}_2$, $m = \frac{1}{2}q^2(q^2 - 1)$ | $m^2$ | |

Of course, our $\boldsymbol{H}$ is one of those $K$, in fact, it may only be the normalizer in $G$ of $\operatorname{Soc} K$. The first row of the table is dismissed by showing in Figure 5.D the interval $[K \div \mathrm{S}_8]$ for a primitive $K \cong \mathrm{PSL}(2,7)$. This row leads to an $\mathcal{M}_2$ in $\mathrm{A}_8$.

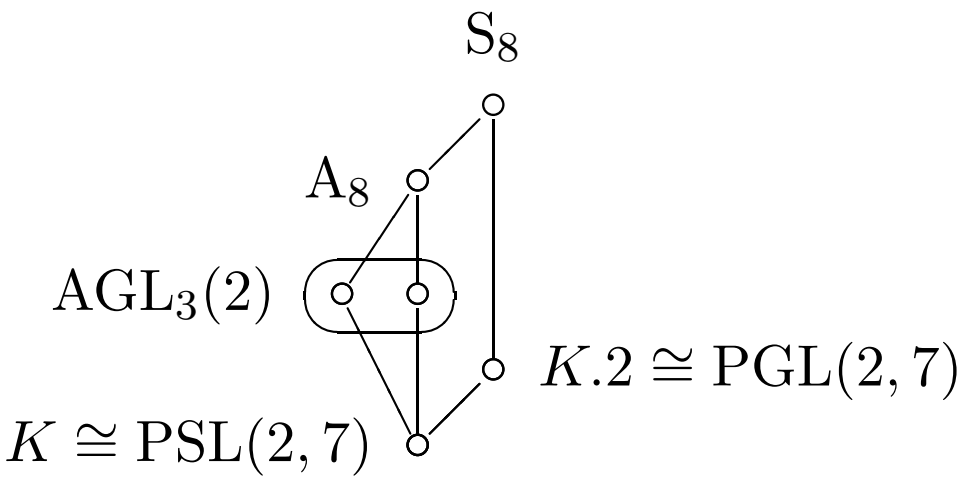

FIGURE 5.D. An interval in the subgroup lattice of $\mathrm{S}_8$, $K$ is primitive

The other rows present a similar pattern. Note first that $G$ and hence $\boldsymbol{H}$ are even. We may assume that $K$, being a candidate for $\boldsymbol{H}$, is the normalizer in $G$ of $\operatorname{Soc} K$. However, such a $K$ is still low. We only show this for rows 2 and 4.

Consider row 2 first. To view the socle of $K$ as a subgroup of

$$G = \mathrm{S}_6 \operatorname{wr} \mathrm{S}_2 = (\mathrm{S}_6 \times \mathrm{S}_6) \rtimes \langle \boldsymbol{t} \rangle,$$

($\boldsymbol{t}$ being the automorphism swapping the coordinates), call $\boldsymbol{a}$ an outer automorphism of order 2 of $\mathrm{S}_6$. One may assume that

$$\operatorname{Soc} K = \{ \ (s, s^{\boldsymbol{a}}) \ \mid \ s \in \mathrm{A}_6 \ \}.$$

Since $\boldsymbol{a}$ has order 2, $\boldsymbol{t}$ normalizes $\operatorname{Soc} K$ and so $K \geqslant \langle D, \boldsymbol{t} \rangle = \langle \boldsymbol{t} \rangle D$ where

$$D = \{ \ (s, s^{\boldsymbol{a}}) \ \mid \ s \in \mathrm{S}_6 \ \}.$$

Clearly, $\langle \boldsymbol{t} \rangle D$ has the same order of $\operatorname{Aut} \mathrm{A}_6$. Thus $K = \langle \boldsymbol{t} \rangle D$ and $K$ is low because it is properly contained in $(\mathrm{A}_6 \times \mathrm{A}_6)K$ which is properly contained in $G$.

Now consider row 4. To view the socle of $K$ as a subgroup of $G = (\mathrm{S}_m \times \mathrm{S}_m) \rtimes \langle \boldsymbol{t} \rangle$ ($\boldsymbol{t}$ being the automorphism swapping the coordinates), call $\boldsymbol{b}$ an outer automorphism of $\mathrm{P\Gamma Sp}_4(q)$ whose square is a field automorphism. One may assume that

$$\operatorname{Soc} K = \{ \ (s, s^{\boldsymbol{b}}) \ \mid \ s \in \mathrm{PSp}_4(q) \ \}.$$

Note that $\boldsymbol{g} = (\boldsymbol{b}^2, 1)\boldsymbol{t}$ normalizes $\operatorname{Soc} K$ and so $K \geqslant \langle D, \boldsymbol{g} \rangle$ where

$$D = \{ \ (s, s^{\boldsymbol{b}}) \ \mid \ s \in \mathrm{P\Gamma Sp}_4(q) \ \}.$$

Clearly, $\langle D, \boldsymbol{g} \rangle$ has the same order of $\operatorname{Aut} \mathrm{PSp}_4(q)$. Thus $K = \langle D, \boldsymbol{g} \rangle$ and $K$ is low because it is contained in $(\mathrm{A}_m \times \mathrm{A}_m)K$ which is contained in $G$.

## 5.16. Proof of Theorem B

In view of the O'Nan-Scott Theorem and of Figure 5.B, we only have to deal with second maximal subgroups of type

**PA:** settled in 5.9.3;

**HS,SD:** settled in 5.12.3;

**CD:** settled in 5.13.2;

**HA:** settled in 5.14.2; and

**AS:** settled in §5.15.

This completes the proof of Theorem B and with it, this chapter.

CHAPTER 6

# Almost simple second maximal subgroups

Recall from previous chapter that

$$\boldsymbol{S} \equiv \operatorname{Sym}\Omega \cong \mathrm{S}_n, \quad n \geqslant 5$$

$$\boldsymbol{A} \equiv \operatorname{Alt}\Omega \cong \mathrm{A}_n$$

$$\boldsymbol{U} = \boldsymbol{S} \text{ or } \boldsymbol{A}$$

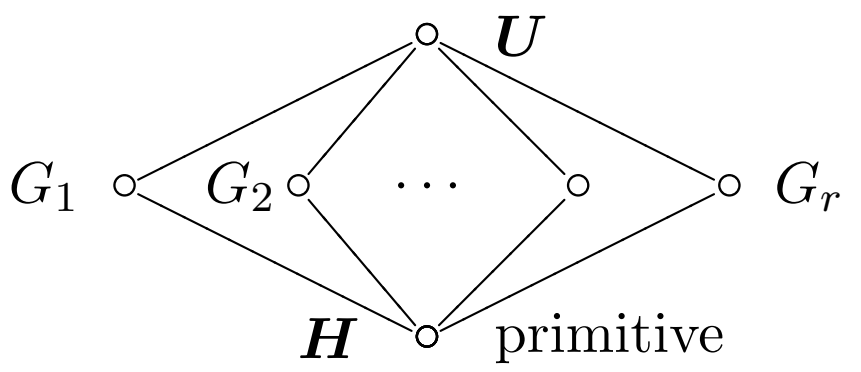

where we analyzed the subgroups $\boldsymbol{H}$ such that either $\boldsymbol{H}$ or one of the $G_i$ is not almost simple. In this chapter we deal with the remaining case and we prove Theorem C.

**Theorem C.** *Suppose that $\boldsymbol{H}$ and all the $G_i$ are almost simple. If $r \geqslant 2$ or if $G_1 \neq \mathrm{N}_{\boldsymbol{U}} \operatorname{Soc} \boldsymbol{H}$, then $\boldsymbol{H}$ and all the $G_i$ can be recovered at once from a same bunch of Table* 6.B*. In particular, $r \leqslant 3$.*

The key is that if $G$ is one of the $G_i$ and $K$ is a point stabilizer in $G$, then by the Frattini Argument $G = \boldsymbol{H}K$. Since both $K$ and $\boldsymbol{H}$ are maximal in $G$, this is a maximal factorization of the almost simple group $G$. The maximal factorizations of the almost simple groups has been classified in [**LPS90**]. An exhaustive analysis is then all what is needed to obtain our result. We begin with a little discussion on group factorizations.

## 6.1. Group factorizations

Given two subsets $\mathsf{H}$ and $\mathsf{K}$ of a group, we denote their product by

$$\mathsf{HK} = \left\{ hk \;\middle|\; h \in \mathsf{H},\, k \in \mathsf{K} \right\}.$$

The product of two subsets may not be commutative, also, the product of two subgroups may not be a subgroup. However, the product of two subgroups $H$ and $K$ is a subgroup if and only if they commute, that is, $HK = KH$.

A finite group $G$ is said *factorizable* if $G = HK$ for some proper subgroups $H$ and $K$ of $G$. Then the expression $HK$ is called a *factorization* of $G$ and it is a *maximal factorization* if both $H$ and $K$ are maximal subgroups of $G$. Any factorization of $G$ is "contained" in at least one maximal factorization. Recall that when $K$ is subgroup of $G$ we denote the action of $G$ on the coset space $\binom{G}{K\diamond}$ of right cosets of $K$ in $G$ by

$$/\binom{G}{K\diamond} : \star / : G \longrightarrow \operatorname{Sym}\binom{G}{K\diamond}, \qquad g \mapsto /\binom{G}{K\diamond} : g/.$$

The following conditions on a finite group $G$ are equivalent:

- $HK$ is a factorization of $G$;
- $/\binom{G}{K\diamond} : H/$ is transitive;
- $/\binom{G}{H\diamond} : K/$ is transitive;
- the permutation characters $\chi_1$ of $G$ on $\binom{G}{H\diamond}$ and $\chi_2$ of $G$ on $\binom{G}{K\diamond}$ satisfy $(\chi_1, \chi_2) = 1$.

Note that $G$ acts faithfully on $\binom{G}{K\diamond}$ if and only if $K$ is *core-free* in $G$, that is, $\operatorname{Core}_G K = \{1\}$, where $\operatorname{Core}_G K = \bigcap_{g \in G} K^g$. Consequently, one looks first for the *core-free factorization*s $HK$ of $G$ where both $H$ and $K$ are core-free in $G$. However, not all the core-free factorizations of a group $G$ are "contained" in core-free and maximal factorizations of $G$.

Another observation is that if $G$ acts transitively on a set $\Omega$, and $H$ is a transitive subgroup of $G$, then for $\omega \in \Omega$, $G = HG_\omega$. This last observation is known as the Frattini Argument.

When the attention focus not so much on given factorizations of a group but on "factorizations of groups" as mathematical objects, the notation presented so far may lead to ambiguities. For example, if $HK$ is a factorization of $G$, then the equation $HK = KH$ holds but the two expressions $HK$ and $KH$ are visibly different. Also, the expression $HK$ is the same used for any product of subgroups and does not encode the essential information that this product is a group, in fact the whole group. For these reasons we occasionally adopt a wrapper as follows.

A *group factorization* is an ordered triple $(G, H, K)$ such that $HK$ is a factorization of $G$. Note that the group factorizations $(G, H, K)$ and $(G, K, H)$ are now undoubtedly different; we say that they are one the *dual* of the other. The attributes "maximal" and "core-free" extend obviously from factorizations of groups to group factorizations. Also, now that the order of the factors is given, we call *Frattini Factorization* a group factorization $(U, H, K)$ such that $U$ is a symmetric or alternating group and $K$ is one of its point stabilizers. Similarly, the *homogeneous factorizations* are the $(U, H, K)$ where $U$ is as above a symmetric or alternating group of some degree, say $t$, and $K$ is the stabilizer in $U$ of a subset of some order $k$, $1 \leqslant k < t$. In [**LPS90**] it is shown that all the core-free and maximal factorizations of the symmetric or alternating groups of degree larger than 10 are either homogeneous factorizations or their duals.

Two group factorizations $(G, H, K)$ and $(E, P, Q)$ are isomorphic if there is a group isomorphism $\boldsymbol{f} : G \longrightarrow E$ such that $H\boldsymbol{f} = P$ and $K\boldsymbol{f} = Q$. Of course, the interest focus on the isomorphism type of a group factorization $(G, H, K)$ which we denote by $[G, H, K]$. Note that if $g \in G$, then

$$[G, H, K^g] = [G, H, K] = [G, H^g, K].$$

To prove the equation on the right for example, note that $g = hk$ for convenient $h \in H$ and $k \in K$, so, $(G, H^g, K) = (G^k, H^k, K^k)$. The following result is also obvious and requires no proof.

6.1.1. PROPOSITION. *If $(G, H, K)$ and $(E, P, Q)$ are isomorphic group factorizations, then $/\binom{G}{K\diamond} : H/$ and $/\binom{E}{Q\diamond} : P/$ are permutationally isomorphic and the intervals $\left[/\binom{G}{K\diamond} : H/ \div \operatorname{Sym}\binom{G}{K\diamond}\right]$ and $\left[/\binom{E}{Q\diamond} : P/ \div \operatorname{Sym}\binom{E}{Q\diamond}\right]$ are lattice isomorphic.*

Note that if $(G, H, K)$ is a Frattini Factorization, then $H$ and $/\binom{G}{K\diamond} : H/$ are permutationally isomorphic.

Suppose now that $HK$ is a core-free and maximal factorization of the almost simple group $G$ of socle $L$. We may assume that $G$ is an automorphic extension of $L$, that is, $G$ is a group of automorphisms of $L$ containing $\operatorname{Inn} L$, the group of the inner automorphisms of $L$. Now, identify $L$ with $\operatorname{Inn} L$; put $X = L \cap H$ and $Y = L \cap K$. The ordered triple $(L, X, Y)$ is then called the *skeleton* of $(G, H, K)$. Of course, $H$ and $K$ may be recovered from $G$, $X$, $Y$ as $H = \mathrm{N}_G X$ and $K = \mathrm{N}_G Y$. Indeed, for each automorphic extension $E$ of $L$, with only occasional few exceptions, $(\mathrm{N}_E X)(\mathrm{N}_E Y)$ is a core-free and maximal factorization of $E$. For this reason, a classification of the core-free and maximal factorizations of the almost simple groups (except the Frattini Factorizations) is given by listing the skeletons $(L, X, Y)$ of the related group factorizations together with notes ruling out the automorphic extensions $E$ of $L$ such that $(\mathrm{N}_E X)(\mathrm{N}_E Y)$ is not a maximal factorization of $E$ (some authors occasionally use smaller subgroups in the roles of our $X$ and $Y$, but we shall keep to the convention just stated).

Since for each automorphism $\boldsymbol{\lambda}$ of $L$ and factorization $HK$ of the automorphic extension $E$ of $L$, $[E, H, K] = [E^{\boldsymbol{\lambda}}, H^{\boldsymbol{\lambda}}, K^{\boldsymbol{\lambda}}]$, one may describe the skeletons $(L, X, Y)$ up to automorphisms of $L$. This is carried out in [**LPS90**], while a classification of the remaining maximal but non core-free factorizations of the almost simple groups (except the Frattini Factorizations) appears in [**LPS96**].

## 6.2. Proof of Theorem C

Suppose that $\boldsymbol{H} \lessdot E \lessdot \boldsymbol{U}$ where both $\boldsymbol{H}$ and $E$ are almost simple. Put $M = \operatorname{Soc} \boldsymbol{H}$ and $L = \operatorname{Soc} E$.

6.2.1. Lemma. *Either $E = \mathrm{N}_{\boldsymbol{U}} M$ or $M < L$.*

Proof. If $M = L$, then $E$ is contained in $\mathrm{N}_{\boldsymbol{U}} M$ and so, by maximality of $E$ in $\boldsymbol{U}$, $E = \mathrm{N}_{\boldsymbol{U}} M$. Otherwise, $M \neq L$ and we now show that in fact $M < L$. If this fails, then $M \cap L = \{1\}$ because $M \cap L$ is normal in $M$. Then $E/L$ has a subgroup isomorphic to the non abelian simple group $M$, which is impossible by the Schreier Conjecture. This contradiction ensures that $M < L$ when they are different. $\square$

6.2.2. Lemma. *If $M < L$, then there is a subgroup $Q$ of $E$ such that $(E, \boldsymbol{H}, Q)$ is a non Frattini, core-free and maximal group factorization such that $\boldsymbol{H} \cap Q \lessdot \boldsymbol{H}$. Moreover, if $(L, X, Y)$ is the skeleton of $(E, \boldsymbol{H}, Q)$, then $M$ and $/\binom{L}{Y\diamond} : M/$ are permutationally isomorphic.*

Proof. Call $Q$ the stabilizer of a point in $E$. By Frattini Argument $E = \boldsymbol{H}Q$. Since $\boldsymbol{H}$ is primitive, $\boldsymbol{H} \cap Q \lessdot \boldsymbol{H}$. Since $E$ is primitive, $Q \lessdot E$. By hypothesis $\boldsymbol{H} \lessdot E$, so $(E, \boldsymbol{H}, Q)$ is a maximal group factorization such that $\boldsymbol{H} \cap Q \lessdot \boldsymbol{H}$. It is not a Frattini Factorization because $E \lessdot \boldsymbol{U}$. Clearly, $Q$ being a point stabilizer is core-free in $E$. Regarding $\boldsymbol{H}$, observe that $G \trianglelefteq E$ implies $L \leqslant G$. Thus either $\operatorname{Core}_E \boldsymbol{H} = \{1\}$ or $L \leqslant \boldsymbol{H}$. But $L \leqslant \boldsymbol{H}$ implies $M \trianglelefteq L$ and hence $M = L$ which is against our assumptions. Therefore both $\boldsymbol{H}$ and $Q$ are core-free in $E$. To finish observe that $Y = L \cap Q$ is the stabilizer of a point in $L$, thus

$$/\tbinom{L}{Y\diamond} : \star / : L \longrightarrow \operatorname{Sym} \tbinom{L}{Y\diamond}$$

is a permutational isomorphism. Therefore $M$ and $/\binom{L}{Y\diamond} : M/$ are permutationally isomorphic. $\square$

Now, if $M < L$, say $(L, X, Y)$ the skeleton of the maximal group factorization $(E, \boldsymbol{H}, Q)$ of the proof above. Put

$$G = \mathrm{N}_{\operatorname{Sym}\binom{L}{Y\diamond}} /\tbinom{L}{Y\diamond} : L/, \qquad N = \mathrm{N}_{\operatorname{Sym}\binom{L}{Y\diamond}} /\tbinom{L}{Y\diamond} : M/, \qquad \text{and } H = G \cap N.$$

Provided $\boldsymbol{H}$ is second maximal in $\boldsymbol{U}$, either

(1) $E \cong G$ and $\boldsymbol{H} \cong H$, or
(2) $E \cong G^{\mathrm{e}}$ and $\boldsymbol{H} \cong H^{\mathrm{e}}$.

Of course, at least one of $H$ and $H^{\mathrm{e}}$ is high in this case. Note that if $M = L$ and $M = \operatorname{Soc} X$ for an $X$ appearing in one of the skeletons $(L, X, Y)$ as above, then with the same notation we have either $E \cong N$ or $E \cong N^{\mathrm{e}}$.

This reduces the classification up to permutational isomorphism of the $\boldsymbol{H}$ as in Theorem C (and of its overgroups) to the investigation of the core-free and maximal group factorization $(E, P, Q)$ such that both $E$ and $P$ are almost simple and $P \cap Q \lessdot P$ (except the Frattini Factorizations). These are described in Table 6.A by listing their skeletons $(L, X, Y)$. In fact, Table 6.A is a repetition of the list given in [**LPS87**] enriched with the normalizers above. Table 6.B is a copy of Table 6.A but stripped of the skeletons for which we know that both $H$ and $H^{\mathrm{e}}$ are low. Moreover, Table 6.B is sorted by the permutational isomorphism type of $\operatorname{Soc} X$, so, one easily recovers up to permutational isomorphism all the subgroups of $\boldsymbol{U}$ containing a given $\boldsymbol{H}$, by choosing the bunch of consecutive lines of the table having $\operatorname{Soc} X$ isomorphic to $\operatorname{Soc} \boldsymbol{H}$. In fact, we carefully subdivided the skeletons in a way that if a $\boldsymbol{H}$ matches two different skeletons, then they must appear in the same bunch. However, this does not necessarily mean that skeletons of a same bunch induce permutationally isomorphic $\boldsymbol{H}$.

The number of permutationally isomorphic subgroups of $\boldsymbol{A}$ containing a given $H^{\mathrm{e}}$ for example, is given by $\operatorname{hm}(G^{\mathrm{e}}, H^{\mathrm{e}})$, which appear in column "$r_e$", so the sum of column "$r_e$", restricted to the bunch related to $H$, gives an upper bound for the number of maximal subgroups of $\boldsymbol{A}$ which contain $H^{\mathrm{e}}$ and do not normalize $\operatorname{Soc} H$. However, by lemma 6.2.1, there is at most 1 maximal subgroup of $\boldsymbol{A}$ normalizing $\operatorname{Soc} H$. Similarly, using column "$r_o$", one finds an upper bound for the number of maximal subgroups of $\boldsymbol{S}$ containing an odd $H$. This completes the proof of Theorem C.

## 6.3. Notation in the tables

Tables list the skeletons $(L, X, Y)$ (up to automorphisms of $L$) of the core-free and maximal group factorizations $(E, P, Q)$ such that both $E$ and $P$ are almost simple and $P \cap Q \lessdot P$, except the Frattini Factorizations.

Each skeleton (or parametrised family of skeletons) extends on one or more physical rows and has a number which appears in first column. For example, we will refer to the skeletons having number 17 as the skeletons in $\sharp 17$.

There are three more reference columns. Columns denoted by "[**LPS87**]", "[**LPS90**]" report the reference to the correspondent group factorizations as given in those works. For example, a "III.17" in column "[**LPS87**]" means that this is the skeleton of the group factorization at line 17 of Table III of [**LPS87**]. We tried to use the notations of those works whenever possible. In particular, [**LPS90**] describes certain subgroups $P_i$, $N_i^{\pm}$ of the classical groups and we use these names where it is appropriate to do so. The last reference column is denoted by "!" and points to a part of this work where the reader may find some interesting comments or explanations. These explanations are given following the order given in Table 6.A, so, in reading about $\sharp n$, we assume that the reader is aware of the explanations given for $\sharp m$ for any $m < n$. Accordingly, a missing entry in column "!" implies that entries of the other columns were determined in an entirely similar way to the correspondent entries of some previous case.

Note that each skeleton $(L, X, Y)$ is paired to a skeleton $(L, Y, X)$ unless the last is the skeleton of a Frattini Factorization or $Y$ is not primitive of type almost simple. This correspondence is made explicit in column "↻" where we write "np" if $Y$ is not primitive, "nas" if $Y$ is not almost simple and "Ff" if $(L, Y, X)$ is the skeleton of a Frattini Factorization.

The column denoted by "deg" report the degree of the permutation group $/\binom{L}{Y\diamond} : L/$, that is, the index of $Y$ in $L$. We may omit this entry if the degree is too large. Other times, we may write something like "$2^n\star$" meaning that $2^n$ is the highest power of 2 dividing that degree.

Both $X$ and $P$ are almost simple with the same socle and there is a column for this socle which we call $M$. The difficulty here is that isomorphic simple groups may have different names for example $\mathrm{A}_6$ and $\mathrm{PSL}(2,9)$ or $\Omega^-(4,q)$ and $\mathrm{PSL}(2,q^2)$ just to name a few. As we want to sort Table 6.A according to this column (producing Table 6.B), it is sensible to have unique names for each isomorphism type of simple groups. We obtain unique name by choosing representatives in this order:

(1) alternating groups,
(2) linear groups,
(3) unitary groups,
(4) symplectic groups,
(5) orthogonal groups,
(6) Suzuki and Ree groups,
(7) remaining exceptional simple groups of Lie type,
(8) sporadic simple groups.

And if there are isomorphisms between groups of the same class as we distinguished above, we choose in the given order the name correspondent to

(1) smallest dimension,
(2) smallest characteristic of the field.

For example, the name chosen for $\mathrm{SL}(3,2) \cong \mathrm{PSL}(2,7)$ is $\mathrm{PSL}(2,7)$. All this has an important consequence: parametrised families of skeletons are partitioned according to our choice of names. For example the parametrised family of skeletons $\big(\mathrm{A}_{q+1}, \mathrm{PSL}(2,q), \mathrm{A}_{q+1} \cap (\mathrm{S}_2 \times \mathrm{S}_{q-1})\big)$ is split in the family with all those skeletons but the one for $q = 9$, plus the skeleton $\big(\mathrm{A}_{10}, \mathrm{A}_6, \mathrm{A}_{10} \cap (\mathrm{S}_2 \times \mathrm{S}_8)\big)$ in view of the isomorphism $\mathrm{PSL}(2,9) \cong \mathrm{A}_6$. Of course, this results in a refinement of the original list given in [**LPS87**].

Then we have columns for $G = \mathrm{N}_{\boldsymbol{S}}\, L$ where $\boldsymbol{S}$ is the symmetric group on $\binom{L}{Y\diamond}$ and $L$ is identified with $/\binom{L}{Y\diamond} : L/$, for $N = \mathrm{N}_{\boldsymbol{S}}\, M$ and for $H = \mathrm{N}_G\, M$. We may omit to write in these columns when the normalizers are not needed to determine $r_o$ and $r_e$ where

$$r_o = \begin{cases} \mathrm{hm}(G, H) & \text{if } H \text{ is odd and second maximal in } \boldsymbol{S}, \\ 0 & \text{otherwise} \end{cases}$$

and

$$r_e = \begin{cases} \mathrm{hm}(G^{\mathrm{e}}, H^{\mathrm{e}}) & \text{if } H^{\mathrm{e}} \text{ is second maximal in } \boldsymbol{A}, \\ 0 & \text{otherwise.} \end{cases}$$

Of course, $\boldsymbol{A}$ is the alternating group related to $\boldsymbol{S}$, that is, $\mathrm{Alt}\binom{L}{Y\diamond}$. It is understood that a missing or non numerical entry in column "$r_o$" implies $r_o = 0$. The same holds for column "$r_e$".

TABLE 6.A. Skeletons $(L, X, Y)$ of the core-free and maximal group factorizations $(E, P, Q)$ such that both $E$ and $P$ are almost simple and $P \cap Q \lessdot P$, except the Frattini Factorizations

| ♯ | $L$ | $X$ | $Y$ | $\circlearrowright$ | deg | [**LPS87**] | [**LPS90**] |
|---|---|---|---|---|---|---|---|
| 1 | $A_8$ | $A_7$ | $AGL_3(2)$ | nas | 15 | III.1 | Th. D |
| 2 | $A_8$ | $A_8 \cap (S_2 \times S_6)$ | $AGL_3(2)$ | nas | 15 | III.2 | Th. D |
| 3 | $A_7$ | $A_6$ | $GL_3(2)$ | Ff | 15 | III.3 | Th. D |
| 4 | $A_9$ | $A_9 \cap (S_2 \times S_7)$ | $P\Gamma L_2(8)$ | 31 | 120 | III.4 | Th. D |
| 5 | $A_{12}$ | $A_{11}$ | $M_{12}$ | Ff | 2520 | III.5 | Th. D |
| 6 | $A_{12}$ | $A_{12} \cap (S_2 \times S_{10})$ | $M_{12}$ | 20 | 2520 | III.6 | Th. D |
| 7 | $A_{11}$ | $A_{10}$ | $M_{11}$ | Ff | 2520 | III.7 | Th. D |
| 8 | $A_{24}$ | $A_{23}$ | $M_{24}$ | Ff | $2^{11}\star$ | III.8 | Th. D |
| 9 | $A_{24}$ | $A_{24} \cap (S_2 \times S_{22})$ | $M_{24}$ | 26 | $2^{11}\star$ | III.9 | Th. D |
| 10 | $A_{23}$ | $A_{22}$ | $M_{23}$ | Ff | $2^{11}\star$ | III.10 | Th. D |
| 11 | $A_{176}$ | $A_{175}$ | HS | Ff | $2^{163}\star$ | III.11 | Th. D |
| 12 | $A_{276}$ | $A_{275}$ | $Co_3$ | Ff | $2^{262}\star$ | III.12 | Th. D |
| 13 | $A_c$ | $A_{c-1}$ | $Sp_{2d}(2)$ | Ff | | III.13 | Th. D |
| | $c = 2^{2d-1} \pm 2^{d-1}$, $d \geqslant 3$. | | | | | | |
| 14 | $A_8$ | $A_7$ | $PSL_2(7)$ | Ff | 120 | III.14 | Th. D |
| 15 | $A_{p+1}$ | $A_p$ | $PSL_2(p)$ | Ff | $(p-2)!$ | III.14 | Th. D |
| | $p$ prime, $p > 7$; for $p = 7$ see ♯14. | | | | | | |
| 16 | $A_{2^d}$ | $A_{2^d-1}$ | $AGL_d(2)$ | Ff | | III.15 | Th. D |
| | $d > 3$; for $d = 3$ see ♯1. | | | | | | |
| 17 | $A_{2l}$, $l \geqslant 3$ | $A_{2l-1}$ | $A_{2l} \cap \lfloor S_l \operatorname{wr} S_2 \rfloor$ | Ff | $\frac{1}{2}\binom{2l}{l}$ | III.16 | Th. D |
| 18 | $A_{11}$ | $M_{11}$ | $A_{11} \cap (S_2 \times S_9)$ | np | 55 | III.17 | Th. D |
| 19 | $A_{11}$ | $M_{11}$ | $A_{11} \cap (S_3 \times S_8)$ | np | 165 | III.17 | Th. D |
| 20 | $A_{12}$ | $M_{12}$ | $A_{12} \cap (S_2 \times S_{10})$ | 6 | 66 | III.18 | Th. D |
| 21 | $A_{12}$ | $M_{12}$ | $A_{12} \cap (S_3 \times S_9)$ | np | 220 | III.18 | Th. D |
| 22 | $A_{12}$ | $M_{12}$ | $A_{12} \cap (S_4 \times S_8)$ | np | 495 | III.18 | Th. D |
| 23 | $A_{22}$ | $M_{22}$ | $A_{22} \cap (S_2 \times S_{20})$ | np | 231 | III.19 | Th. D |

| $M$ | $G$ | $N$ | $H$ | $r_o$ | $r_e$ | ! | ♯ |
|---|---|---|---|---|---|---|---|
| $A_7$ | $L$ | $M$ | $M$ | | 1 | 6.4.1 | 1 |
| $A_6$ | $L$ | $X$ | $X$ | | 1 | 6.4.1 | 2 |
| $A_6$ | $L$ | $A_8 \cap (S_2 \times S_6)$ | $M$ | $H < L < A_8 < \boldsymbol{A}$ | | | 3 |
| $A_7$ | $L$ | $X$ | $X$ | | 1 | 6.4.1 | 4 |
| $A_{11}$ | $L$ | $M$ | $M$ | | 1 | 6.4.1 | 5 |
| $A_{10}$ | $L$ | $X$ | $X$ | | 1 | 6.4.1 | 6 |
| $A_{10}$ | $L$ | $A_{12} \cap (S_2 \times S_{10})$ | $M$ | $H < L < A_{12} < \boldsymbol{A}$ | | | 7 |
| $A_{23}$ | $L$ | $M$ | $M$ | | 1 | 6.4.2 | 8 |
| $A_{22}$ | $L$ | $X$ | $X$ | | 1 | 6.4.2 | 9 |
| $A_{22}$ | $L$ | $A_{24} \cap (S_2 \times S_{22})$ | $M$ | $H < L < A_{24} < \boldsymbol{A}$ | | | 10 |
| $A_{175}$ | $L$ | $M$ | $M$ | | 1 | 6.4.3 | 11 |
| $A_{275}$ | $L$ | $M$ | $M$ | | 1 | 6.4.3 | 12 |
| $A_{c-1}$ | $L$ | $M$ | $M$ | | 1 | 6.4.4 | 13 |
| $A_7$ | $S_8 < \boldsymbol{A}$ | $S_7$ | $N$ | | 1 | 6.4.5 | 14 |
| $A_p$ | $S_{p+1}$ | $S_p$ | $N$ | $\leqslant 1$ | 1 | 6.4.5 | 15 |
| $A_{2^d-1}$ | $L$ | $M$ | $M$ | | 1 | 6.4.6 | 16 |
| $A_{2l-1}$ | $S_{2l} < \boldsymbol{A}$ | $S_{2l-1}$ | $N$ | | 1 | 6.4.7 | 17 |
| $M_{11}$ | $S_{11}$ | $M$ | $M$ | | 1 | | 18 |
| $M_{11}$ | $S_{11} < \boldsymbol{A}$ | $M$ | $M$ | $H < L < G < \boldsymbol{A}$ | | | 19 |
| $M_{12}$ | $S_{12} < \boldsymbol{A}$ | $M$ | $M$ | $H < L < G < \boldsymbol{A}$ | | | 20 |
| $M_{12}$ | $S_{12}$ | $M$ | $M$ | | 1 | | 21 |
| $M_{12}$ | $S_{12} < \boldsymbol{A}$ | $M.2$ | $M$ | $H < L < G < \boldsymbol{A}$ | | | 22 |
| $M_{22}$ | $S_{22} < \boldsymbol{A}$ | $M.2$ | $N$ | | 1 | | 23 |

**Table 6.A continued**

| ♯ | $L$ | $X$ | $Y$ | $\circlearrowright$ | deg | [**LPS87**] | [**LPS90**] |
|---|---|---|---|---|---|---|---|
| 24 | $\mathrm{A}_{23}$ | $\mathrm{M}_{23}$ | $\mathrm{A}_{23} \cap (\mathrm{S}_2 \times \mathrm{S}_{21})$ | np | 253 | III.20 | Th. D |
| 25 | $\mathrm{A}_{23}$ | $\mathrm{M}_{23}$ | $\mathrm{A}_{23} \cap (\mathrm{S}_3 \times \mathrm{S}_{20})$ | np | 1771 | III.20 | Th. D |
| 26 | $\mathrm{A}_{24}$ | $\mathrm{M}_{24}$ | $\mathrm{A}_{24} \cap (\mathrm{S}_2 \times \mathrm{S}_{22})$ | 9 | 276 | III.21 | Th. D |
| 27 | $\mathrm{A}_{24}$ | $\mathrm{M}_{24}$ | $\mathrm{A}_{24} \cap (\mathrm{S}_3 \times \mathrm{S}_{21})$ | np | 2024 | III.21 | Th. D |
| 28 | $\mathrm{A}_{176}$ | HS | $\mathrm{A}_{176} \cap (\mathrm{S}_2 \times \mathrm{S}_{174})$ | np | 15400 | III.22 | Th. D |
| 29 | $\mathrm{A}_{276}$ | $\mathrm{Co}_3$ | $\mathrm{A}_{276} \cap (\mathrm{S}_2 \times \mathrm{S}_{274})$ | np | 37950 | III.23 | Th. D |
| 30 | $\mathrm{A}_{10}$ | $\mathrm{A}_6.2_3$ | $\mathrm{A}_{10} \cap (\mathrm{S}_2 \times \mathrm{S}_8)$ | np | 45 | III.24 | Th. D |
| 31 | $\mathrm{A}_{q+1}$ | $\mathrm{N}_L(\mathrm{PSL}_2\, q)$ | $\mathrm{A}_{q+1} \cap (\mathrm{S}_2 \times \mathrm{S}_{q-1})$ | np, 4 | $\frac{q(q+1)}{2}$ | III.24 | Th. D |
| | $q > 5$, $q \neq 9$; for $q = 9$ see ♯30. | | | | | | |
| 32 | $\mathrm{A}_{q^2+1}$ | $\mathrm{Aut}\,\mathrm{Sz}(q)$ | $\mathrm{A}_{q^2+1} \cap (\mathrm{S}_2 \times \mathrm{S}_{q^2-1})$ | np | $\frac{q^2(q^2+1)}{2}$ | III.25 | Th. D |
| | $q = 2^{2a+1} \geqslant 8$. | | | | | | |
| 33 | $\mathrm{M}_{11}$ | $\mathrm{PSL}_2(11)$ | $\mathrm{M}_{10}$ | np | 11 | IV.1 | 6.1 |
| 34 | $\mathrm{M}_{11}$ | $\mathrm{PSL}_2(11)$ | $\mathrm{M}_9 : 2$ | np | 55 | IV.2 | 6.1 |
| 35 | $\mathrm{M}_{12}$ | $\mathrm{M}_{11}$ | $\mathrm{M}_{11}$ | 35 | 12 | IV.3 | 6.2 |
| 36 | $\mathrm{M}_{12}$ | $\mathrm{PSL}_2(11)$ | $\mathrm{M}_{11}$ | np | 12 | IV.4 | 6.2 |
| 37 | $\mathrm{M}_{12}$ | $\mathrm{M}_{11}$ | $\mathrm{M}_{10} : 2$ | np | 66 | IV.5 | 6.2 |
| 38 | $\mathrm{M}_{24}$ | $\mathrm{M}_{23}$ | $\mathrm{M}_{12} : 2$ | np | 1288 | IV.6 | 6.5 |
| 39 | $\mathrm{M}_{24}$ | $\mathrm{M}_{23}$ | $2^6 : 3.\mathrm{S}_6$ | np | 1771 | IV.7 | 6.5 |
| 40 | $\mathrm{M}_{24}$ | $\mathrm{M}_{23}$ | $\mathrm{PSL}_2(23)$ | 41 | 40320 | IV.8 | 6.5 |
| 41 | $\mathrm{M}_{24}$ | $\mathrm{PSL}_2(23)$ | $\mathrm{M}_{23}$ | 40 | 24 | IV.9 | 6.5 |
| 42 | $\mathrm{M}_{24}$ | $\mathrm{PSL}_2(23)$ | $\mathrm{M}_{22} : 2$ | np | 276 | IV.10 | 6.5 |
| 43 | HS | $\mathrm{M}_{22}$ | $\mathrm{PSU}_3(5) : 2$ | np | 176 | IV.11 | 6.7 |
| 44 | $\mathrm{PSL}_{2m}(q)$ | $\mathrm{N}_L\, \mathrm{PSp}_{2m}(q)$ | $P_1$ | nas | $\frac{q^{2m}-1}{q-1}$ | VI.1 | 3.1.1 |
| | $m \geqslant 2$ but $q \neq 2$ when $m = 2$. | | | | | | |
| 45 | $\mathrm{Sp}_{2m}(q)$ | $\mathrm{Sp}_{2a}(q^b).b$ | $\mathrm{O}^{\pm}_{2m}(q)$ | 47, 48, np | $\frac{q^m(q^m \pm 1)}{2}$ | VI.2 | 3.2.1d |
| | $q$ even; $m = ab > b$, $b$ prime. | | | | | | |
| 46 | $\mathrm{Sp}_{2m}(q)$ | $\mathrm{Sp}_2(q^m).m$ | $\mathrm{O}^{\pm}_{2m}(q)$ | 49, np | $\frac{q^m(q^m \pm 1)}{2}$ | VI.2 | 3.2.1d |
| | $q$ even; $m$ an odd prime. | | | | | | |
| 47 | $\mathrm{Sp}_{2m}(q)$ | $\mathrm{O}^{\pm}_{2m}(q)$ | $\mathrm{Sp}_{2a}(q^b).b$ | 45 | | VI.3 | 3.2.1d |
| | $q$ even; $m = ab > b > 2$, $b$ prime. | | | | | | |
| 48 | $\mathrm{Sp}_{4a}(q)$ | $\mathrm{O}^{-}_{4a}(q)$ | $\mathrm{Sp}_{2a}(q^2).2$ | 45 | | VI.3 | 3.2.1d |
| | $q$ even; $a > 1$. | | | | | | |
| 49 | $\mathrm{Sp}_4(q)$ | $\mathrm{O}^{-}_4(q)$ | $\mathrm{Sp}_2(q^2).2$ | 46 | $\frac{q^2(q^2-1)}{2}$ | VI.3 | 3.2.1d |
| | $q$ even; $q \neq 2$. | | | | | | |

| $M$ | $G$ | $N$ | $H$ | $r_o$ | $r_e$ | ! | ♯ |
|---|---|---|---|---|---|---|---|
| $M_{23}$ | $S_{23}$ | $M$ | $M$ | | 1 | | 24 |
| $M_{23}$ | $S_{23} < \boldsymbol{A}$ | $M$ | $M$ | $H < L < G < \boldsymbol{A}$ | | | 25 |
| $M_{24}$ | $S_{24} < \boldsymbol{A}$ | $M$ | $M$ | $H < L < G < \boldsymbol{A}$ | | | 26 |
| $M_{24}$ | $S_{24}$ | $M$ | $M$ | | 1 | | 27 |
| HS | $S_{176} < \boldsymbol{A}$ | $M$ | $M$ | $H < L < G < \boldsymbol{A}$ | | | 28 |
| $Co_3$ | $S_{276} < \boldsymbol{A}$ | $M$ | $M$ | $H < L < G < \boldsymbol{A}$ | | | 29 |
| $A_6$ | $S_{10} < \boldsymbol{A}$ | $M.2^2$ | $N$ | | 1 | 6.4.8 | 30 |
| $PSL_2(q)$ | $S_{q+1}$ | $P\Gamma L_2(q)$ | $N$ | | 1 | 6.4.8 | 31 |
| $Sz(q)$ | $S_{q^2+1}$ | $X$ | $X$ | | 1 | 3.10 | 32 |
| $PSL_2(11)$ | $L$ | $M$ | $M$ | | 1 | | 33 |
| $PSL_2(11)$ | $L$ | $M.2$ | $M$ | | 2 | | 34 |
| $M_{11}$ | $L$ | $M$ | $M$ | | 1 | 6.4.9 | 35 |
| $PSL_2(11)$ | $L$ | $M.2$ | $M$ | | 2 | | 36 |
| $M_{11}$ | $L$ | $M$ | $M$ | | 1 | 6.4.9 | 37 |
| $M_{23}$ | $L$ | $M$ | $M$ | | 1 | | 38 |
| $M_{23}$ | $L$ | $M$ | $M$ | $H = H^{\mathrm{e}}$ low, see ♯25 | | | 39 |
| $M_{23}$ | $L$ | $M$ | $M$ | | 1 | | 40 |
| $PSL_2(23)$ | $L$ | $M.2$ | $M$ | | 2 | | 41 |
| $PSL_2(23)$ | $L$ | $M.2$ | $M$ | $H < L < A_{24} < \boldsymbol{A}$ | | 6.4.10 | 42 |
| $M_{22}$ | $L$ | $M$ | $M$ | | 1 | | 43 |
| $PSp_{2m}(q)$ | $P\Gamma L_{2m}(q)$ | $P\Gamma Sp_{2m}(q)$ | $N$ | $\leqslant 1$ | $\leqslant 2$ | 6.4.11 | 44 |
| $PSp_{2a}(q^b)$ | $P\Gamma Sp_{2m}(q)$ | $P\Gamma Sp_{2a}(q^b)$ | $N$ | $\leqslant 1$ | 1 | 6.4.12 | 45 |
| $PSL_2(q^m)$ | $P\Gamma Sp_{2m}(q)$ | $P\Gamma L_2(q^m)$ | $N$ | $\leqslant 1$ | 1 | 6.4.12 | 46 |
| $\Omega^{\pm}_{2m}(q)$ | $P\Gamma Sp_{2m}(q)$ | $P\Gamma O^{\pm}_{2m}(q)$ | $N$ | $\leqslant 1$ | 1 | 6.4.13 | 47 |
| $\Omega^{-}_{4a}(q)$ | $P\Gamma Sp_{4a}(q)$ | $P\Gamma O^{-}_{4a}(q)$ | $N$ | $\leqslant 1$ | 1 | 6.4.13 | 48 |
| $PSL_2(q^2)$ | $P\Gamma Sp_4(q)$ | $P\Gamma O^{-}_4(q)$ | $N$ | $\leqslant 1$ | 1 | 6.4.13 | 49 |

**Table 6.A continued**

| ♯ | $L$ | $X$ | $Y$ | ↻ | deg | [**LPS87**] | [**LPS90**] |
|---|---|---|---|---|---|---|---|
| 50 | $\mathrm{Sp}_{2m}(q)$ | $\mathrm{O}^-_{2m}(q)$ | $\mathrm{Sp}_m(q) \,\mathrm{wr}\, \mathrm{S}_2$ | np | | VI.4 | 3.2.4b |
| | $q$ even, $m$ even, $m \neq 2$. | | | | | | |
| 51 | $\mathrm{Sp}_4(q)$ | $\mathrm{O}^-_4(q)$ | $\mathrm{Sp}_2(q) \,\mathrm{wr}\, \mathrm{S}_2$ | np | $\frac{q^2(q^2+1)}{2}$ | VI.4 | 3.2.4b |
| | $q$ even, $q \neq 2$. | | | | | | |
| 52 | $\mathrm{Sp}_{2m}(2)$ | $\mathrm{O}^{\pm}_{2m}(2)$ | $\mathrm{O}^{\mp}_{2m}(2)$ | 52 | $2^{m-1}(2^m \mp 1)$ | VI.5 | 3.2.4e |
| | $m > 3$. Note that $X \cap Y$ is a subgroup of type $N_1$ of $X$. | | | | | | |
| 53 | $\mathrm{Sp}_6(2)$ | $\mathrm{O}^+_6(2)$ | $\mathrm{O}^-_6(2)$ | 54 | 28 | VI.5 | 3.2.4e |
| 54 | $\mathrm{Sp}_6(2)$ | $\mathrm{O}^-_6(2)$ | $\mathrm{O}^+_6(2)$ | 53 | 36 | VI.5 | 3.2.4e |
| 55 | $\mathrm{P}\Omega^+_{2m}(q)$ | $\mathrm{St}_L\{\mathsf{E},\mathsf{F}\}$ | $N_1$ | np | $\frac{q^{m-1}(q^m-1)}{(2,q-1)}$ | VI.6 | 3.6.1b |
| | $q = 2, 3; m \geqslant 5; \mathsf{V} = \mathsf{E} \oplus \mathsf{F}$ with $\mathsf{E}, \mathsf{F}$ tot. singular. | | | | | | |
| 56 | $\mathrm{P}\Omega^+_{2m}(q)$ | $N_1$ | $P_m, P_{m-1}$ | np | $\prod_{i=1}^{m-1}(q^i+1)$ | VI.7 | 3.6.1a |
| | $q$ odd; $m \geqslant 4$. Note that if $m = 4$, then $Y$ is conjugate in $\mathrm{Aut}\, L$ to $P_1$. | | | | | | |
| 57 | $\Omega^+_{2m}(q)$ | $N_1$ | $P_m, P_{m-1}$ | np | $\prod_{i=1}^{m-1}(q^i+1)$ | VI.7 | 3.6.1a |
| | $q$ even; $m \geqslant 4$. Note that if $m = 4$, then $Y$ is conjugate in $\mathrm{Aut}\, L$ to $P_1$. | | | | | | |
| 58 | $\mathrm{P}\Omega^+_8(q)$ | $N_1$ | $X^\tau$ | 58 | $\frac{q^3(q^4-1)}{2}$ | VI.8 | 5.1.15 |
| | $q$ odd; $\tau$ is a triality automorphism. | | | | | | |
| 59 | $\Omega^+_8(q)$ | $N_1$ | $X^\tau$ | 59 | $q^3(q^4-1)$ | VI.8 | 5.1.15 |
| | $q$ even; $\tau$ is a triality automorphism. | | | | | | |
| 60 | $\Omega_7(q)$ | $\mathrm{G}_2(q)$ | $N_1^+, N_1^-$ | np, 77 | $\frac{q^3(q^3\pm1)}{2}$ | VI.9 | 5.1.14 |
| | $q$ odd. | | | | | | |
| 61 | $\mathrm{Sp}_6(q)$ | $\mathrm{G}_2(q)$ | $\mathrm{O}^{\pm}_6(q)$ | np, 81 | $\frac{q^3(q^3\pm1)}{2}$ | VI.9 | 5.2.3b |
| | $q$ even. | | | | | | |
| 62 | $\Omega_7(q)$ | $\mathrm{G}_2(q)$ | $P_1$ | np | $\frac{q^6-1}{q-1}$ | VI.10 | 5.1.14 |
| | $q$ odd. | | | | | | |
| 63 | $\mathrm{Sp}_6(q)$ | $\mathrm{G}_2(q)$ | $P_1$ | np | $\frac{q^6-1}{q-1}$ | VI.10 | 5.2.3b |
| | $q$ even. | | | | | | |
| 64 | $\mathrm{Sp}_6(q)$ | $\mathrm{G}_2(q)$ | $N_2$ | np | $\frac{q^4(q^6-1)}{q^2-1}$ | VI.11 | 5.2.3b |
| | $q$ even, $q \neq 2$. | | | | | | |
| 65 | $\mathrm{Sp}_4(q)$ | $\mathrm{Sz}(q)$ | $\mathrm{O}^+_4(q)$ | np | $\frac{q^2(q^2+1)}{2}$ | VI.12 | 5.1.7b |
| | $q = 2^{2a+1} \geqslant 8$. | | | | | | |

| $M$ | $G$ | $N$ | $H$ | $r_o$ | $r_e$ | ! | ♯ |
|---|---|---|---|---|---|---|---|
| $\Omega^-_{2m}(q)$ | $\mathrm{P\Gamma Sp}_{2m}(q)$ | $\mathrm{P\Gamma O}^-_{2m}(q)$ | $N$ | $\leqslant 1$ | 1 | 6.4.14 | 50 |
| $\mathrm{PSL}_2(q^2)$ | $\mathrm{P\Gamma Sp}_4(q)$ | $\mathrm{P\Gamma O}^-_4(q)$ | $N$ | $\leqslant 1$ | 1 | 6.4.14 | 51 |
| $\Omega^{\pm}_{2m}(2)$ | $L$ | $X$ | $X$ | | 1 | | 52 |
| $\mathrm{A}_8$ | $L$ | $X$ | $X$ | | 1 | | 53 |
| $\mathrm{PSp}_4(3)$ | $L$ | $X$ | $X$ | | 1 | | 54 |
| $\mathrm{PSL}_m(q)$ | $\mathrm{PO}^+_{2m}(q)$ | $\mathrm{Aut}\,M$ | $N$ | $\leqslant 2$ | $\leqslant 2$ | 6.4.15 | 55 |
| $\Omega_{2m-1}(q)$ | $\langle \mathrm{PSO}, \bar\delta, \bar\phi \rangle$ | $\mathrm{Aut}\,M$ | $\langle X, \bar\phi \rangle$ | | $\leqslant 1$ | 6.4.16 | 56 |
| $\mathrm{Sp}_{2m-2}(q)$ | $\langle L, \bar\phi \rangle$ | $\mathrm{Aut}\,M$ | $N$ | $\leqslant 1$ | $\leqslant 1$ | 6.4.16 | 57 |
| $\Omega_7(q)$ | $\langle \mathrm{PO}, \bar\phi \rangle$ | | | both $H$ and $H^{\mathrm{e}}$ low | | 6.4.17 | 58 |
| $\mathrm{Sp}_6(q)$ | $\mathrm{P\Gamma O}^+_8(q)$ | $\mathrm{P\Gamma Sp}_6(q)$ | $N$ | | $\leqslant 1$ | 6.4.17 | 59 |
| $\mathrm{G}_2(q)$ | $\mathrm{P\Gamma O}_7(q)$ | $\langle M, \bar\phi \rangle$ | $N$ | | $\leqslant 1$ | 6.4.18 | 60 |
| $\mathrm{G}_2(q)$ | $\mathrm{P\Gamma Sp}_6(q)$ | $\langle M, \bar\phi \rangle$ | $N$ | $\leqslant 1$ | 1 | | 61 |
| $\mathrm{G}_2(q)$ | $\mathrm{P\Gamma O}_7(q)$ | $\langle M, \bar\phi \rangle$ | $N$ | | $\leqslant 1$ | 6.4.18 | 62 |
| $\mathrm{G}_2(q)$ | $\mathrm{P\Gamma Sp}_6(q)$ | $\langle M, \bar\phi \rangle$ | $N$ | $1, q = 4$ | 1 | | 63 |
| $\mathrm{G}_2(q)$ | $\mathrm{P\Gamma Sp}_6(q)$ | $\mathrm{Aut}\,M$ | $N$ | $\leqslant 1$ | 1 | | 64 |
| $\mathrm{Sz}(q)$ | $\mathrm{P\Gamma Sp}_4(q)$ | $\mathrm{Aut}\,\mathrm{Sz}(q)$ | $N$ | | 1 | 3.10 | 65 |

**Table 6.A continued**

| ♮ | $L$ | $X$ | $Y$ | $\circlearrowright$ | deg | [**LPS87**] | [**LPS90**] |
|---|---|---|---|---|---|---|---|
| 66 | $\Omega_{24}^{+}(2)$ | $\mathrm{Co}_1$ | $\mathrm{Sp}_{22}(2)$ | 67 | $2^{11}\star$ | VI.13 | 4.4.2 |
| 67 | $\Omega_{24}^{+}(2)$ | $\mathrm{Sp}_{22}(2)$ | $\mathrm{Co}_1$ | 66 | $2^{111}\star$ | VI.14 | 4.4.2 |
| 68 | $\Omega_{10}^{-}(2)$ | $\mathrm{A}_{12}$ | $P_1$ | np | 495 | VI.15 | 5.2.16 |
| 69 | $\mathrm{PSU}_9(2)$ | $\mathrm{J}_3$ | $P_1$ | np | 43605 | VI.16 | 5.2.12 |
| 70 | $\mathrm{P}\Omega_8^{+}(3)$ | $\Omega_8^{+}(2)$ | $P_4$ | np | 1120 | VI.17 | 5.1.15 |
| | Note that an automorphism of $L$ fixing $[X]_L$ sends $Y$ to $P_3$ and then to $P_1$. | | | | | | |
| 71 | $\mathrm{P}\Omega_8^{+}(3)$ | $\Omega_7(3)$ | $\Omega_8^{+}(2)$ | np | 28431 | VI.18 | 5.1.15 |
| | Note that $X$ is paired to three of the six $L$-conjugacy classes of $[Y]_{\mathrm{Aut}\,L}$. | | | | | | |
| 72 | $\mathrm{Sp}_8(2)$ | $\mathrm{PSL}_2(17)$ | $\mathrm{O}_8^{+}(2)$ | np | 136 | VI.19 | 5.1.9 |
| 73 | $\mathrm{Sp}_8(2)$ | $\mathrm{S}_{10}$ | $\mathrm{O}_8^{-}(2)$ | np | 120 | VI.20 | 5.1.9 |
| 74 | $\Omega_8^{+}(2)$ | $\mathrm{A}_9$ | $\mathrm{Sp}_6(2)$ | 75 | 120 | VI.21 | 5.1.15 |
| | $X$ is paired to only two of the three $L$-conjugacy classes of $[Y]_{\mathrm{Aut}\,L}$. | | | | | | |
| 75 | $\Omega_8^{+}(2)$ | $\mathrm{Sp}_6(2)$ | $\mathrm{A}_9$ | 74 | 960 | VI.22 | 5.1.15 |
| | $X$ is paired to only two of the three $L$-conjugacy classes of $[Y]_{\mathrm{Aut}\,L}$. | | | | | | |
| 76 | $\Omega_7(3)$ | $\mathrm{G}_2(3)$ | $\mathrm{Sp}_6(2)$ | | 3159 | VI.23 | 5.1.14a |
| | $X$ is paired to just one of the two $L$-conjugacy classes of $[Y]_{\mathrm{Aut}\,L}$. | | | | | | |
| 77 | $\Omega_7(3)$ | $N_1^{+}$ | $\mathrm{G}_2(3)$ | 60 | 1080 | VI.24 | 5.1.14 |
| 78 | $\mathrm{PSU}_6(2)$ | $\mathrm{M}_{22}$ | $\mathrm{PSU}_5(2)$ | 79 | 672 | VI.25 | 5.1.13 |
| 79 | $\mathrm{PSU}_6(2)$ | $\mathrm{PSU}_5(2)$ | $\mathrm{M}_{22}$ | 78 | 20736 | VI.26 | 5.1.13 |
| 80 | $\mathrm{PSU}_6(2)$ | $\mathrm{PSU}_5(2)$ | $\mathrm{PSU}_4(3).2$ | np | 1408 | VI.27 | 5.1.13 |
| 81 | $\mathrm{Sp}_6(2)$ | $\mathrm{O}_6^{+}(2)$ | $\mathrm{G}_2(2)$ | 61 | 120 | VI.28 | 5.1.8 |
| 82 | $\mathrm{PSU}_4(3)$ | $\mathrm{PSL}_3(4)$ | $P_1$ | np | 280 | VI.29 | 5.2.7 |
| | $H$ and $H^{\mathrm{e}}$ are low unless $H$ is even and $G^{\mathrm{e}} = L.(2^2)_{133}$. | | | | | | |
| 83 | $\mathrm{PSU}_3(3)$ | $\mathrm{PSL}_2(7)$ | $P_1$ | np | 28 | VI.30 | 5.1.10 |
| | Check as in 3.4.1 that $N$ is even. If $G$ were odd then $N < G^{\mathrm{e}} = L$, which is not. | | | | | | |
| 84 | $\mathrm{G}_2(4)$ | $\mathrm{PSU}_3(4).2$ | $\mathrm{J}_2$ | 85 | 416 | V.1 | Th. B |
| 85 | $\mathrm{G}_2(4)$ | $\mathrm{J}_2$ | $\mathrm{PSU}_3(4).2$ | 84 | 2016 | V.2 | Th. B |
| 86 | $\mathrm{F}_4(q)$ | ${}^3\mathrm{D}_4(q).3$ | $\mathrm{Sp}_8(q)$ | np | $q^8(q^8+q^4+1)$ | V.3 | Th. B |
| | $q = 2^f$. See also [**Kle88**] and [**LS87**]. | | | | | | |

| $M$ | $G$ | $N$ | $H$ | $r_o$ | $r_e$ | ! | ♯ |
|---|---|---|---|---|---|---|---|
| $\mathrm{Co}_1$ | $L$ | $M$ | $M$ | | 1 | | 66 |
| $\mathrm{Sp}_{22}(2)$ | $L$ | $M$ | $M$ | | 1 | | 67 |
| $\mathrm{A}_{12}$ | $\mathrm{O}_{10}^-(2)$ | $\mathrm{S}_{12}$ | $N$ | $\leqslant 1$ | $\leqslant 1$ | | 68 |
| $\mathrm{J}_3$ | $L.\mathrm{S}_3$ | $M.2$ | $N$ | both $H$ and $H^{\mathrm{e}}$ low | | | 69 |
| $\Omega_8^+(2)$ | $L.\mathrm{D}_8$ | $M.2$ | $N$ | both $H$ and $H^{\mathrm{e}}$ low | | | 70 |
| $\Omega_7(3)$ | $L.\mathrm{S}_3$ | $M.2$ | $N$ | both $H$ and $H^{\mathrm{e}}$ low | | | 71 |
| $\mathrm{PSL}_2(17)$ | $L$ | $M$ | $M$ | | 1 | | 72 |
| $\mathrm{A}_{10}$ | $L$ | $X$ | $X$ | | 1 | | 73 |
| $\mathrm{A}_9$ | $L.2$ | $M$ | $M$ | | $H = H^{\mathrm{e}}$ low | | 74 |
| $\mathrm{Sp}_6(2)$ | $L.2$ | $M$ | $M$ | | $H = H^{\mathrm{e}}$ low | | 75 |
| $\mathrm{G}_2(3)$ | $L$ | $M.2$ | $M$ | | 2 | | 76 |
| $\mathrm{PSL}_4(3)$ | $L$ | $M.2^2$ | $X = M.2_2$ | | 2 | | 77 |
| $\mathrm{M}_{22}$ | $L.\mathrm{S}_3$ | $M.2$ | $N$ | both $H$ and $H^{\mathrm{e}}$ low | | | 78 |
| $\mathrm{PSU}_5(2)$ | $L.2$ | $M.2$ | $N$ | 1 | 1 | | 79 |
| $\mathrm{PSU}_5(2)$ | $L$ | $M.2$ | $M$ | | 2 | | 80 |
| $\mathrm{A}_8$ | $L$ | $X$ | $N$ | | 1 | | 81 |
| $\mathrm{PSL}_3(4)$ | $L.\mathrm{D}_8$ | $M.\mathrm{D}_{12}$ | $M.2^2$ | | $\leqslant 3$ | | 82 |
| $\mathrm{PSL}_2(7)$ | $\mathrm{PGU}_3(3)$ | $\mathrm{PGL}_2(7)$ | $N$ | | 1 | | 83 |
| $\mathrm{PSU}_3(4)$ | $\mathrm{Aut}\,L$ | $\mathrm{Aut}\,M$ | $N$ | $\leqslant 1$ | 1 | | 84 |
| $\mathrm{J}_2$ | $\mathrm{Aut}\,L$ | $\mathrm{Aut}\,M$ | $N$ | $\leqslant 1$ | 1 | | 85 |
| ${}^3\mathrm{D}_4(q)$ | $L.f$ | $X.f$ | $N$ | $\leqslant 1$ | 1 | [**HLS87**] | 86 |

**end of Table 6.A**

TABLE 6.B. Certain primitive almost simple second maximal subgroups

| ♯ | $L$ | $X$ | $Y$ | deg | $M$ | $G$ | $N$ | $H$ | $r_o$ | $r_e$ | ! |
|---|---|---|---|---|---|---|---|---|---|---|---|
| 2 | $A_8$ | $A_8 \cap (S_2 \times S_6)$ | $AGL_3(2)$ | 15 | $A_6$ | $L$ | $X$ | $X$ | | 1 | 6.4.1 |
| 30 | $A_{10}$ | $A_6.2_3$ | $A_{10} \cap (S_2 \times S_8)$ | 45 | $A_6$ | $S_{10} < \boldsymbol{A}$ | $M.2^2$ | $N$ | | 1 | 6.4.8 |
| 1 | $A_8$ | $A_7$ | $AGL_3(2)$ | 15 | $A_7$ | $L$ | $M$ | $M$ | | 1 | 6.4.1 |
| 4 | $A_9$ | $A_9 \cap (S_2 \times S_7)$ | $P\Gamma L_2(8)$ | 120 | $A_7$ | $L$ | $X$ | $X$ | | 1 | 6.4.1 |
| 14 | $A_8$ | $A_7$ | $PSL_2(7)$ | 120 | $A_7$ | $S_8 < \boldsymbol{A}$ | $S_7$ | $N$ | | 1 | 6.4.5 |
| 53 | $Sp_6(2)$ | $O_6^+(2)$ | $O_6^-(2)$ | 28 | $A_8$ | $L$ | $X$ | $X$ | | 1 | |
| 81 | $Sp_6(2)$ | $O_6^+(2)$ | $G_2(2)$ | 120 | $A_8$ | $L$ | $X$ | $N$ | | 1 | |
| 73 | $Sp_8(2)$ | $S_{10}$ | $O_8^-(2)$ | 120 | $A_{10}$ | $L$ | $X$ | $X$ | | 1 | |
| 6 | $A_{12}$ | $A_{12} \cap (S_2 \times S_{10})$ | $M_{12}$ | 2520 | $A_{10}$ | $L$ | $X$ | $X$ | | 1 | 6.4.1 |
| 5 | $A_{12}$ | $A_{11}$ | $M_{12}$ | 2520 | $A_{11}$ | $L$ | $M$ | $M$ | | 1 | 6.4.1 |
| 68 | $\Omega_{10}^-(2)$ | $A_{12}$ | $P_1$ | 495 | $A_{12}$ | $O_{10}^-(2)$ | $S_{12}$ | $N$ | $\leqslant 1$ | $\leqslant 1$ | |
| 9 | $A_{24}$ | $A_{24} \cap (S_2 \times S_{22})$ | $M_{24}$ | $2^{11}\star$ | $A_{22}$ | $L$ | $X$ | $X$ | | 1 | 6.4.2 |
| 8 | $A_{24}$ | $A_{23}$ | $M_{24}$ | $2^{11}\star$ | $A_{23}$ | $L$ | $M$ | $M$ | | 1 | 6.4.2 |
| 11 | $A_{176}$ | $A_{175}$ | HS | $2^{163}\star$ | $A_{175}$ | $L$ | $M$ | $M$ | | 1 | 6.4.3 |
| 12 | $A_{276}$ | $A_{275}$ | $Co_3$ | $2^{262}\star$ | $A_{275}$ | $L$ | $M$ | $M$ | | 1 | 6.4.3 |
| 13 | $A_c$ | $A_{c-1}$ | $Sp_{2d}(2)$ | | $A_{c-1}$ | $L$ | $M$ | $M$ | | 1 | 6.4.4 |
| | $c = 2^{2d-1} \pm 2^{d-1}$, $d \geqslant 3$. | | | | | | | | | | |
| 15 | $A_{p+1}$ | $A_p$ | $PSL_2(p)$ | $(p-2)!$ | $A_p$ | $S_{p+1}$ | $S_p$ | $N$ | $\leqslant 1$ | 1 | 6.4.5 |
| | $p$ prime, $p > 7$; for $p = 7$ see ♯14. | | | | | | | | | | |
| 16 | $A_{2^d}$ | $A_{2^d-1}$ | $AGL_d(2)$ | | $A_{2^d-1}$ | $L$ | $M$ | $M$ | | 1 | 6.4.6 |
| | $d > 3$; for $d = 3$ see ♯1. | | | | | | | | | | |
| 17 | $A_{2l}$, $l \geqslant 3$ | $A_{2l-1}$ | $A_{2l} \cap \lfloor S_l \operatorname{wr} S_2 \rfloor$ | $\frac{1}{2}\binom{2l}{l}$ | $A_{2l-1}$ | $S_{2l} < \boldsymbol{A}$ | $S_{2l-1}$ | $N$ | | 1 | 6.4.7 |

**Table 6.B continued**

| ♯ | $L$ | $X$ | $Y$ | deg | $M$ | $G$ | $N$ | $H$ | $r_o$ | $r_e$ | ! |
|---|---|---|---|---|---|---|---|---|---|---|---|
| 83 | $\mathrm{PSU}_3(3)$ | $\mathrm{PSL}_2(7)$ | $P_1$ | 28 | $\mathrm{PSL}_2(7)$ | $\mathrm{PGU}_3(3)$ | $\mathrm{PGL}_2(7)$ | $N$ | | 1 | |
| | Check as in 3.4.1 that $N$ is even. If $G$ were odd then $N < G^{\mathrm{e}} = L$, which is not. | | | | | | | | | | |
| 36 | $\mathrm{M}_{12}$ | $\mathrm{PSL}_2(11)$ | $\mathrm{M}_{11}$ | 12 | $\mathrm{PSL}_2(11)$ | $L$ | $M.2$ | $M$ | | 2 | |
| 33 | $\mathrm{M}_{11}$ | $\mathrm{PSL}_2(11)$ | $\mathrm{M}_{10}$ | 11 | $\mathrm{PSL}_2(11)$ | $L$ | $M$ | $M$ | | 1 | |
| 34 | $\mathrm{M}_{11}$ | $\mathrm{PSL}_2(11)$ | $\mathrm{M}_9:2$ | 55 | $\mathrm{PSL}_2(11)$ | $L$ | $M.2$ | $M$ | | 2 | |
| 72 | $\mathrm{Sp}_8(2)$ | $\mathrm{PSL}_2(17)$ | $\mathrm{O}_8^+(2)$ | 136 | $\mathrm{PSL}_2(17)$ | $L$ | $M$ | $M$ | | 1 | |
| 41 | $\mathrm{M}_{24}$ | $\mathrm{PSL}_2(23)$ | $\mathrm{M}_{23}$ | 24 | $\mathrm{PSL}_2(23)$ | $L$ | $M.2$ | $M$ | | 2 | |
| 31 | $\mathrm{A}_{t+1}$ | $\mathrm{N}_L(\mathrm{PSL}_2\,t)$ | $\mathrm{A}_{t+1}\cap(\mathrm{S}_2\times\mathrm{S}_{t-1})$ | $\frac{t(t+1)}{2}$ | $\mathrm{PSL}_2(t)$ | $\mathrm{S}_{t+1}$ | $\mathrm{P\Gamma L}_2(t)$ | $N$ | | 1 | 6.4.8 |
| | $t>5$, $t\neq 9$. | | | | | | | | | | |
| 51 | $\mathrm{Sp}_4(q)$ | $\mathrm{O}_4^-(q)$ | $\mathrm{Sp}_2(q)\,\mathrm{wr}\,\mathrm{S}_2$ | $\frac{q^2(q^2+1)}{2}$ | $\mathrm{PSL}_2(q^2)$ | $\mathrm{P\Gamma Sp}_4(q)$ | $\mathrm{P\Gamma O}_4^-(q)$ | $N$ | $\leqslant 1$ | 1 | 6.4.14 |
| | $q$ even; $q\neq 2$. | | | | | | | | | | |
| 46 | $\mathrm{Sp}_{2m}(r)$ | $\mathrm{Sp}_2(r^m).m$ | $\mathrm{O}_{2m}^+(r)$ | $\frac{r^m(r^m+1)}{2}$ | $\mathrm{PSL}_2(r^m)$ | $\mathrm{P\Gamma Sp}_{2m}(r)$ | $\mathrm{P\Gamma L}_2(r^m)$ | $N$ | $\leqslant 1$ | 1 | 6.4.12 |
| | $r$ even; $m$ an odd prime. | | | | | | | | | | |
| 49 | $\mathrm{Sp}_4(q)$ | $\mathrm{O}_4^-(q)$ | $\mathrm{Sp}_2(q^2).2$ | $\frac{q^2(q^2-1)}{2}$ | $\mathrm{PSL}_2(q^2)$ | $\mathrm{P\Gamma Sp}_4(q)$ | $\mathrm{P\Gamma O}_4^-(q)$ | $N$ | $\leqslant 1$ | 1 | 6.4.13 |
| | $q$ even; $q\neq 2$. | | | | | | | | | | |
| 46 | $\mathrm{Sp}_{2m}(r)$ | $\mathrm{Sp}_2(r^m).m$ | $\mathrm{O}_{2m}^-(q)$ | $\frac{r^m(r^m-1)}{2}$ | $\mathrm{PSL}_2(r^m)$ | $\mathrm{P\Gamma Sp}_{2m}(r)$ | $\mathrm{P\Gamma L}_2(r^m)$ | $N$ | $\leqslant 1$ | 1 | 6.4.12 |
| | $r$ even; $m$ an odd prime. | | | | | | | | | | |
| 82 | $\mathrm{PSU}_4(3)$ | $\mathrm{PSL}_3(4)$ | $P_1$ | 280 | $\mathrm{PSL}_3(4)$ | $L.\mathrm{D}_8$ | $M.\mathrm{D}_{12}$ | $M.2^2$ | | $\leqslant 3$ | |
| | $H$ and $H^{\mathrm{e}}$ are low unless $H$ is even and $G^{\mathrm{e}} = L.(2^2)_{133}$. | | | | | | | | | | |
| 77 | $\Omega_7(3)$ | $N_1^+$ | $\mathrm{G}_2(3)$ | 1080 | $\mathrm{PSL}_4(3)$ | $L$ | $M.2^2$ | $X = M.2_2$ | | 2 | |
| 55 | $\mathrm{P\Omega}_{2m}^+(q)$ | $\mathrm{St}_L\{\mathsf{E},\mathsf{F}\}$ | $N_1$ | $\frac{q^{m-1}(q^m-1)}{(2,q-1)}$ | $\mathrm{PSL}_m(q)$ | $\mathrm{PO}_{2m}^+(q)$ | $\mathrm{Aut}\,M$ | $N$ | $\leqslant 2$ | $\leqslant 2$ | 6.4.15 |
| | $q=2,3; m\geqslant 5; \mathsf{V}=\mathsf{E}\oplus\mathsf{F}$ with $\mathsf{E},\mathsf{F}$ tot. singular. | | | | | | | | | | |
| 84 | $\mathrm{G}_2(4)$ | $\mathrm{PSU}_3(4).2$ | $\mathrm{J}_2$ | 416 | $\mathrm{PSU}_3(4)$ | $\mathrm{Aut}\,L$ | $\mathrm{Aut}\,M$ | $N$ | $\leqslant 1$ | 1 | |
| 80 | $\mathrm{PSU}_6(2)$ | $\mathrm{PSU}_5(2)$ | $\mathrm{PSU}_4(3).2$ | 1408 | $\mathrm{PSU}_5(2)$ | $L$ | $M.2$ | $M$ | | 2 | |
| 79 | $\mathrm{PSU}_6(2)$ | $\mathrm{PSU}_5(2)$ | $\mathrm{M}_{22}$ | 20736 | $\mathrm{PSU}_5(2)$ | $L.2$ | $M.2$ | $N$ | 1 | 1 | |

**Table 6.B continued**

| ♯ | $L$ | $X$ | $Y$ | deg | $M$ | $G$ | $N$ | $H$ | $r_o$ | $r_e$ | ! |
|---|---|---|---|---|---|---|---|---|---|---|---|
| 54 | $\mathrm{Sp}_6(2)$ | $\mathrm{O}_6^-(2)$ | $\mathrm{O}_6^+(2)$ | 36 | $\mathrm{PSp}_4(3)$ | $L$ | $X$ | $X$ | | 1 | |
| 59 | $\Omega_8^+(q)$ | $N_1$ | $X^\tau$ | $q^3(q^4-1)$ | $\mathrm{Sp}_6(q)$ | $\mathrm{P\Gamma O}_8^+(q)$ | $\mathrm{P\Gamma Sp}_6(q)$ | $N$ | | $\leqslant 1$ | 6.4.17 |
| | $q$ even; $\tau$ is a triality automorphism. Note that $\mathrm{Sp}_6(2)$ has degree 120. | | | | | | | | | | |
| 67 | $\Omega_{24}^+(2)$ | $\mathrm{Sp}_{22}(2)$ | $\mathrm{Co}_1$ | $2^{111}\star$ | $\mathrm{Sp}_{22}(2)$ | $L$ | $M$ | $M$ | | 1 | |
| 57 | $\Omega_{2m}^+(q)$ | $N_1$ | $P_m, P_{m-1}$ | $\prod_{i=1}^{m-1}(q^i+1)$ | $\mathrm{Sp}_{2m-2}(q)$ | $\langle L, \bar\phi\rangle$ | $\mathrm{Aut}\,M$ | $N$ | $\leqslant 1$ | $\leqslant 1$ | 6.4.16 |
| | $q$ even; $m \geqslant 4$. Note that if $m=4$, then $Y$ is conjugate in $\mathrm{Aut}\,L$ to $P_1$. Note that these are groups of odd degree. | | | | | | | | | | |
| 44 | $\mathrm{PSL}_{2m}(q)$ | $\mathrm{N}_L\,\mathrm{PSp}_{2m}(q)$ | $P_1$ | $\frac{q^{2m}-1}{q-1}$ | $\mathrm{PSp}_{2m}(q)$ | $\mathrm{P\Gamma L}_{2m}(q)$ | $\mathrm{P\Gamma Sp}_{2m}(q)$ | $N$ | $\leqslant 1$ | $\leqslant 2$ | 6.4.11 |
| | $m \geqslant 2$ but $q \neq 2$ when $m=2$. | | | | | | | | | | |
| 45 | $\mathrm{Sp}_{2m}(q)$ | $\mathrm{Sp}_{2a}(q^b).b$ | $\mathrm{O}_{2m}^\pm(q)$ | $\frac{q^m(q^m\pm1)}{2}$ | $\mathrm{PSp}_{2a}(q^b)$ | $\mathrm{P\Gamma Sp}_{2m}(q)$ | $\mathrm{P\Gamma Sp}_{2a}(q^b)$ | $N$ | $\leqslant 1$ | 1 | 6.4.12 |
| | $q$ even; $m = ab > b$, $b$ prime. | | | | | | | | | | |
| 56 | $\mathrm{P\Omega}_{2m}^+(q)$ | $N_1$ | $P_m, P_{m-1}$ | $\prod_{i=1}^{m-1}(q^i+1)$ | $\Omega_{2m-1}(q)$ | $\langle \mathrm{PSO}, \bar\delta, \bar\phi\rangle$ | $\mathrm{Aut}\,M$ | $\langle X, \bar\phi\rangle$ | | $\leqslant 1$ | 6.4.16 |
| | $q$ odd; $m \geqslant 4$. Note that if $m=4$, then $Y$ is conjugate in $\mathrm{Aut}\,L$ to $P_1$. | | | | | | | | | | |
| 52 | $\mathrm{Sp}_{2m}(2)$ | $\mathrm{O}_{2m}^\pm(2)$ | $\mathrm{O}_{2m}^\mp(2)$ | $2^{m-1}(2^m \mp 1)$ | $\Omega_{2m}^\pm(2)$ | $L$ | $X$ | $X$ | | 1 | |
| | $m > 3$. Note that $X \cap Y$ is a subgroup of type $N_1$ of $X$. | | | | | | | | | | |
| 48 | $\mathrm{Sp}_{4a}(q)$ | $\mathrm{O}_{4a}^-(q)$ | $\mathrm{Sp}_{2a}(q^2).2$ | | $\Omega_{4a}^-(q)$ | $\mathrm{P\Gamma Sp}_{4a}(q)$ | $\mathrm{P\Gamma O}_{4a}^-(q)$ | $N$ | $\leqslant 1$ | 1 | 6.4.13 |
| | $q$ even; $a > 1$. Note that $X \cap Y$ lie in the class $\mathcal{C}_3$ of subgroups of $X$. | | | | | | | | | | |
| 50 | $\mathrm{Sp}_{2m}(q)$ | $\mathrm{O}_{2m}^-(q)$ | $\mathrm{Sp}_m(q)\,\mathrm{wr}\,\mathrm{S}_2$ | | $\Omega_{2m}^-(q)$ | $\mathrm{P\Gamma Sp}_{2m}(q)$ | $\mathrm{P\Gamma O}_{2m}^-(q)$ | $N$ | $\leqslant 1$ | 1 | 6.4.14 |
| | $q$ even, $m$ even, $m \neq 2$. Note that $X \cap Y$ is a subgroup of type $N_m$ of $X$. | | | | | | | | | | |
| 47 | $\mathrm{Sp}_{2m}(q)$ | $\mathrm{O}_{2m}^\pm(q)$ | $\mathrm{Sp}_{2a}(q^b).b$ | | $\Omega_{2m}^\pm(q)$ | $\mathrm{P\Gamma Sp}_{2m}(q)$ | $\mathrm{P\Gamma O}_{2m}^\pm(q)$ | $N$ | $\leqslant 1$ | 1 | 6.4.13 |
| | $q$ even; $m = ab > b > 2$, $b$ prime. Note that $X \cap Y$ lie in the class $\mathcal{C}_3$ of subgroups of $X$. | | | | | | | | | | |

**Table 6.B continued**

| ♯ | $L$ | $X$ | $Y$ | deg | $M$ | $G$ | $N$ | $H$ | $r_o$ | $r_e$ | ! |
|---|---|---|---|---|---|---|---|---|---|---|---|
| 32 | $\mathrm{A}_{q^2+1}$ | $\mathrm{Aut\,Sz}(q)$ | $\mathrm{A}_{q^2+1} \cap (\mathrm{S}_2 \times \mathrm{S}_{q^2-1})$ | $\frac{q^2(q^2+1)}{2}$ | $\mathrm{Sz}(q)$ | $\mathrm{S}_{q^2+1}$ | $X$ | $X$ | | 1 | 3.10 |
| 65 | $\mathrm{Sp}_4(q)$ | $\mathrm{Sz}(q)$ | $\mathrm{O}_4^+(q)$ | $\frac{q^2(q^2+1)}{2}$ | $\mathrm{Sz}(q)$ | $\mathrm{P\Gamma Sp}_4(q)$ | $\mathrm{Aut\,Sz}(q)$ | $N$ | | 1 | 3.10 |
| | $q = 2^{2a+1} \geqslant 8$. | | | | | | | | | | |
| 61 | $\mathrm{Sp}_6(q)$, $q$ even | $\mathrm{G}_2(q)$ | $\mathrm{O}_6^{\pm}(q)$ | $\frac{q^3(q^3\pm1)}{2}$ | $\mathrm{G}_2(q)$ | $\mathrm{P\Gamma Sp}_6(q)$ | $\langle M, \bar{\phi} \rangle$ | $N$ | $\leqslant 1$ | 1 | |
| 63 | $\mathrm{Sp}_6(q)$, $q$ even | $\mathrm{G}_2(q)$ | $P_1$ | $\frac{q^6-1}{q-1}$ | $\mathrm{G}_2(q)$ | $\mathrm{P\Gamma Sp}_6(q)$ | $\langle M, \bar{\phi} \rangle$ | $N$ | $1, q=4$ | 1 | |
| 64 | $\mathrm{Sp}_6(q)$ | $\mathrm{G}_2(q)$ | $N_2$ | $\frac{q^4(q^6-1)}{q^2-1}$ | $\mathrm{G}_2(q)$ | $\mathrm{P\Gamma Sp}_6(q)$ | $\mathrm{Aut}\,M$ | $N$ | $\leqslant 1$ | 1 | |
| | $q$ even, $q \neq 2$. | | | | | | | | | | |
| 76 | $\Omega_7(3)$ | $\mathrm{G}_2(3)$ | $\mathrm{Sp}_6(2)$ | 3159 | $\mathrm{G}_2(3)$ | $L$ | $M.2$ | $M$ | | 2 | |
| | $X$ is paired to just one of the two $L$-conjugacy classes of $[Y]_{\mathrm{Aut}\,L}$. | | | | | | | | | | |
| 60 | $\Omega_7(q)$, $q$ odd | $\mathrm{G}_2(q)$ | $N_1^+, N_1^-$ | $\frac{q^3(q^3\pm1)}{2}$ | $\mathrm{G}_2(q)$ | $\mathrm{P\Gamma O}_7(q)$ | $\langle M, \bar{\phi} \rangle$ | $N$ | | $\leqslant 1$ | 6.4.18 |
| 62 | $\Omega_7(q)$, $q$ odd | $\mathrm{G}_2(q)$ | $P_1$ | $\frac{q^6-1}{q-1}$ | $\mathrm{G}_2(q)$ | $\mathrm{P\Gamma O}_7(q)$ | $\langle M, \bar{\phi} \rangle$ | $N$ | | $\leqslant 1$ | 6.4.18 |
| 86 | $\mathrm{F}_4(q)$ | ${}^3\mathrm{D}_4(q).3$ | $\mathrm{Sp}_8(q)$ | $q^8(q^8+q^4+1)$ | ${}^3\mathrm{D}_4(q)$ | $L.f$ | $X.f$ | $N$ | $\leqslant 1$ | 1 | [**HLS87**] |
| | $q = 2^f$. See also [**Kle88**] and [**LS87**]. | | | | | | | | | | |
| 66 | $\Omega_{24}^+(2)$ | $\mathrm{Co}_1$ | $\mathrm{Sp}_{22}(2)$ | $2^{11}\star$ | $\mathrm{Co}_1$ | $L$ | $M$ | $M$ | | 1 | |
| 85 | $\mathrm{G}_2(4)$ | $\mathrm{J}_2$ | $\mathrm{PSU}_3(4).2$ | 2016 | $\mathrm{J}_2$ | $\mathrm{Aut}\,L$ | $\mathrm{Aut}\,M$ | $N$ | $\leqslant 1$ | 1 | |
| 35 | $\mathrm{M}_{12}$ | $\mathrm{M}_{11}$ | $\mathrm{M}_{11}$ | 12 | $\mathrm{M}_{11}$ | $L$ | $M$ | $M$ | | 1 | 6.4.9 |
| 18 | $\mathrm{A}_{11}$ | $\mathrm{M}_{11}$ | $\mathrm{A}_{11} \cap (\mathrm{S}_2 \times \mathrm{S}_9)$ | 55 | $\mathrm{M}_{11}$ | $\mathrm{S}_{11}$ | $M$ | $M$ | | 1 | |
| 37 | $\mathrm{M}_{12}$ | $\mathrm{M}_{11}$ | $\mathrm{M}_{10} : 2$ | 66 | $\mathrm{M}_{11}$ | $L$ | $M$ | $M$ | | 1 | 6.4.9 |
| 21 | $\mathrm{A}_{12}$ | $\mathrm{M}_{12}$ | $\mathrm{A}_{12} \cap (\mathrm{S}_3 \times \mathrm{S}_9)$ | 220 | $\mathrm{M}_{12}$ | $\mathrm{S}_{12}$ | $M$ | $M$ | | 1 | |
| 43 | HS | $\mathrm{M}_{22}$ | $\mathrm{PSU}_3(5) : 2$ | 176 | $\mathrm{M}_{22}$ | $L$ | $M$ | $M$ | | 1 | |
| 23 | $\mathrm{A}_{22}$ | $\mathrm{M}_{22}$ | $\mathrm{A}_{22} \cap (\mathrm{S}_2 \times \mathrm{S}_{20})$ | 231 | $\mathrm{M}_{22}$ | $\mathrm{S}_{22} < \boldsymbol{A}$ | $M.2$ | $N$ | | 1 | |
| 24 | $\mathrm{A}_{23}$ | $\mathrm{M}_{23}$ | $\mathrm{A}_{23} \cap (\mathrm{S}_2 \times \mathrm{S}_{21})$ | 253 | $\mathrm{M}_{23}$ | $\mathrm{S}_{23}$ | $M$ | $M$ | | 1 | |
| 38 | $\mathrm{M}_{24}$ | $\mathrm{M}_{23}$ | $\mathrm{M}_{12} : 2$ | 1288 | $\mathrm{M}_{23}$ | $L$ | $M$ | $M$ | | 1 | |
| 40 | $\mathrm{M}_{24}$ | $\mathrm{M}_{23}$ | $\mathrm{PSL}_2(23)$ | 40320 | $\mathrm{M}_{23}$ | $L$ | $M$ | $M$ | | 1 | |
| 27 | $\mathrm{A}_{24}$ | $\mathrm{M}_{24}$ | $\mathrm{A}_{24} \cap (\mathrm{S}_3 \times \mathrm{S}_{21})$ | 2024 | $\mathrm{M}_{24}$ | $\mathrm{S}_{24}$ | $M$ | $M$ | | 1 | |

**end of Table 6.B**

## 6.4. Explanation of the tables

Throughout this section we refer to the group factorizations listed in our tables and consequently we use the same notations. This comprises names $L$, $X$, $Y$, $\boldsymbol{S} = \mathrm{Sym}\binom{L}{Y^\diamond}$, $\boldsymbol{A} = \mathrm{Alt}\binom{L}{Y^\diamond}$, $M = \mathrm{Soc}\, X$, $G = \mathrm{N}_{\boldsymbol{S}}(L)$, $N = \mathrm{N}_{\boldsymbol{S}}(M)$, $H = \mathrm{N}_G(M)$ and finally

$$r_o = \begin{cases} \mathrm{hm}(G, H) & \text{if } H \text{ is odd and second maximal in } \boldsymbol{S}, \\ 0 & \text{otherwise} \end{cases}$$

and

$$r_e = \begin{cases} \mathrm{hm}(G^{\mathrm{e}}, H^{\mathrm{e}}) & \text{if } H^{\mathrm{e}} \text{ is second maximal in } \boldsymbol{A}, \\ 0 & \text{otherwise.} \end{cases}$$

Then we write $\boldsymbol{H}$ for an almost simple second maximal subgroup of $\boldsymbol{S}$ or $\boldsymbol{A}$ which matches with the skeleton $(L, X, Y)$. In particular $\mathrm{Soc}\, \boldsymbol{H} = M$ and either $\boldsymbol{H} = H$ or $\boldsymbol{H} = H^{\mathrm{e}}$. We set as we did in 2.1, $\boldsymbol{U} = \boldsymbol{H}\boldsymbol{A}$.

Explanations are given following the order given in Table 6.A, so, in reading about $\sharp n$, we assume that the reader is aware of the explanations given for $\sharp m$ for any $m < n$. In particular, a missing explanation implies that normalizers, $r_o$ and $r_e$ were determined in an entirely similar way to some previous case.

We illustrate in great detail how to determine $r_o$ and $r_e$ up to $\sharp 17$ which is explained in 6.4.7. After this, one should be familiar with the techniques used and so we will only add comments to explain the less obvious cases.

Occasionally, we talk of the $p$-part of a positive integer, meaning the highest power of $p$ dividing it. For example the $p$-part of $n$ is $p^e$ where

$$e = \lfloor \log_p n \rfloor$$

and $\lfloor x \rfloor$ is the unique integer $z$ such that $z \leqslant x < z+1$. We also denote the exponent of the $p$-part of $n!$ by $\lfloor \log_p n! \rfloor$. Therefore,

$$\lfloor \log_p \binom{n}{k} \rfloor = \lfloor \log_p n! \rfloor - \lfloor \log_p k! \rfloor - \lfloor \log_p (n-k)! \rfloor$$

and if $n = a_0 + a_1 p + a_2 p^2 + \ldots + a_r p^r$ with $0 \leqslant a_i < p$ for all $i$, then

$$\lfloor \log_p n! \rfloor = \sum_{i=1}^{r} a_i \frac{p^i - 1}{p - 1} < \sum_{i=1}^{r} a_i p^i \leqslant n. \tag{6.4.A}$$

In particular, $\lfloor \log_2 n! \rfloor = n - c$ where $c$ is the number of 1 appearing in the binary expansion of $n$.

**6.4.1. Duals of the homogeneous factorizations.** These are the group factorizations $(U, P, Q)$ where $U$ is a symmetric or alternating group in its action of smallest degree and $P$ is the stabilizer in $U$ of a set of points.

We claim that if the degree is not 6, then

$$\mathrm{hm}_S(U, P) \leqslant \left|\mathrm{N}_S(P) : P\right|,$$

where $S = \mathrm{Sym}\binom{U}{Q\diamond}$ and $U$, $P$ are identified with their images under $/\binom{U}{Q\diamond} : \star/$. By 2.8.1 it is enough to show that

$$\mathrm{hm}_S(P, U) = \left|U : P\right|.$$

But this is well known (see 2.7.1).

**6.4.2. The coset space of right cosets of $\mathrm{M}_{24}$ in $\mathrm{A}_{24}$.** Here $\boldsymbol{S}$ is the symmetric group on the set of the right cosets $\binom{\mathrm{A}_{24}}{\mathrm{M}_{24}\diamond}$. A proper subgroup of $\mathrm{M}_{24}$ has index at least 24 so that it is not possible to embed $\mathrm{M}_{24}$ in a symmetric group of degree smaller than 24. Consequently, all the subgroups of $\mathrm{S}_{24}$ isomorphic to $\mathrm{M}_{24}$ are transitive. Since $\mathrm{M}_{24}$ has only one conjugation class of subgroups of index 24, all the subgroups of $\mathrm{S}_{24}$ isomorphic to $\mathrm{M}_{24}$ are conjugate in $\mathrm{S}_{24}$ and primitive. Of course, being simple, they all lie in $\mathrm{A}_{24}$. This proves the uniqueness up to permutational isomorphism of $/\binom{\mathrm{A}_{24}}{\mathrm{M}_{24}\diamond} : \mathrm{A}_{24}/$. Note that $\mathrm{M}_{24}$ has trivial outer automorphism group. In particular, it is a self-normalizing subgroup of $\mathrm{S}_{24}$. If there were a subgroup $T$ of $\boldsymbol{S}$ isomorphic to $\mathrm{S}_{24}$ and containing $/\binom{\mathrm{A}_{24}}{\mathrm{M}_{24}\diamond} : \mathrm{A}_{24}/$, then the stabilizer of a point in $T$ would have a $\mathrm{M}_{24}$ as normal subgroup of index 2 which we have seen it is impossible. Therefore $/\binom{\mathrm{A}_{24}}{\mathrm{M}_{24}\diamond} : \mathrm{A}_{24}/$ is self-normalizing in $\boldsymbol{S}$.

Similarly, $/\binom{\mathrm{A}_{23}}{\mathrm{M}_{23}\diamond} : \mathrm{A}_{23}/$ is self-normalizing in $\mathrm{Sym}\binom{\mathrm{A}_{23}}{\mathrm{M}_{23}\diamond}$. Moreover, because $/\binom{\mathrm{A}_{23}}{\mathrm{M}_{23}\diamond} : \mathrm{A}_{23}/$ and $/\binom{\mathrm{A}_{24}}{\mathrm{M}_{24}\diamond} : \mathrm{A}_{23}/$ are permutationally isomorphic, the latter must be self-normalizing in $\boldsymbol{S}$. This is enough to show that in both ♯8 and ♯9, $G$ is even and so $r_o = 0$. By 6.4.1, $r_e = 1$.

**6.4.3.** Consider ♯11 where

$$(L, X, Y) = (\mathrm{A}_{176}, \mathrm{A}_{175}, \mathrm{HS}).$$

$G = L$ because $\mathrm{Aut}\,\mathrm{HS}$ has no subgroup of index 176. Note that $\mathrm{A}_{175} \cap \mathrm{HS}$ is isomorphic to $\mathrm{PSU}_3(5).2$ which is primitive of type almost simple of degree 175. This is a self-normalizing subgroup of $\mathrm{S}_{175}$. It follows that $\mathrm{A}_{175}$ is self-normalizing in $\mathrm{Sym}\binom{\mathrm{A}_{175}}{(\mathrm{A}_{175}\cap\mathrm{HS})\diamond}$ (see the argument with $\mathrm{M}_{24}$ and $\mathrm{A}_{24}$ in 6.4.2 for example) and hence it is self-normalizing in $\mathrm{Sym}\binom{\mathrm{A}_{176}}{\mathrm{HS}\diamond}$ as well. This proves that $N = M$.

Now consider ♯12 where

$$(L, X, Y) = (\mathrm{A}_{276}, \mathrm{A}_{275}, \mathrm{Co}_3).$$

$G = L$ because $\operatorname{Aut} L \cong L$. Note that $A_{275} \cap Co_3$ is isomorphic to $McL.2$ which is primitive of type almost simple of degree 275. This is a self-normalizing subgroup of $S_{275}$. One shows as above that $N = M$.

**6.4.4. Two-transitive symplectic groups as stabilizers.** We are concerned with ♯13 where the skeletons are as follows

$$(L, X, Y) = \big(A_c, A_{c-1}, Sp_{2d}(2)\big), \qquad c = 2^{d-1}(2^d \pm 1), \qquad d \geqslant 3.$$

For $d \geqslant 3$ put $\mathsf{Q}^+ = \binom{Sp_{2d}(2)}{O^+_{2d}(2)\diamond}$ and $\mathsf{Q}^- = \binom{Sp_{2d}(2)}{O^-_{2d}(2)\diamond}$. Then $|\mathsf{Q}^+| = 2^{d-1}(2^d+1)$ and $|\mathsf{Q}^-| = 2^{d-1}(2^d-1)$. By $Sp_{2d}(2)$ we mean the subgroups $/\binom{Sp_{2d}(2)}{O^+_{2d}(2)\diamond} : Sp_{2d}(2)/$ of $\operatorname{Sym}\mathsf{Q}^+$ and $/\binom{Sp_{2d}(2)}{O^-_{2d}(2)\diamond} : Sp_{2d}(2)/$ of $\operatorname{Sym}\mathsf{Q}^-$ respectively. These are the well known 2-transitive symplectic groups (see a description in [**DM96**, 7.7] for example). Since

$$PSp_{2d}(2) = P\Gamma Sp_{2d}(2) \cong \Gamma Sp_{2d}(2) = Sp_{2d}(2) \cong \operatorname{Aut} Sp_{2d}(2),$$

$Sp_{2d}(2)$ is self-normalizing in $\operatorname{Sym}\mathsf{Q}^+$ and $\operatorname{Sym}\mathsf{Q}^-$. In particular (see the argument with $M_{24}$ and $A_{24}$ in 6.4.2 for example), $L$ is self-normalizing in $\boldsymbol{S}$ and so $r_o = 0$. By definition

$$A_{c-1} \cap Sp_{2d}(2) = O^{\pm}_{2d}(2).$$

This is a transitive subgroup of degree $c-1$ (because $Sp_{2d}(2)$ is 2-transitive) with trivial outer automorphism group. Therefore it is a self-normalizing subgroup of $S_{c-1}$. It follows that $A_{c-1}$ is self-normalizing in $\boldsymbol{S}$; then $N = M$ and $r_e = 1$.

Note that the smallest $M$ are $A_{27}$, $A_{35}$, $A_{119}$, $A_{135}$, $A_{495}$ and $A_{527}$ so that none of the alternating groups appearing in column $M$ of Table 6.B above ♯13 may be abstractly isomorphic to an $M$ of ♯13.

**6.4.5.** Here we refer to ♯14, ♯15 where $(L, X, Y) = \big(A_{p+1}, A_p, PSL_2(p)\big)$, $p > 5$ and $p$ is prime. Clearly, $G = S_{p+1}$ and $H = N = S_p$. Both $r_o$ and $r_e$ are at most 1 because of 6.4.1. Note however, that in ♯14, $G$ is even. In fact, $S_7$ has only one conjugacy class of maximal subgroups of index 120 and hence $H$ corresponds to the $X$ of the skeleton in ♯4. In particular, $H$ is even because it is contained in the simple group $A_9$. If $G$ were odd, then $G^{\mathrm{e}} = A_8$ and so

$$S_7 = H < G^{\mathrm{e}} = A_8$$

which is impossible.

Finally, we prove that if $\boldsymbol{H}$ matches with ♯15, then it may not match with any of the skeletons which are shown in Table 6.B above ♯15.

First of all, it may not match ♯5 with $p = 11$ or ♯8 with $p = 23$ because comparing the degrees

$$9! \neq 2520 \qquad \text{and} \qquad 21! \neq |A_{24} : M_{24}|.$$

More effort is required to show that $\boldsymbol{H}$ may not match ♯13. If there were a match, then

$$p = c - 1 \qquad \text{and } (p-2)! = \frac{c!}{2\left|\mathrm{Sp}_{2d}(2)\right|} \qquad \text{with } c = 2^{d-1}(2^d \pm 1).$$

This would imply $2\left|\mathrm{Sp}_{2d}(2)\right| = c(c-1)(c-2)$. One easily extracts the 2-part of both members obtaining

$$2^{d^2+1} = 2^d$$

which has no solutions.

**6.4.6.** We are now dealing with ♯1 and ♯16 when

$$(L, X, Y) = \big(\mathrm{A}_{2^d}, \mathrm{A}_{2^d-1}, \mathrm{AGL}_d(2)\big), \qquad d \geqslant 3.$$

Because $\mathrm{AGL}_d(2)$ is self-normalizing in $\boldsymbol{S}$, $L$ is self-normalizing too and so $G = L$. Observe that $\mathrm{A}_{2^d-1} \cap \mathrm{AGL}_d(2)$ is a primitive almost simple group of degree $2^d - 1$ which is isomorphic to its automorphism group. Therefore it is a self-normalizing subgroup of $\mathrm{S}_{2^d-1}$ and hence $\mathrm{A}_{2^d-1}$ is self-normalizing in $\boldsymbol{S}$. This proves that $N = M$ which yields $r_e = 1$.

These factorizations arise for $d \geqslant 3$. However, the skeleton for $d = 3$ appears already in ♯1. There are not other permutational isomorphisms between $M$ of ♯16 and the ones related to skeletons appearing above in Table 6.B. This is obvious up to ♯12 (following the order of Table 6.B) and also for the family in ♯13 because $c = 2^{e-1}(2^e \pm 1)$ is never a power of 2. For the family in ♯15 note that if $p = 2^d - 1$, then the degree of those $M$ is $(2^d-3)!$ while the degree of these $M$ is $(2^d)!/2\left|\mathrm{AGL}_d(2)\right|$. If they were equal then

$$2^d(2^d-1)(2^d-2) = 2^{d+\binom{d}{2}+1}\prod_{i=1}^{d}(2^i-1)$$

which is readily impossible (compare the 2-part of both members for example).

**6.4.7. The set of the equipartitions in two subsets.** Consider ♯17 when

$$(L, X) = (\mathrm{A}_{2l}, \mathrm{A}_{2l-1}),$$

$Y$ is the stabilizer of an equipartition in two subsets and $l \geqslant 3$.

Recall 3.6.1; $G$, is even because $\binom{2x}{x}$ is even for all $x$. Then $r_o = 0$ and arguing as in 6.4.1 one sees that $r_e = 1$.

We now prove that if $\boldsymbol{H}$ matches with ♯17, then it may not match with any of the skeletons which are shown in Table 6.B above ♯17.

This is easy for the skeletons in ♯1, ♯4, ♯14 ($l = 4$), ♯74 ($l = 5$) and ♯5 ($l = 6$): one just compares the degrees.

For skeletons in ♯8, ♯11, ♯12, ♯13, ♯15 and ♯16 proceed as follows. Assume that there were a match and call $Z$ what is there called $Y$. $Z$ is 2-transitive of degree $2l$ and so $\mathrm{A}_{2l-1} \cap Z$ is transitive of degree $2l-1$. However, $\mathrm{A}_{2l-1} \cap Y$ is not transitive of degree $2l-1$. Thus

$$/\binom{\mathrm{A}_{2l-1}}{(\mathrm{A}_{2l-1} \cap Y)_\diamond} : \mathrm{A}_{2l-1}/$$

and

$$/\binom{\mathrm{A}_{2l-1}}{(\mathrm{A}_{2l-1} \cap Z)_\diamond} : \mathrm{A}_{2l-1}/$$

may not be permutationally isomorphic.

**6.4.8. The set of subsets of order two of a projective line.** Here we refer to ♯30, ♯31 where

$$(L, X, Y) = (\mathrm{A}_{q+1}, \mathrm{N}_L(\mathrm{PSL}_2\, q), \mathrm{A}_{q+1} \cap (\mathrm{S}_2 \times \mathrm{S}_{q-1})), \qquad q > 5.$$

In §3.9 we collected conclusive results about the parity of the projective groups acting on the related projective space. In particular, Table 3.C shows that if $q$ is a power of the prime $p$, then

$$X = \begin{cases} \mathrm{P\Gamma L}_2(q) & \text{if } q \text{ is even,} \\ \mathrm{P\Sigma L}_2(q) & \text{if } p \equiv 1 \mod 4 \text{ or } \log_p q \text{ is odd,} \\ \langle \mathrm{PSL}, \mathrm{PG}(\boldsymbol{\pi}^2), \mathrm{G}_1(\boldsymbol{u})\,\mathrm{PG}(\boldsymbol{\pi}) \rangle & \text{otherwise.} \end{cases}$$

The subgroup $G$ is isomorphic to $\mathrm{S}_{q+1}$. By 3.6.1, $G$ is odd in $\boldsymbol{S}$ if and only if $q$ is even.

Assume $q$ even first. Then $H$ may not be bigger than $X$ and so it is even. Consequently, $r_o = 0$ in this case. To determine $r_e$ note that by 3.9.2, $\mathrm{hm}_{\boldsymbol{S}}(H, L) = 2\left|L : H\right|$ but because $\left|G : L\right| = 2$, $r_e = 1$.

Now, assume $q$ odd so $G$ is even. Again, $H$ is isomorphic to $\mathrm{P\Gamma L}$ because this is the normalizer of $X$ in $G$. By 3.9.2, $\mathrm{hm}_{\boldsymbol{S}}(H, G) = \left|G : H\right|$ and thus $r_e = 1$ even if $q$ is odd.

**6.4.9.** Consider ♯35 first. $L = \mathrm{M}_{12}$ has two conjugacy classes of maximal subgroups isomorphic to $\mathrm{M}_{11}$, fused in $\mathrm{Aut}\, L$. Obviously, factorizations $L = XY$ with $X \cong \mathrm{M}_{11} \cong Y$ arise only for $X$ and $Y$ lying in different classes. Suppose now that for some $\boldsymbol{s} \in \boldsymbol{S}$, $X^{\boldsymbol{s}} < L$. Of course, $X^{\boldsymbol{s}}$ is transitive as a subgroup of $\boldsymbol{S}$. By the Frattini Argument, $L = X^{\boldsymbol{s}} Y$, therefore $X^{\boldsymbol{s}}$ is conjugate to $X$ in $L$. This proves that $\mathrm{hm}_{\boldsymbol{S}}(X, L) = \left|L : X\right|$ and thus that

$$r_e = \mathrm{hm}_{\boldsymbol{S}}(X, L) = \frac{\left|\mathrm{N}_{\boldsymbol{S}}\, X : X\right|}{\left|\mathrm{N}_{\boldsymbol{S}}\, L : L\right|} = 1.$$

Now consider ♯37. Looking at the character table of $L$ in the ATLAS, one sees that there are also factorizations $L = XZ$ where $X \cong \mathrm{M}_{11}$ as above but $Z \cong \mathrm{M}_{10} : 2$. However, for each of these $Z$, only the $\mathrm{M}_{11}$ conjugate to $X$ in $L$ give a factorization. In particular, $r_e = 1$ again.

**6.4.10. Degree 276.** We are particularly interested to the group factorization in ♯42:

$$(L, X, Y) = (\mathrm{M}_{24}, \mathrm{PSL}_2(23), \mathrm{M}_{22}\,.2).$$

However, this is strictly related to the group factorization obtained by putting $q = 23$ in ♯31 (see 6.4.8). Note how the different groups relate to each other.

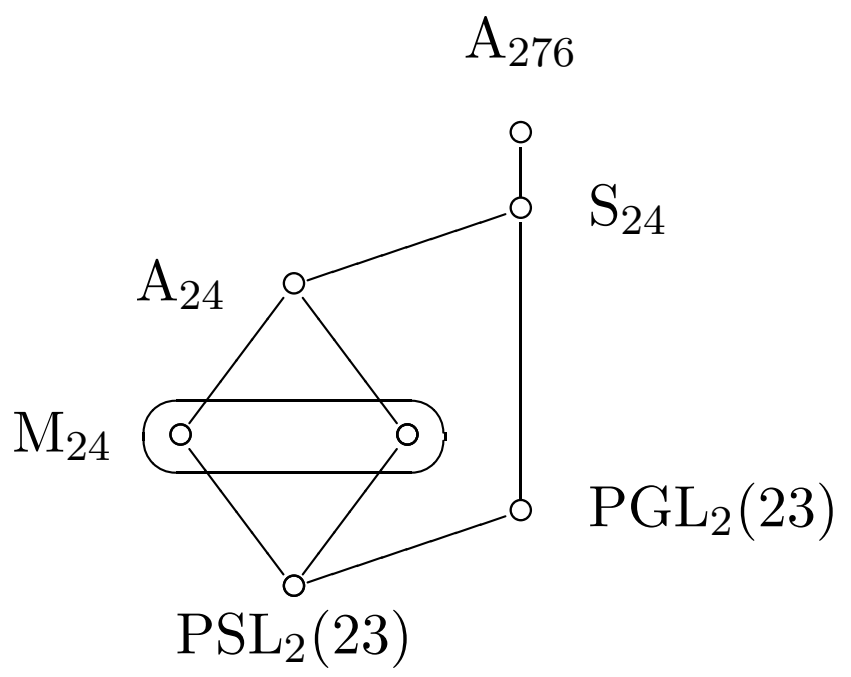

**6.4.11.** Consider ♯44. Here $L$ is a projective special linear group $\mathrm{PSL}(2m, q)$ where $q = p^f$, $p$ is prime. $L$ has 2 conjugacy classes of maximal subgroups of type $P_1$ in class $\mathcal{C}_1$: the class of the stabilizers of points and the class of the stabilizers of hyperplanes. Also, $L$ has $c = (q-1, m)$ conjugacy classes of maximal subgroups of symplectic type in class $\mathcal{C}_8$. For each maximal subgroup $X$ of symplectic type, class $\mathcal{C}_8$, and for each maximal subgroup $Y$ of type $P_1$, class $\mathcal{C}_1$, there is a factorization $L = XY$. $X$ is almost simple with socle $M$ isomorphic to $\mathrm{PSp}(2m, q)$ while $Y$ has a non trivial normal $p$-subgroup [**KL90**, 4.1.13] and therefore may not be almost simple.

We may assume that $Y$ is the stabilizer of a point which makes $\boldsymbol{S}$ the symmetric group on $\mathrm{PG}(2m, q)$. Then $G$, the normalizer of $L$ in $\boldsymbol{S}$, may not be bigger than $\mathrm{P\Gamma L}$ and $N$, the normalizer of $M$ in $\boldsymbol{S}$, may not be bigger than $\mathrm{P\Gamma Sp}$ (not even when $m = 2$ and $M$ has exceptional outer automorphisms). In fact, equality holds in both cases.

As usual we compute $r_o$ with Pálfy's formulae 2.8.1 in the assumption that $[H \div \boldsymbol{S}] \cong \mathcal{M}_r$ for some $r$; if not so, $r_o = 0$ by definition. To compute $\mathrm{hm}_{\boldsymbol{S}}(H, G)$ suppose $H^{\boldsymbol{s}} < G$ for some $\boldsymbol{s} \in \boldsymbol{S}$. Then $H^{\boldsymbol{s}}$ must be maximal in $G$ and primitive. This shows that the skeleton of $(G, H^{\boldsymbol{s}}, \mathrm{N}_G\, Y)$ must appear in our tables and hence that $H^{\boldsymbol{s}}$ is of symplectic type in class $\mathcal{C}_8$. There is one and only one conjugacy class of these subgroups in $G$ therefore $\mathrm{hm}_{\boldsymbol{S}}(H,G) = \left|G : H\right|$ and $r_o = 1$.

Although the method is the same, some more effort is needed to compute $r_e$. Start first with $e = \mathrm{hm}_{\boldsymbol{S}}(H^{\mathrm{e}}, G^{\mathrm{e}})/\left|G^{\mathrm{e}} : H^{\mathrm{e}}\right|$ which is the number of $G^{\mathrm{e}}$ conjugacy classes of maximal subgroups of $L$ of symplectic type, class $\mathcal{C}_8$. The action induced by the conjugation in $G$ on the set of the $c$ $L$-classes is described in [**KL90**]. It is a transitive action and the stabilizer of an $L$-class with notations set in 3.9 is

$$\left\langle \mathrm{PSL}, \mathrm{G}_1(\boldsymbol{u})^c, \mathrm{PG}(\boldsymbol{\pi})\right\rangle$$

Note that $\mathrm{PGL}\,/\,\mathrm{PSL} \cong \left\langle \mathrm{G}_1(\boldsymbol{u})\right\rangle$ is a cyclic group of order

$$(q-1, 2m) = \begin{cases} 2c & \text{if } q \text{ is odd,} \\ \mathrm{c} & \text{otherwise.} \end{cases}$$

The possibilities for $G^{\mathrm{e}}$ are given in table 3.C. In particular one checks that $G^{\mathrm{e}}$ has at most 2 conjugacy classes and hence $r_e \leqslant 2$. The exact value of $r_e$ is obtained by establishing $s = \left|H : H^{\mathrm{e}}\right|/\left|G : G^{\mathrm{e}}\right|$ that depends on the parity of $H$. To determine the parity of $H$, we may assume that PGSp is generated by PSp together with $\mathrm{G}_1(\boldsymbol{m})$ where $\boldsymbol{m}$ is defined by

$$\begin{aligned} e_i\boldsymbol{m} &= \mathfrak{u} e_i \qquad \text{for all } i = 0\ldots m-1, \\ e_i\boldsymbol{m} &= e_i \qquad \text{for all } i = m\ldots 2m-1. \end{aligned}$$

This $\mathrm{G}_1(\boldsymbol{m})$ moves $\frac{(q^m-1)(q^m-1)}{q-1}$ points in orbits of length $q-1$ because a point is moved if and only if at least one of the first $m$ coordinates and at least one of the second $m$ coordinates is not 0. Of course, this permutation is even if $q$ is even. On the other hand, if $q$ is odd, then

$$\mathrm{Par}\,\mathrm{G}_1(\boldsymbol{m}) = \left(\frac{q^m-1}{q-1}\right)^2 \equiv (1+q+\ldots+q^{m-1})^2 \equiv m \mod 2.$$

Therefore $\mathrm{PGSp}(2m,q)$ is even if and only if $mq$ is even. From this computation one also sees that

$$X \cong \begin{cases} \mathrm{PGSp} & \text{if } m \text{ is even,} \\ \mathrm{PSp} & \text{otherwise.} \end{cases}$$

We remind the reader that according to 3.9.5, if $p \equiv 3 \mod 4$ and $f$ is even, then $\mathrm{PG}\,\boldsymbol{\pi}$ is even if and only if $m$ is even. This is enough to write table 6.C from which one extracts the exact value of $r_e$.

We now have to show that if $\boldsymbol{H}$ matches with ♯44, then it may not match with ♯59, ♯67 or ♯57. For the first two, this is straightforward because when $q$ is even ♯44 where we set $\frac{1}{2}\,\mathrm{P\Gamma L} := \langle \mathrm{PSL}, \mathrm{PG}(\boldsymbol{\pi}^2), \mathrm{G}_1(\boldsymbol{u})\,\mathrm{PG}(\boldsymbol{\pi})\rangle$ and
$\frac{1}{2}\,\mathrm{P\Gamma Sp} := \langle \mathrm{PSp}, \mathrm{PG}(\boldsymbol{\pi}^2), \mathrm{G}_1(\boldsymbol{m})\,\mathrm{PG}(\boldsymbol{\pi})\rangle$.

TABLE 6.C. Even part of $\mathrm{P\Gamma Sp}(2m, p^f)$

| $p$ | $f$ | $m$ | $G^{\mathrm{e}}$ | $H^{\mathrm{e}}$ | $e$ | $s$ | $r_e$ |
|---|---|---|---|---|---|---|---|
| 2 | 2 | all | PGL | PGSp | 1 | 1 | 1 |
| 2 | $\neq 2$ | all | PΓL | PΓSp | 1 | 1 | 1 |
| 1 mod 4 | all | even | PΣL | PΓSp | 2 | $\frac{1}{2}$ | 1 |
| | | odd | PΣL | PΣSp | 2 | 1 | 2 |
| 3 mod 4 | odd | even | PΣL | PΓSp | 2 | $\frac{1}{2}$ | 1 |
| | | odd | PΣL | PΣSp | 2 | 1 | 2 |
| 3 mod 4 | even | even | PΣL | PΓSp | 2 | $\frac{1}{2}$ | 1 |
| | | odd | $\frac{1}{2}$ PΓL | $\frac{1}{2}$ PΓSp | 1 | 1 | 1 |

has odd degree, while in the other cases the degree is even. For ♯57 this reduces to prove that

$$(1+q^{m-1})(1+q+\ldots+q^{m-1}) < (1+q^{m-1})(1+q^{m-2})\cdots(1+q)$$

which is readily cheked.

**6.4.12.** The parametrized family of skeletons

$$(L,X,Y) = \big(\mathrm{Sp}_{2m}(q),\, \mathrm{Sp}_{2a}(q^b).b,\, \mathrm{O}^{\pm}_{2m}(q)\big), \qquad m = ab, \qquad b \text{ prime}$$

is split in two subfamilies according to $a > 1$ or $a = 1$ because $\mathrm{Sp}_{2a}(q^b) \cong \mathrm{SL}_{2a}(q^b)$ when $a = 1$. However, the analysis is not affected considerably by this isomorphism. These group factorizations have a subgroup of class $\mathcal{C}_3$ as middle term and a subgroup of type $\mathrm{O}^+$, $\mathrm{O}^-$ in class $\mathcal{C}_8$ as last term. Note that when $m = 2$ these skeletons are isomorphic to the ones in ♯49, ♯51 in view of the exceptional outer automorphism of $\mathrm{Sp}_4(q)$ which swaps the $\mathrm{O}^+_4(q)$ with the subgroups in class $\mathcal{C}_2$ and the $\mathrm{O}^-_4(q)$ with the subgroups in class $\mathcal{C}_3$. In order to avoid crossings in Table 6.B, we subdivide the subfamily for $a = 1$ into ♯46 (where we assume $m > 2$), ♯49 and ♯51 while the subfamily for $a > 1$ corresponds to ♯45. We now proceed with $L$, $X$, $Y$ as above.

First of all, $G$ contains and in fact equals $\mathrm{P\Gamma Sp}_{2m}(q)$ because that is the largest subgroup of $\mathrm{Aut}\,L$ which preserves the $L$-conjugacy class of $\mathrm{O}^{\pm}_{2m}(q)$ (this is true even when $m = 2$).

By [**LPS90**, 3.2.1d] stabilizers of points in $M$ are subgroups of class $\mathcal{C}_8$ when $a > 1$ and subgroups of class $\mathcal{C}_2$ or $\mathcal{C}_3$ when $a = 1$. As the $\mathrm{P\Gamma Sp}_{2a}(q^b)$-conjugacy classes are $M$-conjugacy classes in each of those cases, $N = \mathrm{P\Gamma Sp}_{2a}(q^b)$. Since $N$ is contained in $G$, $H = N$.

Provided $\boldsymbol{H} \in \{H, H^{\mathrm{e}}\}$,

$$\frac{\left|\mathrm{N}(\boldsymbol{H}) : \boldsymbol{H}\right|}{\left|\mathrm{N}(G \cap \boldsymbol{U}) : G \cap \boldsymbol{U}\right|} = 1.$$

We claim that $\mathrm{hm}(\boldsymbol{H}, G \cap \boldsymbol{U}) = \left|G \cap \boldsymbol{U} : \boldsymbol{H}\right|$ yielding both $r_o \leqslant 1$ and $r_e = 1$.

In fact, if $\boldsymbol{H}^{\boldsymbol{s}} < G$ with $\boldsymbol{s} \in \boldsymbol{S}$, then $\boldsymbol{H}^{\boldsymbol{s}} \lessdot G \cap \boldsymbol{U}$ and $\boldsymbol{H}^{\boldsymbol{s}}$ is primitive of type almost simple. Therefore the socle of $\boldsymbol{H}^{\boldsymbol{s}}$ appears as middle term of a skeleton in our tables such that the first term is a symplectic group and the last term is one of its subgroups of class $\mathcal{C}_8$. A quick view leaves us with $\boldsymbol{H}^{\boldsymbol{s}}$ in classes $\mathcal{C}_3$, $\mathcal{C}_8$ (when $q = 2 \neq m$), or class $\mathcal{S}$. But none of the shown subgroups in class $\mathcal{C}_8$ or in class $\mathcal{S}$ is isomorphic to $\boldsymbol{H}$, therefore $\boldsymbol{H}^{\boldsymbol{s}}$ is in class $\mathcal{C}_3$. Then $\boldsymbol{H}^{\boldsymbol{s}}$ is conjugated to $\boldsymbol{H}$ in $L$ and *a fortiori* in $G \cap \boldsymbol{U}$. This completes the proof because $\boldsymbol{H}$ is selfnormalizing in $G \cap U$, so, the number of conjugates of $\boldsymbol{H}$ is the index of $\boldsymbol{H}$.

Suppose now that $\boldsymbol{H}$ matches ♯45. We claim that $\boldsymbol{H}$ does not match any of ♯75, ♯59, ♯67, ♯57 and ♯44. Because $b$ is prime, $q^b \neq 2$ and this rules out ♯75 and ♯67. Moreover, the degree is even and this rules out ♯57 and ♯44. If $\boldsymbol{H}$ matches ♯59, then

$$\frac{q^{3b}(q^{3b} \pm 1)}{2} = q^{3b}(q^{4b} - 1)$$

which is impossible.

**6.4.13.** The parametrised family of skeletons

$$\left(\mathrm{Sp}_{2m}(q),\, \mathrm{O}^{\pm}_{2m}(q),\, \mathrm{Sp}_{2a}(q^b).b\right)$$

is split in the three subfamilies ♯47, ♯48 and ♯49 because when $b = 2$ only type $\mathrm{O}^-$ appears as middle term and when $(a, b) = (1, 2)$, $\mathrm{O}^-_4(q) \cong \mathrm{PSL}_2(q^2)$. In all cases, both in $L$ and in $M$, stabilizers of points are subgroups of class $\mathcal{C}_3$. There is one and only one conjugacy class of those subgroups; This class is preserved by the related automorphism group. Therefore $G$ and $N$ are isomorphic to the automorphism groups of $L$, $M$ respectively. A similar argument to the one at the end of 6.4.12 shows that $r_o \leqslant 1$ and $r_e = 1$.

**6.4.14.** Consider ♯50, ♯51 where

$$(L, X, Y) = (\mathrm{Sp}_{2m}(q), \mathrm{O}^-_{2m}(q), \mathrm{Sp}_m(q) \,\mathrm{wr}\, \mathrm{S}_2).$$

From [**LPS90**, 3.2.4b]

$$X \cap Y \cong \mathrm{O}^-_m(q) \times \mathrm{O}^+_m(q).$$

This is a subgroup of class $\mathcal{C}_1$ of $\mathrm{O}_{2m}^-(q)$. Note that if $m = 2$, then

$$\mathrm{O}_4^-(q) \cong \mathrm{PSL}_2(q^2).2$$

and $X \cap Y$ corresponds to a subgroup of class $\mathcal{C}_2$ of $\mathrm{PSL}_2(q^2).2$. More details about this special case are given in 6.4.12

**6.4.15.** In ♯55 $L = \mathrm{P}\Omega_{2m}^+(q)$; adopting notations of [**KL90**] $X$ is a subgroup of type $\mathrm{GL}_m(q).2$ class $\mathcal{C}_2$, $Y$ is a subgroup of type $\mathrm{O}_1 \perp \mathrm{O}_{2m-1}$ (or $\mathrm{Sp}_{2m-2}(q)$ when $q = 2$) class $\mathcal{C}_1$. We may assume that $L = \mathrm{P}\Omega(\mathsf{V}, Q)$, $X = \mathrm{St}_L\{\mathsf{E}, \mathsf{F}\}$, $Y = \mathrm{St}_L \langle \mathsf{v} \rangle$ where

$$\mathsf{E} = \langle \mathsf{e}_0, \ldots, \mathsf{e}_{m-1} \rangle, \qquad \mathsf{F} = \langle \mathsf{f}_0, \ldots, \mathsf{f}_{m-1} \rangle,$$

$\mathsf{v} = \mathsf{e}_0 + \alpha \mathsf{f}_0$ and $\vec{b} = (\mathsf{e}_0, \ldots, \mathsf{e}_{m-1}, \mathsf{f}_0, \ldots, \mathsf{f}_{m-1})$ is a basis of $\mathsf{V}$ such that

$$(\mathsf{u} + \mathsf{v})Q = \mathsf{u}Q + \mathsf{v}Q + \mathsf{u} \circ \mathsf{v},$$
$$\mathsf{e}_i \circ \mathsf{f}_j = \delta_{ij}, \qquad \mathsf{e}_i Q = \mathsf{f}_i Q = 0.$$

Assume first $q = 3$ and $m$ odd.

There are exactly 2 conjugacy classes of subgroups of $L$ of type $\mathrm{O}_1 \perp \mathrm{O}_{2m-1}$ and the largest subgroup of $\mathrm{Aut}\, L$ preserving those 2 classes is $\mathrm{PO}_{2m}^+(3)$.

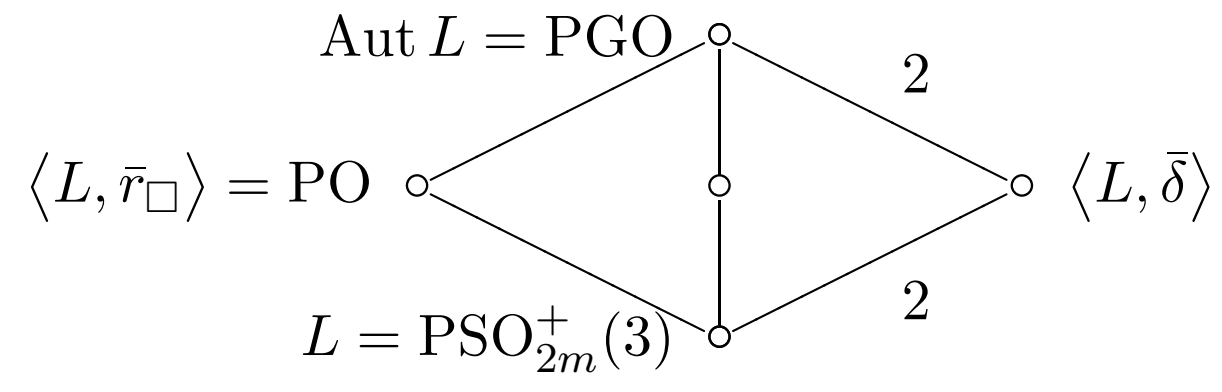

Thus $G = \mathrm{PO}_{2m}^+(3)$. The elements of $X$ are represented with respect to $\vec{b}$ by matrices

$$\begin{bmatrix} \mathsf{A} & 0 \\ 0 & \mathsf{B} \end{bmatrix}$$

where $\mathsf{A} \in \mathrm{GL}_m(3)$ and $\mathsf{B} = (\mathsf{A}^t)^{-1}$. Thus $X \cong \mathrm{PGL}_m(3) = \mathrm{PSL}_m(3)$. We claim that $N$, the full normalizer of $X$, is isomorphic to $\mathrm{Aut}\, X$ and that $N < G$. Note first that $|\mathrm{Aut}\, X : X| = 2$. Now consider $\boldsymbol{g} \in \mathrm{GL}_{2m}(3)$ such that

$$\mathsf{e}_0 \boldsymbol{g} = \alpha \mathsf{f}_0, \qquad (\alpha \mathsf{f}_0)\boldsymbol{g} = \mathsf{e}_0,$$
$$\mathsf{e}_i \boldsymbol{g} = \mathsf{f}_i, \qquad \mathsf{f}_i \boldsymbol{g} = \mathsf{e}_i \quad \forall i > 0.$$

This $\boldsymbol{g}$ preserves $Q$, therefore $\bar{\boldsymbol{g}} \in \mathrm{PO}_{2m}^+(3)$. Moreover, $\boldsymbol{g}$ preserves $\mathsf{v}$ and $\{\mathsf{E}, \mathsf{F}\}$, therefore $\bar{\boldsymbol{g}} \in N$. Since $\bar{\boldsymbol{g}} \notin X$,

$$\mathrm{Aut}\, X \cong \langle X, \bar{\boldsymbol{g}} \rangle < \mathrm{PO}_{2m}^+(3).$$

Regarding $r_o$ and $r_e$, note that there are at most 2 conjugacy classes of subgroups of the same type as $X$ in both $G$ and $G^{\mathrm{e}}$. Therefore by the Pálfy Lemma, both $r_o$ and $r_e$ are less than or equal to 2.

When $q = 3$ and $m$ is even $\operatorname{Aut} L/L$ is isomorphic to the dihedral group of order 8 and $\operatorname{Aut} M/M$ is isomorphic to the elementary abelian group of order 4:

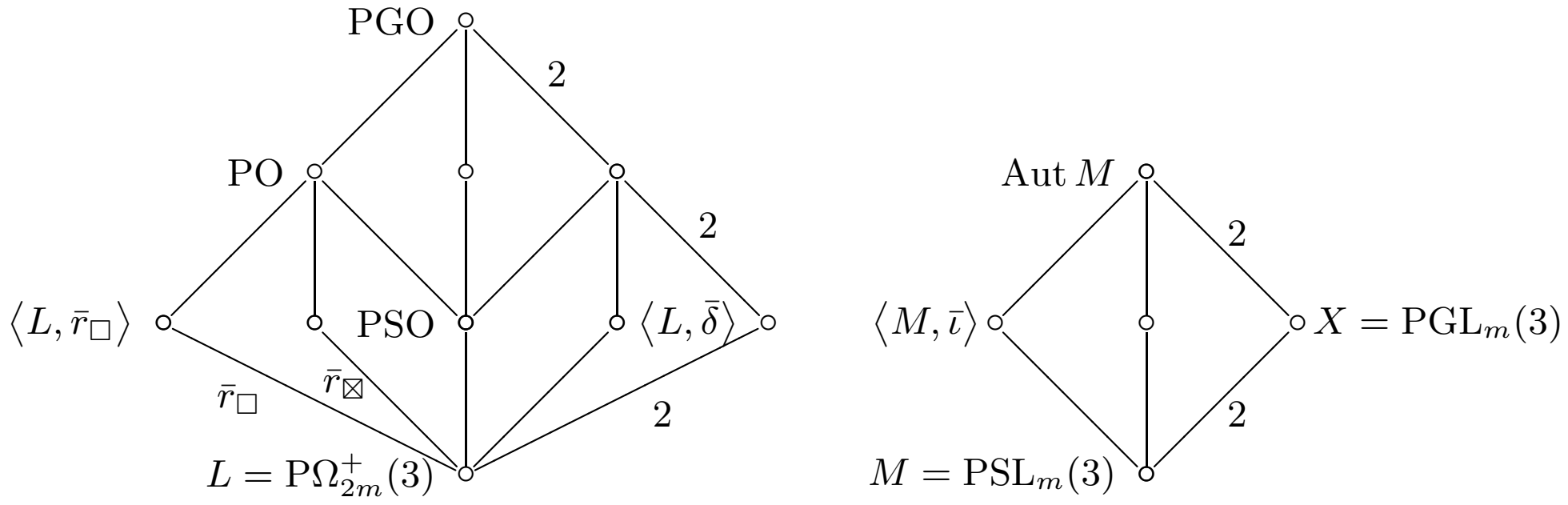

Again, $G = \mathrm{PO}^+_{2m}(3)$ and $\operatorname{Aut} M \cong N < G$. However, because $m$ is even, the $\boldsymbol{g}$ above belongs to $\mathrm{PSO}^+_{2m}(3)$ and so $N < \mathrm{PSO}$. If $H$ (which is equal to $N$) is even, then it is low. Otherwise $H^{\mathrm{e}} = X$ because $X < L$ which is simple. Therefore $H^{\mathrm{e}}$ is low again and $r_o = r_e = 0$ in this case.

As far as the case $q = 2$ is concerned, there is only one conjugacy class of subgroups of the same type as $Y$, therefore the full normalizer of $L$ is isomorphic to $\operatorname{Aut} L$. But now

$$\mathrm{PGO}^+_{2m}(2) = \mathrm{PO}^+_{2m}(2)$$

and so $G = \mathrm{PO}^+_{2m}(2)$ holds. As in the case $q = 3$ but $m$ odd, $X = M$ when $q = 2$ and $|\operatorname{Aut} M : M| = 2$. Arguing as above, one finds that again $\operatorname{Aut} M \cong N < G$ and that both $r_o$ and $r_e$ are less than or equal to 2.

**6.4.16.** We now deal with skeletons

$$(L, X, Y) = (\mathrm{P}\Omega^+_{2m}(q), N_1, P_i)$$

where $P_i$ is the stabilizer of a totally singular subspace of dimension $m$ (when $i \in \{m, m-1\}$) or 1 (only but not necessarily if $m = 4$) and the middle term, $N_1$, is the stabilizer of a non singular point. Referring to notations and results of both [**KL90**] and [**Kle87**] we have

$$X \cong \begin{cases} \mathrm{Sp}_{2m-2}(q) & \text{if } q \text{ is even,} \\ \Omega_{2m-1}(q) & \text{if } q \text{ is odd and } D = \square, \\ \Omega_{2m-1}(q).2 & \text{if } q \text{ is odd and } D = \boxtimes. \end{cases}$$

To discuss $G$, the full normalizer of $L$, it is convenient to view the automorphism group of $L$ first. When $m \neq 4$ this is isomorphic to $\mathrm{P\Gamma O}$ which is given by $\mathrm{PGO} \rtimes \langle \mathrm{PG}(d, q, \boldsymbol{\pi}) \rangle$ where $\boldsymbol{\pi}$ is the Frobenius automorphism. However, to match notations of [**KL90**], for the remainder of this section we put

$$\bar{\phi} = \mathrm{PG}(d, q, \boldsymbol{\pi}).$$

The structure of $\mathrm{PGO} / \mathrm{P}\Omega$ is known

$$\mathrm{PGO} / \mathrm{P}\Omega \cong \begin{cases} 2 & \text{if } q \text{ is even,} \\ \mathrm{D}_8 & \text{if } q \text{ is odd and } D = \square, \\ 2 \times 2 & \text{if } q \text{ is odd and } D = \boxtimes. \end{cases}$$

the subgroup lattice for $q$ odd matches the one given for $q = 3$ in 6.4.15. When $m = 4$ the automorphism group of $L$ is isomorphic to $\mathrm{S}_3$ or $\mathrm{S}_4$ depending on whether $q$ is even or odd. Note that if $q$ is odd and $m = 4$, then $D = \square$ and so $\mathrm{PGO} / \mathrm{P}\Omega$ is a dihedral group of order 8. The full normalizer of $L$ is given by [**KL90**, 2.7.4(iii)] or by the Result Matrix of [**Kle87**] as the largest subgroup of $\mathrm{Aut}\, L$ preserving the $L$-conjugacy class of $Y$. It is a subgroup of index 2 of $\mathrm{P\Gamma O}$ which contains $\bar{\phi}$ but avoids PO.

The intersection $X \cap Y$ is a subgroup of type $P_{m-1}$ class $\mathcal{C}_1$ of $X$, of which there is one and only one $M$-conjugacy class. Therefore

$$N \cong \mathrm{Aut}\, M \cong (M.2) \rtimes \langle \bar{\phi} \rangle.$$

Clearly, $\langle X, \bar{\phi} \rangle$ is contained in $H$, the intersection of $N$ with $G$. In fact, $\langle X, \bar{\phi} \rangle = N$ when $q$ is even or $q$ is odd but $D = \boxtimes$. When $q$ is odd and $D = \square$, note that $N \cap \mathrm{PGO}$ is contained in the subgroup generated by $L$ and the reflection in the non singular point stabilized by $X$. Therefore $N \cap \mathrm{PGO} \cap G$ is contained in $L$, that is,

$$N \cap \mathrm{PGO} \cap G = X$$

and hence $H = \langle X, \bar{\phi} \rangle$ again.

Since $X$ is even (being contained in the simple group $L$), the parity of $H$ is equal to the parity of $\bar{\phi}$. It follows that if both $H$ and $q$ are odd, then both $H$ and $H^{\mathrm{e}}$ are low. We leave it to the reader to check that in all other cases both $r_o$ and $r_e$, provided it makes sense to consider them, are at most 1.

**6.4.17.** Consider ♯58 first, so $q$ is odd. Following [**Kle87**] say $R_{+1}$, $K_1^1$, $K_1^2$, $R_{-1}$, $K_1^3$, $K_1^4$ the 6 $L$-conjugacy classes whose union is the $(\mathrm{Aut}\, L)$-conjugacy class of $Y$. We may assume that $Y \in R_{+1}$ ($R_{+1}$, $R_{-1}$ are the two conjugacy classes of stabilizers of non singular points). In [**Kle87**] it is shown that $G = \langle \mathrm{PO}, \bar{\phi} \rangle$. This $G$ fixes $R_{+1}$, $R_{-1}$ and $G/L$ acts regularly on the set of the other 4 conjugacy classes $K_1^i$. By 3.3.2, $X$ may not lift to a maximal subgroup of $G$. Moreover, even if $G$ is odd, $G^{\mathrm{e}}$ still moves each $K_1^i$ and so $X$ may not lift to a maximal subgroup of $G^{\mathrm{e}}$ either. Therefore $H$ and $H^{\mathrm{e}}$ are low subgroups in this case.

We appeal to a similar argument for ♯59 to show that $H$ is low. However, if $H$ is even but $G$ is odd, $H$ is maximal in $G^{\mathrm{e}}$. A conjugate of $H$ by an element of $\boldsymbol{S}$ may only lie in $K_1^1$ or in $K_1^2$, so, $\mathrm{hm}(H, G^{\mathrm{e}}) = 2\big|G^{\mathrm{e}} : H\big|$. But $\big|G : G^{\mathrm{e}}\big| = 2$ and hence $r_e = \mathrm{hm}(G^{\mathrm{e}}, H) = 1$. This was only in the assumption that $H$ is even and $G$ is odd which I do not investigate whether it does actually happen or not.

**6.4.18. Skeleton $(L, X, Y) = (\Omega_7(q), \mathrm{G}_2(q), N_1^{\pm})$, $q$ odd.** Note that there are two $L$-conjugacy classes of $\mathrm{G}_2(q)$ in $L$ and that for each of those subgroups there are both factorizations

$$L = \mathrm{G}_2(q)N_1^{+} \qquad \text{and} \qquad L = \mathrm{G}_2(q)N_1^{-}.$$

Clearly, the full normalizer of $L$ is

$$G = \mathrm{P\Gamma O}_7(q) \cong L.(2 \times f)$$

where $q = p^f$ and $p$ is prime. Observe that adopting [**Asc87**] notations $X \cap Y$ is a subgroup of $X$ of type $\mathcal{M}_4$ (when $Y = N_1^{+}$) or $\mathcal{M}_5$ (when $Y = N_1^{-}$). Then by [**Asc87**, 17.3, 17.4] the full normalizer of $X$, $N$, is $\langle X, \bar{\phi}\rangle$ where $\bar{\phi}$ is the Frobenius automorphism. Therefore $N$ is contained in $\langle L, \bar{\phi}\rangle$ which is contained in $G$ and so $H = N$ and $H$ is low. If $\bar{\phi}$ is odd, again $H^{\mathrm{e}} = N^{\mathrm{e}}$ and $H^{\mathrm{e}}$ is low. Otherwise, $H$ is even and by [**KL90**, 2.6.3] the parity of $G$ is the parity of $\bar{r}_{\square}\bar{r}_{\boxtimes}$ where $\bar{r}_{\square}$, $\bar{r}_{\boxtimes}$ are reflections in vectors whose norm is a square or a non square. One can ascertain the parity of those reflections and see that the parity of $G$ depends on $q \mod 8$. We prefer to skip these details observing that even when $H$ is even and $G$ is odd, $\big|G : G^{\mathrm{e}}\big| = 2$ and $\mathrm{hm}(H, G^{\mathrm{e}}) = 2\big|G^{\mathrm{e}} : H\big|$, so, by Pálfy Formula $r_e = 1$.

A similar argument leads to the same conclusion when $Y$ is the stabilizer of a totally singular point. In this case $X \cap Y$ is a maximal parabolic subgroup of $X$ of type $\mathcal{M}_1$.

APPENDIX A

# Intervals in the subgroup lattices of direct products

## Introduction

The subgroup lattice of a group is the set of all its subgroups endowed with the partial order given by the inclusion among subgroups. We denote the subgroup lattice of a group $G$ by $\mathcal{S}\mathrm{ub}\, G$. An interval $[H \div K]$ of this subgroup lattice, where $H$ and $K$ are subgroups of $G$ such that $H$ is also a subgroup of $K$, is simply the sublattice made of all the subgroups of $K$ containing $H$.

The subgroups of a direct product of two groups have been known for more than a century now but we understand that very little has been written on how these subgroups relate to each other and in particular on the intervals of the related subgroup lattice. In this paper we address the following issues.

(1) Which intervals of $\mathcal{S}\mathrm{ub}(L \times R)$ are not intervals of $\mathcal{S}\mathrm{ub}\, L$ or $\mathcal{S}\mathrm{ub}\, R$?
(2) Which intervals of $\mathcal{S}\mathrm{ub}(L \times R)$ are not intervals of $\mathcal{S}\mathrm{ub}(L^2)$ or $\mathcal{S}\mathrm{ub}(R^2)$?
(3) Describing the maximal inclusions among subgroups of $L \times R$.
(4) Describing the *shortcut intervals* of $\mathcal{S}\mathrm{ub}(L \times R)$, that is, the intervals $[H \div K]$ such that for some subgroup $G$, $H \lessdot G \lessdot K$.

A solution to problem (3) for example, was given in [**Fla94**] for the case where both $L$ and $R$ are abelian, and in [**Thé97**] for the case where $K = L \times R$.

In this paper we provide the tools to deal with the previous problems in the general case of groups with operators, in fact, of algebras satisfying A.1.6. These tools are four canonical homomorphisms which are defined for any given interval of $\mathcal{S}\mathrm{ub}(L \times R)$. Our main theorem, Theorem A, says that each interval is a composition of at most four intervals whose canonical homomorphisms are either isomorphisms or trivial in pairs. An immediate consequence of this theorem is the description of the maximal inclusions in the subgroup lattice of a direct product.

Section A.4 is devoted to the description of the shortcut intervals. Theorem B says that if these intervals are not already intervals of $\mathcal{S}\mathrm{ub}\, L$, $\mathcal{S}\mathrm{ub}\, R$, then they are particular instances of the intervals represented in Figure A.A, A.F, A.G.

In the next section we introduce the terminology and notations used in this paper. Most of them are standard and we tried to emphasize the new ones appropriately.

## A.1. Preliminaries and Notations

**Binary relations.** By *binary relation* on a set $\mathsf{X}$ we mean a subset of $\mathsf{X}^2$. If $\pitchfork$ is a binary relation and if $(\mathsf{x}_1, \mathsf{x}_2) \in \pitchfork$, then we write $\mathsf{x}_1 \pitchfork \mathsf{x}_2$. If $\mathsf{A}$ is a subset of $\mathsf{X}$, then the image of $\mathsf{A}$ by $\pitchfork$ is

$$\mathsf{A} \pitchfork := \left\{ \mathsf{x} \in \mathsf{X} \;\middle|\; \exists \mathsf{a} \in \mathsf{A} \ \mathsf{a} \pitchfork \mathsf{x} \right\},$$

while the restriction of $\pitchfork$ to $\mathsf{A}$ is

$$\pitchfork^{\mathsf{A}} := \mathsf{A}^2 \cap \pitchfork . \qquad \text{(A.1.A)}$$

We say that a binary relation is contained in another one if it is contained as a subset. Given two binary relations $\pitchfork$, $\between$, there is a composition defined by

$$\mathsf{x}_1 \pitchfork \between \mathsf{x}_2 \qquad \text{if and only if} \qquad \exists \mathsf{x}_3 \quad \mathsf{x}_1 \pitchfork \mathsf{x}_3 \between \mathsf{x}_2 .$$

We say that the binary relations $\pitchfork$, $\between$ are *permutable* when $\pitchfork \between = \between \pitchfork$.

**Posets.** A *partially ordered set* or *poset* is a set $\mathcal{P}$ together with a binary relation $\leqslant$ such that the following conditions are satisfied for all $\mathsf{x}, \mathsf{y}, \mathsf{z} \in \mathcal{P}$:

$$\mathsf{x} \leqslant \mathsf{x},$$

$$\text{if } \mathsf{x} \leqslant \mathsf{y} \text{ and } \mathsf{y} \leqslant \mathsf{x} \text{ then } \mathsf{x} = \mathsf{y},$$

$$\text{if } \mathsf{x} \leqslant \mathsf{y} \text{ and } \mathsf{y} \leqslant \mathsf{z} \text{ then } \mathsf{x} \leqslant \mathsf{z}.$$

Then we write $\mathsf{x} < \mathsf{y}$ when $\mathsf{x} \leqslant \mathsf{y}$ and $\mathsf{x} \neq \mathsf{y}$. Two elements $\mathsf{x}$, $\mathsf{y}$ are said to be *comparable* when

$$\mathsf{x} \leqslant \mathsf{y} \qquad \text{or} \qquad \mathsf{y} \leqslant \mathsf{x}.$$

A *chain* is a subset of $\mathcal{P}$ whose elements are pairwise comparable. If $\mathsf{x} < \mathsf{y}$ and there is no element $\mathsf{z} \in \mathcal{P}$ such that $\mathsf{x} < \mathsf{z} < \mathsf{y}$, then we say that $\mathsf{x}$ is maximal in $\mathsf{y}$ and we write $\mathsf{x} \lessdot \mathsf{y}$.

A homomorphism of posets is a map between posets which preserves the order. The inverse of a bijective homomorphism needs not to be a homomorphism, when it is, the homomorphism is said to be an isomorphism. Given a subset $\mathcal{S}$ of $\mathcal{P}$, a *lower bound* of $\mathcal{S}$, if there is one, is an element $\mathsf{b} \in \mathcal{P}$ such that $\mathsf{b} \leqslant \mathsf{s}$ for all $\mathsf{s} \in \mathcal{S}$. A *minimum* of $\mathcal{S}$, if there is one, is a lower bound of $\mathcal{S}$ which lies in $\mathcal{S}$. If there is a minimum of $\mathcal{S}$ then it is unique and we write $\mathsf{m} = \min \mathcal{S}$. Similarly we define *upper bound*s and *maximum* of $\mathcal{S}$ ($\max \mathcal{S}$). $\mathcal{S}$ is *bounded* if it has a minimum and a maximum. We denote the *greatest lower bound* of $\mathcal{S}$ by $\bigwedge \mathcal{S}$; this is the maximum of the lower bounds of $\mathcal{S}$ and is sometime called the *meet* of $\mathcal{S}$. Similarly the *least upper bound* or *join* of $\mathcal{S}$ is denoted by $\bigvee \mathcal{S}$. An element $\mathsf{x} \in \mathcal{S}$ is called *minimal* (*maximal*) if there is no $\mathsf{s} \in \mathcal{S}$ such that $\mathsf{s} < \mathsf{x}$ ($\mathsf{x} < \mathsf{s}$). However, by abuse of language, when $\mathcal{S}$ is bounded the minimal elements or *atoms* of $\mathcal{S}$ are the minimal elements of $\mathcal{S}-\{\min\mathcal{S}\}$ and the maximal elements of $\mathcal{S}$ are the maximal elements of $\mathcal{S}-\{\max\mathcal{S}\}$. The *length* of $\mathcal{S}$ is the number of proper inclusions in a maximal chain of $\mathcal{S}$. It is easy to check that for each positive integer $r$, there is only one bounded poset with length 2 and order $r+2$ which we denote by $\mathcal{M}_r$, see Figure A.A.

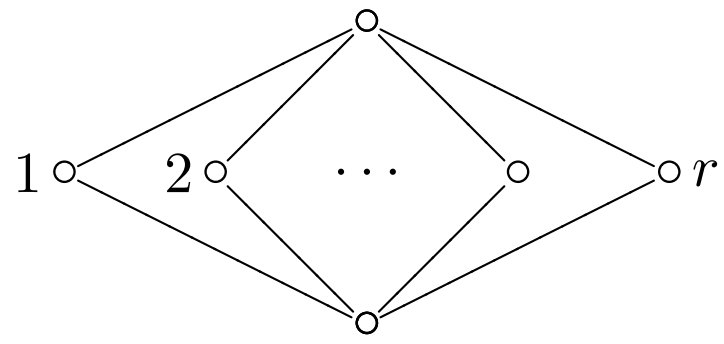

FIGURE A.A. $\mathcal{M}_r$ has length 2 and order $r+2$

**Intervals.** Given $\mathsf{a} \leqslant \mathsf{b}$ in a poset $\mathcal{P}$ we define the *interval* $[\mathsf{a} \div \mathsf{b}]$ as the set of the elements $\mathsf{c} \in \mathcal{P}$ such that $\mathsf{a} \leqslant \mathsf{c} \leqslant \mathsf{b}$. In particular, an interval is a non empty bounded subset. An interval $[\mathsf{a} \div \mathsf{b}]$ is said to be *trivial* if $\mathsf{a} = \mathsf{b}$, *simple* if $\mathsf{a} \lessdot \mathsf{b}$, and *composed* otherwise. If $\mathcal{I}_1$, $\mathcal{I}_2$ are intervals such that $\max\mathcal{I}_1 = \min\mathcal{I}_2$ then $\mathcal{I}_1 \cup \mathcal{I}_2$ is usually far from being an interval. However, we define a composition

$$\mathcal{I}_1 \rtimes \mathcal{I}_2 = [\min\mathcal{I}_1 \div \max\mathcal{I}_2]. \tag{A.1.B}$$

A.1.1. *Definition*. We say that $\mathcal{I}$ has a *shortcut* or even that $\mathcal{I}$ is a *shortcut interval*, if $\mathcal{I} = \mathcal{I}_1 \rtimes \mathcal{I}_2$ for convenient simple intervals $\mathcal{I}_1$, $\mathcal{I}_2$.

**Equivalences.** An *equivalence* on $\mathsf{X}$ is a binary relation $\varepsilon$ on $\mathsf{X}$ such that the following conditions are satisfied for all $\mathsf{x},\mathsf{y},\mathsf{z} \in \mathsf{X}$:

$$\mathsf{x}\ \varepsilon\ \mathsf{x},$$

$$\text{if } \mathsf{x}\ \varepsilon\ \mathsf{y} \text{ then } \mathsf{y}\ \varepsilon\ \mathsf{x},$$

$$\text{if } \mathsf{x}\ \varepsilon\ \mathsf{y} \text{ and } \mathsf{y}\ \varepsilon\ \mathsf{z} \text{ then } \mathsf{x}\ \varepsilon\ \mathsf{z}.$$

There are always two trivial equivalences which are $\mathsf{X}^2$ and its diagonal

$$\boldsymbol{\Delta}_{\mathsf{X}} = \{\ (\mathsf{x},\mathsf{x})\ \mid\ \mathsf{x} \in \mathsf{X}\ \}. \tag{A.1.C}$$

If $\varepsilon$ is an equivalence and $\mathsf{x} \in \mathsf{X}$ then the class of $\mathsf{x}$ is

$$[\mathsf{x}]_\varepsilon := \{\ \mathsf{y} \in \mathsf{X}\ \mid\ \mathsf{x}\ \varepsilon \mathsf{y}\ \}.$$

The *quotient* of $\mathsf{X}$ by $\varepsilon$ is $\mathsf{X}/\varepsilon := \{\ [\mathsf{x}]_\varepsilon\ \mid\ \mathsf{x} \in \mathsf{X}\ \}$, this is always a partition of $\mathsf{X}$, that is, either $[\mathsf{x}]_\varepsilon \cap [\mathsf{y}]_\varepsilon = \emptyset$ or $[\mathsf{x}]_\varepsilon = [\mathsf{y}]_\varepsilon$.

If $\boldsymbol{f} : \mathsf{X} \longrightarrow \mathsf{Y}$, then the binary relation

$$\mathsf{x}\ \kappa_{\boldsymbol{f}}\ \mathsf{y} \qquad \text{iff} \qquad \mathsf{x}\boldsymbol{f} = \mathsf{y}\boldsymbol{f} \tag{A.1.D}$$

is an equivalence and all the equivalences arise in this way. The classes $[\mathsf{x}]_{\kappa_{\boldsymbol{f}}}$ are also called *fibers* of $\boldsymbol{f}$. The arbitrary intersection of equivalences (as subsets of $\mathsf{X}^2$) is again an equivalence, hence any set of equivalences has a meet. Moreover, there is the join which is not the union (as it would be if they were just binary relations) but the intersection of all the equivalences which contain the union. The following lemma is well known and plays a key role in §A.2:

A.1.2. Lemma. *If $\varepsilon$, $\vartheta$ are equivalences then the following are equivalent:*

(a) *$\varepsilon\,\vartheta$ is an equivalence,*
(b) *$\varepsilon\,\vartheta = \vartheta\,\varepsilon$,*
(c) *$\varepsilon\,\vartheta = join\{\varepsilon, \vartheta\}$.*

**Algebras.** We usually refer to the terminology used in [**Grä79**]. However, roughly speaking, by algebra we mean a set with operations. By subalgebra we mean a subset which is closed under these operations and we write $\mathfrak{A} \leqslant \mathfrak{B}$ when $\mathfrak{A}$ is a subalgebra of $\mathfrak{B}$. If the symbol $*$ denotes an operation of arity $m$ of the algebras $\mathfrak{A}$, $\mathfrak{B}$, then by $*$-homomorphism we mean a map $\boldsymbol{f} : \mathfrak{A} \to \mathfrak{B}$ such that

$$(\mathfrak{a}_1, \ldots, \mathfrak{a}_m) * \boldsymbol{f} = (\mathfrak{a}_1\boldsymbol{f}, \ldots, \mathfrak{a}_m\boldsymbol{f}) * .$$

If no ambiguity arises, it is common practice to leave the specifications of the operations preserved by a homomorphism to the context. The *congruences* of an algebra are the equivalences induced by homomorphisms. Then, the quotient sets inherit from the algebra the same operations which were preserved by the homomorphisms, in the same way in which a factor group inherits the group operation from the original group. If $\vartheta$ is a congruence of $\mathfrak{B}$ and if $\mathfrak{A} \leqslant \mathfrak{B}$, then $\mathfrak{A}\,\vartheta$ is a subalgebra of $\mathfrak{B}$ and $\mathfrak{A} \leqslant \mathfrak{A}\,\vartheta$. We say that $\mathfrak{A}$ is $\vartheta$-closed when $\mathfrak{A} = \mathfrak{A}\,\vartheta$. Finally, by *section* of $\mathfrak{B}$ we mean a quotient $\mathfrak{A}/\,\alpha$ where $\alpha$ is a congruence of some subalgebra $\mathfrak{A}$ of $\mathfrak{B}$.

**Lattices.** A *lattice* is a poset $\mathcal{L}$ such that for every pair of elements $x, y \in \mathcal{L}$, there is the greatest lower bound $x \wedge y$ and the least upper bound $x \vee y$. Thus a lattice $\mathcal{L}$ is an algebra $\langle \mathfrak{L}; \wedge, \vee \rangle$ such that the following laws hold:

$$\begin{aligned} x \wedge y &= y \wedge x, & x \vee y &= y \vee x, \\ (x \wedge y) \wedge z &= x \wedge (y \wedge z), & (x \vee y) \vee z &= x \vee (y \vee z), \\ x \wedge (x \vee y) &= x, & x \vee (x \wedge y) &= x. \end{aligned}$$

Conversely, any such algebra can be regarded as a poset by defining

$$x \leqslant y \qquad \text{iff} \qquad x = x \wedge y.$$

A lattice is said to be *meet complete, join complete*, or *complete* when for any subset $\mathsf{S}$ there is $\bigwedge \mathsf{S}$, $\bigvee \mathsf{S}$ or both respectively. Note that a meet complete bounded lattice is a complete lattice because for any subset $\mathsf{S}$

$$\bigvee \mathsf{S} = \bigwedge \{\, x \in \mathfrak{L} \;\mid\; \forall s \in \mathsf{S},\; s \leqslant x \,\}.$$

Example of complete bounded lattices are

$$\begin{aligned} \mathcal{S}\mathrm{ub}\,\mathfrak{A} &\equiv \text{subalgebras of } \mathfrak{A}, \\ \mathcal{R}\mathrm{el}\,\mathfrak{A} &\equiv \text{binary relations on } \mathfrak{A}, \\ \mathcal{E}\mathrm{quiv}\,\mathfrak{A} &\equiv \text{equivalences on } \mathfrak{A}, \\ \mathcal{C}\mathrm{on}\,\mathfrak{A} &\equiv \text{congruences of } \mathfrak{A}, \\ \mathcal{S}\mathrm{ub}_{\vartheta}\,\mathfrak{A} &\equiv \vartheta\text{-closed subalgebras of } \mathfrak{A},\ \vartheta \in \mathcal{C}\mathrm{on}\,\mathfrak{A}. \end{aligned}$$

By sublattice, meet-sublattice and join-sublattice we mean subalgebra for the operations $\{\wedge, \vee\}$, $\{\wedge\}$ and $\{\vee\}$ respectively. For example $\mathcal{S}\mathrm{ub}_{\vartheta}\,\mathfrak{A}$ is a sublattice of $\mathcal{S}\mathrm{ub}\,\mathfrak{A}$ and $\mathcal{C}\mathrm{on}\,\mathfrak{A}$ is a complete meet-sublattice of $\mathcal{E}\mathrm{quiv}\,\mathfrak{A}$ which is a complete meet-sublattice of $\mathcal{R}\mathrm{el}\,\mathfrak{A}$. In fact, $\mathcal{C}\mathrm{on}\,\mathfrak{A}$ is a sublattice of $\mathcal{E}\mathrm{quiv}\,\mathfrak{A}$ [**Grä79**, I, §10] but $\mathcal{E}\mathrm{quiv}\,\mathfrak{A}$ is not necessarily a sublattice of $\mathcal{R}\mathrm{el}\,\mathfrak{A}$. Similarly, by homomorphism, meet-homomorphism, join-homomorphism we mean a homomorphism for the operations $\{\wedge, \vee\}$, $\{\wedge\}$ and $\{\vee\}$ respectively. For example any algebra homomorphism from $\mathfrak{A}$ to $\mathfrak{B}$ induces a join-homomorphism from $\mathcal{S}\mathrm{ub}\,\mathfrak{A}$ to $\mathcal{S}\mathrm{ub}\,\mathfrak{B}$. It is well known that a bijective map between lattices is an isomorphism iff it is a meet-homomorphism, iff it is a join-homomorphism, iff it is an isomorphism of posets.

**Subcongruence Lattice.** By *subcongruence* we mean a congruence of a subalgebra. If we regard a subcongruence $\gamma$ of a subalgebra $\mathfrak{C}$ of $\mathfrak{A}$ as a binary relation on $\mathfrak{A}$, the subalgebra $\mathfrak{C}$ is recovered as

(A.1.E) $$\mathfrak{C} = \operatorname{Supp}\gamma := \{\, \mathfrak{a} \in \mathfrak{A} \;\mid\; \mathfrak{a}\ \gamma\ \mathfrak{a} \,\}.$$

For this reason, subcongruences on $\mathfrak{A}$ and sections of $\mathfrak{A}$ are one and the same. However, for the reminder of this work, we pursue the subcongruence point of view.

The *subcongruence lattice* of $\mathfrak{A}$, $\mathcal{SC}\mathrm{on}\,\mathfrak{A}$, is the set of the subcongruences of $\mathfrak{A}$ endowed with the order given by the inclusion among subsets of $\mathfrak{A}^2$. It is easy to see that if $\vartheta_i$ is a congruence of a subalgebra $\mathfrak{B}_i$ for $i = 1, 2$, then $\vartheta_1 \cap \vartheta_2$ is a congruence of $\mathfrak{B}_1 \cap \mathfrak{B}_2$. While the join of $\vartheta_1$, $\vartheta_2$ will be the intersection of all the congruences $\vartheta$ of $\mathfrak{B}_1 \vee \mathfrak{B}_2$ such that $\vartheta_1 \leqslant \vartheta$ and $\vartheta_2 \leqslant \vartheta$. In particular $\mathcal{SC}\mathrm{on}\,\mathfrak{A}$ is is a meet-sublattice of $\mathcal{R}\mathrm{el}\,\mathfrak{A}$.

Let $\gamma, \delta \in \mathcal{SC}\mathrm{on}\,\mathfrak{A}$, put $\mathfrak{C} = \operatorname{Supp}\gamma$ and $\mathfrak{D} = \operatorname{Supp}\delta$. If we regard $\gamma$, $\delta$ as binary relations on $\mathfrak{A}$, their product $\gamma\,\delta$ is a binary relation on $\mathfrak{A}$ which is really a subcartesian product of $\mathfrak{D}\,\gamma \times \mathfrak{C}\,\delta$. We define $\gamma * \delta$ as the least subcongruence of $\mathfrak{A}$ which contains $\gamma\,\delta$. It is in fact the least congruence of $\mathfrak{D}\,\gamma \vee \mathfrak{C}\,\delta$ which contains $\gamma\,\delta$ and $\delta\,\gamma$, hence:

$$\gamma * \delta = \delta * \gamma\,.$$

Furthermore, one checks that

$$\text{(A.1.F)} \qquad \text{if } \gamma \leqslant \delta, \quad \text{then } \operatorname{Supp}(\gamma * \delta) = \mathfrak{C}\,\delta \text{ and } \gamma * \delta = \delta^{\mathfrak{C}\,\delta}\,.$$

We say that $\delta$ is $\gamma$-invariant when $\gamma * \delta = \delta$. When $\gamma \leqslant \delta$, this happens if and only if $\operatorname{Supp}\delta = \mathfrak{C}\,\delta$. If $\delta_1$, $\delta_2$ are $\gamma$-invariant, the interval of the $\gamma$-invariant subcongruences from $\delta_1$ to $\delta_2$ is

$$\text{(A.1.G)} \qquad [\delta_1 \div \delta_2]_\gamma := \left\{\ \delta \ \middle|\ \delta_1 \leqslant \delta \leqslant \delta_2 \text{ and } \gamma * \delta = \delta\ \right\}.$$

The $\gamma$-invariant intervals for which $\gamma = \delta_1$ will occur frequently in §A.3. Therefore we give them a name:

A.1.3. *Definition*. By *bottom invariant* intervals of $\mathcal{SC}\mathrm{on}\,\mathfrak{A}$, we mean the $[\gamma \div \varepsilon]_\gamma$ where $\gamma \leqslant \varepsilon$ and $\varepsilon$ is $\gamma$-invariant.

If $\mathsf{S}$ is a set of $\gamma$-invariant subcongruences containing $\gamma$ and $\mathfrak{C} = \operatorname{Supp}\gamma$, then

$$\operatorname{Supp}(\bigvee \mathsf{S}) = \bigvee\{\operatorname{Supp}\delta \mid \delta \in \mathsf{S}\} = \bigvee\{\mathfrak{C}\,\delta \mid \delta \in \mathsf{S}\} \leqslant \mathfrak{C}(\bigvee \mathsf{S}).$$

Beside, $\mathfrak{C}(\bigvee \mathsf{S}) \leqslant \operatorname{Supp}(\bigvee \mathsf{S})$, therefore $\bigvee \mathsf{S}$ is $\gamma$-invariant. Then, for $\gamma \leqslant \varepsilon$ we define $\varepsilon_\gamma$ to be the largest subcongruence which is $\gamma$-invariant and contained in $\varepsilon$:

$$\text{(A.1.H)} \qquad \varepsilon_\gamma := \bigvee\{\delta \leqslant \varepsilon \mid \gamma * \delta = \delta\}.$$

Now to each interval $\mathcal{I} = [\gamma \div \varepsilon]$ of $\mathcal{SC}\mathrm{on}\,\mathfrak{A}$ we can associate a poset homomorphism $\delta \mapsto \delta_\gamma$ from $\mathcal{I}$ onto a bottom invariant interval

$$\text{(A.1.I)} \qquad \mathcal{I}_\bullet := [\gamma \div \varepsilon_\gamma]_\gamma.$$

We observe that $\mathcal{I}_\bullet$ is a complete bounded join sublattice of $\mathcal{I}$, hence a lattice.

**Canonical homomorphisms.** Let $\mathfrak{A}$ be an algebra, $\vartheta \in \mathcal{C}\mathrm{on}\,\mathfrak{A}$. The canonical projection of $\mathfrak{A}$ onto $\mathfrak{A}/\,\vartheta$ is the homomorphism $\mathfrak{a} \mapsto [\mathfrak{a}]_\vartheta$. Other canonical homomorphisms arise as follows.

Homomorphism Theorem. *Let* $\boldsymbol{f} : \mathfrak{A} \to \mathfrak{B}$ *be an algebra homomorphism,* $\kappa = \kappa_{\boldsymbol{f}}$ *and* $\boldsymbol{f}^\flat : \mathfrak{A} \to \mathfrak{A}/\,\kappa$ *the canonical projection; then there is one and only one homomorphism* $\boldsymbol{f}^\sharp : \mathfrak{A}/\,\kappa \to \mathfrak{A}\boldsymbol{f}$ *such that* $\boldsymbol{f} = \boldsymbol{f}^\flat \boldsymbol{f}^\sharp$*. In fact,* $\boldsymbol{f}^\sharp$ *is an isomorphism.*

As the image by $\boldsymbol{f}$ of subalgebras are subalgebras, $\boldsymbol{f}$ induces a lattice isomorphism of $\mathcal{S}\mathrm{ub}_{\kappa_{\boldsymbol{f}}}\mathfrak{A}$ onto $\mathcal{S}\mathrm{ub}(\mathfrak{A}\boldsymbol{f})$. In particular, when $\mathfrak{A} = \mathfrak{F}/\varphi$ and $\mathfrak{C}$ is a $\varphi$-closed subalgebra of $\mathfrak{F}$ then $\mathfrak{C}$ corresponds to a unique subalgebra $\bar{\mathfrak{C}}$ of $\mathfrak{A}$ whose image by $\boldsymbol{f}$ is a subalgebra $\mathfrak{D}$ of $\mathfrak{A}\boldsymbol{f}$. Then we define $\boldsymbol{f}$ *below* $\mathfrak{C}$ as

$$\overline{\boldsymbol{f}}^{\mathfrak{C}} : \frac{\mathfrak{C}}{\varphi^{\mathfrak{C}}} \longrightarrow \mathfrak{D}, \qquad [\mathfrak{c}]_{(\varphi^{\mathfrak{C}})} \mapsto [\mathfrak{c}]_{\varphi}\boldsymbol{f}. \tag{A.1.J}$$

For each congruence $\vartheta$ containing $\kappa$ there is a congruence of $\mathfrak{A}\boldsymbol{f}$ defined by

$$\vartheta(\boldsymbol{f},\boldsymbol{f}) = \left\{ (\mathfrak{a}_1\boldsymbol{f}, \mathfrak{a}_2\boldsymbol{f}) \;\middle|\; (\mathfrak{a}_1, \mathfrak{a}_2) \in \vartheta \right\}.$$

Clearly, $(\boldsymbol{f},\boldsymbol{f})$ induces a lattice isomorphism of the interval $\left[\kappa \div \mathfrak{A}^2\right]$ of $\mathcal{C}\mathrm{on}\,\mathfrak{A}$ onto $\mathcal{C}\mathrm{on}(\mathfrak{A}\boldsymbol{f})$. We define $\boldsymbol{f}$ *above* $\vartheta$ as

$$\underline{\boldsymbol{f}}_{\vartheta} : \frac{\mathfrak{A}}{\vartheta} \longrightarrow \frac{\mathfrak{A}\boldsymbol{f}}{\vartheta(\boldsymbol{f},\boldsymbol{f})}, \qquad [\mathfrak{a}]_{\vartheta} \mapsto [\mathfrak{a}\boldsymbol{f}]_{\vartheta(\boldsymbol{f},\boldsymbol{f})}\,. \tag{A.1.K}$$

As a first consequence of the Homomorphism Theorem we have:

ISOMORPHISM THEOREM. *Let* $\mathfrak{C} \leqslant \mathfrak{A}$, $\vartheta \in \mathcal{C}\mathrm{on}\,\mathfrak{A}$, *then there is an isomorphism*

$$\boldsymbol{\chi}_{\mathfrak{C},\vartheta} : \frac{\mathfrak{C}}{\vartheta^{\mathfrak{C}}} \longrightarrow \frac{\mathfrak{C}\,\vartheta}{\vartheta^{\mathfrak{C}\,\vartheta}}, \qquad [\mathfrak{c}]_{(\vartheta^{\mathfrak{C}})} \mapsto [\mathfrak{c}]_{(\vartheta^{\mathfrak{C}\,\vartheta})}\,.$$

In the Isomorphism Theorem, if we call $\mathfrak{F} = \mathfrak{C}\,\vartheta$, $\varphi = \vartheta^{\mathfrak{F}}$ and $\gamma = \vartheta^{\mathfrak{C}}$, then we get $\mathfrak{C}/\gamma \cong \mathfrak{F}/\varphi$. But we also have $[\mathfrak{C} \div \mathfrak{F}] \cong [\gamma \div \varphi]_{\gamma}$, as we now show.

A.1.4. LEMMA. *Let* $\mathfrak{C} \leqslant \mathfrak{F}$, $\gamma \in \mathcal{C}\mathrm{on}\,\mathfrak{C}$, $\varphi \in \mathcal{C}\mathrm{on}\,\mathfrak{F}$ *such that* $\mathfrak{C}\,\varphi = \mathfrak{F}$ *and* $\gamma = \varphi^{\mathfrak{C}}$, *then* $[\mathfrak{C} \div \mathfrak{F}] \cong [\gamma \div \varphi]_{\gamma}$.

PROOF. We consider the map $\boldsymbol{\zeta} : \mathfrak{M} \mapsto \varphi^{\mathfrak{M}}$ from $[\mathfrak{C} \div \mathfrak{F}]$ to $[\gamma \div \varphi]_{\gamma}$. It is a lattice homomorphism because $\varphi^{\mathfrak{M}} = \varphi \wedge \mathfrak{M}^2$. From $\mathfrak{C}\,\varphi = \mathfrak{F}$ we get $\mathfrak{M} = \mathfrak{C}\,\varphi^{\mathfrak{M}}$ therefore $\boldsymbol{\zeta}$ is injective. But if $\mu \in [\gamma \div \varphi]_{\gamma}$ and $\mathfrak{M} = \mathrm{Supp}\,\mu$, we claim that $\mu = \varphi^{\mathfrak{M}}$. By definition we have $\mu \leqslant \varphi^{\mathfrak{M}}$. And if $\mathfrak{m}_1 \;\varphi\; \mathfrak{m}_2$ with $\mathfrak{m}_i \in \mathfrak{M}$, then let $\mathfrak{c}_i \in \mathfrak{C}$ such that $\mathfrak{c}_i \;\mu\; \mathfrak{m}_i$. As $\mu \leqslant \varphi$, we get, by the transitivity of $\varphi$, $\mathfrak{c}_1 \;\varphi\; \mathfrak{c}_2$. Since $\gamma = \varphi^{\mathfrak{C}}$, we have $\mathfrak{c}_1 \;\gamma\; \mathfrak{c}_2$ and hence $\mathfrak{c}_1 \;\mu\; \mathfrak{c}_2$. Therefore, by transitivity of $\mu$, $\mathfrak{m}_1 \;\mu\; \mathfrak{m}_2$ as we claimed. □

A direct product of two algebras, $\mathfrak{L} \times \mathfrak{R}$, is given together with two canonical homomorphisms called *projections*:

$$\begin{aligned} \boldsymbol{\lambda} : \mathfrak{L} \times \mathfrak{R} &\longrightarrow \mathfrak{L}, & \qquad \boldsymbol{\rho} : \mathfrak{L} \times \mathfrak{R} &\longrightarrow \mathfrak{R} \\ (\mathfrak{l}, \mathfrak{r}) &\mapsto \mathfrak{l} & (\mathfrak{l}, \mathfrak{r}) &\mapsto \mathfrak{r}. \end{aligned} \tag{A.1.L}$$

These induce join-homomorphisms from $\mathcal{S}\mathrm{ub}(\mathfrak{L} \times \mathfrak{R})$ on to $\mathcal{S}\mathrm{ub}\,\mathfrak{L}$ and on to $\mathcal{S}\mathrm{ub}\,\mathfrak{R}$ which we call again $\boldsymbol{\lambda}$ and $\boldsymbol{\rho}$. If $\mathfrak{H} \leqslant \mathfrak{L} \times \mathfrak{R}$, then $\kappa_{\boldsymbol{\lambda}}^{\mathfrak{H}} \vee \kappa_{\boldsymbol{\rho}}^{\mathfrak{H}}$ is a congruence of $\mathfrak{H}$ which corresponds by $\boldsymbol{\lambda}$ to the congruence

$$\mathfrak{H}\tilde{\boldsymbol{\lambda}} := \left(\kappa_{\boldsymbol{\lambda}}^{\mathfrak{H}} \vee \kappa_{\boldsymbol{\rho}}^{\mathfrak{H}}\right)(\boldsymbol{\lambda}, \boldsymbol{\lambda}) \tag{A.1.M}$$

of $\mathfrak{H}\boldsymbol{\lambda}$. We have

A.1.5. Lemma. *If $\kappa_{\boldsymbol{\lambda}}^{\mathfrak{H}}$, $\kappa_{\boldsymbol{\rho}}^{\mathfrak{H}}$ permute then the following are equivalent:*

(a) $\mathfrak{l}_1 \ \left(\mathfrak{H}\tilde{\boldsymbol{\lambda}}\right) \ \mathfrak{l}_2$.
(b) *There is* $\mathfrak{r} \in \mathfrak{R}$ *such that* $(\mathfrak{l}_1, \mathfrak{r}) \in \mathfrak{H}$, $(\mathfrak{l}_2, \mathfrak{r}) \in \mathfrak{H}$.
(c) $\mathfrak{l}_1, \mathfrak{l}_2 \in \mathfrak{H}\boldsymbol{\lambda}$ *and for all* $\mathfrak{r} \in \mathfrak{R}$ $\big((\mathfrak{l}_1, \mathfrak{r}) \in \mathfrak{H} \ \Leftrightarrow \ (\mathfrak{l}_2, \mathfrak{r}) \in \mathfrak{H}\big)$.

Proof. If $x_1 \ \left(\mathfrak{H}\tilde{\boldsymbol{\lambda}}\right) \ x_2$ then there are $y_1, y_2 \in \mathfrak{R}$ such that

$$(x_1, y_1) \ \kappa_{\boldsymbol{\lambda}}^{\mathfrak{H}} \kappa_{\boldsymbol{\rho}}^{\mathfrak{H}} \, (x_2, y_2).$$

Hence $(x_2, y_2) \in \mathfrak{H}$ and there is $(x, y) \in \mathfrak{H}$ such that

$$(x_1, y_1) \ \kappa_{\boldsymbol{\lambda}} \ (x, y) \ \kappa_{\boldsymbol{\rho}} \ (x_2, y_2).$$

Thus $x_1 = x$, $y = y_2$ and $(x_1, y_2) \in \mathfrak{H}$. Therefore (b) holds. Now, if $(x_1, y) \in \mathfrak{H}$ then

$$(x_1, y) \ \kappa_{\boldsymbol{\lambda}}^{\mathfrak{H}} \ (x_1, y_2) \ \kappa_{\boldsymbol{\rho}}^{\mathfrak{H}} \ (x_2, y_2)$$

therefore there is $(x_3, y_3) \in \mathfrak{H}$ such that

$$(x_1, y) \ \kappa_{\boldsymbol{\rho}} \ (x_3, y_3) \ \kappa_{\boldsymbol{\lambda}} \ (x_2, y_2).$$

This implies $y = y_3$, $x_3 = x_2$ and $(x_2, y) \in \mathfrak{H}$. Similarly, one proves that $(x_2, y) \in \mathfrak{H}$ implies $(x_1, y) \in \mathfrak{H}$ so that (c) holds. And of course (c) implies (a). □

A.1.6. Assumption. **For the reminder of this work we shall only consider classes of algebras closed for the formation of subalgebras, homomorphic images and direct products and such that for each algebra $\mathfrak{A}$ of such classes**

$$\vartheta, \varepsilon \in \mathcal{C}\mathrm{on}\,\mathfrak{A} \quad \Longrightarrow \quad \vartheta\,\varepsilon = \varepsilon\,\vartheta, \tag{A.1.N}$$

**and**

$$\vartheta \in \mathcal{C}\mathrm{on}\,\mathfrak{A}, \ \mathfrak{A} \geqslant \mathfrak{C} \geqslant \mathfrak{B} = \mathfrak{B}\,\vartheta \quad \Longrightarrow \quad \mathfrak{C} = \mathfrak{C}\,\vartheta\,. \tag{A.1.O}$$

For example the varieties of groups, groups with operators, hence modules, rings (hence boolean algebras) or algebras over fields enjoy (A.1.N) and (A.1.O). However, we will not make use of (A.1.O) before §A.3.

A.1.7. COROLLARY. *For each interval $\mathcal{I}$ of $\mathcal{S}\mathrm{ub}(\mathfrak{L}\times\mathfrak{R})$, there are four (poset) homomorphisms:*

$$\boldsymbol{\lambda}:\mathcal{I}\to\mathcal{S}\mathrm{ub}\,\mathfrak{L},\qquad \tilde{\boldsymbol{\lambda}}_{\bullet}:\mathcal{I}\to\mathcal{S}\mathcal{C}\mathrm{on}\,\mathfrak{L},$$

$$\boldsymbol{\rho}:\mathcal{I}\to\mathcal{S}\mathrm{ub}\,\mathfrak{R},\qquad \tilde{\boldsymbol{\rho}}_{\bullet}:\mathcal{I}\to\mathcal{S}\mathcal{C}\mathrm{on}\,\mathfrak{R}.$$

PROOF. $\boldsymbol{\lambda}$, $\boldsymbol{\rho}$ are just the restrictions of the canonical projections, thus they are join-homomorphisms. Given $\mathfrak{H}_1,\mathfrak{H}_2\in\mathcal{I}$, apply (a)⇒(c) of A.1.5 (once with $\mathfrak{H}=\mathfrak{H}_1$, once with $\mathfrak{H}=\mathfrak{H}_2$), to see that

$$(\mathfrak{H}_1\wedge\mathfrak{H}_2)\tilde{\boldsymbol{\lambda}}=\big(\mathfrak{H}_1\tilde{\boldsymbol{\lambda}}\big)\wedge\big(\mathfrak{H}_2\tilde{\boldsymbol{\lambda}}\big)\wedge\big((\mathfrak{H}_1\wedge\mathfrak{H}_2)\boldsymbol{\lambda}\big)^2. \tag{A.1.P}$$

Therefore provided $\gamma=(\min\mathcal{I})\tilde{\boldsymbol{\lambda}}$ and because of the comment following (A.1.H), the map

$$\tilde{\boldsymbol{\lambda}}_{\bullet}:\mathcal{I}\to\mathcal{S}\mathcal{C}\mathrm{on}\,\mathfrak{L};\qquad \mathfrak{H}\mapsto\big(\mathfrak{H}\tilde{\boldsymbol{\lambda}}\big)_{\gamma} \tag{A.1.Q}$$

is a poset homomorphism. Analogously, there are poset homomorphisms $\tilde{\boldsymbol{\rho}}$ and $\tilde{\boldsymbol{\rho}}_{\bullet}$ on the right. □

A.1.8. REMARK. *$\tilde{\boldsymbol{\lambda}}_{\bullet}$, $\tilde{\boldsymbol{\rho}}_{\bullet}$ are meet-homomorphisms from $\mathcal{I}$ to $\big(\mathcal{I}\tilde{\boldsymbol{\lambda}}\big)_{\bullet}$ and from $\mathcal{I}$ to $\big(\mathcal{I}\tilde{\boldsymbol{\rho}}\big)_{\bullet}$ respectively.*

PROOF. Put $\gamma=(\min\mathcal{I})\tilde{\boldsymbol{\lambda}}$; it is enough to show that

$$(\mathfrak{H}_1\wedge\mathfrak{H}_2)\tilde{\boldsymbol{\lambda}}_{\bullet}=\Big(\big(\mathfrak{H}_1\tilde{\boldsymbol{\lambda}}\big)\wedge\big(\mathfrak{H}_2\tilde{\boldsymbol{\lambda}}\big)\Big)_{\gamma} \tag{A.1.R}$$

and hence, in view of (A.1.P), that $\mathfrak{C}\vartheta\leqslant(\mathfrak{H}_1\wedge\mathfrak{H}_2)\boldsymbol{\lambda}$ where $\mathfrak{C}=(\min\mathcal{I})\boldsymbol{\lambda}$ and $\vartheta=\big(\mathfrak{H}_1\tilde{\boldsymbol{\lambda}}\big)\wedge\big(\mathfrak{H}_2\tilde{\boldsymbol{\lambda}}\big)$. In fact, if $\big((\mathfrak{h}\boldsymbol{\lambda})\ \vartheta\ \mathfrak{l}\big)$ with $\mathfrak{h}\in\mathfrak{H}_1\wedge\mathfrak{H}_2$, then by (c) of A.1.5 we have $(\mathfrak{l},\mathfrak{h}\boldsymbol{\rho})\in\mathfrak{H}_1\wedge\mathfrak{H}_2$ which in turn implies $\mathfrak{l}\in(\mathfrak{H}_1\wedge\mathfrak{H}_2)\boldsymbol{\lambda}$. □

## A.2. The subalgebras of the direct product of two algebras

The subgroups of the direct product of two groups are well known and the first detailed and satisfactory description, although hard to read nowadays, seems to be the one given by Goursat in 1889 [**Gou89**, pp. 43–48]. Goursat description was rediscovered by Remak [**Rem30**] and then again, in a more general setting, by Fuchs [**Fuc52**] who showed that this description is meaningful even for classes of algebras like groups with operators, rings, modules, algebras over fields and boolean algebras. The content of his work is essentially summarised by saying that any subalgebra of a direct product of two algebras is a pull-back. Two years later, Fleischer [**Fle55**] remarked that Fuchs theorems hold in any class of algebras (closed for direct products, subalgebras and homomorphic images) such that for each algebra its congruences are permutable. This is the widest possible generalisation of Goursat's work. In fact, if there are two non permutable congruences $\vartheta_1$, $\vartheta_2$ of a given algebra $\mathfrak{A}$, then there are subalgebras of $(\mathfrak{A}/\vartheta_1) \times (\mathfrak{A}/\vartheta_2)$ which cannot be described in terms of pull-backs. In this section we present Goursat's description at the light of Fleischer's comments without avoiding to emphasize a more group theoretical point of view at the end.

Among the subalgebras of a direct product $\mathfrak{L} \times \mathfrak{R}$ there are the *cartesian* subalgebras. They are the subalgebras of the form $\mathfrak{C} \times \mathfrak{D}$ for subalgebras $\mathfrak{C}$ of $\mathfrak{L}$ and $\mathfrak{D}$ of $\mathfrak{R}$. If $|\mathfrak{L}| > 1$ then an example of a non cartesian subalgebra is the diagonal $\mathbf{\Delta}_{\mathfrak{L}}$ of $\mathfrak{L}^2$. If $\boldsymbol{\alpha} : \mathfrak{L} \to \mathfrak{R}$ is an isomorphism then we generalize this example by defining the diagonal subalgebra induced by $\boldsymbol{\alpha}$ as

$$\left\{\ (\mathfrak{l}, \mathfrak{r}) \in \mathfrak{L} \times \mathfrak{R} \ \middle|\ \mathfrak{l}\boldsymbol{\alpha} = \mathfrak{r}\ \right\}.$$

This is a subalgebra even if $\boldsymbol{\alpha}$ is a homomorphism defined on a subalgebra $\mathfrak{C}$ of $\mathfrak{L}$. As the roles of $\mathfrak{L}$, $\mathfrak{R}$ could be interchanged, we say that $\mathfrak{H} \leqslant \mathfrak{L} \times \mathfrak{R}$ is a *Goursat subalgebra* when there are $\mathfrak{C} \leqslant \mathfrak{L}$, $\mathfrak{D} \leqslant \mathfrak{R}$ and homomorphisms $\boldsymbol{\gamma} : \mathfrak{C} \to \mathfrak{A}$, $\boldsymbol{\delta} : \mathfrak{D} \to \mathfrak{A}$ such that the following diagram commutes:

$$\begin{array}{ccc} \mathfrak{H} & \xrightarrow{\rho} & \mathfrak{D} \\ {\scriptstyle\lambda}\downarrow & & \downarrow{\scriptstyle\delta} \\ \mathfrak{C} & \xrightarrow{\gamma} & \mathfrak{A} \end{array}$$

For example, if $\gamma \in \mathcal{C}\mathrm{on}\,\mathfrak{C}$, $\delta \in \mathcal{C}\mathrm{on}\,\mathfrak{D}$ and $\boldsymbol{\alpha} : \mathfrak{C}/\gamma \to \mathfrak{D}/\delta$ is an isomorphism, then

$$\text{(A.2.A)} \qquad \mathbf{G}_{\boldsymbol{\alpha}} := \{(\mathfrak{c}, \mathfrak{d}) \in \mathfrak{C} \times \mathfrak{D} \mid [\mathfrak{c}]_\gamma \boldsymbol{\alpha} = [\mathfrak{d}]_\delta\}$$

is a Goursat subalgebra. It is straightforward to check that all the Goursat subalgebras arise in this way, in fact, if the previous diagram commutes, then $\mathfrak{H} = \mathbf{G}_{\gamma^\sharp(\delta^\sharp)^{-1}}$.

A.2.1. LEMMA. *For an $\boldsymbol{\alpha}$ defined as above, $\mathbf{G}_{\boldsymbol{\alpha}}\boldsymbol{\lambda} = \mathfrak{C}$ and $\mathbf{G}_{\boldsymbol{\alpha}}\tilde{\boldsymbol{\lambda}} = \gamma$.*

PROOF. The first equality is obvious, for the second one observe that by (a)⇒(b) of A.1.5 we have $\mathfrak{c}_1 \ \big(\mathbf{G}_{\boldsymbol{\alpha}}\tilde{\boldsymbol{\lambda}}\big)\ \mathfrak{c}_2$ if and only if there is $\mathfrak{d}$ such that

$$[\mathfrak{c}_1]_\gamma \boldsymbol{\alpha} = [\mathfrak{d}]_\delta = [\mathfrak{c}_2]_\gamma \boldsymbol{\alpha},$$

which happens if and only if $\mathfrak{c}_1 \ \gamma \ \mathfrak{c}_2$. □

Goursat subalgebras are actually all the subalgebras of $\mathfrak{L} \times \mathfrak{R}$.

FLEISCHER THEOREM. *Let $\mathfrak{H}$ be a subalgebra of $\mathfrak{L} \times \mathfrak{R}$, put*

$$\mathfrak{C} = \mathfrak{H}\boldsymbol{\lambda}, \qquad \gamma = \mathfrak{H}\tilde{\boldsymbol{\lambda}},$$
$$\mathfrak{D} = \mathfrak{H}\boldsymbol{\rho}, \qquad \delta = \mathfrak{H}\tilde{\boldsymbol{\rho}}.$$

*Then* $\gamma \in \mathcal{C}\mathrm{on}\,\mathfrak{C}$, $\delta \in \mathcal{C}\mathrm{on}\,\mathfrak{D}$ *and the projections* $\boldsymbol{\lambda}$, $\boldsymbol{\rho}$ *restricted to* $\mathfrak{H}$ *induce algebra isomorphisms*

$$\boldsymbol{\gamma} : \frac{\mathfrak{H}}{\kappa^{\mathfrak{H}}_{\boldsymbol{\lambda}}\,\kappa^{\mathfrak{H}}_{\boldsymbol{\rho}}} \longrightarrow \frac{\mathfrak{C}}{\gamma}, \qquad \boldsymbol{\delta} : \frac{\mathfrak{H}}{\kappa^{\mathfrak{H}}_{\boldsymbol{\lambda}}\,\kappa^{\mathfrak{H}}_{\boldsymbol{\rho}}} \longrightarrow \frac{\mathfrak{D}}{\delta}$$

*such that* $\mathfrak{H} = \mathbf{G}_{\boldsymbol{\gamma}^{-1}\boldsymbol{\delta}}$ *is the Goursat subalgebra induced by* $\boldsymbol{\gamma}^{-1}\boldsymbol{\delta}$.

Proof. We only need to prove that $\mathbf{G}_{\boldsymbol{\gamma}^{-1}\boldsymbol{\delta}} \leqslant \mathfrak{H}$ because it is readily seen that $\mathfrak{H} \leqslant \mathbf{G}_{\boldsymbol{\gamma}^{-1}\boldsymbol{\delta}}$. To this end, let us say $(\mathfrak{c}, \mathfrak{d}) \in \mathbf{G}_{\boldsymbol{\gamma}^{-1}\boldsymbol{\delta}}$ and $[\mathfrak{c}]_\gamma \boldsymbol{\gamma}^{-1} = [(\mathfrak{c}, \mathfrak{d}_1)]_{\kappa^{\mathfrak{H}}_{\boldsymbol{\lambda}}\kappa^{\mathfrak{H}}_{\boldsymbol{\rho}}}$. Then $[\mathfrak{d}]_\delta = [\mathfrak{c}]_\gamma \boldsymbol{\gamma}^{-1}\boldsymbol{\delta} = [\mathfrak{d}_1]_\delta$. By A.1.5 there is $\mathfrak{c}_1$ such that $(\mathfrak{c}_1, \mathfrak{d}) \in \mathfrak{H}$ and $(\mathfrak{c}_1, \mathfrak{d}_1) \in \mathfrak{H}$. Thus

$$(\mathfrak{c}, \mathfrak{d}_1)\ \kappa^{\mathfrak{H}}_{\boldsymbol{\rho}}\ (\mathfrak{c}_1, \mathfrak{d}_1)\ \kappa^{\mathfrak{H}}_{\boldsymbol{\lambda}}\ (\mathfrak{c}_1, \mathfrak{d}).$$

By permutability there is $(\mathfrak{c}_2, \mathfrak{d}_2) \in \mathfrak{H}$ such that

$$(\mathfrak{c}, \mathfrak{d}_1)\ \kappa^{\mathfrak{H}}_{\boldsymbol{\lambda}}\ (\mathfrak{c}_2, \mathfrak{d}_2)\ \kappa^{\mathfrak{H}}_{\boldsymbol{\rho}}\ (\mathfrak{c}_1, \mathfrak{d}).$$

Therefore $\mathfrak{c} = \mathfrak{c}_2$, $\mathfrak{d}_2 = \mathfrak{d}$, and so $(\mathfrak{c}, \mathfrak{d}) \in \mathfrak{H}$. □

For example the cartesian subalgebra $\mathfrak{C} \times \mathfrak{D}$ corresponds to the isomorphism of trivial sections $\mathfrak{C}/\mathfrak{C}^2 \to \mathfrak{D}/\mathfrak{D}^2$. If we apply Fleischer theorem to $\mathbf{G}_{\boldsymbol{\alpha}}$ then $\boldsymbol{\gamma}^{-1}\boldsymbol{\delta} = \boldsymbol{\alpha}$, in particular we have

A.2.2. Corollary. *The set of subalgebras of* $\mathfrak{L} \times \mathfrak{R}$ *is in bijection with the set of isomorphisms between sections of* $\mathfrak{L}$ *and* $\mathfrak{R}$.

A.2.3. Corollary. $\mathcal{S}\mathrm{ub}(\mathfrak{L} \times \mathfrak{R}) = \mathcal{S}\mathrm{ub}(\mathfrak{L}) \times \mathcal{S}\mathrm{ub}(\mathfrak{R})$ *if and only if the sections of* $\mathfrak{L}$ *isomorphic to sections of* $\mathfrak{R}$ *are trivial.*

Note that if there are non trivial isomorphic sections, then $\mathcal{S}\mathrm{ub}(\mathfrak{L} \times \mathfrak{R})$ contains the congruence lattices of these isomorphic sections as intervals.

A.2.4. Proposition. *Let* $\boldsymbol{\alpha} : \mathfrak{C}/\,\gamma \to \mathfrak{D}/\,\delta$ *be an isomorphism, then the interval* $[\mathbf{G}_{\boldsymbol{\alpha}} \div \mathfrak{C} \times \mathfrak{D}]$ *is canonically lattice isomorphic to* $\mathcal{C}\mathrm{on}(\mathfrak{C}/\,\gamma)$.

Proof. We only need to show that $\tilde{\boldsymbol{\lambda}}_\bullet = \tilde{\boldsymbol{\lambda}} : [\mathbf{G}_{\boldsymbol{\alpha}} \div \mathfrak{C} \times \mathfrak{D}] \to \left[\gamma \div \mathfrak{C}^2\right]$ has an inverse. In fact, for $\vartheta \in \mathcal{C}\mathrm{on}\,\mathfrak{C}$ such that $\gamma \leqslant \vartheta \leqslant \mathfrak{C}^2$, the Goursat algebra induced by $\underline{\boldsymbol{\alpha}}_\vartheta$ has image $\vartheta$ by $\tilde{\boldsymbol{\lambda}}_\bullet$. Beside, if $\mathbf{G}_{\boldsymbol{\alpha}} \leqslant \mathbf{G}_{\boldsymbol{\vartheta}} \leqslant \mathfrak{C} \times \mathfrak{D}$, then $\boldsymbol{\vartheta} : \mathfrak{C}/\,\varepsilon \to \mathfrak{D}/\,\varphi$ for suitable $\varepsilon$, $\varphi$. Next lemma shows that $\varphi = \varepsilon(\boldsymbol{\alpha}, \boldsymbol{\alpha})$ and this implies $\boldsymbol{\vartheta} = \underline{\boldsymbol{\alpha}}_\varepsilon$ because if $(\mathfrak{c}, \mathfrak{d}) \in \mathbf{G}_{\boldsymbol{\alpha}} \leqslant \mathbf{G}_{\boldsymbol{\vartheta}}$, then $[\mathfrak{c}]_\varepsilon \boldsymbol{\vartheta} = [\mathfrak{d}]_\varphi = [\mathfrak{d}]_{\varepsilon(\boldsymbol{\alpha},\boldsymbol{\alpha})} = [\mathfrak{c}]_\varepsilon \underline{\boldsymbol{\alpha}}_\vartheta$. □

A.2.5. Lemma. *Let* $\boldsymbol{\alpha} : \frac{\mathfrak{C}}{\gamma} \to \frac{\mathfrak{D}}{\delta}$ *and* $\boldsymbol{\beta} : \frac{\mathfrak{C}}{\varepsilon} \to \frac{\mathfrak{F}}{\varphi}$ *be isomorphisms. If* $\mathbf{G}_{\boldsymbol{\alpha}} \leqslant \mathbf{G}_{\boldsymbol{\beta}}$ *then* $\varepsilon^{\mathfrak{C}}(\boldsymbol{\alpha}, \boldsymbol{\alpha}) = \varphi^{\mathfrak{D}}$ *and* $(\mathfrak{C}\,\varepsilon)\boldsymbol{\beta} = \mathfrak{D}\,\varphi$.

PROOF. let $\mathfrak{d}_1, \mathfrak{d}_2 \in \mathfrak{D}$ and $\mathfrak{c}_1, \mathfrak{c}_2 \in \mathfrak{C}$ such that $[\mathfrak{c}_i]_\gamma \boldsymbol{\alpha} = [\mathfrak{d}_i]_\delta$. This implies $[\mathfrak{c}_i]_\varepsilon \boldsymbol{\beta} = [\mathfrak{d}_i]_\varphi$ and since $\boldsymbol{\beta}$ is an isomorphism we have

$$\mathfrak{d}_1 \ \varphi \ \mathfrak{d}_2 \quad \Longleftrightarrow \quad \mathfrak{c}_1 \ \varepsilon \ \mathfrak{c}_2,$$

that is exactly $\varphi^{\mathfrak{D}} = \varepsilon^{\mathfrak{C}}(\boldsymbol{\alpha\alpha})$. Now, assume $\mathfrak{e} \in \mathfrak{C}\,\varepsilon$, then there is $\mathfrak{c} \in \mathfrak{C}$ such that $(\mathfrak{c} \ \varepsilon \ \mathfrak{e})$. There is $(\mathfrak{c}, \mathfrak{d}) \in \mathbf{G}_{\boldsymbol{\alpha}} \leqslant \mathbf{G}_{\boldsymbol{\beta}}$. Therefore $[\mathfrak{e}]_\varepsilon \boldsymbol{\beta} = [\mathfrak{c}]_\varepsilon \boldsymbol{\beta} = [\mathfrak{d}]_\varphi$. This proves that $(\mathfrak{C}\,\varepsilon)\boldsymbol{\beta} \leqslant \mathfrak{D}\,\varphi$ and similarly that $(\mathfrak{D}\,\varphi)\boldsymbol{\beta}^{-1} \leqslant \mathfrak{C}\,\varepsilon$, that is $\mathfrak{D}\,\varphi \leqslant (\mathfrak{C}\,\varepsilon)\boldsymbol{\beta}$. □

**Example: the variety of Groups.** In the variety of groups we agree to call "1" the identity element. If $G$ is a group and $\vartheta \in \mathcal{C}\mathrm{on}\, G$, then $[1]_\vartheta$ is a normal subgroup of $G$ and this correspondence induces a lattice isomorphism between $\mathcal{C}\mathrm{on}\, G$ and the lattice of the normal subgroups of $G$, which we write $\mathcal{NS}\mathrm{ub}\, G$. We will always identify $\mathcal{C}\mathrm{on}\, G$ with $\mathcal{NS}\mathrm{ub}\, G$ via this isomorphism. Similarly, we will always identify $\mathcal{SC}\mathrm{on}\, G$ with $\mathrm{Sec}\, G = \{\ (\mathcal{A}, A) \ \mid \ \mathcal{A} \trianglelefteq A \leqslant G\ \}$ via the bijection $\vartheta \mapsto ([1]_\vartheta, \mathrm{Supp}\,\vartheta)$. We observe that if $(\mathcal{A}, A)$ is $(\mathcal{B}, B)$-invariant and $B \leqslant A$, then $\mathcal{A}$ is $B$-invariant (by conjugation). Conversely, if $\mathcal{A}$ is $B$-invariant, then $B\mathcal{A}$ is a group and $(\mathcal{A}, B\mathcal{A})$ is a $(\mathcal{B}, B)$-invariant section for each $\mathcal{B} \leqslant B \cap \mathcal{A}$ which is normal in $B$. Therefore Lemma A.1.4 may be translated as follows.

A.2.6. LEMMA. *Let $A \leqslant B$, $\mathcal{A} \trianglelefteq A$, $\mathcal{B} \trianglelefteq B$ such that $A\mathcal{B} = B$ and $\mathcal{A} = \mathcal{B} \cap A$. Then $[A \div B] \cong [\mathcal{A} \div \mathcal{B}]_A$. where the second is the interval of the $A$-invariant (by conjugation) subgroups between $\mathcal{A}$ and $\mathcal{B}$.*

By identifying $L$ with $L \times \{1\}$ and $R$ with $\{1\} \times R$, we will often regard the direct product $L \times R$ as the internal product $LR$. The first component of

$$\tilde{\boldsymbol{\lambda}} : \mathcal{S}\mathrm{ub}(LR) \longrightarrow \mathrm{Sec}\, L; \qquad\qquad H \mapsto \big([1]_{H\tilde{\boldsymbol{\lambda}}}, H\boldsymbol{\lambda}\big),$$

regarded as a map from $\mathcal{S}\mathrm{ub}(LR)$ to $\mathcal{S}\mathrm{ub}\, L$, is a meet-homomorphism because $[1]_{H\tilde{\boldsymbol{\lambda}}} = L \wedge H$. Therefore, rather than $\tilde{\boldsymbol{\lambda}}_\bullet$, $\tilde{\boldsymbol{\rho}}_\bullet$ which are only defined on intervals, we will use the global

$$\text{(A.2.B)} \qquad \tilde{\boldsymbol{\lambda}}_1 : \mathcal{S}\mathrm{ub}(LR) \longrightarrow \mathcal{S}\mathrm{ub}\, L; \qquad\qquad H \mapsto L \wedge H,$$

$$\text{(A.2.C)} \qquad \tilde{\boldsymbol{\rho}}_1 : \mathcal{S}\mathrm{ub}(LR) \longrightarrow \mathcal{S}\mathrm{ub}\, R; \qquad\qquad H \mapsto H \wedge R.$$

If $\boldsymbol{\alpha} : A/\mathcal{A} \to B/\mathcal{B}$ is an isomorphism of sections of $L$ and $R$, then

$$\big\langle A, \mathbf{G}_{\boldsymbol{\alpha}} \big\rangle = A\mathbf{G}_{\boldsymbol{\alpha}} = AB \qquad \text{and} \qquad \mathbf{G}_{\boldsymbol{\alpha}} \cap L = \mathcal{A}.$$

The isomorphism A.2.4 for groups is well known, (see [**Ros65**] for example). In more details, $\tilde{\boldsymbol{\lambda}}_1$ induces a lattice isomorphism of $[\mathbf{G}_\alpha \div A\mathcal{B}]$ on to $[\mathcal{A} \div A]_A$, the interval of the $A$-invariant subgroups of $A$ containing $\mathcal{A}$; the inverse being induced by $\mathcal{C} \mapsto \mathcal{C}\mathbf{G}_\alpha$ ($\mathcal{C} \trianglelefteq LR$ therefore $\mathcal{C}\mathbf{G}_\alpha \leqslant AB$).

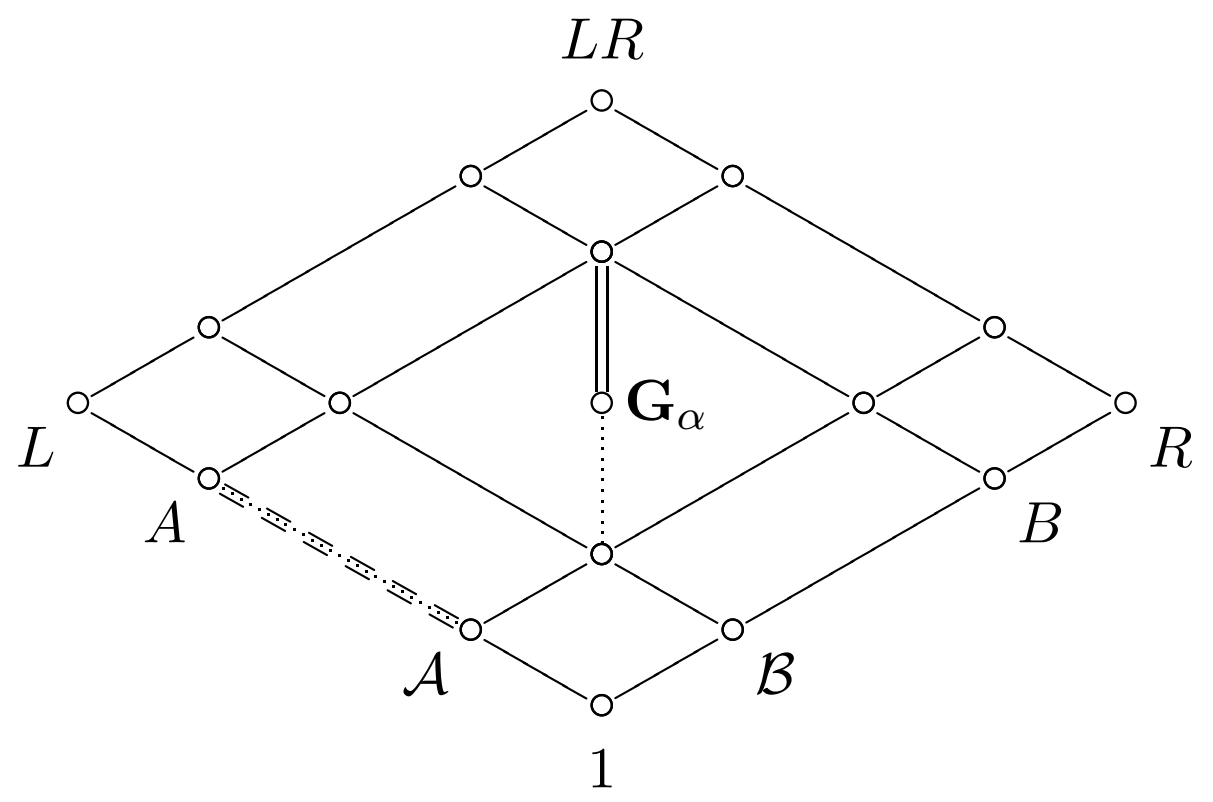

$$[\mathcal{A}\mathcal{B} \div \mathbf{G}_\alpha] \overset{\boldsymbol{\lambda}}{\cong} [\mathcal{A} \div A]$$

$$[\mathbf{G}_\alpha \div AB] \overset{\tilde{\boldsymbol{\lambda}}_1}{\cong} [\mathcal{A} \div A]_A$$

FIGURE A.B. A Goursat subgroup of $L \times R$, (some edges may collapse)

## A.3. Intervals

The aim of this section is to exploit the four canonical homomorphisms of A.1.7 to describe the simple intervals of $\mathcal{S}\text{ub}(\mathfrak{L} \times \mathfrak{R})$. To this end we will first describe a class of intervals whose homomorphisms behave in the best possible way.

*Definition.* An interval $\mathcal{I}$ of $\mathcal{S}\text{ub}(\mathfrak{L} \times \mathfrak{R})$ is said to be *elementary* if at least two of the four canonical homomorphisms of $\mathcal{I}$ are injective.

Note that by A.3.11 below, if three canonical homomorphisms of $\mathcal{I}$ are isomorphisms, then $\mathcal{I}$ is trivial. Also, recall that $\boldsymbol{\lambda}$, $\boldsymbol{\rho}$, being induced by algebra homomorphisms, send intervals to intervals. We shall show that $\tilde{\boldsymbol{\lambda}}_\bullet$, $\tilde{\boldsymbol{\rho}}_\bullet$, send intervals to bottom invariant intervals. As a consequence, an elementary interval $\mathcal{I}$ is always canonically isomorphic to an interval of $\mathcal{S}\text{ub}\,\mathfrak{L}$ or $\mathcal{S}\text{ub}\,\mathfrak{R}$, or to a bottom invariant interval of $\mathcal{S}\mathcal{C}\text{on}\,\mathfrak{L}$ or $\mathcal{S}\mathcal{C}\text{on}\,\mathfrak{R}$.

Examples of elementary intervals are given by describing six subclasses or types which are closed under the composition of intervals whenever it makes sense.

**Types 2L, 2R.** Let $\mathfrak{A} < \mathfrak{C} \leqslant \mathfrak{L}$ and $\mathfrak{B} \leqslant \mathfrak{R}$. The intervals of type **2L** are the

$$\mathcal{I} = [\mathfrak{A} \times \mathfrak{B} \div \mathfrak{C} \times \mathfrak{B}].$$

We notice that while $\mathcal{I}\boldsymbol{\rho}$, $\mathcal{I}\tilde{\boldsymbol{\rho}}_\bullet$ are trivial intervals, $\boldsymbol{\lambda}$ and $\tilde{\boldsymbol{\lambda}}_\bullet$ are isomorphisms onto $[\mathfrak{A} \div \mathfrak{C}]$, $\left[\mathfrak{A}^2 \div \mathfrak{C}^2\right]_\bullet$ respectively (see A.1.4). The intervals of type **2R** are defined in a similar way on the right.

**Type 3A.** These are the intervals of type

$$\mathcal{I} = \left[\mathbf{G}_\varphi \div \mathbf{G}_{\underline{\varphi}_\gamma}\right]$$

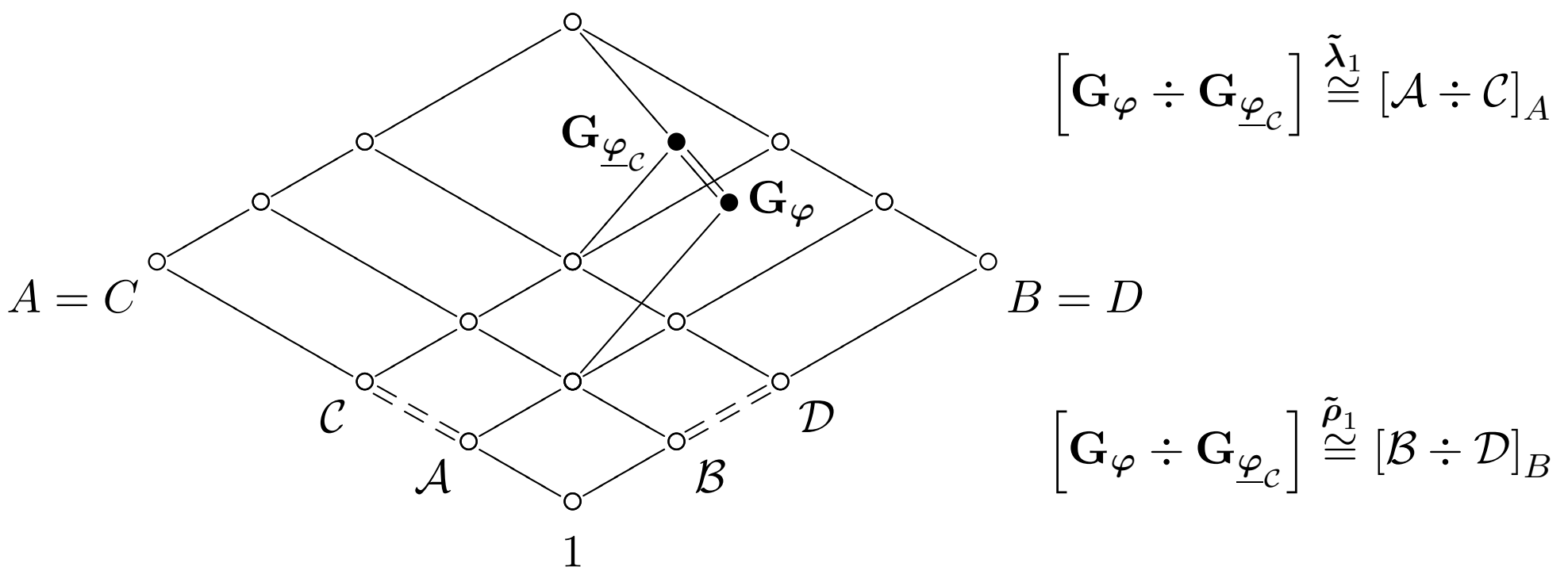

FIGURE A.C. Groups: intervals of type **3A**

where $\varphi : \mathfrak{A}/\alpha \to \mathfrak{B}/\beta$ is an isomorphism, $\alpha < \gamma \in \mathcal{C}\mathrm{on}\,\mathfrak{A}$ and $\underline{\varphi}_{\gamma}$ is defined as in (A.1.K). In this case, while $\mathcal{I}\boldsymbol{\lambda}$, $\mathcal{I}\boldsymbol{\rho}$ are trivial, $\tilde{\boldsymbol{\lambda}}_{\bullet}$ and $\tilde{\boldsymbol{\rho}}_{\bullet}$ are injective. In fact, by A.2.4, they are isomorphisms on to $[\alpha \div \gamma]$ and $[\beta \div \gamma(\varphi, \varphi)]$ respectively.

A.3.1. LEMMA. *Let* $\varphi : \frac{\mathfrak{A}}{\alpha} \to \frac{\mathfrak{B}}{\beta}$, $\psi : \frac{\mathfrak{C}}{\gamma} \to \frac{\mathfrak{D}}{\delta}$ *be isomorphisms. If*

$$\mathfrak{A} = \mathfrak{C} \quad \textit{and} \quad \mathfrak{B} = \mathfrak{B}\delta,$$

*then* $\mathbf{G}_{\varphi} \leqslant \mathbf{G}_{\psi}$ *if and only if* $\psi = \underline{\varphi}_{\gamma}$. *Furthermore, If* $\mathfrak{C}$, $\mathfrak{D}$ *are groups and* $\mathcal{C} = [1]_{\gamma}$, *then they are also equivalent to* $\mathbf{G}_{\psi} = \mathcal{C}\mathbf{G}_{\varphi}$.

PROOF. Assume $\mathbf{G}_{\varphi} \leqslant \mathbf{G}_{\psi}$, then by A.2.5 we have

$$\mathfrak{B} = \mathfrak{B}\delta = (\mathfrak{A}\gamma)\psi = \mathfrak{C}\psi = \mathfrak{D}.$$

Then $\mathbf{G}_{\varphi} \leqslant \mathbf{G}_{\psi} \leqslant \mathfrak{A} \times \mathfrak{B}$ and by A.2.4, $\psi = \underline{\varphi}_{\gamma}$. If also $\mathfrak{A}$, $\mathfrak{B}$ are groups and $\mathcal{C} = [1]_{\gamma}$, then $\mathcal{C}\mathbf{G}_{\varphi} \subseteq \mathbf{G}_{\psi}$. But if $[\mathfrak{a}]_{\gamma}\psi = [\mathfrak{b}]_{\delta}$, then there is $\mathfrak{a}_1 \in \mathfrak{A}$ such that $[\mathfrak{a}_1]_{\alpha}\varphi = [\mathfrak{b}]_{\beta}$, thus $(\mathfrak{a}_1, \mathfrak{b}) \in \mathbf{G}_{\varphi}$ and $(\mathfrak{a}, \mathfrak{b}) = (\mathfrak{a}\mathfrak{a}_1^{-1}, 1)(\mathfrak{a}_1, \mathfrak{b}) \in \mathcal{C}\mathbf{G}_{\varphi}$. □

**Type 3B.** These are the intervals of type

$$\mathcal{I} = \left[\mathbf{G}_{\overline{\psi}^{\mathfrak{A}}} \div \mathbf{G}_{\psi}\right]$$

where $\psi : \mathfrak{C}/\gamma \to \mathfrak{D}/\delta$ is an isomorphism, $\mathfrak{A}\gamma = \mathfrak{A} < \mathfrak{C}$ and $\overline{\psi}^{\mathfrak{A}}$ is defined as in (A.1.J). It is easy to check that as $\mathfrak{A}$ is $\gamma$-closed, $\mathcal{I}\tilde{\boldsymbol{\lambda}}_{\bullet}$ and $\mathcal{I}\tilde{\boldsymbol{\rho}}_{\bullet}$ are trivial. If we call $\mathfrak{H} = \mathbf{G}_{\overline{\psi}^{\mathfrak{A}}}$ and $\vartheta = \kappa_{\boldsymbol{\lambda}}^{\mathbf{G}_{\psi}} \kappa_{\boldsymbol{\rho}}^{\mathbf{G}_{\psi}}$ for short, we show that $\mathfrak{H}$ is $\vartheta$-closed too. To this end let $(\mathfrak{c}, \mathfrak{d}) \in \mathbf{G}_{\psi}$, such that $(\mathfrak{c}, \mathfrak{d})\ \vartheta\ (\mathfrak{a}, \mathfrak{b})$ for some $(\mathfrak{a}, \mathfrak{b}) \in \mathfrak{H}$. Then

$$(\mathfrak{c}, \mathfrak{d})\ \kappa_{\boldsymbol{\lambda}}\ (\mathfrak{c}_1, \mathfrak{d}_1)\ \kappa_{\boldsymbol{\rho}}\ (\mathfrak{a}, \mathfrak{b})$$

for some $(\mathfrak{c}_1, \mathfrak{d}_1) \in \mathbf{G}_{\psi}$. But this forces $\mathfrak{c}_1 = \mathfrak{c}$, $\mathfrak{d}_1 = \mathfrak{b}$ and thus $(\mathfrak{c}, \mathfrak{b}) \in \mathbf{G}_{\psi}$. In particular $[\mathfrak{c}]_{\gamma}\psi = [\mathfrak{b}]_{\delta} = [\mathfrak{a}]_{\gamma}$ and by $\mathfrak{A} = \mathfrak{A}\gamma$ we conclude $\mathfrak{c} \in \mathfrak{A}$ which proves that

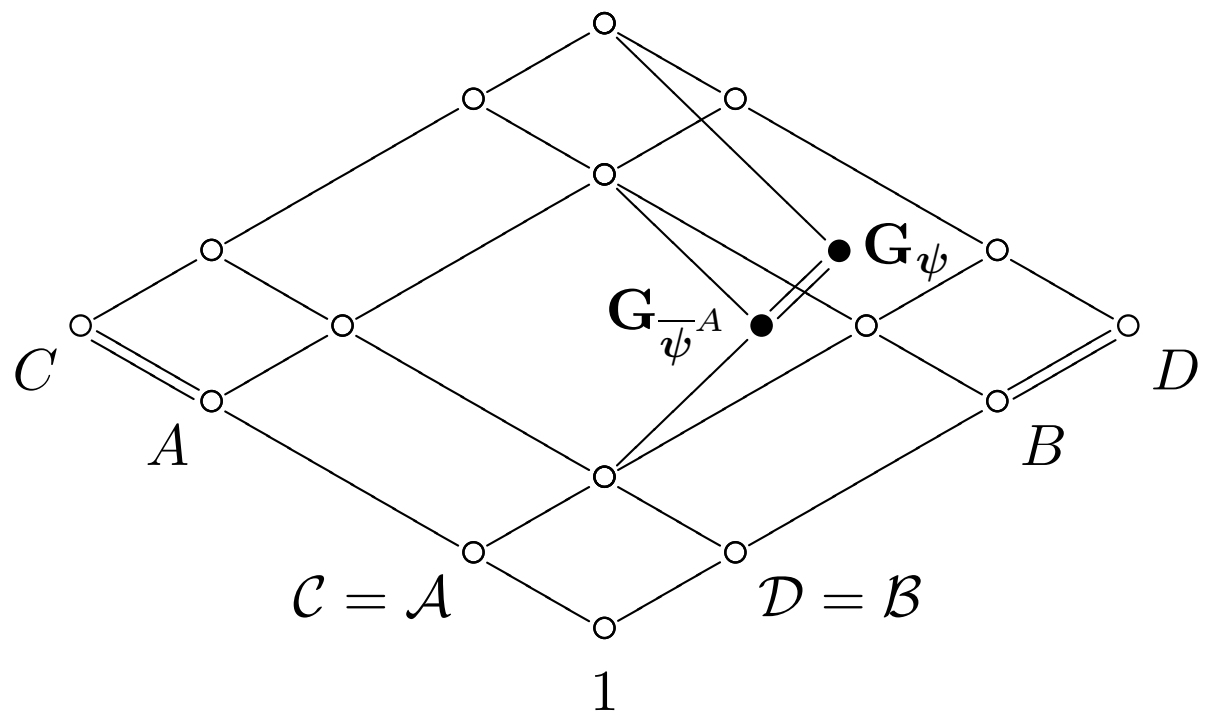

$$\left[\mathbf{G}_{\overline{\boldsymbol{\psi}}^A} \div \mathbf{G}_{\boldsymbol{\psi}}\right] \overset{\boldsymbol{\lambda}}{\cong} [A \div C]$$

$$\left[\mathbf{G}_{\overline{\boldsymbol{\psi}}^A} \div \mathbf{G}_{\boldsymbol{\psi}}\right] \overset{\boldsymbol{\rho}}{\cong} [B \div D]$$

FIGURE A.D. Groups: intervals of type **3B**

$(\mathfrak{c}, \mathfrak{d}) \in \mathfrak{H}$. Now, by definition of $\mathbf{G}_{\boldsymbol{\psi}}\tilde{\boldsymbol{\lambda}} = \gamma$, $\boldsymbol{\lambda}$ induces an algebra isomorphisms from $\mathbf{G}_{\boldsymbol{\psi}}/\,\vartheta$ on to $\mathfrak{C}/\gamma$, therefore from (A.1.O) and the closure of $\mathfrak{H}$, $\mathfrak{A}$, we see that

$$\boldsymbol{\lambda} : \mathcal{I} \longrightarrow [\mathfrak{A} \div \mathfrak{C}] \qquad \text{and} \qquad \boldsymbol{\rho} : \mathcal{I} \longrightarrow [\mathfrak{A}\boldsymbol{\psi} \div \mathfrak{D}]$$

are isomorphisms.

A.3.2. LEMMA. *Let* $\boldsymbol{\varphi} : \frac{\mathfrak{A}}{\alpha} \to \frac{\mathfrak{B}}{\beta}$, $\boldsymbol{\psi} : \frac{\mathfrak{C}}{\gamma} \to \frac{\mathfrak{D}}{\delta}$ *be isomorphisms. If*

$$\mathfrak{A} = \mathfrak{A}\gamma \quad \textit{and} \quad \mathfrak{B} = \mathfrak{B}\delta,$$

*then the following are equivalent:*

(a) $\mathbf{G}_{\boldsymbol{\varphi}} \leqslant \mathbf{G}_{\boldsymbol{\psi}}$.
(b) $\underline{\boldsymbol{\varphi}}_{(\gamma^{\mathfrak{A}})} = \overline{\boldsymbol{\psi}}^{\mathfrak{A}}$.
(c) $\mathbf{G}_{\underline{\boldsymbol{\varphi}}_{(\gamma^{\mathfrak{A}})}} = \mathbf{G}_{\boldsymbol{\psi}} \wedge (\mathfrak{A} \times \mathfrak{D})$.

PROOF. We only prove that (a)⇒(b), the other implications being trivial. If $\mathbf{G}_{\boldsymbol{\varphi}} \leqslant \mathbf{G}_{\boldsymbol{\psi}}$, then by A.2.5 $\gamma^{\mathfrak{A}}(\boldsymbol{\varphi}\boldsymbol{\varphi}) = \delta^{\mathfrak{B}}$. Then b holds because if $[\mathfrak{a}]_\alpha \boldsymbol{\varphi} = [\mathfrak{b}]_\beta$, then

$$[\mathfrak{a}]_{(\gamma^{\mathfrak{A}})} \underline{\boldsymbol{\varphi}}_{(\gamma^{\mathfrak{A}})} = [\mathfrak{b}]_{(\gamma^{\mathfrak{A}})(\boldsymbol{\varphi}\boldsymbol{\varphi})} = [\mathfrak{b}]_{(\delta^{\mathfrak{B}})} = [\mathfrak{b}]_\delta = [\mathfrak{a}]_\gamma \boldsymbol{\psi} = [\mathfrak{a}]_{(\gamma^{\mathfrak{A}})} \overline{\boldsymbol{\psi}}^{\mathfrak{A}}.$$

□

**Types 4L, 4R.** Let $\mathfrak{A} < \mathfrak{C}$, $\alpha = \gamma^{\mathfrak{A}}$ for some $\gamma \in \mathcal{C}\mathrm{on}\,\mathfrak{C}$ such that $\mathfrak{A}\gamma = \mathfrak{C}$. The intervals of type **4L** are the

$$\textbf{(4L)} \qquad \mathcal{I} = \left[\mathbf{G}_{\boldsymbol{\chi}_{\mathfrak{A},\gamma}\boldsymbol{\psi}} \div \mathbf{G}_{\boldsymbol{\psi}}\right]$$

where $\boldsymbol{\psi} : \mathfrak{C}/\gamma \to \mathfrak{D}/\delta$ is an isomorphism and $\boldsymbol{\chi}_{\mathfrak{A},\gamma}$ is given by the Isomorphism Theorem. It is readily seen that $\mathcal{I}\boldsymbol{\rho}$, $\mathcal{I}\tilde{\boldsymbol{\rho}}$ are trivial.

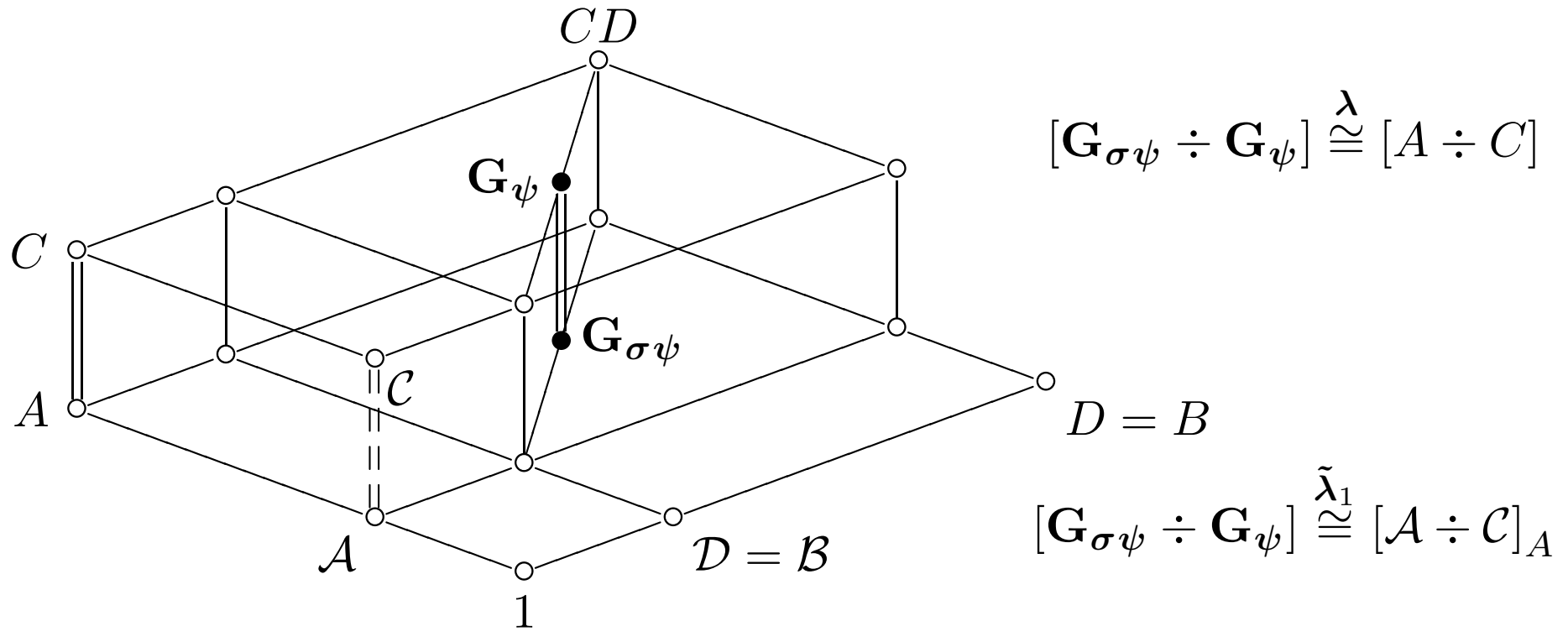

FIGURE A.E. Groups: intervals of type **4L**

A.3.3. LEMMA. *Let* $\boldsymbol{\varphi} : \frac{\mathfrak{A}}{\alpha} \to \frac{\mathfrak{B}}{\beta}$, $\boldsymbol{\psi} : \frac{\mathfrak{C}}{\gamma} \to \frac{\mathfrak{D}}{\delta}$ *be isomorphisms. If* $\mathfrak{B} = \mathfrak{D}$ *then the following are equivalent:*

(a) $\mathbf{G}_{\boldsymbol{\varphi}} \leqslant \mathbf{G}_{\boldsymbol{\psi}}$.

(b) $\underline{\boldsymbol{\varphi}}_{(\gamma^{\mathfrak{A}})} = \boldsymbol{\chi}_{\mathfrak{A},\gamma}\boldsymbol{\psi}$.

(c) $\mathbf{G}_{\underline{\boldsymbol{\varphi}}_{(\gamma^{\mathfrak{A}})}} = \mathbf{G}_{\boldsymbol{\psi}} \wedge (\mathfrak{A} \times \mathfrak{D})$.

*Furthermore, if* $\mathfrak{C}$, $\mathfrak{D}$ *are groups and* $\mathcal{C} = [1]_\gamma$, *then they are also equivalent to*

(d) $\mathbf{G}_{\boldsymbol{\psi}} = \mathcal{C}\mathbf{G}_{\boldsymbol{\varphi}}$.

PROOF. We only prove that (a)⇒(b), the other implications being trivial ((a)⇒(d) as in A.3.1). If $\mathbf{G}_{\boldsymbol{\varphi}} \leqslant \mathbf{G}_{\boldsymbol{\psi}}$, then by A.2.5 $\gamma^{\mathfrak{A}}(\boldsymbol{\varphi}, \boldsymbol{\varphi}) = \delta$ and since $\mathfrak{B} = \mathfrak{D}$, $\mathfrak{A}\gamma = \mathfrak{C}$. Then b holds because if $[\mathfrak{a}]_\alpha \boldsymbol{\varphi} = [\mathfrak{b}]_\beta$, then

$$[\mathfrak{a}]_{(\gamma^{\mathfrak{A}})}\underline{\boldsymbol{\varphi}}_{(\gamma^{\mathfrak{A}})} = [\mathfrak{b}]_{(\gamma^{\mathfrak{A}})(\boldsymbol{\varphi}\boldsymbol{\varphi})} = [\mathfrak{b}]_\delta = [\mathfrak{a}]_\gamma \boldsymbol{\psi} = [\mathfrak{a}]_{(\gamma^{\mathfrak{A}})}\boldsymbol{\chi}_{\mathfrak{A},\gamma}\boldsymbol{\psi}.$$

□

A.3.4. COROLLARY. *Let* $\boldsymbol{\varphi} : \frac{\mathfrak{A}}{\alpha} \to \frac{\mathfrak{B}}{\beta}$, $\boldsymbol{\psi} : \frac{\mathfrak{C}}{\gamma} \to \frac{\mathfrak{B}}{\delta}$ *be isomorphisms such that* $\mathbf{G}_{\boldsymbol{\varphi}} \leqslant \mathbf{G}_{\boldsymbol{\psi}}$. *Then* $\tilde{\boldsymbol{\lambda}} : [\mathbf{G}_{\boldsymbol{\varphi}} \div \mathbf{G}_{\boldsymbol{\psi}}] \to [\alpha \div \gamma]_\alpha$ *is an isomorphism.*

PROOF. if $\mathbf{G}_{\boldsymbol{\vartheta}} \in [\mathbf{G}_{\boldsymbol{\varphi}} \div \mathbf{G}_{\boldsymbol{\psi}}]$ then there are $\mathfrak{M} = \mathbf{G}_{\boldsymbol{\vartheta}}\boldsymbol{\lambda}$, $\mu = \mathbf{G}_{\boldsymbol{\vartheta}}\tilde{\boldsymbol{\lambda}}$, $\nu = \mathbf{G}_{\boldsymbol{\vartheta}}\tilde{\boldsymbol{\rho}}$ such that $\boldsymbol{\vartheta} : \mathfrak{M}/\mu \to \mathfrak{B}/\nu$. By A.3.3 $\mu$ is $\alpha$-invariant and $\underline{\boldsymbol{\varphi}}_{(\mu^{\mathfrak{A}})} = \boldsymbol{\chi}_{\mathfrak{A},\mu}\boldsymbol{\vartheta}$. In particular $\boldsymbol{\vartheta} = \boldsymbol{\chi}_{\mathfrak{A},\mu}^{-1}\underline{\boldsymbol{\varphi}}_{(\mu^{\mathfrak{A}})}$ so that $\tilde{\boldsymbol{\lambda}}$ is injective. But $\tilde{\boldsymbol{\lambda}}$ is also surjective because if $\mu \in [\alpha \div \gamma]_\alpha$ and $\boldsymbol{\vartheta} = \boldsymbol{\chi}_{\mathfrak{A},\mu}^{-1}\underline{\boldsymbol{\varphi}}_{(\mu^{\mathfrak{A}})}$, then $\underline{\boldsymbol{\vartheta}}_{(\gamma^{\mathfrak{M}})} = \boldsymbol{\chi}_{\mathfrak{A},\gamma^{\mathfrak{M}}}^{-1}\underline{\boldsymbol{\varphi}}_{(\gamma^{\mathfrak{A}})}$ hence $\mathbf{G}_{\boldsymbol{\vartheta}} \leqslant \mathbf{G}_{\boldsymbol{\psi}}$ because by A.3.3 again, applied to $\mathbf{G}_{\boldsymbol{\varphi}} \leqslant \mathbf{G}_{\boldsymbol{\psi}}$, we have $\underline{\boldsymbol{\varphi}}_{(\gamma^{\mathfrak{A}})} = \boldsymbol{\chi}_{\mathfrak{A},\gamma}\boldsymbol{\psi}$ and therefore $\underline{\boldsymbol{\vartheta}}_{(\gamma^{\mathfrak{M}})} = \boldsymbol{\chi}_{\mathfrak{M},\gamma}\boldsymbol{\psi}$. □

A.3.5. COROLLARY. *Let $\mathcal{I}$ be as in* (**4L**). *The following diagram commutes and all the arrows are isomorphisms.*

$$\begin{array}{ccc} \mathcal{I} & \xrightarrow{\boldsymbol{\lambda}} & [\mathfrak{A} \div \mathfrak{C}] \\ {\scriptstyle \tilde{\boldsymbol{\lambda}}}\big\downarrow & \swarrow {\scriptstyle \boldsymbol{\zeta}} & \\ [\alpha \div \gamma]_{\alpha} & & \end{array}$$

PROOF. By A.3.4 and A.1.4, $\tilde{\boldsymbol{\lambda}}$ and $\boldsymbol{\zeta}$ are isomorphisms. It is enough to prove that $\boldsymbol{\lambda} = \tilde{\boldsymbol{\lambda}}\boldsymbol{\zeta}^{-1}$. In fact, if $\mathbf{G}_{\boldsymbol{\vartheta}} \in \mathcal{I}$ then there are $\mathfrak{M} \in [\mathfrak{A} \div \mathfrak{C}]$, $\mu \in \mathcal{C}\mathrm{on}\,\mathfrak{M}$ such that $\boldsymbol{\vartheta} : \mathfrak{M}/\mu \to \mathfrak{D}/\delta$. Since $\mu$ is $\alpha$-invariant,

$$\mathbf{G}_{\boldsymbol{\vartheta}}\boldsymbol{\lambda} = \mathfrak{M} = \mathfrak{A}\,\mu = \mu\,\boldsymbol{\zeta}^{-1} = \mathbf{G}_{\boldsymbol{\vartheta}}\tilde{\boldsymbol{\lambda}}\boldsymbol{\zeta}^{-1}.$$

□

The intervals of type **4R** are defined in a similar way on the right.

———— ○ ————

From now on, let $\mathcal{I} = [\mathbf{G}_{\boldsymbol{\varphi}} \div \mathbf{G}_{\boldsymbol{\psi}}]$ where $\boldsymbol{\varphi} : \frac{\mathfrak{A}}{\alpha} \to \frac{\mathfrak{B}}{\beta}$ and $\boldsymbol{\psi} : \frac{\mathfrak{C}}{\gamma} \to \frac{\mathfrak{D}}{\delta}$ are isomorphisms. The left and right skeletons of $\mathcal{I}$ are

(A.3.A) $$(\alpha, \mathfrak{A}^2, \gamma^{\mathfrak{A}\gamma}, \mathfrak{C}^2) \quad \text{and} \quad (\beta, \mathfrak{B}^2, \delta^{\mathfrak{B}\delta}, \mathfrak{D}^2).$$

Their orders are $l = \big|\{\alpha, \mathfrak{A}^2, \gamma^{\mathfrak{A}\gamma}, \mathfrak{C}^2\}\big|$ and $r = \big|\{\beta, \mathfrak{B}^2, \delta^{\mathfrak{B}\delta}, \mathfrak{D}^2\}\big|$. The six given types of elementary intervals are determined by the left skeleton and the order of the right skeleton. In fact,

A.3.6. PROPOSITION. *$\mathcal{I}$ is of a given type if and only if $l$, $r$, and its left skeleton are as in Table A.A.*

PROOF. The proposition is trivial for types **2L**, **2R**. Since the "only if" is trivial, we just prove the "if".

Assume $l = 4$, $r = 2$. Since $\alpha < \mathfrak{A}^2$ and $\gamma < \mathfrak{C}^2$, we also have $\beta < \mathfrak{B}^2$ and $\delta^{\mathfrak{B}\delta} < \mathfrak{D}^2$. Since $r = 2$, it must be $\beta = \delta^{\mathfrak{B}\delta}$, $\mathfrak{B} = \mathfrak{D}$. By A.3.3 we have $\underline{\boldsymbol{\varphi}}_{(\gamma^{\mathfrak{A}})} = \boldsymbol{\chi}_{\mathfrak{A},\gamma}\boldsymbol{\psi}$. By A.2.5 we have $\delta^{\mathfrak{B}} = (\gamma^{\mathfrak{A}})(\boldsymbol{\varphi}, \boldsymbol{\varphi})$, hence $\gamma^{\mathfrak{A}} = \alpha$ because $\delta^{\mathfrak{B}} = \beta = \alpha(\boldsymbol{\varphi}, \boldsymbol{\varphi})$. This proves that $\mathcal{I}$ is of type **4L**. The case **4R** is similar therefore we now assume $1 < l = r \leqslant 3$.

If $\alpha < \gamma^{\mathfrak{A}\gamma} \leqslant \mathfrak{A}^2 = \mathfrak{C}^2$, then $\mathfrak{C} = \mathfrak{A}$ and by A.2.5, $\mathfrak{B}\delta = (\mathfrak{A}\gamma)\boldsymbol{\psi} = \mathfrak{C}\boldsymbol{\psi} = \mathfrak{B}$. By A.3.1 we have $\underline{\boldsymbol{\varphi}}_{\gamma} = \boldsymbol{\psi}$ hence $\mathcal{I}$ is of type **3A**.

If $\alpha = \gamma^{\mathfrak{A}\gamma} \leqslant \mathfrak{A}^2 < \mathfrak{C}^2$, then $\mathfrak{A} = \mathfrak{A}\gamma < \mathfrak{C}$ and by A.2.5, $\beta = \delta^{\mathfrak{B}}$, $\mathfrak{B}\delta < \mathfrak{D}$. If also $\mathfrak{B}\delta \neq \mathfrak{B}$, then $\beta < \delta^{\mathfrak{B}\delta} < \mathfrak{D}^2$ hence $l = r = 3$, that is, $\beta = \mathfrak{B}^2$ and $\alpha = \gamma^{\mathfrak{A}\gamma} < \mathfrak{A}^2$. Then $\beta = \delta^{\mathfrak{B}} = (\gamma^{\mathfrak{A}})(\boldsymbol{\varphi}, \boldsymbol{\varphi}) < \mathfrak{B}^2$, a contradiction. Therefore it must be $\mathfrak{B}\delta = \mathfrak{B}$ and by A.3.2, $\mathcal{I}$ is of type **3B**. □

TABLE A.A. Elementary intervals

| TYPE | $l$ | $r$ | | TYPE | $l$ | $r$ | |
|---|---|---|---|---|---|---|---|
| **2L** | 2 | 1 | | **2R** | 1 | 2 | |
| **4L** | 4 | 2 | | **4R** | 2 | 4 | |
| **3A** | 2, 3 | $l$ | $\alpha < \gamma^{\mathfrak{A}\gamma} \leqslant \mathfrak{A}^2 = \mathfrak{C}^2$ | **3B** | 2, 3 | $l$ | $\alpha = \gamma^{\mathfrak{A}\gamma} \leqslant \mathfrak{A}^2 < \mathfrak{C}^2$ |

A.3.7. REMARK. *There are non elementary intervals with* $1 < l = r \leqslant 3$.

In fact, in the variety of groups, If $A \lhd C$ and $\boldsymbol{\psi}$ is the identity of $C/A$, then the following intervals are all non elementary:

| interval | $l$ | $r$ | |
|---|---|---|---|
| $[\mathbf{G}_{\boldsymbol{\psi}} \div \mathbf{G}_{\boldsymbol{\psi}}]$ | 2 | 2 | $\mathcal{A} = \mathcal{C}^{AC} < A^2 = C^2$ |
| $[A \times A \div C \times C]$ | 2 | 2 | $\mathcal{A} = A^2 < \mathcal{C}^{AC} = C^2$ |
| $[\boldsymbol{\Delta}_A \div C \times C]$ | 3 | 3 | $\mathcal{A} < A^2 < \mathcal{C}^{AC} = C^2$ |
| $[\boldsymbol{\Delta}_A \div \mathbf{G}_{\boldsymbol{\psi}}]$ | 3 | 3 | $\mathcal{A} = A^2 = \mathcal{C}^{AC} < C^2$ |

We now show that the six given types of elementary intervals arise in a natural way.

A.3.8. PROPOSITION. $\mathbf{G}_{\boldsymbol{\varphi}} \leqslant \mathbf{G}_{\boldsymbol{\psi}}$ *if and only if the following diagram makes sense and commutes:*

(A.3.B)
$$\begin{array}{ccc}
\frac{\mathfrak{A}\gamma}{\gamma^{\mathfrak{A}\gamma}} & \xrightarrow{\overline{\boldsymbol{\psi}}^{(\mathfrak{A}\gamma)}} & \frac{\mathfrak{B}\delta}{\delta^{\mathfrak{B}\delta}} \\
{\scriptstyle \chi_{\mathfrak{A},\gamma}}\uparrow & & \uparrow{\scriptstyle \chi_{\mathfrak{B},\delta}} \\
\frac{\mathfrak{A}}{\gamma^{\mathfrak{A}}} & \xrightarrow[\underline{\boldsymbol{\varphi}}_{(\gamma^{\mathfrak{A}})}]{} & \frac{\mathfrak{B}}{\delta^{\mathfrak{B}}}
\end{array}$$

PROOF. For short, put $\underline{\boldsymbol{\varphi}} = \underline{\boldsymbol{\varphi}}_{(\gamma^{\mathfrak{A}})}$, $\overline{\boldsymbol{\psi}} = \overline{\boldsymbol{\psi}}^{(\mathfrak{A}\gamma)}$, $\boldsymbol{\sigma} = \chi_{\mathfrak{A},\gamma}$ and $\boldsymbol{\tau} = \chi_{\mathfrak{B},\delta}$. The diagram makes sense when $\alpha \leqslant \gamma$ (in particular $\mathfrak{A} \leqslant \mathfrak{C}$), $\beta \leqslant \delta$ (in particular $\mathfrak{B} \leqslant \mathfrak{D}$), $(\mathfrak{A}\gamma)\boldsymbol{\psi} = \mathfrak{B}\delta$ and $\gamma^{\mathfrak{A}}(\boldsymbol{\varphi}, \boldsymbol{\varphi}) = \delta^{\mathfrak{B}}$. By A.2.5 this conditions are all satisfied when $\mathbf{G}_{\boldsymbol{\varphi}} \leqslant \mathbf{G}_{\boldsymbol{\psi}}$, and it is routine to check that in fact the diagram commutes. Conversely, if the diagram makes sense, then

$$\mathbf{G}_{\boldsymbol{\varphi}} \leqslant \mathbf{G}_{\underline{\boldsymbol{\varphi}}} \qquad \text{by A.3.1, interval trivial or } \mathbf{3A},$$
$$\mathbf{G}_{\underline{\boldsymbol{\varphi}}} \leqslant \mathbf{G}_{\boldsymbol{\sigma}^{-1}\underline{\boldsymbol{\varphi}}} \qquad \text{by A.3.3, interval trivial or } \mathbf{xL},$$
$$\mathbf{G}_{\boldsymbol{\sigma}^{-1}\underline{\boldsymbol{\varphi}}} \leqslant \mathbf{G}_{\boldsymbol{\sigma}^{-1}\underline{\boldsymbol{\varphi}}\boldsymbol{\tau}} \qquad \text{by A.3.3, interval trivial or } \mathbf{xR},$$

if also the diagram commutes then $\boldsymbol{\sigma}^{-1}\underline{\boldsymbol{\varphi}}\boldsymbol{\tau} = \overline{\boldsymbol{\psi}}$ and

$$\mathbf{G}_{\boldsymbol{\sigma}^{-1}\underline{\boldsymbol{\varphi}}\boldsymbol{\tau}} = \mathbf{G}_{\overline{\boldsymbol{\psi}}} \leqslant \mathbf{G}_{\boldsymbol{\psi}} \qquad \text{by A.3.2, interval trivial or } \mathbf{3B}.$$

□

This also proves our main theorem:

**Theorem A.** *If $\mathcal{I}$ is an interval of $\mathcal{S}\mathrm{ub}(\mathfrak{L} \times \mathfrak{R})$ then*

$$\mathcal{I} = \mathcal{I}_1 \rtimes \mathcal{I}_2 \rtimes \mathcal{I}_3 \rtimes \mathcal{I}_4$$

*where $\mathcal{I}_i$ are trivial or elementary of type* **3A**, **xL**, **xR**, **3B** *respectively.*

This theorem is of little use for "large" intervals. For example, when $L$, $R$ are groups, it says that

$$[1 \div L \times R] = [1 \div L \times 1] \rtimes [L \times 1 \div L \times R],$$

which is nothing new. But Theorem A.3 can be of some help with "smaller" intervals. For example we have straightaway:

A.3.9. Corollary. *Simple intervals are elementary.*

This is conclusive, because

A.3.10. Corollary. *The elementary intervals are exactly the ones listed in Table A.A.*

Proof. Assume for example that $\boldsymbol{\lambda}$ and $\tilde{\boldsymbol{\lambda}}_\bullet$ of $\mathcal{I}$ are injective. Let

$$\mathcal{I} = \mathcal{I}_1 \rtimes \mathcal{I}_2 \rtimes \mathcal{I}_3 \rtimes \mathcal{I}_4$$

as in Theorem A.3. Since $\mathcal{I}_1\boldsymbol{\lambda}$, $\mathcal{I}_3\boldsymbol{\lambda}$, $\mathcal{I}_4\tilde{\boldsymbol{\lambda}}_\bullet$ are trivial, so are $\mathcal{I}_1$, $\mathcal{I}_2$, $\mathcal{I}_4$, therefore $\mathcal{I} = \mathcal{I}_2$ is elementary of type **xL**. One proves the other cases in the same way. $\square$

The homomorphisms $\boldsymbol{\rho}$ and $\tilde{\boldsymbol{\lambda}}_\bullet$ can never be both injective. This follows from a more general fact.

A.3.11. Proposition. *Let $\mathcal{I}$ be an interval and let $\boldsymbol{\mu}, \boldsymbol{\nu} \in \{\boldsymbol{\rho}, \tilde{\boldsymbol{\lambda}}_\bullet\}$ with $\boldsymbol{\mu} \neq \boldsymbol{\nu}$. Then $\mathcal{I}\boldsymbol{\mu}$ is trivial if and only if $\mathcal{I} \overset{\boldsymbol{\nu}}{\cong} \mathcal{I}\boldsymbol{\nu}$.*

Proof. Let $\mathcal{I} = \mathcal{I}_1 \rtimes \mathcal{I}_2 \rtimes \mathcal{I}_3 \rtimes \mathcal{I}_4$ as in Theorem A.3. If $\mathcal{I} \overset{\boldsymbol{\rho}}{\cong} \mathcal{I}\boldsymbol{\rho}$, then $\mathcal{I}_1$, $\mathcal{I}_2$ are trivial hence $\mathcal{I}\tilde{\boldsymbol{\lambda}}_\bullet \subseteq \mathcal{I}_3\tilde{\boldsymbol{\lambda}}_\bullet \rtimes \mathcal{I}_4\tilde{\boldsymbol{\lambda}}_\bullet$ is trivial. Similarly, if $\mathcal{I} \overset{\tilde{\boldsymbol{\lambda}}_\bullet}{\cong} \mathcal{I}\tilde{\boldsymbol{\lambda}}_\bullet$, then $\mathcal{I}_3$, $\mathcal{I}_4$ are trivial and so is $\mathcal{I}\boldsymbol{\rho}$. Assume now $\mathcal{I}\tilde{\boldsymbol{\lambda}}_\bullet$ trivial and assume $\mathfrak{H}_1\boldsymbol{\rho} = \mathfrak{H}_2\boldsymbol{\rho}$ with $\mathfrak{H}_i \in \mathcal{I}$. Since $\boldsymbol{\rho}$ is a join-hom, we have $(\mathfrak{H}_1 \vee \mathfrak{H}_2)\boldsymbol{\rho} = \mathfrak{H}_i\boldsymbol{\rho}$, therefore $\tilde{\boldsymbol{\lambda}}_\bullet$, $\boldsymbol{\rho}$ are trivial on the intervals $[\mathfrak{H}_i \div \mathfrak{H}_1 \vee \mathfrak{H}_2]$ which therefore need to be trivial. This forces $\mathfrak{H}_1 = \mathfrak{H}_1 \vee \mathfrak{H}_2 = \mathfrak{H}_2$, thus $\mathcal{I} \overset{\boldsymbol{\rho}}{\cong} \mathcal{I}\boldsymbol{\rho}$. Finally, if $\mathcal{I}\boldsymbol{\rho}$ is trivial, then $\mathcal{I} = \mathcal{I}_1 \rtimes \mathcal{I}_2$ and the result follows by A.3.4. $\square$

A.3.12. Corollary. *$\mathcal{I}$ is elementary if and only if its four canonical homomorphisms are injective and trivial in pairs.*

The intervals of $\mathcal{S}\mathrm{ub}(\mathfrak{L} \times \mathfrak{R})$ which are not already intervals of $\mathcal{S}\mathrm{ub}\,\mathfrak{L}$, $\mathcal{S}\mathrm{ub}\,\mathfrak{R}$, or bottom invariant intervals of $\mathcal{SC}\mathrm{on}\,\mathfrak{L}$, $\mathcal{SC}\mathrm{on}\,\mathfrak{R}$, are called *novelty.* We are usually interested in intervals which are novelty, hence in intervals $\mathcal{I}$ such that $\left|\mathcal{I}\boldsymbol{\mu}\right| > 1$ for all $\boldsymbol{\mu} \in \{\boldsymbol{\lambda}, \tilde{\boldsymbol{\lambda}}_{\bullet}, \boldsymbol{\rho}, \tilde{\boldsymbol{\rho}}_{\bullet}\}$. We end this section providing a necessary condition for a shortcut interval to be a novelty.

A.3.13. *Definition.* We say that $\mathcal{I}$ has a shortcut **xX**–**yY**, when $\mathcal{I} = \mathcal{I}_1 \rtimes \mathcal{I}_2$ with $\mathcal{I}_1$ simple of type **xX** and $\mathcal{I}_2$ simple of type **yY**.

A.3.14. PROPOSITION. *If a shortcut interval $\mathcal{I}$ is a novelty, then $\mathcal{I}$ has one of the following shortcuts:*

$$\mathbf{2L}\text{–}\mathbf{2R},\ \mathbf{3B}\text{–}\mathbf{3A},\ \mathbf{3A}\text{–}\mathbf{3B},\ \mathbf{4L}\text{–}\mathbf{4R}.$$

PROOF. Elementary types are closed by composition thus no **xX**–**xX** is novelty. Also no interval with shortcuts **xL**–**3B**or **3B**–**xL** is a novelty because its $\tilde{\boldsymbol{\rho}}_{\bullet}$ is trivial. For a similar reason no interval with shortcuts **3A**–**xR**, **xR**–**3A**, **3A**–**xL**, **xL**–**3A**, **xR**–**3B**or **3B**–**xR** is a novelty. There are no intervals at all with shortcuts **2X**–**4Y** or **4Y**–**2X** because this would require a cartesian subalgebra which is not cartesian! If $\mathcal{I}$ has a shortcut **2R**–**2L**, then it also has a shortcut **2L**–**2R**. Similarly, if $\mathcal{I}$ has a shortcut **4R**–**4L**, then it also has a shortcut **4L**–**4R**. Of the 36 possible formal combination (but a few of them could not be realised anyway), we remain with just

$$36 - 6 - 4 - 12 - 8 - 2 = 4$$

shortcuts, which are the one listed. □

## A.4. Shortcut intervals in the subgroup lattice of direct products

We are going to apply the results of the previous section to the variety of groups. For the reminder of this section $L$, $R$ are groups. Say $\mathbf{G}_{\varphi}$ and $\mathbf{G}_{\psi}$ two Goursat subgroups of $L \times R$, where

$$\varphi : \frac{A}{\mathcal{A}} \longrightarrow \frac{B}{\mathcal{B}} \qquad \text{and} \qquad \psi : \frac{C}{\mathcal{C}} \longrightarrow \frac{D}{\mathcal{D}}$$

are isomorphisms of sections of $L$, $R$ respectively. Applying A.3.9, A.3.10 we get sufficient and necessary conditions for $\mathbf{G}_{\varphi}$ being maximal in $\mathbf{G}_{\psi}$.

A.4.1. THEOREM. $\mathbf{G}_{\varphi} \lessdot \mathbf{G}_{\psi}$ *if and only if one of the following conditions holds:*

(a) $\mathcal{C} \leqslant A \lessdot C$ *and* $\varphi = \overline{\psi}^{A}$,

(b) $\mathcal{C} \nleqslant A \lessdot C$ *and* $\varphi = \chi_{A,\mathcal{C}}\psi$,

(c) $\mathcal{D} \nleqslant B \lessdot D$ *and* $\varphi\chi_{B,\mathcal{D}} = \psi$,

(d) $[\mathcal{A} \div \mathcal{C}]_{A}$ *is simple and* $\psi = \underline{\varphi}_{\mathcal{C}}$ *(hence* $\mathcal{A} \lhd C$ *and* $\mathcal{C}/\mathcal{A} \underset{min}{\lhd} C/\mathcal{A}$*).*

A.4.2. Corollary. *A subgroup $G$ is maximal in $L \times R$ if and only if one of the following conditions holds:*

(a) $G = A \times R$ *with* $A \lessdot L$,

(b) $G = L \times B$ *with* $B \lessdot R$,

(c) $G = \mathbf{G}_{\varphi}$ *where* $\varphi$ *is a group isomorphism* $\varphi : \frac{L}{\mathcal{A}} \longrightarrow \frac{R}{\mathcal{B}}$ *and* $L/\mathcal{A}$ *is simple, that is,* $\mathcal{A}$ *is maximal among the normal subgroups of* $L$.

A.4.3. Corollary. *Let* $\varphi : \frac{L}{\mathcal{A}} \longrightarrow \frac{R}{\mathcal{B}}$ *be an isomorphism where* $\mathcal{A}$ *is a maximal normal subgroup of* $L$. *A subgroup* $G$ *is maximal in* $\mathbf{G}_{\varphi}$ *if and only if one of the following conditions holds:*

(a) $G = \mathbf{G}_{\psi}$ *where* $\psi : C/(C \cap \mathcal{A}) \longrightarrow R/\mathcal{B}$ *is the canonical isomorphism induced by* $\varphi$, *being* $C \lessdot L$ *and* $\mathcal{A} \nleqslant C$;

(b) $G = \mathbf{G}_{\psi}$ *where* $\psi : L/\mathcal{A} \longrightarrow D/(D \cap \mathcal{B})$ *is the canonical isomorphism induced by* $\varphi$, *being* $D \lessdot R$ *and* $\mathcal{B} \nleqslant D$;

(c) $G = \mathbf{G}_{\overline{\varphi}^A}$ *with* $A \lessdot L$ *and* $\mathcal{A} \leqslant A$.

As we often refer to the bottom invariant interval of $\mathcal{SC}\mathrm{on}\, L$, $\mathcal{SC}\mathrm{on}\, R$, that is, to interval of subgroups which are normalised by some other subgroup, we feel obliged to state the following.

A.4.4. Proposition. $\mathbf{G}_{\varphi}$ *is* $\mathbf{G}_{\psi}$*–invariant (by conjugation) if and only if the following diagram makes sense and commutes.*

$$\begin{array}{ccc} \frac{C}{\mathcal{C}} & \xrightarrow{\ \boldsymbol{t}\ } & \mathrm{Aut}\,\frac{A}{\mathcal{A}} \\ {\scriptstyle \psi}\downarrow & & \downarrow{\scriptstyle \boldsymbol{t}_{\varphi}} \\ \frac{D}{\mathcal{D}} & \xrightarrow{\ \boldsymbol{t}\ } & \mathrm{Aut}\,\frac{B}{\mathcal{B}} \end{array}$$

Explanation. The horizontal arrows are given by the conjugation. For example

$$(\mathcal{C}c)\boldsymbol{t} : \mathcal{A}a \mapsto \mathcal{A}(c^{-1}ac).$$

While $\boldsymbol{t}_{\varphi}$ is defined by $\boldsymbol{f} \mapsto \varphi^{-1}\boldsymbol{f}\varphi$. The diagram makes sense when the arrows are well defined, that is,

(A.4.A) $$A^C = A, \quad [A, \mathcal{C}] \leqslant \mathcal{A} \qquad \text{and} \qquad B^D = D, \quad [B, \mathcal{D}] \leqslant \mathcal{B}.$$

In particular $\mathbf{G}_{\varphi} \trianglelefteq L \times R$ if and only if $[A, L] \leqslant \mathcal{A}$ and $[B, R] \leqslant \mathcal{B}$.

**Three shortcut intervals.** We give examples of intervals with a shortcut **3B**–**3A** arising from three different settlements. In all the examples we consider the interval $\mathcal{I} = [\mathbf{\Delta}_A \div (\mathcal{C} \times 1)\mathbf{\Delta}_C]$ where $A < C$ are finite groups and $[A \div C]$, $[1 \div \mathcal{C}]_C$ are simple. In particular $\mathcal{C} \cong S_1 \times \cdots \times S_n$ for convenient isomorphic finite simple groups $S_i$. We also set, for short, $\boldsymbol{\psi}$ to be the identity of $C/\mathcal{C}$ and $\boldsymbol{\varphi}$ to be the identity of $A$, so that $\mathcal{I} = [\mathbf{G}_{\varphi} \div \mathbf{G}_{\psi}]$.

A.4.5. Example. Case $1 = A \cap \mathcal{C} \leqslant A \lessdot A\mathcal{C} = C$.

We prove that $\mathcal{I} \cong \mathcal{M}_{p^r+1}$ for some prime $p$ and $r \geqslant 0$ and that any such $\mathcal{M}_{p^r+1}$ can be represented in this way.

$\mathcal{I}$ has also shortcuts **xL**–**xR**, **xR**–**xL**, in fact, in $\mathcal{I}$ there are

$$(1 \times \mathcal{C})\boldsymbol{\Delta}_A = (\mathcal{C} \times 1)\boldsymbol{\Delta}_C \cap (A \times C), \tag{A.4.B}$$

$$(\mathcal{C} \times 1)\boldsymbol{\Delta}_A = (\mathcal{C} \times 1)\boldsymbol{\Delta}_C \cap (C \times A) \tag{A.4.C}$$

which are both atoms and maximal elements of $\mathcal{I}$ because of A.3.3. Moreover, if $\boldsymbol{\vartheta} \in \mathrm{N}_{\mathrm{Aut}\,C}\,\mathcal{C}$, the subgroup of the automorphisms of $C$ which fix $\mathcal{C}$, then

$$\boldsymbol{\chi}_{A,\mathcal{C}}\underline{\boldsymbol{\vartheta}}_{\mathcal{C}} = \overline{\boldsymbol{\vartheta}}^{A}\boldsymbol{\chi}_{A,\mathcal{C}}. \tag{A.4.D}$$

If also $\boldsymbol{\vartheta} \in \mathrm{C}_{\mathrm{Aut}\,C}\,A$, the subgroup of the automorphisms of $C$ whose restriction to $A$ is $\boldsymbol{\varphi}$, then $\boldsymbol{\varphi} = \overline{\boldsymbol{\vartheta}}^{A}$ and $\underline{\boldsymbol{\vartheta}}_{\mathcal{C}} = \boldsymbol{\chi}_{A,\mathcal{C}}^{-1}\boldsymbol{\varphi}\boldsymbol{\chi}_{A,\mathcal{C}} = \boldsymbol{\psi}$, hence

$$\boldsymbol{\Delta}_A = \mathbf{G}_{\varphi} \lessdot \mathbf{G}_{\vartheta} \lessdot \mathbf{G}_{\psi} = (\mathcal{C} \times 1)\boldsymbol{\Delta}_{C}.$$

Beside, if $\mathbf{G}_{\varphi} < \mathbf{G}_{\vartheta} < \mathbf{G}_{\psi}$ and $\mathbf{G}_{\vartheta}$ is not contained in (A.4.B), (A.4.C), then it has to be a subcartesian product of $C^2$ and by A.3.1, $\underline{\boldsymbol{\vartheta}}_{\mathcal{C}} = \boldsymbol{\psi}$. In particular $\mathcal{C}\boldsymbol{\vartheta} = \mathcal{C}$. Since $\mathbf{G}_{\vartheta} < \mathbf{G}_{\psi}$, $\mathcal{C}_1 = \mathbf{G}_{\vartheta}\tilde{\boldsymbol{\lambda}}_1$ and $\mathcal{C}_2 = \mathbf{G}_{\vartheta}\tilde{\boldsymbol{\rho}}_1$ are less than $\mathcal{C}$, hence $\mathcal{C}_1 = 1 = \mathcal{C}_2$ and $\boldsymbol{\vartheta} \in \mathrm{Aut}\,C$. Thus, by (A.4.D), $\boldsymbol{\vartheta} \in \mathrm{N}_{\mathrm{Aut}\,C}\,\mathcal{C} \cap \mathrm{C}_{\mathrm{Aut}\,C}\,A$ which is canonically isomorphic to $\mathrm{Aut}_A\,\mathcal{C}$, the group of the automorphisms of $\mathcal{C}$ which commute with those induced by conjugation by $A$. It is known[1] that the order of this group must be $p^r - 1$ for some prime $p$, because $\mathcal{C} \times 1$, $1 \times \mathcal{C}$ are two distinct normal subgroups of $\mathbf{G}_{\psi}$. However, we now show this by distinguishing two cases:

(1) $\mathcal{C}$ is non abelian. By [**Bad93**, 3.2][2] and the maximality of $A$, for each $s \in S_1$, there is $a \in A$ such that $\boldsymbol{t}_s = \boldsymbol{t}_a$ where $\boldsymbol{t}_s$ is the conjugation by $s$ in $S_1$ and $\boldsymbol{t}_a$ is the restriction of the conjugation by $a$ in $C$ to $S_1$. This implies $\mathrm{C}_{\mathrm{Aut}\,C}\,A \leqslant \mathrm{C}_{\mathrm{Aut}\,C}\,S_1 = 1$ because $C = \langle S_1, A \rangle$.

(2) $\mathcal{C}$ is elementary abelian of exponent $p$. Since $\mathcal{C}$ is a simple $A$-module, by Schur Lemma, $\mathrm{End}_A\,\mathcal{C}$ is a division ring (of exponent $p$), hence $\left|\mathrm{End}_A\,\mathcal{C}\right|$ is a power of $p$. But $\left|\mathrm{Aut}_A\,\mathcal{C}\right| = \left|\mathrm{End}_A\,\mathcal{C}\right| - 1$.

Finally, let $\mathbb{F}$ be a finite field of order $p^r$, $\mathcal{C} = \mathbb{F}$, $A = \mathbb{F}^{\times}$ the multiplicative group acting on $\mathcal{C}$ by multiplication on the right and $C = \mathcal{C} \rtimes A$. Then $\mathrm{Aut}_A\,\mathcal{C} = A$ has

[1] see for example [**Pál95**].

[2] depend on "Schreier conjecture".

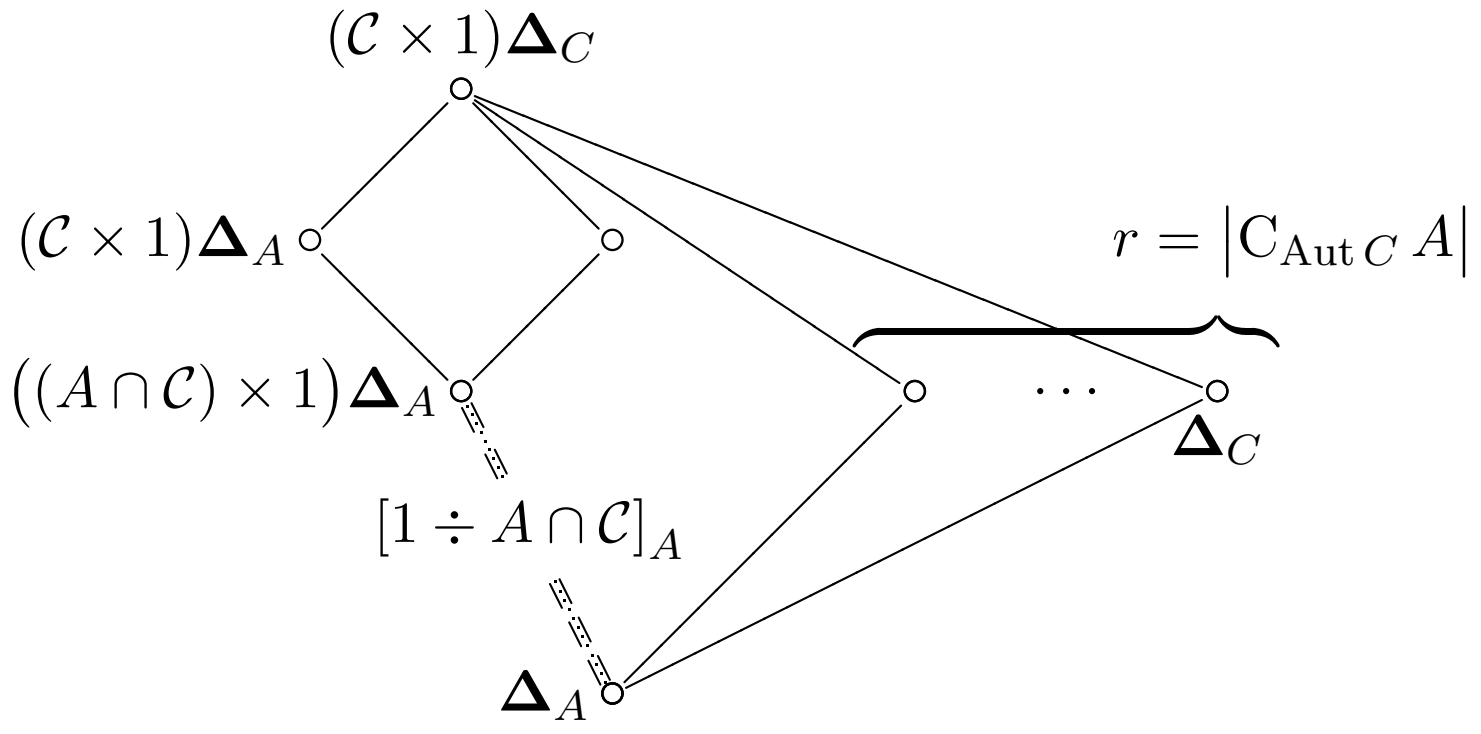

$1 < A \cap \mathcal{C} \leqslant A \lessdot A\mathcal{C} = C, \qquad [1 \div \mathcal{C}]_C$ simple

FIGURE A.F. A shortcut interval

order $p^r - 1$. It is not difficult to check that this is in disguise the example 3.1 of [**Pál95**].

A.4.6. EXAMPLE. Case $1 < A \cap \mathcal{C} \leqslant A \lessdot A\mathcal{C} = C$, see Figure A.F.

Here $\mathcal{I} = \mathcal{I}_1 \rtimes \mathcal{I}_2$ where $\max \mathcal{I}_1 = \mathbf{G}_{\underline{\varphi}_{(A \cap \mathcal{C})}} = \big((A \cap \mathcal{C}) \times 1\big)\boldsymbol{\Delta}_A$. By A.3.1, $\mathcal{I}_1 \cong [1 \div A \cap \mathcal{C}]_A$. Assume now $\mathbf{G}_{\boldsymbol{\vartheta}} \in \mathcal{I} - \mathcal{I}_1$, $\boldsymbol{\vartheta} : M/\mathsf{M} \to N/\mathsf{N}$. Intersecting (A.4.B) with $A \times A$, one sees that $\max \mathcal{I}_1 = \mathbf{G}_{\boldsymbol{\psi}} \cap (A \times A)$ therefore $M$ and $N$ cannot be both equal to $A$. Then we distinguish three cases:

(1) $M = A$, $N = C$. We apply A.3.3 to $\mathbf{G}_{\boldsymbol{\vartheta}} \leqslant \mathbf{G}_{\boldsymbol{\psi}}$ to see that $\underline{\boldsymbol{\vartheta}}_{(A \cap \mathcal{C})} = \boldsymbol{\chi}_{A,\mathcal{C}}\boldsymbol{\psi}$. In particular $(A \cap \mathcal{C})\boldsymbol{\vartheta} = \mathcal{C}$. If $\mathsf{M} < A \cap \mathcal{C}$, then $\mathsf{N} < \mathcal{C}$ hence $\mathsf{N} = 1$ because $\mathsf{N} \trianglelefteq C$. We apply A.3.2 to $\mathbf{G}_{\boldsymbol{\varphi}} \leqslant \mathbf{G}_{\boldsymbol{\vartheta}}$ to see that $\underline{\boldsymbol{\varphi}}_{\mathsf{M}} = \overline{\boldsymbol{\vartheta}}^A = \boldsymbol{\vartheta}$, which contradicts $N = C$. Therefore $\mathsf{M} = A \cap \mathcal{C}$ and $\boldsymbol{\vartheta} = \boldsymbol{\chi}_{A,\mathcal{C}}\boldsymbol{\psi}$, that is, $\mathbf{G}_{\boldsymbol{\vartheta}}$ is equal to (A.4.B).
(2) $M = C$, $N = A$. In an entirely similar way one proves that $\mathbf{G}_{\boldsymbol{\vartheta}}$ is equal to (A.4.C).
(3) $M = C = N$. In particular either $\mathsf{M} = 1$ or $\mathsf{M} = \mathcal{C}$. We apply A.3.1 to $\mathbf{G}_{\boldsymbol{\vartheta}} \leqslant \mathbf{G}_{\boldsymbol{\psi}}$ to see that $\underline{\boldsymbol{\vartheta}}_{\mathcal{C}} = \boldsymbol{\psi}$. If $\mathsf{M} = \mathcal{C}$, then $\boldsymbol{\vartheta} = \boldsymbol{\psi}$, otherwise $\mathcal{C} = \mathcal{C}\boldsymbol{\vartheta} > \mathsf{M}\boldsymbol{\vartheta} = \mathsf{N}$, hence $\boldsymbol{\vartheta} \in \mathrm{Aut}\, C$. Applying A.3.2 to $\mathbf{G}_{\boldsymbol{\varphi}} \leqslant \mathbf{G}_{\boldsymbol{\vartheta}}$, one has $\boldsymbol{\vartheta} \in \mathrm{C}_{\mathrm{Aut}\,C}\, A$. Beside, if $\boldsymbol{\vartheta} \in \mathrm{C}_{\mathrm{Aut}\,C}\, A$, then $\mathbf{G}_{\boldsymbol{\varphi}} \lessdot \mathbf{G}_{\boldsymbol{\vartheta}} \lessdot \mathbf{G}_{\boldsymbol{\psi}}$.

Thus $\mathcal{I} - \mathcal{I}_1 - \mathcal{I}_2$ is in bijection with $\mathrm{C}_{\mathrm{Aut}\,C}\, A$, $\mathcal{I}_2 \cong \mathcal{M}_2$ and $\mathcal{I}$ is as in Figure A.F.

A.4.7. EXAMPLE. Case $1 < \mathcal{C} \leqslant A \lessdot C$, see Figure A.G.

Here $\mathcal{I} = \mathcal{I}_1 \rtimes \mathcal{I}_2$ where

$$\mathcal{I}_1 = [\boldsymbol{\Delta}_A \div (\mathcal{C} \times 1)\boldsymbol{\Delta}_A] \cong [1 \div \mathcal{C}]_A \qquad \text{and} \qquad \mathcal{I}_2 = [(\mathcal{C} \times 1)\boldsymbol{\Delta}_A \div (\mathcal{C} \times 1)\boldsymbol{\Delta}_C].$$

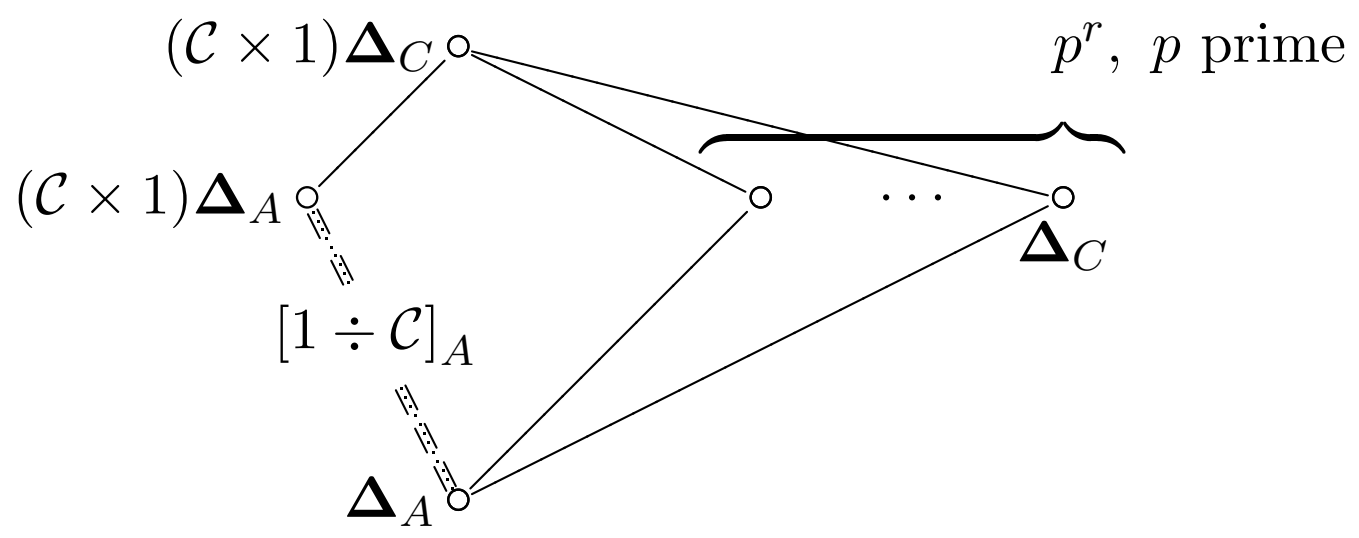

$1 < \mathcal{C} \leqslant A \lessdot C, \qquad [1 \div \mathcal{C}]_C$ simple

FIGURE A.G. Another shortcut interval

Arguing as in the previous example, one sees that $\mathbf{G}_{\vartheta} \in \mathcal{I} - \mathcal{I}_1 - \mathcal{I}_2$ if and only if $\boldsymbol{\vartheta} \in \operatorname{Aut} C$, $\boldsymbol{\varphi} = \overline{\boldsymbol{\vartheta}}^A$ and $\boldsymbol{\psi} = \underline{\boldsymbol{\vartheta}}_{\mathcal{C}}$. By [**Rob82**, §4], the set of these $\boldsymbol{\vartheta}$ is in bijection with a subgroup of $\operatorname{Der}(C/\mathcal{C}, \mathrm{Z}(\mathcal{C}))$ which has prime power order because either $\mathcal{C}$ is non abelian and $\mathrm{Z}(\mathcal{C}) = 1$, or $\mathcal{C}$ is elementary abelian and so is $\operatorname{Der}(C/\mathcal{C}, \mathrm{Z}(\mathcal{C}))$.

———— ∘ ————

A.4.8. REMARK. *In example A.4.6, $\mathrm{C}_{\operatorname{Aut} C} A$ may be non trivial and in example A.4.7, $[1 \div \mathcal{C}]_A$ may be a chain of any arbitrary length.*

In fact, if $C$ is the alternating group on $n$ symbols, then the centralizer in $\operatorname{Aut} C$ of a transposition is a maximal subgroup $A$ of $C$ so that $\mathrm{C}_{\operatorname{Aut} C} A$ has order 2.

And for each $m \geqslant 3$, let $p$ be a prime larger than $m$. Let $V_m$ be the subspace of the homogeneous component of degree $m-1$ of $\mathbb{F}_p[X, Y]$, the commutative polynomial ring in two variables $X$, $Y$, on the prime field $\mathbb{F}_p$. Let $Q = \mathrm{SL}(2, p)$ and $P = Q \cap \mathrm{T}(2, p)$ where $\mathrm{T}(2, p)$ are the upper triangular $2 \times 2$ matrixes. Then $Q$ acts irreducibly on the left of $V_m$ and if we put $C = V_m \rtimes Q$, $A = V_m \rtimes P$ and $\mathcal{C} = V_m$, then we are in the assumptions of example A.4.7 and $[1 \div \mathcal{C}]_A$ is a chain of length $m$ [**Alp86**, pp. 15–16].

We are now ready to state our second main result:

**Theorem B.** *If a shortcut interval $\mathcal{I}$ of $\mathcal{S}\mathrm{ub}(L \times R)$ is a novelty, then either $\mathcal{I} \cong \mathcal{M}_1$ or $\mathcal{I}$ is isomorphic to one of the examples A.4.5–A.4.7.*

PROOF. Assume that $\mathcal{I} = [\mathbf{G}_{\varphi} \div \mathbf{G}_{\psi}]$.

If $\mathcal{I}$ has a shortcut **3B**–**3A**, then there must be an isomorphism $\boldsymbol{\vartheta} : C/\mathcal{A} \to D/\mathcal{B}$. Thus we may assume, up to algebra isomorphisms, that $\mathcal{A} = 1 = \mathcal{B}$, $C = D$ and $\boldsymbol{\vartheta}$ is the identity of $C$. Therefore $\mathcal{I}$ must be as in one of the examples.

If $\mathcal{I}$ has a shortcut **3A**–**3B** and $\mathcal{I} \neq \mathcal{M}_1$, let $\mathbf{G}_{\vartheta} \in \mathcal{I} - \{\mathbf{G}_{\varphi}, \mathbf{G}_{\psi}, \mathbf{G}_{(\overline{\psi}^A)}\}$ with $\vartheta : M/\mathsf{M} \to N/\mathsf{N}$. Since

$$\mathbf{G}_{(\overline{\psi}^A)} = (\mathcal{C} \times 1)\mathbf{G}_{\varphi} = (1 \times \mathcal{C})\mathbf{G}_{\varphi} = (A \times D) \cap \mathbf{G}_{\psi} = (D \times A) \cap \mathbf{G}_{\psi},$$

it must be $M = C = N$ and $\mathsf{M} = \mathcal{A} = \mathsf{N}$. Therefore $\mathcal{I}$ has a shortcut **3B**–**3A** too and we conclude as before.

If $\mathcal{I}$ has a shortcut **xL**–**xR** and $\mathcal{I} \neq \mathcal{M}_2$, then there has to be $\mathbf{G}_{\vartheta} \in \mathcal{I}$ such that $\vartheta : C/\mathsf{M} \to D/\mathsf{N}$ and $M < \mathcal{C}$, $\mathsf{N} < \mathcal{D}$. As $[\mathcal{A} \div \mathcal{C}]_A$ is simple, it must be $\mathcal{A} = \mathsf{M} \lhd C$ and similarly $\mathcal{B} = \mathsf{N} \lhd D$. Therefore again, $\mathcal{I}$ has a shortcut **3B**–**3A**.

By A.3.14 the proof is complete □

An elementary interval of type **3A** is a modular lattice. Consequently, if it also has a shortcut, it must be isomorphic to $\mathcal{M}_r$ for some $r$. It is not difficult to prove the same result for each interval which has a shortcut **3A**–**xL** or **3A**–**xR**. However, in each case the maximum of the interval has two distinct normal subgroups and as a consequence $r - 1$ must be a power of prime (see [**Pál95**]). This completes the reduction of the search of the shortcut intervals of $\mathcal{S}\mathrm{ub}(L \times R)$ to the shortcut intervals of $\mathcal{S}\mathrm{ub}\, L$, $\mathcal{S}\mathrm{ub}\, R$.

Finally, we end this section with the following corollary of Theorem A.4.

A.4.9. Corollary. *If a shortcut interval of* $\mathcal{S}\mathrm{ub}(L \times R)$ *has more than* 4 *points, then it is already an interval of* $\mathcal{S}\mathrm{ub}(L^2)$ *or* $\mathcal{S}\mathrm{ub}(R^2)$.

Proof. If such an interval, say $\mathcal{I}$, is a novelty then this follows by the theorem. Otherwise $\mathcal{I}$ is an interval of $\mathcal{S}\mathrm{ub}\, L$, or $\mathcal{S}\mathrm{ub}\, R$, or a bottom invariant interval of $\mathcal{SC}\mathrm{on}\, L$ or $\mathcal{SC}\mathrm{on}\, R$. But the bottom invariant intervals of $\mathcal{SC}\mathrm{on}\, L$, for example, are intervals of $\mathcal{S}\mathrm{ub}(L^2)$ for A.3.4. □

# Index